\documentclass[10pt,english]{article}
\usepackage{ae,aecompl}
\usepackage[T1]{fontenc}
\usepackage[latin9]{inputenc}
\usepackage[a4paper]{geometry}
\usepackage[active]{srcltx}
\usepackage{babel}
\usepackage{array}
\usepackage{rotating}
\usepackage{multirow}
\usepackage{amsmath}
\usepackage{amssymb}
\usepackage{cancel}
\usepackage{graphicx}
\usepackage[unicode=true,pdfusetitle,
 bookmarks=true,bookmarksnumbered=false,bookmarksopen=false,
 breaklinks=false,pdfborder={0 0 0},pdfborderstyle={},backref=false,colorlinks=false]
 {hyperref}
\usepackage{breakurl}

\makeatletter

\providecommand{\tabularnewline}{\\}
\newcommand{\lyxdot}{.}

\numberwithin{equation}{section}
\numberwithin{figure}{section}

\date{}
\usepackage{abstract}

\usepackage{eurosym}
\title{Nitsche's Method for resolving boundary conditions on embedded interfaces using XFEM in Code Aster}
\author{Nanda Gopala Kilingar}
\makeatother

\begin{document}
\newcommand{\horrule}[1]{\rule{\linewidth}{#1}} 	% Horizontal rule 

\title{\usefont{OT1}{qhv}{b}{n}\normalfont \LARGE \textsc{Master of Science Thesis\\
Ecole Centrale de Nantes}
\vspace*{1cm}\\
[0.4cm]\huge\textbf{Nitsche's Method for resolving boundary conditions on embedded interfaces using XFEM in Code Aster}\\
\vspace{4cm}}
\maketitle
\begin{center}
\LARGE \textbf{Nanda Gopala Kilingar}
\par\end{center}

\begin{center}
\vspace{0.5cm}\date{August 2016}\\
\horrule{0.5pt}\\
\par\end{center}

\begin{center}
{\large{}Submitted in fulfillment of the requirements for the degree
of }\\
{\large{}Master of Science in Computational Mechanics}{\large\par}
\par\end{center}
\title*{Nitsche's Method for resolving boundary conditions on embedded interfaces using XFEM in Code Aster}

\begin{eqnarray*}
\mbox{Supervisors:} &  & \mbox{\textbf{Patrick Massin}, Director, IMSIA, EDF R\&D Paris-Saclay}\\
 &  & \mbox{\textbf{Alexandre Martin}, Research Engineer, IMSIA CNRS}
\end{eqnarray*}

\begin{center}
\includegraphics[width=5cm]{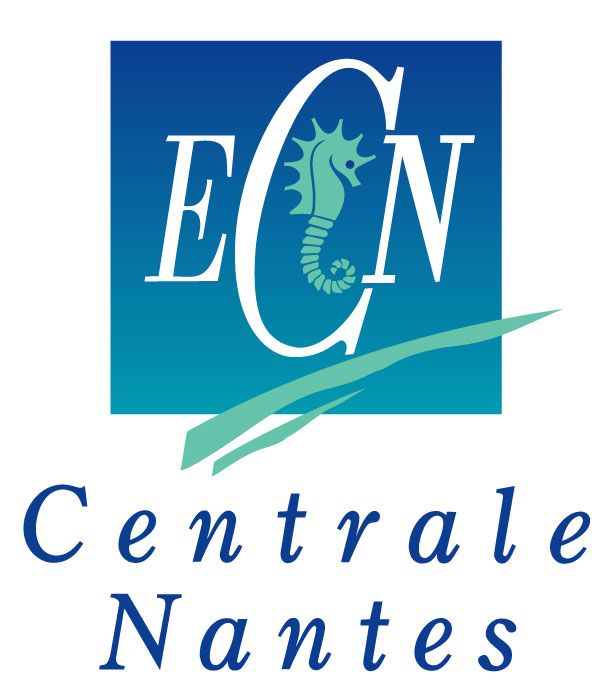}
\par\end{center}

\thispagestyle{empty}

\newpage{}

\vspace*{2cm}
\begin{abstract}
{\large{}As X-FEM approximation doesn't need meshing of the crack,
the method has garnered a lot of attention from industrial point of
view. This thesis report summarises some of the concepts involved
in Nitsche's approach for resolving boundary conditions in embedded
interfaces using XFEM. We consider here cases in which the jump of
a field across the interface is given, as well as cases in which the
primary field on the interface is given. We will first derive the
basics of Nitsche's method and then discretize it with X-FEM using
shifted basis enrichment. We will then implement this on an open source
platform, Code-Aster.}{\large\par}
\end{abstract}
\pagenumbering{roman}

\newpage{}

\subsection*{Acknowledgments}

I would like to thank Parick Massin for providing me the opportunity
to conduct an internship at EDF R\&D and his constant guidance during
my time working in a state-of-art facility with the best infrastructure
possible. I thank Alexandre Martin for his supervision and his inputs
that led me to the successful completion of the internship. I would
like to thank Marcel Ndeffo for his help with various parts of the
subject and also my fellow interns and colleagues of IMSIA at EDF
R\&D.

I would like to thank my professors and lecturers at Ecole Centrale
Nantes for all the knowledge they provided in the preceding semesters
that led me to undergo this internship. Last, but not the least, I
would like to express my gratitude to ECN for giving me an opportunity
to pursue my master's studies and providing a platform for me to expand
my scope of understanding.

\subsection*{\newpage{}}

\tableofcontents{}\newpage{}

\listoffigures

\listoftables
\newpage{}

\pagenumbering{arabic}
\setcounter{page}{1}

\section{Presentation}

\subsection{EDF R\&D}

As the leading global electricity provider, EDF operates in every
energy business line, from generation to customer offer, from transmission
to distribution and from research to innovation. With sales of upto
\euro73 billion, about 55\% of it comes from France while the rest
is generated by international sales and other activities. Generating
about 623.5TWh with around 76.6\% of it coming from nuclear sources,
EDF is almost 87\% $\mbox{CO}_{2}$ free. It invests about \euro650
million research and development alone.\cite{key-12}

\subsubsection*{Strategy}

As today's increasingly digital world dramatically changes the way
we produce and consume, research into electricity generation, transmission
and consumption is of decisive importance. To succeed in the energy
transition, the 2,100 EDF\textquoteright s R\&D division staff (representing
29 nationalities) are currently working on many different projects
designed simultaneously to deliver low-carbon power generation, smarter
energy transmission grids and more responsible energy consumption.
The missions of EDF\textquoteright s R\&D are structured around 3
key priorities.\cite{key-13}
\begin{description}
\item [{Priority}] 1: consolidating and developing competitive, low-carbon
energy generation mixes: One of the major challenges presented by
the energy transition is to ensure the efficient coexistence of traditional
generating methods \textendash{} particularly in terms of improving
nuclear plant safety, efficiency and operating life even further \textendash{}
with the development of renewables.
\item [{Priority}] 2: developing new energy services for customers: Responding
to customer expectations means thinking about new solutions that respond
effectively to variable energy demand while also limiting carbon emissions.
This involves:

\begin{itemize}
\item promoting new ways of using electricity more efficiently (heat pumps,
electric mobility, etc.) 
\item developing digital energy services (real-time consumption control,
smart load balancing, etc.) 
\item developing solutions that encourage energy savings (insulation, appliances,
etc.) 
\item supporting local authorities in their energy plans for sustainable
cities and regions 
\end{itemize}
\item [{Priority}] 3: preparing the electrical systems of tomorrow: This
involves developing smart management tools that will make electrical
systems more flexible and adaptable, encouraging the injection of
intermittent energy sources into the grid, and designing new sustainable
energy solutions at local and regional level.
\end{description}

\subsection{IMSIA}

The IMSIA, Institute of Mechanical Sciences and its Industrial Applications
is a mixed EDF-CNRS research unit created in january 2004. The laboratory
is part of the research facilities of EDF. Its human resources come
from three thematic research departments of EDF R\&D (Mechanical Analyses
and Acoustics (AMA), Material and Mechanics of Components (MMC), Neutronic
Simulation, Information Technology and Scientific Computation (SINETICS)).
Mechanical resistance of structures confronted to ageing problems,
under the constraints of maintained safety and economical performance,
constitutes an important matter for a society facing decisive economic
choices and requiring at the same time an improved safety with respect
to industrial risks. In that perspective, increasing the lifetime
of installations, following and validating maintenance repairs or
structural modifications, monitoring their real behaviour with respect
to design specifications and the need of in service lifetime monitoring,
constitute the key issues that need to be associated to sustainable
development and that require numerous multidisciplinary scientific
progresses. These societal issues are shared with the Engineering
Department of the CNRS and are beyond the sole preoccupations of EDF.
The laboratory is devoted to three main research operations : 
\begin{itemize}
\item Damage and rupture of structures (metallic and civil engineering ones)
; 
\item Data identification, assimilation, exploitation and reduction (loadings,
material properties) and coupled problems involving structures ; 
\item Computational Mechanics : methods, formulations and algorithms for
non linear structural calculations. 
\end{itemize}
The IMSIA relies mainly on Code\_Aster libre, free software under
GNU General Public Licence. It contributes to its evolution in collaboration
with the development team of the software at EDF R\&D and its industrial
an academic partners. The IMSIA is part of the Parisian Federation
for Mechanics Fédération de Recherche Francilienne en Mécanique des
Matériaux, Structures et Procédés (F2M2SP).

\subsection{Code Aster}

Code\_Aster offers a full range of multiphysical analysis and modelling
methods that go well beyond the standard functions of a thermomechanical
calculation code: from seismic analysis to porous media via acoustics,
fatigue, stochastic dynamics, etc. Its modelling, algorithms and solvers
are constantly under construction to improve and complete them (1,200,000
lines of code, 200 operators). Resolutely open, it is linked, coupled
and encapsulated in numerous ways.\cite{key-15}

With the Code\_Aster\textquoteright s architecture, advanced users
can easily work on the code, partly thanks to PYTHON, in order to
write professional applications, introduce finite elements and constitutive
laws or define new exchange formats. The Code\_Aster user describes
the parameters and progression of the survey in a command file. The
grammar and vocabulary of this language, which is specific to Code\_Aster
and written in the PYTHON language, are described in catalogues. This
structuring of the information makes it possible to enhance the language
with new commands at lesser cost or to encapsulate recurring calculation
sequences into macrocommands. A more advanced use enables users to
introduce programming in their datasets: from basic ones (check structures,
loop and tests) to more complex ones using all the richness of PYTHON
(methods, classes, importing graphics or mathematical calculation
modules, etc.) Here is a first basic example: Optimising a pipe bendradius.
Any calculation result can be uploaded in the PYTHON space. Here we
use an indicator for maximal stress in the elbow in order to repeat
the mesh, calculation and postprocessing tasks, thus optimizing the
pipe bend-radius. Another example: with the MEIDEE macro-command,
it is possible to launch calculations for stress identification on
wire structure. Using graphics modules provides an intuitive interface
that helps proceeding to the identification. By encapsulating it into
a macro-command it becomes a professional tool that make the methodology
reliable and durable. 

\subsection{Salome Meca}

The Salome-Meca platform offers a unique environment for the various
phases of a study: 
\begin{itemize}
\item Creating the CAD geometry 
\item Free or structured mesh 
\item Converting to physical data 
\item Launching the Code\_Aster calculation case (ASTK) 
\item Post-processing results 
\end{itemize}
\begin{center}
\begin{figure}
\begin{centering}
\includegraphics[height=10cm]{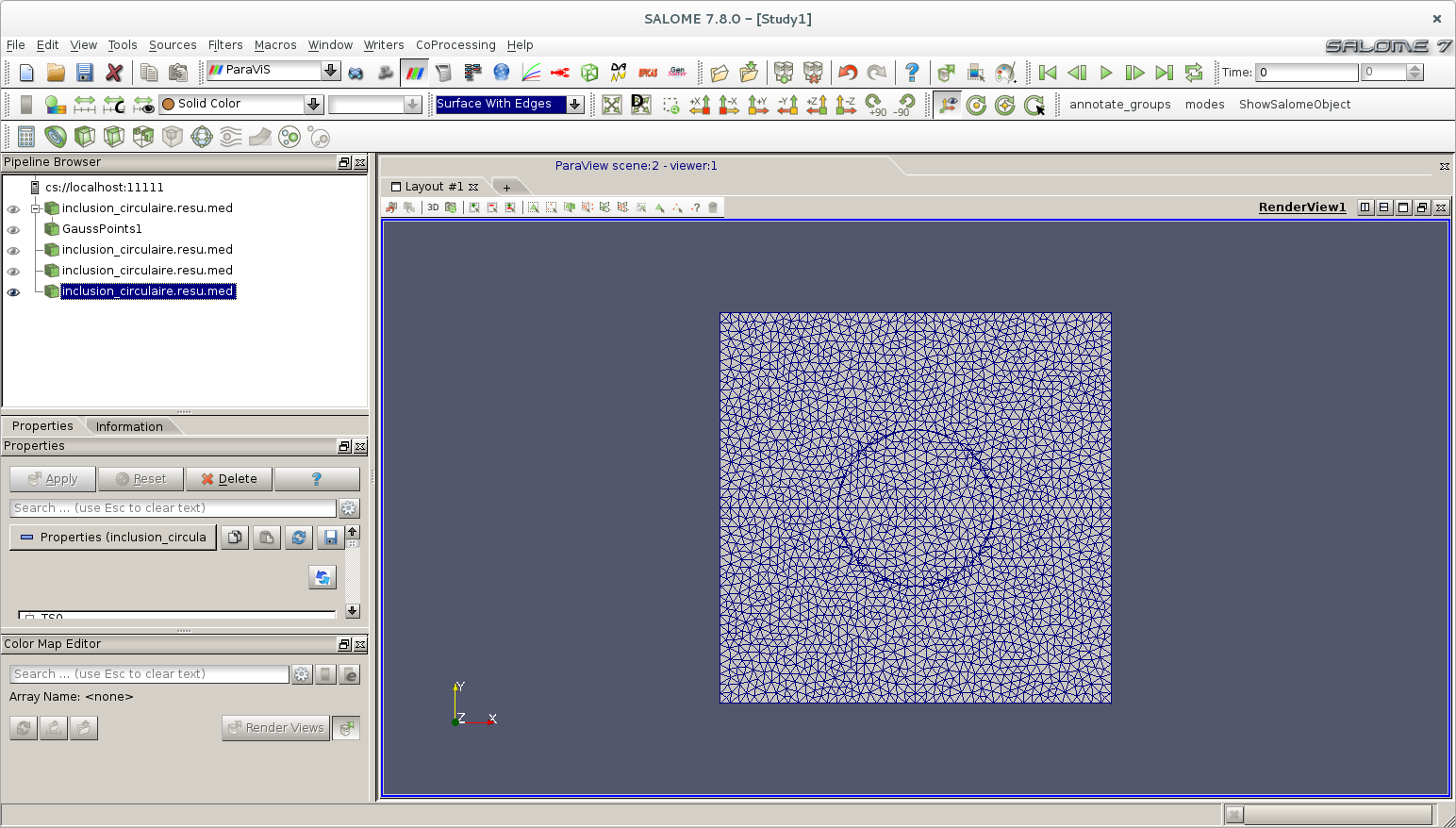}
\par\end{centering}
\caption{GUI of Salome Meca 7.8.0 displaying a 2D surface with mesh during
postprocessing}
\end{figure}
\par\end{center}

\subsection{ASTK}

The provision of a multi-platform, multi-version IT tool that is used
and co-developed by various teams has to be done through a Study and
Developments Manager. This is ASTK\textquoteright s aim: selecting
the code version, defining the files comprised in a study, creating
an overloaded version and accessing configuration management tools
for developers. This interface uses network protocols for transferring
files between clients and server, or for starting remote commands,
including over the Internet. Users can easily distribute their data
files and results to different machines as the interface ensures the
transfer of files, including compressed ones, over the network.

\begin{center}
\begin{figure}
\begin{centering}
\includegraphics[height=8cm]{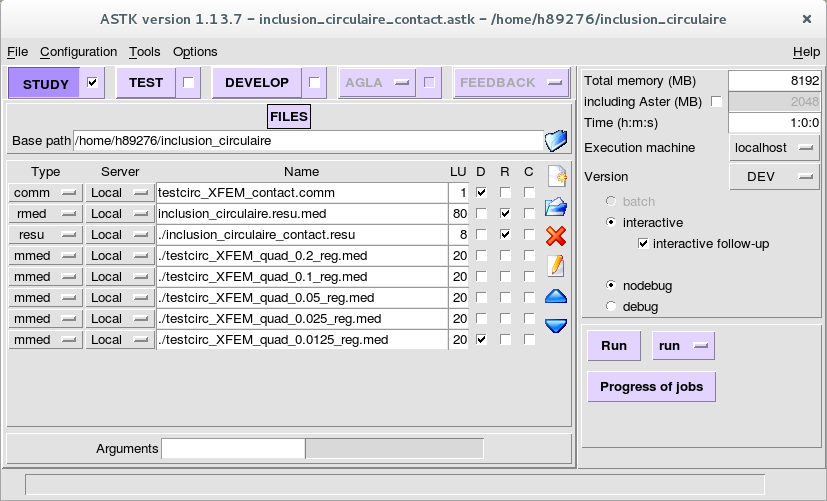}
\par\end{centering}
\caption{ASTK GUI letting the user select the mesh file and the output file
along with the command file}

\end{figure}
\par\end{center}

\newpage{}

\section{Introduction to the Work}

Due to a lot of attention focusing on the development of finite element
methods for embedded interfaces, currently Nitsche's method has been
brought forward to enforce constraints, closed form analytical expressions
for interfacial stabilization terms, and simple flux evaluation by
the help of works done by Dolbow and others. By embedded, we refer
to those methods in which the finite element mesh is not aligned with
the interface geometry (for example X-FEM). Because the interfacial
geometry can be arbitrary with respect to mesh, the robust enforcement
of nonlinear constitutive laws(such as frictional contact) on embedded
interfaces is a challenge that can be tackled with the help of Nitsche's
method.

Focusing on steady problems, but also presenting a case towards time
dependent problems, the method has been implemented in the open source
software Code-Aster. Two main cases of interface problems have been
presented here. 'Jump' problems, which deal with the interfacial problems
concerning those in which the jump in the bulk primary (e.g. displacement,
temperature) and/or secondary (e.g. traction, heat flux) field across
the interface is known or given. The second class of problems are
those in which primary field on the interface is given, referred to
as 'Dirichlet' problems. 

The challenge with 'Jump' type of problems, like a perfectly bonded
material interface in composites, where both jump in displacement
and traction across the interface vanish, the issue often amounts
to the capturing the presence of slope discontinuities that arise
due to the mismatch in material properties. With Nitsche's method,
more general case of non-zero jump can also be considered, which is
more efficient than the enrichment with 'ridge' function and simpler
to implement than the blending of ramp function. 

Gibbs-Thomson conditions arising in crystal growth and solidification
problems, where the interfacial temperature is a function of the interfacial
velocity and curvature, is an example of 'Dirichlet' type of problems.
The problems associated with such type of situations are usually due
to unstable Lagrange multipliers even with the most convenient choice.
Techniques such as penalty methods that may be adequate for enforcing
constraints on stationary interfaces often prove to be lacking when
it comes to yielding accurate, consistent flux quantities. 

We also focus on developing the method taking into consideration some
numerical issues like high sensitivity of normal flux, mild oscillations
and non convergence issues. To balance this, we propose to calculate
a modified numerical flux based on a weighted form. The advantages
of this approach lies firstly in that it is a primal method that does
not introduce additional degrees of freedom at the embedded interface.
Secondly, we obtain stabilization parameters that are based on interfacial
quantities of interest and not necessarily detrimental 'free' parameter.
Finally on extending the method to problems of contact, Nitsche's
method yields more accurate approximations of interfacial traction
fields. 

We implement Nitsche's method focusing on discretization with a shifted
basis enrichment. We also discuss the possibilities of implementing
it in a non-linear Newton loop and obtain the basic matrix form to
do so. We will look at an analytical solution and see how the method
behaves in a simple problem. Finally, we will test the method in Code
Aster on a circular inclusion problem and compare the results with
those obtained by using Lagrange multipliers.

\subsection{Nitsche's Approach on General boundary condition}

Consider the simple 2D Poisson problem: find $u$ such that

\begin{equation}
-\triangle u=f\,\,\,\,\mathrm{{in}}\,\Omega,\label{eq:}
\end{equation}
\begin{equation}
u=u_{0}\,\,\,\,\mathrm{on}\,\Gamma=\partial\Omega_{d},\label{eq:1}
\end{equation}
\begin{equation}
\frac{\partial u}{\partial n}=g\,\,\,\,\mathrm{on}\,\Gamma=\partial\Omega_{n}.\label{eq:-9}
\end{equation}
where $\Omega$ is a bounded domain with polygonal boundary, $f\in L^{2}(\Omega)$,
$u_{0}\in H^{1/2}(\Gamma)$, $g\in L^{2}(\Gamma)$ and $\epsilon\in\mathrm{R}$,
$0\leq\epsilon\leq\infty$. If we consider the penalty method, by
replacing the Dirichlet condition with:

\begin{eqnarray}
\frac{\partial u}{\partial n} & = & \frac{1}{\epsilon}\left[(u_{0}-u)+\epsilon g\right]\,\,\,\,\mathrm{on}\,\Gamma,\label{eq:-1}\\
 & = & \frac{1}{\epsilon}\left[u_{0}-u\right]+g\nonumber 
\end{eqnarray}
where $\epsilon$ is a small parameter, which is problem dependent,
and n is the outward normal. When $\epsilon\rightarrow0$, the solution
to the continuous problem converges to the solution of the Dirichlet
problem (\ref{eq:1}) and $\epsilon\rightarrow\infty$ gives us the
pure Neumann condition (\ref{eq:-9}).
\begin{eqnarray*}
\epsilon\rightarrow0 & \Rightarrow & u=u_{0}\,\,\,\,\mathrm{on}\,\Gamma\\
u-u_{0} & = & 0\\
\frac{\partial u}{\partial n} & = & 0
\end{eqnarray*}
\[
\epsilon\rightarrow\infty\Rightarrow\frac{\partial u}{\partial n}=g\,\,\,\,\mathrm{on}\,\Gamma
\]

The drawbacks of this method are:
\begin{itemize}
\item nonconformity - the method requires coupling of the penalty parameter
to the mesh size
\item possible ill conditioning of the discrete system when $\epsilon$
is too small 
\end{itemize}
Let us consider the variational form of the above simple problem.
Multiply (\ref{eq:}) with $v\in V_{h}$, integrating over the domain
$\Omega$, and using Green's formula, leads to
\begin{equation}
(\nabla u_{h},\nabla v)_{\Omega}-(\frac{\partial u_{h}}{\partial n},v)_{\Gamma}=(f,v)_{\Omega}\label{eq:-2}
\end{equation}

We now multiply the boundary condition (\ref{eq:-1}) by $v$ and
integrating over the domain which gives,
\begin{equation}
(\frac{\partial u_{h}}{\partial n},v)_{\Gamma}+\frac{1}{\epsilon}(u_{h},v)_{\Gamma}=\frac{1}{\epsilon}(u_{0},v)_{\Gamma}+(g,v)_{\Gamma}\label{eq:-3}
\end{equation}

Adding (\ref{eq:-2}) and (\ref{eq:-3}), the problem is equivalent
to finding $u_{h}\in V_{h}$ such that
\begin{equation}
(\nabla u_{h},\nabla v)_{\Omega}\cancel{-(\frac{\partial u_{h}}{\partial n},v)_{\Gamma}}\cancel{+(\frac{\partial u_{h}}{\partial n},v)_{\Gamma}}+\frac{1}{\epsilon}(u_{h},v)_{\Gamma}=(f,v)_{\Omega}+\frac{1}{\epsilon}(u_{0},v)_{\Gamma}+(g,v)_{\Gamma}\,\,\,\,\forall v\in V_{h}
\end{equation}
\begin{equation}
(\nabla u_{h},\nabla v)_{\Omega}+\frac{1}{\epsilon}(u_{h},v)_{\Gamma}=(f,v)_{\Omega}+\frac{1}{\epsilon}(u_{0},v)_{\Gamma}+(g,v)_{\Gamma}\,\,\,\,\forall v\in V_{h}\label{eq:-8}
\end{equation}

Following Juntenen and Stenberg formulations\cite{key-1}, consider,
for simplicity, a regular shaped finite element partitioning ($\mathcal{T}_{h}$)
of the domain ($\Omega\subset\text{R}^{N}$) into triangles or tetrahedra
have been considered. The induced mesh is denoted by $\mathcal{G}_{h}$,
on the boundary $\Gamma$. $K\in\mathcal{T}_{h}$ denotes the element
of the mesh with diameter $h_{K}$ and $E\in\mathcal{G}_{h}$ denotes
the edge or face with diameter $h_{E}$. Further definition consists
of
\begin{equation}
h:=\text{max}\left\{ h_{K}:K\in\mathcal{T}_{h}\right\} 
\end{equation}
and
\begin{equation}
V_{h}:=\left\{ v\in H^{1}\left(\Omega\right):v_{|K}\in\mathcal{P}_{p}\left(K\right)\,\,\,\,\forall K\in\mathcal{T}_{h}\right\} 
\end{equation}
where $\mathcal{P}_{p}(K)$ is the space of polynomials of degree
$p$. Integrating (\ref{eq:-3}) over an element $E$ now gives us:
\begin{equation}
\epsilon(\frac{\partial u_{h}}{\partial n},v)_{E}+(u_{h},v)_{E}=(u_{0},v)_{E}+\epsilon(g,v)_{E}\label{eq:-4}
\end{equation}
 We can write now (\ref{eq:-4}) as 
\begin{equation}
\sum_{E\in\mathcal{G}_{h}}\frac{1}{\epsilon+\gamma h_{E}}\left\{ \epsilon(\frac{\partial u_{h}}{\partial n},v)_{E}+(u_{h},v)_{E}\right\} =\sum_{E\in\mathcal{G}_{h}}\frac{1}{\epsilon+\gamma h_{E}}\left\{ (u_{0},v)_{E}+\epsilon(g,v)_{E}\right\} \label{eq:-5}
\end{equation}
where $\gamma$ is a positive parameter, known as the stability parameter. 

Also from (\ref{eq:-3}), we can write, by multiplying the condition
with $\frac{\partial v}{\partial n}$:
\begin{equation}
\epsilon(\frac{\partial u_{h}}{\partial n},\frac{\partial v}{\partial n})_{E}+(u_{h},\frac{\partial v}{\partial n})_{E}=(u_{0},\frac{\partial v}{\partial n})_{E}+\epsilon(g,\frac{\partial v}{\partial n})_{E}
\end{equation}

Similarly like (\ref{eq:-5}), we can write:
\begin{equation}
\sum_{E\in\mathcal{G}_{h}}-\frac{\gamma h_{E}}{\epsilon+\gamma h_{E}}\left\{ \epsilon(\frac{\partial u_{h}}{\partial n},\frac{\partial v}{\partial n})_{E}+(u_{h},\frac{\partial v}{\partial n})_{E}\right\} =\sum_{E\in\mathcal{G}_{h}}-\frac{\gamma h_{E}}{\epsilon+\gamma h_{E}}\left\{ (u_{0},\frac{\partial v}{\partial n})_{E}+\epsilon(g,\frac{\partial v}{\partial n})_{E}\right\} \label{eq:-6}
\end{equation}

Adding (\ref{eq:-2}), (\ref{eq:-5}) and (\ref{eq:-6}), we get:
\begin{eqnarray*}
(\nabla u_{h},\nabla v)_{\Omega}-(\frac{\partial u_{h}}{\partial n},v)_{\Gamma} & + & \sum_{E\in\mathcal{G}_{h}}\frac{1}{\epsilon+\gamma h_{E}}\left\{ \epsilon(\frac{\partial u_{h}}{\partial n},v)_{E}+(u_{h},v)_{E}\right\} 
\end{eqnarray*}
\[
+\sum_{E\in\mathcal{G}_{h}}-\frac{\gamma h_{E}}{\epsilon+\gamma h_{E}}\left\{ \epsilon(\frac{\partial u_{h}}{\partial n},\frac{\partial v}{\partial n})_{E}+(u_{h},\frac{\partial v}{\partial n})_{E}\right\} 
\]
\begin{equation}
=(f,v)_{\Omega}+\sum_{E\in\mathcal{G}_{h}}\frac{1}{\epsilon+\gamma h_{E}}\left\{ (u_{0},v)_{E}+\epsilon(g,v)_{E}\right\} +\sum_{E\in\mathcal{G}_{h}}-\frac{\gamma h_{E}}{\epsilon+\gamma h_{E}}\left\{ (u_{0},\frac{\partial v}{\partial n})_{E}+\epsilon(g,\frac{\partial v}{\partial n})_{E}\right\} 
\end{equation}

Finally, we find $u_{h}\in V_{h}$ such that: 
\begin{equation}
\mathcal{B}_{h}(u_{h},v)=\mathcal{F}_{h}(v)\,\,\,\,\forall v\in V_{h}
\end{equation}
with:
\begin{eqnarray}
\mathcal{B}_{h}(u,v) & = & (\nabla u,\nabla v)_{\Omega}\\
 &  & +\sum_{E\in\mathcal{G}_{h}}\left\{ -\frac{\gamma h_{E}}{\epsilon+\gamma h_{E}}\left[\left\langle \frac{\partial u}{\partial n},v\right\rangle _{E}+\left\langle u,\frac{\partial v}{\partial n}\right\rangle _{E}\right]+\frac{1}{\epsilon+\gamma h_{E}}\left\langle u,v\right\rangle _{E}-\frac{\epsilon\gamma h_{E}}{\epsilon+\gamma h_{E}}\left\langle \frac{\partial u}{\partial n},\frac{\partial v}{\partial n}\right\rangle _{E}\right\} \nonumber 
\end{eqnarray}
and:
\begin{eqnarray}
\mathcal{F}_{h}(v) & = & (f,v)_{\Omega}\label{eq:-7}\\
 &  & +\sum_{E\in\mathcal{G}_{h}}\left\{ \frac{1}{\epsilon+\gamma h_{E}}\left\langle u_{0},v\right\rangle _{E}-\frac{\gamma h_{E}}{\epsilon+\gamma h_{E}}\left\langle u_{0},\frac{\partial v}{\partial n}\right\rangle _{E}+\frac{\epsilon}{\epsilon+\gamma h_{E}}\left\langle g,v\right\rangle _{E}-\frac{\epsilon\gamma h_{E}}{\epsilon+\gamma h_{E}}\left\langle g,\frac{\partial v}{\partial n}\right\rangle _{E}\right\} \nonumber 
\end{eqnarray}

We can use Nitsche's technique at the limiting condition $\epsilon=0$.
This method can also be extended to the whole range of boundary conditions,
with $\epsilon\geq0$.

If we put $\gamma=0$ in (\ref{eq:-7}), we can see that
\begin{eqnarray*}
(\nabla u_{h},\nabla v)_{\Omega}-(\frac{\partial u_{h}}{\partial n},v)_{\Gamma} & + & \sum_{E\in\mathcal{G}_{h}}\frac{1}{\epsilon+\cancel{\gamma h_{E}}}\left\{ \epsilon(\frac{\partial u_{h}}{\partial n},v)_{E}+(u_{h},v)_{E}\right\} 
\end{eqnarray*}
\[
\cancel{+\sum_{E\in\mathcal{G}_{h}}-\frac{\gamma h_{E}}{\epsilon+\gamma h_{E}}\left\{ \epsilon(\frac{\partial u_{h}}{\partial n},\frac{\partial v}{\partial n})_{E}+(u_{h},\frac{\partial v}{\partial n})_{E}\right\} }
\]
\begin{equation}
=(f,v)_{\Omega}+\sum_{E\in\mathcal{G}_{h}}\frac{1}{\epsilon+\cancel{\gamma h_{E}}}\left\{ (u_{0},v)_{E}+\epsilon(g,v)_{E}\right\} +\cancel{\sum_{E\in\mathcal{G}_{h}}-\frac{\gamma h_{E}}{\epsilon+\gamma h_{E}}\left\{ (u_{0},\frac{\partial v}{\partial n})_{E}+\epsilon(g,\frac{\partial v}{\partial n})_{E}\right\} }
\end{equation}

\begin{equation}
(\nabla u_{h},\nabla v)_{\Omega}\cancel{-(\frac{\partial u_{h}}{\partial n},v)_{\Gamma}}\cancel{+\sum_{E\in\mathcal{G}_{h}}(\frac{\partial u_{h}}{\partial n},v)_{E}}+\sum_{E\in\mathcal{G}_{h}}\frac{1}{\epsilon}(u_{h},v)_{E}=(f,v)_{\Omega}+\sum_{E\in\mathcal{G}_{h}}\frac{1}{\epsilon}(u_{0},v)_{E}+\sum_{E\in\mathcal{G}_{h}}(g,v)_{E}
\end{equation}
\begin{equation}
(\nabla u_{h},\nabla v)_{\Omega}+\sum_{E\in\mathcal{G}_{h}}\frac{1}{\epsilon}(u_{h},v)_{E}=(f,v)_{\Omega}+\sum_{E\in\mathcal{G}_{h}}\frac{1}{\epsilon}(u_{0},v)_{E}+\sum_{E\in\mathcal{G}_{h}}(g,v)_{E}
\end{equation}
which is the same as the traditional form obtained in (\ref{eq:-8})
with . With $\gamma=0,$ the system may become ill-conditioned at
small $\epsilon>0$ values.

For a stabilized method with $\gamma>0,$ at the limit $\epsilon=0,$
we get the Dirichlet problem with Nitsche's application: find $u_{h}\in V_{h}$
such that
\begin{eqnarray}
(\nabla u,\nabla v)_{\Omega}-\left\langle \frac{\partial u_{h}}{\partial n},v\right\rangle _{\Gamma}-\left\langle u_{h},\frac{\partial v}{\partial n}\right\rangle _{\Gamma} & +\sum_{E\in\mathcal{G}_{h}}\frac{1}{\gamma h_{E}}\left\langle u_{h},v\right\rangle _{E}\nonumber \\
=(f,v)_{\Omega}-\left\langle u_{0},\frac{\partial v}{\partial n}\right\rangle _{\Gamma} & +\sum_{E\in\mathcal{G}_{h}}\frac{1}{\gamma h_{E}}\left\langle u_{0},v\right\rangle _{E} & \forall v\in V_{h}\label{eq:-10}
\end{eqnarray}
and at $\epsilon\rightarrow\infty,$ it is a pure Neumann problem
that needs to be solved: find $u_{h}\in V_{h}$ such that

\begin{equation}
(\nabla u,\nabla v)_{\Omega}+\sum_{E\in\mathcal{G}_{h}}\gamma h_{E}\left\langle \frac{\partial u_{h}}{\partial n},\frac{\partial v}{\partial n}\right\rangle _{E}=(f,v)_{\Omega}+(g,v)_{\Gamma}-\sum_{E\in\mathcal{G}_{h}}\gamma h_{E}\left\langle g,\frac{\partial v}{\partial n}\right\rangle _{E}
\end{equation}
This requires that the data satisfy
\begin{equation}
\left(f,1\right)_{\Omega}+\left\langle g,1\right\rangle _{\Gamma}=0
\end{equation}
and this condition is not violated in the above formulations.

\subsection{Comparison with other approaches}

Consider the governing equation
\begin{equation}
-\nabla\centerdot(\kappa\nabla u)=f\,\,\,\,\mathrm{{in}}\,\Omega,\label{eq:-22}
\end{equation}
\begin{equation}
u=u_{0}\,\,\,\,\mathrm{on}\,\Gamma=\partial\Omega\label{eq:1-1}
\end{equation}
We will now discuss some techniques to weakly impose Dirichlet constraints
on embedded surfaces and try to develop one method from another while
briefly discussing the merits of one over the other. This part is
based on the work done by Sanders, Dolbow and Laursen.\cite{key-3}

We consider $\mathbf{n}$ the unit normal to $\Gamma$ which points
out of $\Omega$. The primal variable of $u$ is defined in $\mathbb{U},$
and its variation $\delta u$ is an element of $\mathbb{U}_{0}$:
\[
\mathbb{U}=\left\{ u\in H^{1}(\Omega),u=u_{0}\,\mathrm{on}\,\Gamma\right\} ,
\]
\[
\mathbb{U}_{0}=\left\{ u\in H^{1}(\Omega),u=0\,\mathrm{on}\,\Gamma\right\} .
\]
The potential energy of such a system can be given by 
\begin{equation}
\Pi(u)=\frac{1}{2}\int_{\Omega}\nabla u\centerdot\kappa\nabla u\mathrm{d\Omega}-\int_{\Omega}fu\mathrm{d}\Omega
\end{equation}

The solution $u$ minimizes this potential energy under Dirichlet
constraints. We can transform this constrained problem into an unconstrained
one by using Lagrange multipliers. For the constraint to be enforced,
let's build the Lagrangian of the system, $\mathcal{L}$, by adding
the work of the Lagrange multipliers, $\lambda$ in $\mathbb{L}=H^{-1/2}(\Gamma^{*})$:
\begin{eqnarray}
\mathcal{L}(u,\lambda) & = & \Pi(u)+\int_{\Gamma}\lambda(u-u_{0})\mathrm{d}\mathrm{\Gamma}\nonumber \\
 & = & \frac{1}{2}\int_{\Omega}\nabla u\centerdot\kappa\nabla u\mathrm{d\Omega}-\int_{\Omega}fu\mathrm{d}\Omega+\int_{\Gamma}\lambda(u-u_{0})\mathrm{d}\mathrm{\Gamma}
\end{eqnarray}
 We get a dual variational formulation due to the stationarity of
$\mathcal{L}$: for all $(\delta u,\delta\lambda)\in\mathbb{U}_{0}\times\mathbb{L}$
find $(u,\lambda)\in\mathbb{U\times L}$, such that:
\begin{equation}
\delta\mathcal{L}=\int_{\Omega}\nabla\delta u\centerdot\kappa\nabla u\mathrm{d\Omega}-\int_{\Omega}\delta uf\mathrm{d\Omega}+\int_{\Gamma}\lambda\delta u\mathrm{d}\mathrm{\Gamma}+\int_{\Gamma}\delta\lambda(u-u_{0})\mathrm{d}\mathrm{\Gamma}=0
\end{equation}

\begin{flushright}
\textit{(Lagrangian method)}
\par\end{flushright}

Here $\lambda$ and $u$ cannot be independently determined. And this
brings about a lot of stability issues. We can solve this by the use
of penalty methods. We can interpret Lagrange multipliers as flux
imposed on the boundary condition. This leads us to establish $\lambda=-\kappa\nabla u\centerdot\mathrm{\mathbf{n}}$,
the flux. Now we assume that the flux can be approximated in a spring
like form $\kappa\nabla u\centerdot\mathrm{\mathbf{n}}\approx-\epsilon(u-u_{0}).$
The penalty potential is now, 
\begin{eqnarray}
\Pi^{\mathrm{pen}}(u) & = & \Pi(u)+\int_{\Gamma}\frac{\epsilon}{2}(u-u_{0})(u-u_{0})\mathrm{d}\mathrm{\Gamma}\nonumber \\
 & = & \Pi(u)+\int_{\Gamma}\frac{\epsilon}{2}(u-u_{0})^{2}\mathrm{d}\mathrm{\Gamma}\nonumber \\
 & = & \frac{1}{2}\int_{\Omega}\nabla u\centerdot\kappa\nabla u\mathrm{d\Omega}-\int_{\Omega}fu\mathrm{d}\Omega+\int_{\Gamma}\frac{\epsilon}{2}(u-u_{0})^{2}\mathrm{d}\mathrm{\Gamma}
\end{eqnarray}

The primal penalty variational form is given by: for all $\delta u\in\mathbb{U}_{0}$
find $u\in\mathbb{U}$, such that
\begin{equation}
\delta\Pi^{\mathrm{pen}}=\int_{\Omega}\nabla\delta u\centerdot\kappa\nabla u\mathrm{d\Omega}-\int_{\Omega}\delta uf\mathrm{d\Omega}+\int_{\Gamma}\delta u\epsilon(u-u_{0})\mathrm{d}\mathrm{\Gamma}=0
\end{equation}

\begin{flushright}
\textit{(Penalty method)}
\par\end{flushright}

This is not variationally consistent, as the desired problem is solved
only in the limiting case when $\epsilon\rightarrow\infty$.

Another standard way to improve the behavior of a Lagrangian method
is to stabilize it with a penalty term. This gives the augmented Lagrangian
approach.

\begin{eqnarray}
\mathcal{L}^{\mathrm{aug}}(u,\lambda) & = & \Pi(u)+\int_{\Gamma}\lambda(u-u_{0})\mathrm{d}\mathrm{\Gamma}+\int_{\Gamma}\frac{\epsilon}{2}(u-u_{0})^{2}\mathrm{d}\mathrm{\Gamma}\nonumber \\
 & = & \frac{1}{2}\int_{\Omega}\nabla u\centerdot\kappa\nabla u\mathrm{d\Omega}-\int_{\Omega}fu\mathrm{d}\Omega+\int_{\Gamma}\lambda(u-u_{0})\mathrm{d}\mathrm{\Gamma}+\int_{\Gamma}\frac{\epsilon}{2}(u-u_{0})^{2}\mathrm{d}\mathrm{\Gamma}
\end{eqnarray}

Here the penalty stiffness can be seen as the stabilization parameter,
which does not need large values of it. The dual variational form:
for all $(\delta u,\delta\lambda)\in\mathbb{U}_{0}\times\mathbb{L}$
find $(u,\lambda)\in\mathbb{U\times L}$, such that:
\begin{eqnarray}
\delta\mathcal{L}^{\mathrm{aug}} & = & \int_{\Omega}\nabla\delta u\centerdot\kappa\nabla u\mathrm{d\Omega}-\int_{\Omega}\delta uf\mathrm{d\Omega}+\int_{\Gamma}\delta u\lambda\mathrm{d}\mathrm{\Gamma}+\int_{\Gamma}\delta\lambda(u-u_{0})\mathrm{d}\mathrm{\Gamma}+\int_{\Gamma}\delta u\epsilon(u-u_{0})\mathrm{d}\mathrm{\Gamma}\nonumber \\
 & = & \int_{\Omega}\nabla\delta u\centerdot\kappa\nabla u\mathrm{d\Omega}-\int_{\Omega}\delta uf\mathrm{d\Omega}+\int_{\Gamma}\delta u(\lambda+\epsilon(u-u_{0}))\mathrm{d}\mathrm{\Gamma}+\int_{\Gamma}\delta\lambda(u-u_{0})\mathrm{d}\mathrm{\Gamma}=0
\end{eqnarray}

\begin{flushright}
\textit{(Augmented Lagrangian method)}
\par\end{flushright}

In the same manner to obtain a penalty variational form from a Lagrangian
one, we can utilize the flux relation $\lambda=-\kappa\nabla u\centerdot\mathrm{\mathbf{n}}$
to obtain the potential function that forms the basis of Nitsche's
method.
\begin{eqnarray}
\Pi^{\mathrm{Nit}}(u) & = & \Pi(u)-\int_{\Gamma}(u-u_{0})\kappa\nabla u\centerdot\mathrm{\mathbf{n}}\mathrm{d}\mathrm{\Gamma}+\int_{\Gamma}\frac{\epsilon}{2}(u-u_{0})^{2}\mathrm{d}\mathrm{\Gamma}\nonumber \\
 & = & \frac{1}{2}\int_{\Omega}\nabla u\centerdot\kappa\nabla u\mathrm{d\Omega}-\int_{\Omega}fu\mathrm{d}\Omega-\int_{\Gamma}(u-u_{0})\kappa\nabla u\centerdot\mathrm{\mathbf{n}}\mathrm{d}\mathrm{\Gamma}+\int_{\Gamma}\frac{\epsilon}{2}(u-u_{0})^{2}\mathrm{d}\mathrm{\Gamma}
\end{eqnarray}

We get one-field symmetric variational formulation: for all $\delta u\in\mathbb{U}_{0}$
find $u\in\mathbb{U}$, such that :

\begin{eqnarray}
\delta\Pi^{\mathrm{Nit}}(u) & = & \int_{\Omega}\nabla\delta u\centerdot\kappa\nabla u\mathrm{d\Omega}-\int_{\Omega}\delta uf\mathrm{d\Omega}-\int_{\Gamma}\delta u\kappa\nabla u\centerdot\mathrm{\mathbf{n}}\mathrm{d}\mathrm{\Gamma}\nonumber \\
 &  & -\int_{\Gamma}(u-u_{0})\kappa\nabla\delta u\centerdot\mathrm{\mathbf{n}}\mathrm{d}\mathrm{\Gamma}+\int_{\Gamma}\delta u\epsilon(u-u_{0})\mathrm{d}\mathrm{\Gamma}=0
\end{eqnarray}

\begin{flushright}
\textit{(Nitsche's method)}
\par\end{flushright}

Rearranging, we can write:
\[
\int_{\Omega}\nabla\delta u\centerdot\kappa\nabla u\mathrm{d\Omega}-\int_{\Gamma}\delta u\kappa\nabla u\centerdot\mathrm{\mathbf{n}}\mathrm{d}\mathrm{\Gamma}-\int_{\Gamma}(u-u_{0})\kappa\nabla\delta u\centerdot\mathrm{\mathbf{n}}\mathrm{d}\mathrm{\Gamma}+\int_{\Gamma}\delta u\epsilon(u-u_{0})\mathrm{d}\mathrm{\Gamma}=\int_{\Omega}\delta uf\mathrm{d\Omega}
\]
or:

\[
\int_{\Omega}\nabla\delta u\centerdot\kappa\nabla u\mathrm{d\Omega}-\int_{\Gamma}\delta u\kappa\nabla u\centerdot\mathrm{\mathbf{n}}\mathrm{d}\mathrm{\Gamma}-\int_{\Gamma}u\kappa\nabla\delta u\centerdot\mathrm{\mathbf{n}}\mathrm{d}\mathrm{\Gamma}+\int_{\Gamma}\delta u\epsilon u\mathrm{d}\mathrm{\Gamma}
\]
\begin{eqnarray}
= & \int_{\Omega}\delta uf\mathrm{d\Omega}- & \int_{\Gamma}u_{0}\kappa\nabla\delta u\centerdot\mathrm{\mathbf{n}}\mathrm{d}\mathrm{\Gamma}+\int_{\Gamma}\delta u\epsilon u_{0}\mathrm{d}\mathrm{\Gamma}\label{eq:-11}
\end{eqnarray}

If we compare this with (\ref{eq:-10}), we can see that $\epsilon$
is similar to the term $\sum_{E\in\mathcal{G}_{h}}\frac{1}{\gamma h_{E}}$.

\subsection{Application to Interfaces\label{sec:Application-on-Interface}}

We can now extend Nitsche's method to an interface made of two materials
or even to a crack that divides the surface. We develop this based
on the works of Dolbow and Harari.\cite{key-2} 

Find $u\in\mathbb{U}$, such that: 
\begin{equation}
a_{b}(v,u)+a_{i}(v,u)=l_{b}(v)+l_{i}(v)\,\,\,\,\forall v\in\mathbb{U}_{0}
\end{equation}
where now:
\[
a_{b}(v,u)=\int_{\Omega}\nabla v\centerdot\kappa\nabla u\mathrm{d\Omega}
\]
and: 
\[
l_{b}(v)=\int_{\Omega}vf\mathrm{d\Omega}
\]
are the standard bulk contributions from (\ref{eq:-11}), with: 
\begin{eqnarray*}
\mathbb{U}_{\,\,} & = & \left\{ u\in H^{1}(\Omega^{-}\cup\Omega^{+}),u=u_{0}\,\mathrm{on}\,\Gamma,\mbox{ may be discontinuous on }\mathcal{S}\right\} \mbox{ and}\\
\mathbb{U}_{0} & = & \left\{ v\in H^{1}(\Omega^{-}\cup\Omega^{+}),v=0\,\mathrm{on}\,\Gamma,\mbox{ may be discontinuous on }\mathcal{S}\right\} .
\end{eqnarray*}

$a_{i}(v,u)$ and $l_{i}(v)$ are the interfacial contributions which
depend on the case being considered. These terms are obtained by utilizing
the boundary conditions in a similar manner to that in (\ref{eq:-11}),
but on the interface, if we can consider the domain $\Omega$ to be
divided by an interface $\mathcal{S}$.

\subparagraph*{Dirichlet condition (ex. crack surface)}

Consider the boundary conditions:

\[
u^{+}=g^{+},\,\,\,\,u^{-}=g^{-}\mbox{\,\,\,\,\ on }\mathcal{\mathcal{S}}
\]
where $g^{+}$ and $g^{-}$ are assumed to be sufficiently smooth
functions of the position on the interface. $u^{+}$ and $u^{-}$
are limiting values of the field $u$ as the interface is approached
from either $\Omega^{+}$ or $\Omega^{-}$, respectively. Approaching
this problem as two one-sided problems, we can write
\begin{eqnarray}
a_{i}(v,u) & = & a_{i}(v,u)^{+}+a_{i}(v,u)^{-}\label{eq:-15}\\
l_{i}(v) & = & l_{i}(v)^{+}+l_{i}(v)^{-}\label{eq:-16}
\end{eqnarray}

\begin{figure}
\begin{centering}
\includegraphics[bb=300bp 600bp 460bp 735bp,clip]{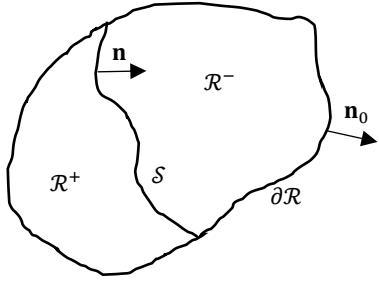}
\par\end{centering}
\caption{\label{fig:Notation-for-two-sided}Notation for two-sided problem}
\end{figure}

We choose the interfacial normal $\mathrm{\mathbf{n}}$ as pointing
outwardly from $\Omega^{+}$. This gives us, from (\ref{eq:-11}),
\begin{eqnarray}
a_{i}(v,u)^{+} & = & -\int_{\mathcal{S}}v^{+}(\kappa^{+}\nabla u^{+}\centerdot\mathrm{\mathbf{n}})\mathrm{d}\mathrm{\Gamma}-\int_{\mathcal{S}}u^{+}(\kappa^{+}\nabla v^{+}\centerdot\mathrm{\mathbf{n}})\mathrm{d}\mathrm{\Gamma}+\int_{\mathcal{S}}v^{+}\alpha^{+}u^{+}\mathrm{d}\mathrm{\Gamma}\label{eq:-17}\\
a_{i}(v,u)^{-} & = & -\int_{\mathcal{S}}v^{-}(\kappa^{-}\nabla u^{-}\centerdot\left(-\mathrm{\mathbf{n}}\right))\mathrm{d}\mathrm{\Gamma}-\int_{\mathcal{S}}u^{-}(\kappa^{-}\nabla v^{-}\centerdot\left(-\mathrm{\mathbf{n}}\right))\mathrm{d}\mathrm{\Gamma}+\int_{\mathcal{S}}v^{-}\alpha^{-}u^{-}\mathrm{d}\mathrm{\Gamma}\nonumber \\
 & = & \int_{\mathcal{S}}v^{-}(\kappa^{-}\nabla u^{-}\centerdot\mathrm{\mathbf{n}})\mathrm{d}\mathrm{\Gamma}+\int_{\mathcal{S}}u^{-}(\kappa^{-}\nabla v^{-}\centerdot\mathrm{\mathbf{n}})\mathrm{d}\mathrm{\Gamma}+\int_{\mathcal{S}}v^{-}\alpha^{-}u^{-}\mathrm{d}\mathrm{\Gamma}\label{eq:-42}
\end{eqnarray}
\begin{eqnarray}
l_{i}(v)^{+} & = & -\int_{\mathcal{S}}g^{+}(\kappa^{+}\nabla v^{+}\centerdot\mathrm{\mathbf{n}})\mathrm{d}\mathrm{\Gamma}+\int_{\mathcal{S}}v^{+}\alpha^{+}g^{+}\mathrm{d}\mathrm{\Gamma}\label{eq:-43}\\
l_{i}(v)^{-} & = & -\int_{\mathcal{S}}g^{-}(\kappa^{-}\nabla v^{-}\centerdot\left(-\mathrm{\mathbf{n}}\right))\mathrm{d}\mathrm{\Gamma}+\int_{\mathcal{S}}v^{-}\alpha^{-}g^{-}\mathrm{d}\mathrm{\Gamma}\nonumber \\
 & = & \int_{\mathcal{S}}g^{-}(\kappa^{-}\nabla v^{-}\centerdot\mathrm{\mathbf{n}})\mathrm{d}\mathrm{\Gamma}+\int_{\mathcal{S}}v^{-}\alpha^{-}g^{-}\mathrm{d}\mathrm{\Gamma}\label{eq:-44}
\end{eqnarray}
where we have chosen $\alpha^{+}$ and $\alpha^{-}$ as the stabilization
parameters for the $\Omega^{+}$ and $\Omega^{-}$ domains respectively
(Figure \ref{fig:Notation-for-two-sided}).

With a Dirichlet condition of this type, a jump $\bar{\mathrm{j}}$
in the flux, if it exists, is unknown and represents a quantity of
interest. This is because we have assumed here that the two domains
are completely different problems. If we consider $\mathcal{L}$ as
a portion of the interface, with $\mathcal{B}=\mathcal{B}^{+}\cup\mathcal{B}^{-}$
being the supports of the weight function $v_{d}$ that covers smoothly,
a point $\mathbf{x}_{\mbox{d}}$ on the interface, and $\mathcal{L}=\mathcal{B}\cap\mathcal{S}$,
we obtain: 
\begin{eqnarray}
\int_{\mathcal{L}}v_{d}\bar{\mathrm{j}}\mbox{d}\Gamma & = & \int_{\mathcal{B}}v_{d}f\,\mbox{d}\Omega-\int_{\mathcal{B}}\nabla v_{d}\kappa\nabla u\,\mbox{d}\Omega\\
 & = & \int_{\mathcal{B}}v_{d}^{+}f\,\mbox{d}\Omega-\int_{\mathcal{B}}\nabla v_{d}^{+}\kappa^{+}\nabla u^{+}\,\mbox{d}\Omega+\int_{\mathcal{B}}v_{d}^{-}f\,\mbox{d}\Omega-\int_{\mathcal{B}}\nabla v_{d}^{-}\kappa^{-}\nabla u^{-}\,\mbox{d}\Omega
\end{eqnarray}
with $\left[\left[\kappa\nabla u\right]\right]\centerdot\mathbf{n}=\left(\kappa\nabla u^{+}-\kappa\nabla u^{-}\right)\centerdot\mathbf{n}=\bar{\mathrm{j}}$
which can be used to approximate the jump in flux across the portion
of interface considered. 

\subparagraph*{Jump Condition (ex. a bimaterial interface)\label{par:Jump-Condition}}

Consider the given problem
\[
-\nabla\centerdot(\kappa^{+}\nabla u)=f\mbox{ in }\Omega^{+}
\]
\[
-\nabla\centerdot(\kappa^{-}\nabla u)=f\mbox{ in }\Omega^{-}
\]
with
\[
u=u_{0}\,\,\,\,\mathrm{on}\,\Gamma=\partial\Omega
\]
\begin{eqnarray*}
\left[\left[u\right]\right] & = & \bar{\mathrm{i}}\,\,\,\,\mathrm{on}\,\mathcal{S}\\
u^{+}-u^{-} & = & \bar{\mathrm{i}}
\end{eqnarray*}
\begin{eqnarray*}
\left[\left[\kappa\nabla u\right]\right]\centerdot\mathbf{n} & = & \bar{\mathrm{j}}\,\,\,\,\mathrm{on}\,\mathcal{S}\\
\left(\kappa^{+}\nabla u^{+}-\kappa^{-}\nabla u^{-}\right).\mathbf{n} & = & \bar{\mathrm{j}}
\end{eqnarray*}
where $\mathbf{n}$ is considered to point outwards from $\Omega^{+}$.
\[
\left(\kappa\nabla u\right).\mathbf{n}_{\text{0}}=0\mbox{ on }\Gamma=\partial\Omega_{n}
\]
Consider the domain $\Omega^{+}$. We can try to obtain the weak Galerkin
formulation.
\begin{equation}
-\int_{\Omega^{+}}v^{+}\nabla\centerdot\left(\kappa^{+}\nabla u^{+}\right)=\int_{\Omega^{+}}v^{+}f\,\mbox{d}\Omega
\end{equation}
We can write:
\[
\int_{\Omega^{+}}\nabla v^{+}.\left(\kappa^{+}\nabla u^{+}\right)-\int_{\Omega^{+}}\nabla.\left(v^{+}\kappa^{+}\nabla u^{+}\right)=\int_{\Omega^{+}}v^{+}f\,\mbox{d}\Omega
\]
from divergence theorem. And also:
\[
\int_{\Omega^{+}}\nabla.\left(v^{+}\kappa^{+}\nabla u^{+}\right)=\int_{S}v^{+}\left(\kappa^{+}\nabla u^{+}.\mathbf{n}\right)+\cancel{\int_{\partial\Omega_{n}}v^{+}\left(\kappa^{+}\nabla u^{+}.\mathbf{n}_{0}\right)}
\]
Thus weak Galerkin formulation of the above problem gives us: find
$u\in\mathbb{U}$ such that 
\begin{equation}
\int_{\Omega^{+}}\nabla v^{+}.\kappa^{+}\nabla u^{+}\mbox{d}\Omega-\int_{\mathcal{S}}v^{+}(\kappa^{+}\nabla u^{+}.\mathbf{n})\mbox{d}\Gamma=\int_{\Omega^{+}}v^{+}f\,\mbox{d}\Omega\mbox{ }\forall v\in\mathbb{U}_{0}
\end{equation}
and similarly on domain $\Omega^{-}$,
\begin{equation}
\int_{\Omega^{-}}\nabla v^{-}.\kappa^{-}\nabla u^{-}\mbox{d}\Omega+\int_{\mathcal{S}}v^{-}(\kappa^{-}\nabla u^{-}.\mathbf{n})\mbox{d}\Gamma=\int_{\Omega^{-}}v^{-}f\,\mbox{d}\Omega\mbox{ }\forall v\in\mathbb{U}_{0}
\end{equation}
We have considered the problem separately in the two domains which
has resulted in the separation of the jump in flux. Adding the two
equations gives:
\begin{equation}
\int_{\Omega}\nabla v.\kappa\nabla u\mbox{d}\Omega-\int_{\mathcal{S}}\left[v^{+}(\kappa^{+}\nabla u^{+}.\mathbf{n})-v^{-}(\kappa^{-}\nabla u^{-}.\mathbf{n})\right]\mbox{d}\Gamma=\int_{\Omega}vf\mbox{d}\Omega
\end{equation}
Also 
\begin{equation}
\kappa^{+}\nabla u^{+}=<\kappa\nabla u>+\frac{1}{2}[[\kappa\nabla u]]
\end{equation}
and
\begin{equation}
\kappa^{-}\nabla u^{-}=<\kappa\nabla u>-\frac{1}{2}[[\kappa\nabla u]]
\end{equation}
Substituting these in the weak formulation gives:

\[
\int_{\Omega}\nabla v.\kappa\nabla u\mbox{d}\Omega-\int_{\mathcal{S}}\left[v^{+}\left(\left(<\kappa\nabla u>+\frac{1}{2}[[\kappa\nabla u]]\right).\mathbf{n}\right)-v^{-}\left(\left(<\kappa\nabla u>-\frac{1}{2}[[\kappa\nabla u]]\right).\mathbf{n}\right)\right]\mbox{d}\Gamma=\int_{\Omega}vf\,\mbox{d}\Omega
\]
\[
\int_{\Omega}\nabla v.\kappa\nabla u\mbox{d}\Omega-\int_{\mathcal{S}}\left[\frac{1}{2}v^{+}[[\kappa\nabla u]].\mathbf{n}+v^{+}<\kappa\nabla u>.\mathbf{n}+\frac{1}{2}v^{-}[[\kappa\nabla u]].\mathbf{n}-v^{-}<\kappa\nabla u>.\mathbf{n}\right]\mbox{d}\Gamma=\int_{\Omega}vf\,\mbox{d}\Omega
\]
\[
\int_{\Omega}\nabla v.\kappa\nabla u\mbox{d}\Omega-\int_{\mathcal{S}}\left[\frac{1}{2}(v^{+}+v^{-})[[\kappa\nabla u]].\mathbf{n}+(v^{+}-v^{-})<\kappa\nabla u>.\mathbf{n}\right]\mbox{d}\Gamma=\int_{\Omega}vf\,\mbox{d}\Omega
\]
\[
\int_{\Omega}\nabla v.\kappa\nabla u\mbox{d}\Omega-\int_{\mathcal{S}}\left[\left\langle v\right\rangle \bar{\mathrm{j}}+\left[\left[v\right]\right]<\kappa\nabla u>.\mathbf{n}\right]\mbox{d}\Gamma=\int_{\Omega}vf\,\mbox{d}\Omega
\]
\begin{equation}
\int_{\Omega}\nabla v.\kappa\nabla u\mbox{d}\Omega-\int_{\mathcal{S}}\left[\left[v\right]\right]<\kappa\nabla u>.\mathbf{n}+\int_{\mathcal{S}}\frac{v^{-}}{2}\left(\kappa^{-}\nabla u^{-}\centerdot\mathrm{\mathbf{n}}\right)\mathrm{d}\mathrm{\Gamma}\mbox{d}\Gamma=\int_{\Omega}f\,\mbox{d}\Omega+\int_{\mathcal{S}}\left\langle v\right\rangle \bar{\mathrm{j}}\mbox{d}\Gamma
\end{equation}

Now we can introduce Nitsche's terms and the corresponding stabilization
terms in this equation to get the variational form and maintain the
symmetric nature of the system with Nitsche's approach:
\begin{eqnarray}
\int_{\Omega}\nabla v.\kappa\nabla u\mbox{d}\Omega-\int_{\mathcal{S}}\left[\left[v\right]\right]<\kappa\nabla u>.\mathbf{n}\mbox{d}\Gamma-\int_{\mathcal{S}}[[u]]<\kappa\nabla v>\centerdot\mathrm{\mathbf{n}}\mathrm{d}\mathrm{\Gamma}+\int_{\mathcal{S}}[[v]]\alpha[[u]]\mathrm{d}\mathrm{\Gamma}\nonumber \\
=\int_{\Omega}vf\,\mbox{d}\Omega+\int_{\mathcal{S}}\left\langle v\right\rangle \bar{\mathrm{j}}\mbox{d}\Gamma-\int_{\mathcal{S}}\bar{\mbox{i}}<\kappa\nabla v>\centerdot\mathrm{\mathbf{n}}\mathrm{d}\mathrm{\Gamma}+\int_{\mathcal{S}}[[v]]\alpha\,\bar{\mbox{i}}\,\mathrm{d}\mathrm{\Gamma}
\end{eqnarray}
This form is both variationally consistent as well as symmetric. Thus,
we can write 

\begin{eqnarray}
a_{i}(v,u) & = & -\int_{\mathcal{S}}[[v]]<\kappa\nabla u>\centerdot\mathrm{\mathbf{n}}\mathrm{d}\mathrm{\Gamma}-\int_{\mathcal{S}}[[u]]<\kappa\nabla v>\centerdot\mathrm{\mathbf{n}}\mathrm{d}\mathrm{\Gamma}+\int_{\mathcal{S}}[[v]]\alpha[[u]]\mathrm{d}\mathrm{\Gamma}\label{eq:-14}\\
l_{i}(v) & = & -\int_{\mathcal{S}}\bar{\mbox{i}}<\kappa\nabla v>\centerdot\mathrm{\mathbf{n}}\mathrm{d}\mathrm{\Gamma}+\int_{\mathcal{S}}[[v]]\alpha\,\bar{\mbox{i}}\,\mathrm{d}\mathrm{\Gamma}+\int_{\mathcal{S}}<v>\bar{\mbox{j}}\mathrm{d}\mathrm{\Gamma}\label{eq:-41}
\end{eqnarray}
$\alpha$ is the integrated stabilizing term for this form here.

\subparagraph*{Comparison}

If we expand the jump condition relations, we have, from (\ref{eq:-14})
and (\ref{eq:-41}), 
\begin{eqnarray*}
a_{i}(v,u) & = & -\int_{\mathcal{S}}\frac{\left(v^{+}-v^{-}\right)}{2}\left(\kappa^{+}\nabla u^{+}\centerdot\mathrm{\mathbf{n}}+\kappa^{-}\nabla u^{-}\centerdot\mathrm{\mathbf{n}}\right)\mathrm{d}\mathrm{\Gamma}\\
 &  & -\int_{\mathcal{S}}\frac{\left(u^{+}-u^{-}\right)}{2}\left(\kappa^{+}\nabla v^{+}\centerdot\mathrm{\mathbf{n}}+\kappa^{-}\nabla v^{-}\centerdot\mathrm{\mathbf{n}}\right)\mathrm{d}\mathrm{\Gamma}+\int_{\mathcal{S}}\left(v^{+}-v^{-}\right)\alpha(u^{+}-u^{-})\mathrm{d}\mathrm{\Gamma}\\
 & = & -\int_{\mathcal{S}}\frac{v^{+}}{2}\left(\kappa^{+}\nabla u^{+}\centerdot\mathrm{\mathbf{n}}\right)\mathrm{d}\mathrm{\Gamma}-\int_{\mathcal{S}}\frac{v^{+}}{2}\left(\kappa^{-}\nabla u^{-}\centerdot\mathrm{\mathbf{n}}\right)\mathrm{d}\mathrm{\Gamma}+\int_{\mathcal{S}}\frac{v^{-}}{2}\left(\kappa^{+}\nabla u^{+}\centerdot\mathrm{\mathbf{n}}\right)\mathrm{d}\mathrm{\Gamma}\\
 &  & -\int_{\mathcal{S}}\frac{u^{+}}{2}\left(\kappa^{+}\nabla v^{+}\centerdot\mathrm{\mathbf{n}}\right)\mathrm{d}\mathrm{\Gamma}-\int_{\mathcal{S}}\frac{u^{+}}{2}\left(\kappa^{-}\nabla v^{-}\centerdot\mathrm{\mathbf{n}}\right)\mathrm{d}\mathrm{\Gamma}+\int_{\mathcal{S}}\frac{u^{-}}{2}\left(\kappa^{+}\nabla v^{+}\centerdot\mathrm{\mathbf{n}}\right)\mathrm{d}\mathrm{\Gamma}\\
 &  & +\int_{\mathcal{S}}\frac{v^{-}}{2}\left(\kappa^{-}\nabla u^{-}\centerdot\mathrm{\mathbf{n}}\right)\mathrm{d}\mathrm{\Gamma}+\int_{\mathcal{S}}\frac{u^{-}}{2}\left(\kappa^{-}\nabla v^{-}\centerdot\mathrm{\mathbf{n}}\right)\mathrm{d}\mathrm{\Gamma}\\
 &  & +\int_{\mathcal{S}}v^{+}\alpha u^{+}\mathrm{d}\mathrm{\Gamma}-\int_{\mathcal{S}}v^{+}\alpha u^{-}\mathrm{d}\mathrm{\Gamma}-\int_{\mathcal{S}}v^{-}\alpha u^{+}\mathrm{d}\mathrm{\Gamma}+\int_{\mathcal{S}}v^{-}\alpha u^{-}\mathrm{d}\mathrm{\Gamma}
\end{eqnarray*}
\begin{eqnarray*}
l_{i}(v) & = & -\int_{\mathcal{S}}\frac{\bar{\mbox{i}}}{2}\left(\kappa^{+}\nabla v^{+}\centerdot\mathrm{\mathbf{n}}+\kappa^{-}\nabla v^{-}\centerdot\mathrm{\mathbf{n}}\right)\mathrm{d}\mathrm{\Gamma}+\int_{\mathcal{S}}\left(v^{+}-v^{-}\right)\alpha\,\bar{\mbox{i}}\,\mathrm{d}\mathrm{\Gamma}+\int_{\mathcal{S}}\frac{1}{2}\left(v^{+}+v^{-}\right)\,\bar{\mbox{j}}\,\mathrm{d}\mathrm{\Gamma}\\
 & = & -\int_{\mathcal{S}}\frac{\bar{\mbox{i}}}{2}\left(\kappa^{+}\nabla v^{+}\centerdot\mathrm{\mathbf{n}}\right)\mathrm{d}\mathrm{\Gamma}-\int_{\mathcal{S}}\frac{\bar{\mbox{i}}}{2}\left(\kappa^{-}\nabla v^{-}\centerdot\mathrm{\mathbf{n}}\right)\mathrm{d}\mathrm{\Gamma}\\
 &  & +\int_{\mathcal{S}}v^{+}\alpha\,\bar{\mbox{i}}\,\mathrm{d}\mathrm{\Gamma}-\int_{\mathcal{S}}v^{-}\alpha\,\bar{\mbox{i}}\,\mathrm{d}\mathrm{\Gamma}+\int_{\mathcal{S}}\frac{v^{+}}{2}\bar{\mbox{j}}\,\mathrm{d}\mathrm{\Gamma}+\int_{\mathcal{S}}\frac{v^{-}}{2}\bar{\mbox{j}}\,\mathrm{d}\mathrm{\Gamma}
\end{eqnarray*}
and from the Dirichlet conditions, by combining the relations for
the two domains (\ref{eq:-17}), (\ref{eq:-42}), (\ref{eq:-43})
and (\ref{eq:-44}), we have:
\[
a_{i}(v,u)=a_{i}(v,u)^{+}+a_{i}(v,u)^{-}
\]
\begin{eqnarray*}
a_{i}(v,u) & = & -\int_{\mathcal{S}}v^{+}\left(\kappa^{+}\nabla u^{+}\centerdot\mathrm{\mathbf{n}}\right)\mathrm{d}\mathrm{\Gamma}-\int_{\mathcal{S}}u^{+}\left(\kappa^{+}\nabla v^{+}\centerdot\mathrm{\mathbf{n}}\right)\mathrm{d}\mathrm{\Gamma}+\int_{\mathcal{S}}v^{+}\alpha^{+}u^{+}\mathrm{d}\mathrm{\Gamma}\\
 &  & +\int_{\mathcal{S}}v^{-}\left(\kappa^{-}\nabla u^{-}\centerdot\mathrm{\mathbf{n}}\right)\mathrm{d}\mathrm{\Gamma}+\int_{\mathcal{S}}u^{-}\left(\kappa^{-}\nabla v^{-}\centerdot\mathrm{\mathbf{n}}\right)\mathrm{d}\mathrm{\Gamma}+\int_{\mathcal{S}}v^{-}\alpha^{-}u^{-}\mathrm{d}\mathrm{\Gamma}
\end{eqnarray*}
\[
l_{i}(v)=l_{i}(v)^{+}+l_{i}(v)^{-}
\]
\begin{eqnarray*}
l_{i}(v) & = & -\int_{\mathcal{S}}g^{+}\left(\kappa^{+}\nabla v^{+}\centerdot\mathrm{\mathbf{n}}\right)\mathrm{d}\mathrm{\Gamma}+\int_{\mathcal{S}}v^{+}\alpha^{+}g^{+}\mathrm{d}\mathrm{\Gamma}+\int_{\mathcal{S}}g^{-}\left(\kappa^{-}\nabla v^{-}\centerdot\mathrm{\mathbf{n}}\right)\mathrm{d}\mathrm{\Gamma}+\int_{\mathcal{S}}v^{-}\alpha^{-}g^{-}\mathrm{d}\mathrm{\Gamma}
\end{eqnarray*}
and
\[
\int_{\mathcal{L}}v_{d}\bar{\mathrm{j}}\mbox{d}\Gamma=\int_{\mathcal{B}}v_{d}^{+}f\mbox{d}\Omega-\int_{\mathcal{B}}\nabla v_{d}^{+}\kappa^{+}\nabla u^{+}\mbox{d}\Omega+\int_{\mathcal{B}}v_{d}^{-}f\mbox{d}\Omega-\int_{\mathcal{B}}\nabla v_{d}^{-}\kappa^{-}\nabla u^{-}\mbox{d}\Omega
\]

\begin{center}
\begin{table}
\begin{centering}
\begin{turn}{90}
\begin{tabular}{|c|c|c|}
\hline 
\noalign{\vskip0.25cm}
Term & Jump Conditions & Dirichlet Condition\tabularnewline[0.25cm]
\hline 
\hline 
\noalign{\vskip0.25cm}
Bulk stiffness contribution terms & $\int_{\Omega}\nabla v\centerdot\kappa\nabla u\mathrm{d\Omega}$ & $\int_{\Omega}\nabla v\centerdot\kappa\nabla u\mathrm{d\Omega}$\tabularnewline[0.25cm]
\hline 
\noalign{\vskip0.25cm}
Bulk force contribution terms & $\int_{\Omega}vf\mathrm{d\Omega}$ & $\int_{\Omega}vf\mathrm{d\Omega}$\tabularnewline[0.25cm]
\hline 
\noalign{\vskip0.25cm}
\multirow{2}{*}{Nitsche's contribution to stiffness} & $-\int_{\mathcal{S}}\frac{v^{+}}{2}(\kappa^{+}\nabla u^{+}\centerdot\mathrm{\mathbf{n}})\mathrm{d}\mathrm{\Gamma}+\int_{\mathcal{S}}\frac{v^{-}}{2}(\kappa^{-}\nabla u^{-}\centerdot\mathrm{\mathbf{n}})\mathrm{d}\mathrm{\Gamma}$ & $-\int_{\mathcal{S}}v^{+}(\kappa^{+}\nabla u^{+}\centerdot\mathrm{\mathbf{n}})\mathrm{d}\mathrm{\Gamma}-\int_{\mathcal{S}}u^{+}(\kappa^{+}\nabla v^{+}\centerdot\mathrm{\mathbf{n}})\mathrm{d}\mathrm{\Gamma}$\tabularnewline[0.25cm]
\noalign{\vskip0.25cm}
 & $-\int_{\mathcal{S}}\frac{u^{+}}{2}(\kappa^{+}\nabla v^{+}\centerdot\mathrm{\mathbf{n}})\mathrm{d}\mathrm{\Gamma}+\int_{\mathcal{S}}\frac{u^{-}}{2}(\kappa^{-}\nabla v^{-}\centerdot\mathrm{\mathbf{n}})\mathrm{d}\mathrm{\Gamma}$ & $+\int_{\mathcal{S}}v^{-}(\kappa^{-}\nabla u^{-}\centerdot\mathrm{\mathbf{n}})\mathrm{d}\mathrm{\Gamma}+\int_{\mathcal{S}}u^{-}(\kappa^{-}\nabla v^{-}\centerdot\mathrm{\mathbf{n}})\mathrm{d}\mathrm{\Gamma}$\tabularnewline[0.25cm]
\hline 
\noalign{\vskip0.25cm}
Nitsche's contribution to force & $-\int_{\mathcal{S}}\frac{\bar{\mbox{i}}}{2}(\kappa^{+}\nabla v^{+}\centerdot\mathrm{\mathbf{n}})\mathrm{d}\mathrm{\Gamma}-\int_{\mathcal{S}}\frac{\bar{\mbox{i}}}{2}(\kappa^{-}\nabla v^{-}\centerdot\mathrm{\mathbf{n}})\mathrm{d}\mathrm{\Gamma}$ & $-\int_{\mathcal{S}}g^{+}(\kappa^{+}\nabla v^{+}\centerdot\mathrm{\mathbf{n}})\mathrm{d}\mathrm{\Gamma}+\int_{\mathcal{S}}g^{-}(\kappa^{-}\nabla v^{-}\centerdot\mathrm{\mathbf{n}})\mathrm{d}\mathrm{\Gamma}$\tabularnewline[0.25cm]
\hline 
\noalign{\vskip0.25cm}
Stabilization contribution to stiffness & $\int_{\mathcal{S}}v^{+}\alpha u^{+}\mathrm{d}\mathrm{\Gamma}+\int_{\mathcal{S}}v^{-}\alpha u^{-}\mathrm{d}\mathrm{\Gamma}$ & $\int_{\mathcal{S}}v^{+}\alpha^{+}u^{+}\mathrm{d}\mathrm{\Gamma}+\int_{\mathcal{S}}v^{-}\alpha^{-}u^{-}\mathrm{d}\mathrm{\Gamma}$\tabularnewline[0.25cm]
\hline 
\noalign{\vskip0.25cm}
Stabilization contribution to force & $\int_{\mathcal{S}}v^{+}\alpha\bar{\mbox{i}}\mathrm{d}\mathrm{\Gamma}-\int_{\mathcal{S}}v^{-}\alpha\bar{\mbox{i}}\mathrm{d}\mathrm{\Gamma}$ & $\int_{\mathcal{S}}v^{+}\alpha^{+}g^{+}\mathrm{d}\mathrm{\Gamma}+\int_{\mathcal{S}}v^{-}\alpha^{-}g^{-}\mathrm{d}\mathrm{\Gamma}$\tabularnewline[0.25cm]
\hline 
\noalign{\vskip0.25cm}
\multirow{2}{*}{Coupling contribution to stiffness} & $-\int_{\mathcal{S}}\frac{v^{+}}{2}(\kappa^{-}\nabla u^{-}\centerdot\mathrm{\mathbf{n}})\mathrm{d}\mathrm{\Gamma}+\int_{\mathcal{S}}\frac{u^{-}}{2}(\kappa^{+}\nabla v^{+}\centerdot\mathrm{\mathbf{n}})\mathrm{d}\mathrm{\Gamma}-\int_{\mathcal{S}}v^{+}\alpha u^{-}\mathrm{d}\mathrm{\Gamma}$ & \multirow{2}{*}{-}\tabularnewline[0.25cm]
\noalign{\vskip0.25cm}
 & $+\int_{\mathcal{S}}\frac{v^{-}}{2}(\kappa^{+}\nabla u^{+}\centerdot\mathrm{\mathbf{n}})\mathrm{d}\mathrm{\Gamma}-\int_{\mathcal{S}}\frac{u^{+}}{2}(\kappa^{-}\nabla v^{-}\centerdot\mathrm{\mathbf{n}})\mathrm{d}\mathrm{\Gamma}-\int_{\mathcal{S}}v^{-}\alpha u^{+}\mathrm{d}\mathrm{\Gamma}$ & \tabularnewline[0.25cm]
\hline 
\noalign{\vskip0.25cm}
Contribution of jump in  & \multirow{2}{*}{$\int_{\mathcal{S}}\frac{v^{+}}{2}\bar{\mbox{j}}\mathrm{d}\mathrm{\Gamma}+\int_{\mathcal{S}}\frac{v^{-}}{2}\bar{\mbox{j}}\mathrm{d}\mathrm{\Gamma}$} & \multirow{2}{*}{-}\tabularnewline[0.25cm]
\noalign{\vskip0.25cm}
interfacial flux to force &  & \tabularnewline[0.25cm]
\hline 
\noalign{\vskip0.25cm}
\multirow{2}{*}{Jump calculation} & \multirow{2}{*}{-} & $\int_{\mathcal{L}}v_{d}\bar{\mathrm{j}}\mbox{d}\Gamma=\int_{\mathcal{B}}v_{d}^{+}f\mbox{d}\Omega+\int_{\mathcal{B}}v_{d}^{-}f\mbox{d}\Omega$\tabularnewline[0.25cm]
\noalign{\vskip0.25cm}
 &  & $-\int_{\mathcal{B}}\nabla v_{d}^{+}\kappa^{+}\nabla u^{+}\mbox{d}\Omega-\int_{\mathcal{B}}\nabla v_{d}^{-}\kappa^{-}\nabla u^{-}\mbox{d}\Omega$\tabularnewline[0.25cm]
\hline 
\end{tabular}
\end{turn}
\par\end{centering}
\caption{Term by term comparison of the variational form of 'Jump' and 'Dirichlet'
problems}
\end{table}
\par\end{center}

\newpage{}

\section{Nitsche's method with Elastostatic (Cauchy Navier) Equations}

Consider the elastostatic equation of the form: 
\begin{equation}
\underbar{\mbox{div}}\underline{\underline{\sigma}}=-\underbar{\ensuremath{f}}
\end{equation}
With $\lambda>0$ being Lame's constant and $\mu>0$ being the shear
modulus of the material we have,
\[
\underline{\underline{\sigma}}=\lambda\mbox{tr}\underline{\underline{\varepsilon}}+2\mu\underline{\underline{\varepsilon}}
\]
where:
\[
\underline{\underbar{\ensuremath{\varepsilon}}}=\frac{1}{2}\left(\underbar{\ensuremath{\underline{\underline{\nabla}}}}^{T}u+\underbar{\ensuremath{\underline{\underbar{\ensuremath{\nabla}}}}}u\right).
\]
Consider the domain to be $\Omega=\Omega^{+}\cup\Omega^{-}$ divided
by the internal interface $\mathcal{S}$. We can write this in simple
terms as,
\begin{equation}
\mbox{div}\sigma^{\pm}=-f\mbox{\,\,\,\,\ in }\Omega^{\pm},\,\,\Omega=\Omega^{+}\cup\Omega^{-}\label{eq:-13}
\end{equation}
with:
\[
\sigma^{\pm}=C^{\pm}:\varepsilon(u)^{\pm}
\]
and with the Dirichlet boundary conditions:
\begin{equation}
u^{\pm}=u_{0}^{\pm}\mathrm{\,\,\,\,on}\,\Gamma_{\mbox{d}}^{\pm},\,\Gamma_{\mbox{d}}=\Gamma_{\mbox{d}}^{+}\cup\Gamma_{\mbox{d}}^{-}
\end{equation}
and Neumann boundary conditions: 
\[
\sigma^{\pm}n^{\pm}=t^{\pm}\mathrm{\,\,\,\,on}\,\Gamma_{\mbox{n}}^{\pm},\,\Gamma_{\mbox{n}}=\Gamma_{\mbox{n}}^{+}\cup\Gamma_{\mbox{n}}^{-}.
\]

The primary unknown, displacement $u_{i}$ over $\Omega$, can be
seen as a collection of displacements over each part, $u_{i}^{+}$
and $u_{i}^{-}$. Thus we can write:
\[
[[u]]=(u^{+}-u^{-})
\]
and: 
\[
<u>=\frac{1}{2}(u^{+}+u^{-})
\]

We can write the variational form of the problem defined by (\ref{eq:-13}):
\begin{equation}
\mbox{Find u\ensuremath{\in\mathbb{U}} , s.t }a_{b}(v,u)+a_{i}(v,u)=l_{b}(v)+l_{i}(v)\,\,\,\,\forall v\in\mathbb{U}_{0}\label{eq:-18}
\end{equation}
where now: 
\begin{eqnarray*}
a_{b}(v,u) & = & \int_{\Omega}\varepsilon(v)\sigma\mathrm{d\Omega}\\
l_{b}(v) & = & \int_{\Omega}vf\,\mathrm{d\Omega}+\int_{\Gamma_{\mbox{n}}}v\,t\,\mathrm{d\Gamma}
\end{eqnarray*}
 are the standard bulk contributions. We take $n^{-}=-n^{+}=n$, and
$\alpha$ as Nitsche's stability parameter. 

\subparagraph*{Jump Condition}

\[
\left[\left[u\right]\right]=\bar{\mathrm{i}},\,\,\,\,\left[\left[\sigma.n\right]\right]=\bar{\mathrm{j}}\,\,\,\,\mathrm{on}\,\mathcal{S}.
\]
 We can write, similarly to (\ref{eq:-14}),
\begin{eqnarray}
a_{i}(v,u) & = & -\int_{\mathcal{S}}[[v]]<\sigma>.\mathbf{n}\mathrm{d\Gamma}-\int_{\mathcal{S}}[[u]]<\sigma(v)>.\mathbf{n}\mathrm{d\Gamma}+\int_{\mathcal{S}}[[v]]\alpha[[u]]\mathrm{d\Gamma}\\
l_{i}(v) & = & -\int_{\mathcal{S}}\bar{\mbox{i}}<\sigma(v)>.\mathbf{n}\mathrm{d}\mathrm{\Gamma}+\int_{\mathcal{S}}[[v]]\alpha\bar{\mbox{i}}\mathrm{d}\mathrm{\Gamma}+\int_{\mathcal{S}}<v>\bar{\mbox{j}}\mathrm{d}\mathrm{\Gamma}
\end{eqnarray}

It is interesting to note here that in the case of internally traction
free problems, we have: 
\[
\sigma^{+}.n=\sigma^{-}.n=0\,\,\,\,\mbox{on }\mathcal{S}
\]
\[
\bar{\mbox{j}}=\sigma^{+}.n-\sigma^{-}.n=0
\]

\subparagraph*{Dirichlet Condition}

\[
u^{+}=g^{+},\,\,\,\,u^{-}=g^{-}\mbox{\,\,\,\,\ on }\mathcal{\mathcal{S}}
\]

We can write
\begin{eqnarray}
a_{i}(v,u) & = & a_{i}(v,u)^{+}+a_{i}(v,u)^{-}\\
l_{i}(v) & = & l_{i}(v)^{+}+l_{i}(v)^{-}
\end{eqnarray}
with the same normal direction as considered previously in (\ref{eq:-17})
and (\ref{eq:-42}): 
\begin{eqnarray}
a_{i}(v,u)^{+} & = & -\int_{\mathcal{S}}v^{+}\left(\sigma^{+}\right).\mathbf{n}\mathrm{d}\mathrm{\Gamma}-\int_{\mathcal{S}}u^{+}\left(\sigma(v)^{+}\right).\mathbf{n}\mathrm{d}\mathrm{\Gamma}+\int_{\mathcal{S}}v^{+}\alpha^{+}u^{+}\mathrm{d}\mathrm{\Gamma}\label{eq:-33}\\
a_{i}(v,u)^{-} & = & +\int_{\mathcal{S}}v^{-}\left(\sigma^{-}\right).\mathbf{n}\mathrm{d}\mathrm{\Gamma}+\int_{\mathcal{S}}u^{-}\left(\sigma(v)^{-}\right).\mathbf{n}\mathrm{d}\mathrm{\Gamma}+\int_{\mathcal{S}}v^{-}\alpha^{-}u^{-}\mathrm{d}\mathrm{\Gamma}
\end{eqnarray}
and
\begin{eqnarray}
l_{i}(v)^{+} & = & -\int_{\mathcal{S}}g^{+}\left(\sigma(v)^{+}\right).\mathbf{n}\mathrm{d}\mathrm{\Gamma}+\int_{\mathcal{S}}v^{+}\alpha^{+}g^{+}\mathrm{d}\mathrm{\Gamma}\\
l_{i}(v)^{-} & = & +\int_{\mathcal{S}}g^{-}\left(\sigma(v)^{-}\right).\mathbf{n}\mathrm{d}\mathrm{\Gamma}+\int_{\mathcal{S}}v^{-}\alpha^{-}g^{-}\mathrm{d}\mathrm{\Gamma}
\end{eqnarray}

This effectively gives us our two 'one-sided' problems.

We can find the approximate flux by doing the same as before in section
(\ref{sec:Application-on-Interface}): 
\begin{eqnarray}
\int_{\mathcal{L}}v_{d}\bar{\mathrm{j}}\mbox{d}\Gamma & = & \int_{\mathcal{B}}v_{d}f\,\mbox{d}\Omega+\int_{\Gamma_{\mbox{n}}}v_{d}t\,\mbox{d}\Gamma-\int_{\mathcal{B}}\varepsilon(v)\sigma\mbox{d}\Omega
\end{eqnarray}

Once the unknown displacement $u$ is obtained, we can calculate the
jump in flux j by simple post processing of the solution.

\subparagraph*{Comparison}

Similarly to what was done in section (\ref{sec:Application-on-Interface}),
we can have a comparison table (Table \ref{Table}) between both types
of problems.

\begin{center}
\begin{table}
\begin{centering}
\begin{turn}{90}
\begin{tabular}{|c|c|c|}
\hline 
\noalign{\vskip0.25cm}
Term & Jump Conditions & Dirichlet Condition\tabularnewline[0.25cm]
\hline 
\hline 
\noalign{\vskip0.25cm}
Bulk stiffness contribution terms & $\int_{\Omega^{-}}\varepsilon(v^{-})\sigma^{-}\mathrm{d\Omega}+\int_{\Omega^{+}}\varepsilon(v^{+})\sigma^{+}\mathrm{d\Omega}$ & $\int_{\Omega^{-}}\varepsilon(v^{-})\sigma^{-}\mathrm{d\Omega}+\int_{\Omega^{+}}\varepsilon(v^{+})\sigma^{+}\mathrm{d\Omega}$\tabularnewline[0.25cm]
\hline 
\noalign{\vskip0.25cm}
Bulk force contribution terms & $\int_{\Omega^{+}}v^{+}f\,\mathrm{d\Omega}+\int_{\Omega^{-}}v^{-}f\mathrm{\,d\Omega}+\int_{\Gamma_{\mbox{n}}^{+}}v^{+}t^{+}\mathrm{d\Gamma}+\int_{\Gamma_{\mbox{n}}^{-}}v^{-}t^{-}\mathrm{d\Gamma}$ & $\int_{\Omega^{+}}v^{+}f\,\mathrm{d\Omega}+\int_{\Omega^{-}}v^{-}f\,\mathrm{d\Omega}+\int_{\Gamma_{\mbox{n}}^{+}}v^{+}t^{+}\mathrm{d\Gamma}+\int_{\Gamma_{\mbox{n}}^{-}}v^{-}t^{-}\mathrm{d\Gamma}$\tabularnewline[0.25cm]
\hline 
\noalign{\vskip0.25cm}
\multirow{2}{*}{Nitsche's contribution to stiffness} & $-\int_{\mathcal{S}}\frac{v^{+}}{2}(\sigma^{+}).\mathbf{n}\mathrm{d}\mathrm{\Gamma}+\int_{\mathcal{S}}\frac{u^{-}}{2}(\sigma(v)^{-}).\mathbf{n}\mathrm{d}\mathrm{\Gamma}$ & $-\int_{\mathcal{S}}v^{+}\left(\sigma^{+}\right).\mathbf{n}\mathrm{d}\mathrm{\Gamma}-\int_{\mathcal{S}}u^{+}\left(\sigma(v)^{+}\right).\mathbf{n}\mathrm{d}\mathrm{\Gamma}$\tabularnewline[0.25cm]
\noalign{\vskip0.25cm}
 & $-\int_{\mathcal{S}}\frac{u^{+}}{2}(\sigma(v)^{+}).\mathbf{n}\mathrm{d}\mathrm{\Gamma}+\int_{\mathcal{S}}\frac{v^{-}}{2}(\sigma^{-}).\mathbf{n}\mathrm{d}\mathrm{\Gamma}$ & $+\int_{\mathcal{S}}v^{-}\left(\sigma^{-}\right).\mathbf{n}\mathrm{d}\mathrm{\Gamma}+\int_{\mathcal{S}}u^{-}\left(\sigma(v)^{-}\right).\mathbf{n}\mathrm{d}\mathrm{\Gamma}$\tabularnewline[0.25cm]
\hline 
\noalign{\vskip0.25cm}
Nitsche's contribution to force & $-\int_{\mathcal{S}}\frac{\bar{\mbox{i}}}{2}(\sigma(v)^{+}).\mathbf{n}\mathrm{d}\mathrm{\Gamma}-\int_{\mathcal{S}}\frac{\bar{\mbox{i}}}{2}(\sigma(v)^{-}).\mathbf{n}\mathrm{d}\mathrm{\Gamma}$ & $-\int_{\mathcal{S}}g^{+}(\sigma(v)^{+}).\mathbf{n}\mathrm{d}\mathrm{\Gamma}+\int_{\mathcal{S}}g^{-}(\sigma(v)^{-}).\mathbf{n}\mathrm{d}\mathrm{\Gamma}$\tabularnewline[0.25cm]
\hline 
\noalign{\vskip0.25cm}
Stabilization contribution to stiffness & $\int_{\mathcal{S}}v^{+}\alpha u^{+}\mathrm{d}\mathrm{\Gamma}+\int_{\mathcal{S}}v^{-}\alpha u^{-}\mathrm{d}\mathrm{\Gamma}$ & $\int_{\mathcal{S}}v^{+}\alpha^{+}u^{+}\mathrm{d}\mathrm{\Gamma}+\int_{\mathcal{S}}v^{-}\alpha^{-}u^{-}\mathrm{d}\mathrm{\Gamma}$\tabularnewline[0.25cm]
\hline 
\noalign{\vskip0.25cm}
Stabilization contribution to force & $\int_{\mathcal{S}}v^{+}\alpha\bar{\mbox{i}}\mathrm{d}\mathrm{\Gamma}-\int_{\mathcal{S}}v^{-}\alpha\bar{\mbox{i}}\mathrm{d}\mathrm{\Gamma}$ & $\int_{\mathcal{S}}v^{+}\alpha^{+}g^{+}\mathrm{d}\mathrm{\Gamma}+\int_{\mathcal{S}}v^{-}\alpha^{-}g^{-}\mathrm{d}\mathrm{\Gamma}$\tabularnewline[0.25cm]
\hline 
\noalign{\vskip0.25cm}
\multirow{3}{*}{Coupling contribution to stiffness} & $-\int_{\mathcal{S}}\frac{v^{+}}{2}(\sigma^{-}).\mathbf{n}\mathrm{d}\mathrm{\Gamma}+\int_{\mathcal{S}}\frac{u^{-}}{2}(\sigma(v)^{+}).\mathbf{n}\mathrm{d}\mathrm{\Gamma}$ & \multirow{3}{*}{-}\tabularnewline[0.25cm]
\noalign{\vskip0.25cm}
 & $+\int_{\mathcal{S}}\frac{v^{-}}{2}(\sigma^{+}).\mathbf{n}\mathrm{d}\mathrm{\Gamma}-\int_{\mathcal{S}}\frac{u^{+}}{2}(\sigma(v)^{-}).\mathbf{n}\mathrm{d}\mathrm{\Gamma}$ & \tabularnewline[0.25cm]
\noalign{\vskip0.25cm}
 & $-\int_{\mathcal{S}}v^{+}\alpha u^{-}\mathrm{d}\mathrm{\Gamma}-\int_{\mathcal{S}}v^{-}\alpha u^{+}\mathrm{d}\mathrm{\Gamma}$ & \tabularnewline[0.25cm]
\hline 
\noalign{\vskip0.25cm}
Contribution of jump in  & \multirow{2}{*}{$\int_{\mathcal{S}}\frac{v^{+}}{2}\bar{\mbox{j}}\mathrm{d}\mathrm{\Gamma}+\int_{\mathcal{S}}\frac{v^{-}}{2}\bar{\mbox{j}}\mathrm{d}\mathrm{\Gamma}$} & \multirow{2}{*}{-}\tabularnewline[0.25cm]
\noalign{\vskip0.25cm}
interfacial flux to force &  & \tabularnewline[0.25cm]
\hline 
\noalign{\vskip0.25cm}
\multirow{1}{*}{Jump calculation} & \multirow{1}{*}{-} & $\int_{\mathcal{L}}v_{d}\bar{\mathrm{j}}\mbox{d}\Gamma=\int_{\mathcal{B}}v_{d}f\mbox{\,\ d}\Omega+\int_{\Gamma_{\mbox{n}}}v_{d}t\,\mbox{d}\Gamma-\int_{\mathcal{B}}\varepsilon(v)\sigma\mbox{d}\Omega$\tabularnewline[0.25cm]
\hline 
\end{tabular}
\end{turn}
\par\end{centering}
\caption{\label{Table}Term by term comparison of the variational form of 'Jump'
and 'Dirichlet' interfacial displacement problems}
\end{table}
\par\end{center}

\subsection{Discretization of the problem}

We consider the XFEM discretization by partitioning the domain into
a set of elements independently of the geometry and of any internal
interface. Near the interface, the enriched approximation of the solution
and its variation over an element take the form
\begin{eqnarray}
\boldsymbol{u}_{h}(\mathbf{x}) & = & \sum_{i\in I}\boldsymbol{u}_{i}N_{i}(\mathbf{x})+\sum_{i\in L}\boldsymbol{a}_{i}N_{i}(\mathbf{x})H(\mathbf{x})\,\,\,\,\,\,\,\,\mathbb{U}_{h}\subset\mathbb{U}\\
\boldsymbol{v}_{h}(\mathbf{x}) & = & \sum_{i\in I}\boldsymbol{v}_{i}N_{i}(\mathbf{x})+\sum_{i\in L}\boldsymbol{b}_{i}N_{i}(\mathbf{x})H(\mathbf{x})\,\,\,\,\,\,\,\,\mathbb{U}_{0h}\subset\mathbb{U}_{0}\label{eq:-12}
\end{eqnarray}
where $I$ is the set of nodes of the mesh, $\boldsymbol{u}_{i}$
is the classical (vectorial) degree of freedom of node $i$ and $N_{i}$
is the shape function associated with that node. $L\subset I$ is
the subset of nodes enriched by the Heaviside function. The corresponding
(vectorial) degrees of freedom are denoted $a_{i}$ . A node belongs
to $L$ if its support is cut in two by the interface. The jump function
$H(\mathbf{x})$ is discontinuous over the interface and constant
on each side.\cite{key-8} With a level set framework, one can define
\begin{equation}
H(\mathbf{x})=\begin{cases}
1 & \mbox{ls}_{n}(\mathbf{x})>0\\
-1 & \mbox{ls}_{n}(\mathbf{x})<0
\end{cases}
\end{equation}

We can now try to obtain the entire problem in matrix form, starting
with one element and assembling the entire system. If we consider
ls$_{n}(\Omega^{+})>0\mbox{ and ls}_{n}(\Omega^{-})<0$, we have:
\[
\boldsymbol{u}_{h}^{+}=\sum_{i\in I}\boldsymbol{u}_{i}N_{i}(\mathbf{x})+\sum_{i\in L}\boldsymbol{a}_{i}(+1)N_{i}(\mathbf{x})
\]
\[
\boldsymbol{u}_{h}^{-}=\sum_{i\in I}\boldsymbol{u}_{i}N_{i}(\mathbf{x})+\sum_{i\in L}\boldsymbol{a}_{i}(-1)N_{i}(\mathbf{x})
\]
The same formulation can also be applied to $\boldsymbol{v}_{h}$.
We consider a quasi-uniform partition $\Omega_{h}$ of the domain
$\Omega$ into non overlapping domains $\Omega_{e}$. Considering
$\mathcal{S}_{h}$ a partition of the interface $\mathcal{S}$ into
a set of non overlapping segments $\mathcal{S}_{e}$. We consider
an 'unfitted' or 'embedded' interface method. From the usual Galerkin
method formulation obtained in (\ref{eq:-18}) in terms of finite-dimensional
solution, we have: 
\begin{equation}
\mbox{Find u\ensuremath{\in\mathbb{U}} , s.t }a_{b}(v_{h},u_{h})+a_{i}(v_{h},u_{h})=l_{b}(v_{h})+l_{i}(v_{h})\,\,\,\,\forall v\in\mathbb{U}_{0}\label{eq:-18-1}
\end{equation}

\subsubsection{Jump Condition}

We can start with the weak formulation we obtained earlier.

\begin{eqnarray}
 &  & \int_{\Omega_{e}}\varepsilon^{T}(\boldsymbol{v}_{h})\sigma(\boldsymbol{u}_{h})\mathrm{d\Omega}-\int_{\mathcal{S}_{e}}[[\boldsymbol{v}_{h}]]^{T}<\sigma(\boldsymbol{u}_{h})>.\mathbf{n}\mathrm{d\Gamma}\label{eq:-24}\\
 &  & -\int_{\mathcal{S}_{e}}(<\sigma(\boldsymbol{v}_{h})>.\mathbf{n})^{T}[[\boldsymbol{u}_{h}]]\mathrm{d\Gamma}+\int_{\mathcal{S}_{e}}[[\boldsymbol{v}_{h}]]^{T}\alpha[[\boldsymbol{u}_{h}]]\mathrm{d\Gamma}\nonumber \\
 & = & \int_{\Omega_{e}}\boldsymbol{v}_{h}^{T}\boldsymbol{f}\mathrm{d\Omega}+\int_{\Gamma_{\mbox{n}_{e}}}\boldsymbol{v}_{h}^{T}\boldsymbol{t}\mathrm{d\Gamma}-\int_{\mathcal{S}_{e}}(<\sigma(\boldsymbol{v}_{h})>.\mathbf{n})^{T}\bar{\mbox{i}}\mathrm{d}\mathrm{\Gamma}\nonumber \\
 &  & +\int_{\mathcal{S}_{e}}[[\boldsymbol{v}_{h}]]^{T}\alpha\bar{\mbox{i}}\mathrm{d}\mathrm{\Gamma}+\int_{\mathcal{S}_{e}}<\boldsymbol{v}_{h}>^{T}\bar{\mbox{j}}\mathrm{d}\mathrm{\Gamma}\nonumber 
\end{eqnarray}
From (\ref{eq:-12}), considering a single element, dropping the subscript
$e$ and using the constitutive law of linear elasticity (Hook's Law),
we can write:
\[
\varepsilon(\boldsymbol{u}_{h})=\boldsymbol{B}_{i}\boldsymbol{u}_{h}
\]
\[
\sigma(\boldsymbol{u}_{h})=\boldsymbol{D}\boldsymbol{B}_{i}\boldsymbol{u}_{h}
\]
\[
\sigma(\boldsymbol{u}_{h}).n=\boldsymbol{n}^{T}\boldsymbol{D}\boldsymbol{B}_{i}\boldsymbol{u}_{h}
\]
where:
\begin{eqnarray*}
\boldsymbol{B}_{i}\boldsymbol{u}_{h} & = & \left[\begin{array}{cccccc}
\varepsilon_{xx} & \varepsilon_{yy} & \varepsilon_{zz} & 2\varepsilon_{yz} & 2\varepsilon_{zx} & 2\varepsilon_{xy}\end{array}\right]^{T}\\
 & = & \nabla\boldsymbol{N}_{i}\boldsymbol{u}_{h},
\end{eqnarray*}
\[
\boldsymbol{n}=\left[\begin{array}{cccccc}
n_{x} & 0 & 0 & 0 & n_{z} & n_{y}\\
0 & n_{y} & 0 & n_{z} & 0 & n_{x}\\
0 & 0 & n_{z} & n_{y} & n_{x} & 0
\end{array}\right]^{T}
\]
and with $\boldsymbol{D}$ being the constitutive relation between
stress and strain in matrix form. We also have:
\begin{eqnarray*}
<\boldsymbol{u}_{h}> & = & \frac{1}{2}(\boldsymbol{u}_{h}^{+}+\boldsymbol{u}_{h}^{-})\\
 & = & \frac{1}{2}(\boldsymbol{u}_{i}+\boldsymbol{a}_{i}+\boldsymbol{u}_{i}-\boldsymbol{a}_{i})\\
 & = & \boldsymbol{u}_{i}
\end{eqnarray*}
\begin{eqnarray*}
<\sigma(\boldsymbol{u}_{h})>.\mathbf{n} & = & \frac{1}{2}(\sigma(\boldsymbol{u}_{h}^{+}).\mathbf{n}+\sigma(\boldsymbol{u}_{h}^{-}).\mathbf{n})\\
 & = & \frac{1}{2}(\sigma(\boldsymbol{u}_{h}^{+}).\mathbf{n}+\sigma(\boldsymbol{u}_{h}^{-}).\mathbf{n})\\
 & = & \frac{1}{2}(\boldsymbol{n}^{T}\boldsymbol{D}^{+}\boldsymbol{B}_{i}\boldsymbol{u}_{i}+\boldsymbol{n}^{T}\boldsymbol{D}^{+}\boldsymbol{B}_{i}\boldsymbol{a}_{i}+\boldsymbol{n}^{T}\boldsymbol{D}^{-}\boldsymbol{B}_{i}\boldsymbol{u}_{i}-\boldsymbol{n}^{T}\boldsymbol{D}^{-}\boldsymbol{B}_{i}\boldsymbol{a}_{i})\\
 & = & \frac{1}{2}(\boldsymbol{n}^{T}\boldsymbol{D}^{+}\boldsymbol{B}_{i}+\boldsymbol{n}^{T}\boldsymbol{D}^{-}\boldsymbol{B}_{i})\boldsymbol{u}_{i}+\frac{1}{2}(\boldsymbol{n}^{T}\boldsymbol{D}^{+}\boldsymbol{B}_{i}-\boldsymbol{n}^{T}\boldsymbol{D}^{-}\boldsymbol{B}_{i})\boldsymbol{a}_{i}
\end{eqnarray*}
and:
\begin{eqnarray*}
[[\boldsymbol{u}_{h}]] & = & \boldsymbol{u}_{h}^{+}-\boldsymbol{u}_{h}^{-}\\
 & = & \boldsymbol{u}_{i}+\boldsymbol{a}_{i}-\boldsymbol{u}_{i}+\boldsymbol{a}_{i}\\
 & = & 2\boldsymbol{a}_{i}
\end{eqnarray*}
 Applying these into the Galerkin form for the element gives us,
\begin{eqnarray}
 &  & \sum_{e}\int_{\Omega_{e}}(\boldsymbol{B}_{i}\boldsymbol{v}_{i}+\boldsymbol{H}\boldsymbol{B}_{i}\boldsymbol{b}_{i})^{T}\sigma(\boldsymbol{u}_{h})\mathrm{d\Omega}-\sum_{e}\int_{\mathcal{S}_{e}}(2\boldsymbol{N}_{i}\boldsymbol{b}_{i})^{T}<\sigma(\boldsymbol{u}_{h})>.\mathbf{n}\mathrm{d\Gamma}\nonumber \\
 &  & -\sum_{e}\int_{\mathcal{S}_{e}}\left[\frac{1}{2}(\boldsymbol{n}^{T}\boldsymbol{D}^{+}\boldsymbol{B}_{i}+\boldsymbol{n}^{T}\boldsymbol{D}^{-}\boldsymbol{B}_{i})\boldsymbol{v}_{i}+\frac{1}{2}(\boldsymbol{n}^{T}\boldsymbol{D}^{+}\boldsymbol{B}_{i}-\boldsymbol{n}^{T}\boldsymbol{D}^{-}\boldsymbol{B}_{i})\boldsymbol{b}_{i}\right]^{T}[[\boldsymbol{u}_{h}]]\mathrm{d\Gamma}\nonumber \\
 &  & +\sum_{e}\int_{\mathcal{S}_{e}}(2\boldsymbol{N}_{i}\boldsymbol{b}_{i})^{T}\alpha[[\boldsymbol{u}_{h}]]\mathrm{d\Gamma}\nonumber \\
 & = & \sum_{e}\int_{\Omega_{e}}(\boldsymbol{N}_{i}\boldsymbol{v}_{i}+\boldsymbol{H}\boldsymbol{N}_{i}\boldsymbol{b}_{i})^{T}\boldsymbol{f}\mathrm{d\Omega}+\sum_{e}\int_{\Gamma_{\mbox{n}_{e}}}(\boldsymbol{N}_{i}\boldsymbol{v}_{i}+\boldsymbol{H}\boldsymbol{N}_{i}\boldsymbol{b}_{i})^{T}\boldsymbol{t}\mathrm{d\Gamma}\nonumber \\
 &  & +\sum_{e}\int_{\mathcal{S}_{e}}(2\boldsymbol{N}_{i}\boldsymbol{b}_{i})^{T}\alpha\bar{\mbox{i}}\mathrm{d}\mathrm{\Gamma}+\sum_{e}\int_{\mathcal{S}_{e}}\frac{1}{2}(2\boldsymbol{N}_{i}\boldsymbol{v}_{i})^{T}\bar{\mbox{j}}\mathrm{d}\mathrm{\Gamma}\nonumber \\
 &  & -\sum_{e}\int_{\mathcal{S}_{e}}\left[\frac{1}{2}(\boldsymbol{n}^{T}\boldsymbol{D}^{+}\boldsymbol{B}_{i}+\boldsymbol{n}^{T}\boldsymbol{D}^{-}\boldsymbol{B}_{i})\boldsymbol{v}_{i}+\frac{1}{2}(\boldsymbol{n}^{T}\boldsymbol{D}^{+}\boldsymbol{B}_{i}-\boldsymbol{n}^{T}\boldsymbol{D}^{-}\boldsymbol{B}_{i})\boldsymbol{b}_{i}\right]^{T}\bar{\mbox{i}}\mathrm{d}\mathrm{\Gamma}
\end{eqnarray}
Expanding the terms inside the brackets, and dropping the summation
over all elements:
\[
\int_{\Omega_{e}}(\boldsymbol{v}_{i}^{T}\boldsymbol{B}_{i}^{T}+\boldsymbol{b}_{i}^{T}\boldsymbol{H}\boldsymbol{B}_{i}{}^{T})\sigma(\boldsymbol{u}_{h})\mathrm{d\Omega}-2\int_{\mathcal{S}_{e}}\boldsymbol{b}_{i}^{T}\boldsymbol{N}_{i}{}^{T}<\sigma(\boldsymbol{u}_{h}).n>\mathrm{d\Gamma}+2\int_{\mathcal{S}_{e}}\boldsymbol{b}_{i}^{T}\boldsymbol{N}_{i}{}^{T}\alpha[[\boldsymbol{u}_{h}]]\mathrm{d\Gamma}
\]
\[
-\frac{1}{2}\int_{\mathcal{S}_{e}}\boldsymbol{v}_{i}^{T}(\boldsymbol{B}_{i}^{T}\left[\boldsymbol{D}^{+}\right]^{T}\boldsymbol{n}+\boldsymbol{B}_{i}^{T}\left[\boldsymbol{D}^{-}\right]^{T}\boldsymbol{n})[[\boldsymbol{u}_{h}]]\mathrm{d\Gamma}-\frac{1}{2}\int_{\mathcal{S}_{e}}\boldsymbol{b}_{i}^{T}(\boldsymbol{B}_{i}^{T}\left[\boldsymbol{D}^{+}\right]^{T}\boldsymbol{n}-\boldsymbol{B}_{i}^{T}\left[\boldsymbol{D}^{-}\right]^{T}\boldsymbol{n}))[[\boldsymbol{u}_{h}]]\mathrm{d\Gamma}
\]
\[
=\int_{\Omega_{e}}(\boldsymbol{v}_{i}^{T}\boldsymbol{N}_{i}^{T}+\boldsymbol{b}_{i}^{T}\boldsymbol{H}\boldsymbol{N}_{i}^{T})\boldsymbol{f}\mathrm{d\Omega}+\int_{\Gamma_{\mbox{n}_{e}}}(\boldsymbol{v}_{i}^{T}\boldsymbol{N}_{i}^{T}+\boldsymbol{b}_{i}^{T}\boldsymbol{H}\boldsymbol{N}_{i}{}^{T})\boldsymbol{t}\mathrm{d\Gamma}+2\int_{\mathcal{S}_{e}}\boldsymbol{b}_{i}^{T}\boldsymbol{N}_{i}{}^{T}\alpha\bar{\mbox{i}}\mathrm{d}\mathrm{\Gamma}+\int_{\mathcal{S}_{e}}\boldsymbol{v}_{i}^{T}\boldsymbol{N}_{i}{}^{T}\bar{\mbox{j}}\mathrm{d}\mathrm{\Gamma}
\]
\begin{equation}
-\frac{1}{2}\int_{\mathcal{S}_{e}}\boldsymbol{v}_{i}^{T}(\boldsymbol{B}_{i}^{T}\left[\boldsymbol{D}^{+}\right]^{T}\boldsymbol{n}+\boldsymbol{B}_{i}^{T}\left[\boldsymbol{D}^{-}\right]^{T}\boldsymbol{n})\bar{\mbox{i}}\mathrm{d}\mathrm{\Gamma}-\frac{1}{2}\int_{\mathcal{S}_{e}}\boldsymbol{b}_{i}^{T}(\boldsymbol{B}_{i}^{T}\left[\boldsymbol{D}^{+}\right]^{T}\boldsymbol{n}-\boldsymbol{B}_{i}^{T}\left[\boldsymbol{D}^{-}\right]^{T}\boldsymbol{n}))\bar{\mbox{i}}\mathrm{d}\mathrm{\Gamma}\label{eq:-25}
\end{equation}
Grouping $\boldsymbol{v}_{i}^{T}$and $\boldsymbol{b}_{i}^{T}$ and
knowing their arbitrariness, we can write (\ref{eq:-25}) in a two
equation form, 
\begin{eqnarray}
 &  & \int_{\Omega_{e}}\boldsymbol{B}_{i}^{T}\sigma(\boldsymbol{u}_{h})\mathrm{d\Omega}-\frac{1}{2}\int_{\mathcal{S}_{e}}(\boldsymbol{B}_{i}^{T}\left[\boldsymbol{D}^{+}\right]^{T}\boldsymbol{n}+\boldsymbol{B}_{i}^{T}\left[\boldsymbol{D}^{-}\right]^{T}\boldsymbol{n})[[\boldsymbol{u}_{h}]]\mathrm{d\Gamma}\nonumber \\
 & = & \int_{\Omega_{e}}\boldsymbol{N}_{i}^{T}\boldsymbol{f}\mathrm{d\Omega}+\int_{\Gamma_{\mbox{n}_{e}}}\boldsymbol{N}_{i}^{T}\boldsymbol{t}\mathrm{d\Gamma}+\int_{\mathcal{S}_{e}}\boldsymbol{N}_{i}{}^{T}\bar{\mbox{j}}\mathrm{d}\mathrm{\Gamma}-\frac{1}{2}\int_{\mathcal{S}_{e}}(\boldsymbol{B}_{i}^{T}\left[\boldsymbol{D}^{+}\right]^{T}\boldsymbol{n}+\boldsymbol{B}_{i}^{T}\left[\boldsymbol{D}^{-}\right]^{T}\boldsymbol{n})\bar{\mbox{i}}\mathrm{d}\mathrm{\Gamma}
\end{eqnarray}

\begin{eqnarray}
 &  & \int_{\Omega_{e}}\boldsymbol{H}\boldsymbol{B}_{i}{}^{T}\sigma(\boldsymbol{u}_{h})\mathrm{d\Omega}-2\int_{\mathcal{S}_{e}}\boldsymbol{N}_{i}{}^{T}<\sigma(\boldsymbol{u}_{h}).n>\mathrm{d\Gamma}\nonumber \\
 &  & -\frac{1}{2}\int_{\mathcal{S}_{e}}(\boldsymbol{B}_{i}^{T}\left[\boldsymbol{D}^{+}\right]^{T}\boldsymbol{n}-\boldsymbol{B}_{i}^{T}\left[\boldsymbol{D}^{-}\right]^{T}\boldsymbol{n}))[[\boldsymbol{u}_{h}]]\mathrm{d\Gamma}+2\int_{\mathcal{S}_{e}}\boldsymbol{N}_{i}{}^{T}\alpha[[\boldsymbol{u}_{h}]]\mathrm{d\Gamma}\nonumber \\
 & = & \int_{\Omega_{e}}\boldsymbol{H}\boldsymbol{N}_{i}^{T}\boldsymbol{f}\mathrm{d\Omega}+\int_{\Gamma_{\mbox{n}_{e}}}\boldsymbol{H}\boldsymbol{N}_{i}{}^{T}\boldsymbol{t}\mathrm{d\Gamma}-\frac{1}{2}\int_{\mathcal{S}_{e}}(\boldsymbol{B}_{i}^{T}\left[\boldsymbol{D}^{+}\right]^{T}\boldsymbol{n}\nonumber \\
 &  & -\boldsymbol{B}_{i}^{T}\left[\boldsymbol{D}^{-}\right]^{T}\boldsymbol{n})\bar{\mbox{i}}\mathrm{d}\mathrm{\Gamma}+2\int_{\mathcal{S}_{e}}\boldsymbol{N}_{i}{}^{T}\alpha\bar{\mbox{i}}\mathrm{d}\mathrm{\Gamma}
\end{eqnarray}
Applying the constitutive law,
\begin{eqnarray}
 &  & \int_{\Omega_{e}}\boldsymbol{B}_{i}^{T}\boldsymbol{D}(\boldsymbol{B}_{j}\boldsymbol{u}_{j}+\boldsymbol{B}_{j}\boldsymbol{H}\boldsymbol{a}_{j})\mathrm{d\Omega}-\int_{\mathcal{S}_{e}}(\boldsymbol{B}_{i}^{T}\left[\boldsymbol{D}^{+}\right]^{T}\boldsymbol{n}+\boldsymbol{B}_{i}^{T}\left[\boldsymbol{D}^{-}\right]^{T}\boldsymbol{n})\boldsymbol{N}_{j}\boldsymbol{a}_{j}\mathrm{d\Gamma}\nonumber \\
 & = & \int_{\Omega_{e}}\boldsymbol{N}_{i}^{T}\boldsymbol{f}\mathrm{d\Omega}+\int_{\Gamma_{\mbox{n}_{e}}}\boldsymbol{N}_{i}^{T}\boldsymbol{t}\mathrm{d\Gamma}+\int_{\mathcal{S}_{e}}\boldsymbol{N}_{i}{}^{T}\bar{\mbox{j}}\mathrm{d}\mathrm{\Gamma}\nonumber \\
 &  & -\frac{1}{2}\int_{\mathcal{S}_{e}}\boldsymbol{B}_{i}^{T}\left[\boldsymbol{D}^{+}\right]^{T}\boldsymbol{n}\bar{\mbox{i}}\mathrm{d}\mathrm{\Gamma}-\frac{1}{2}\int_{\mathcal{S}_{e}}\boldsymbol{B}_{i}^{T}\left[\boldsymbol{D}^{-}\right]^{T}\boldsymbol{n}\bar{\mbox{i}}\mathrm{d}\mathrm{\Gamma}\label{eq:-26}
\end{eqnarray}

\begin{eqnarray}
 &  & \int_{\Omega_{e}}\boldsymbol{H}\boldsymbol{B}_{i}{}^{T}\boldsymbol{D}(\boldsymbol{B}_{j}\boldsymbol{u}_{j}+\boldsymbol{B}_{j}\boldsymbol{H}\boldsymbol{a}_{j})\mathrm{d\Omega}-\int_{\mathcal{S}_{e}}\boldsymbol{N}_{i}{}^{T}(\boldsymbol{n}^{T}\boldsymbol{D}^{+}\boldsymbol{B}_{j}+\boldsymbol{n}^{T}\boldsymbol{D}^{-}\boldsymbol{B}_{j})\boldsymbol{u}_{j}\mathrm{d\Gamma}\nonumber \\
 & - & \int_{\mathcal{S}_{e}}\boldsymbol{N}_{i}{}^{T}(\boldsymbol{n}^{T}\boldsymbol{D}^{+}\boldsymbol{B}_{j}-\boldsymbol{n}^{T}\boldsymbol{D}^{-}\boldsymbol{B}_{j})\boldsymbol{a}_{j}\mathrm{d\Gamma}-\int_{\mathcal{S}_{e}}(\boldsymbol{B}_{i}^{T}\left[\boldsymbol{D}^{+}\right]^{T}\boldsymbol{n}\boldsymbol{N}_{j}-\boldsymbol{B}_{i}^{T}\left[\boldsymbol{D}^{-}\right]^{T}\boldsymbol{n}\boldsymbol{N}_{j})\boldsymbol{a}_{j}\mathrm{d\Gamma}\nonumber \\
 & + & 4\int_{\mathcal{S}_{e}}\boldsymbol{N}_{i}{}^{T}\alpha\boldsymbol{N}_{j}\boldsymbol{a}_{j}\mathrm{d\Gamma}\nonumber \\
= &  & \int_{\Omega_{e}}\boldsymbol{H}\boldsymbol{N}_{i}^{T}\boldsymbol{f}\mathrm{d\Omega}+\int_{\Gamma_{\mbox{n}_{e}}}\boldsymbol{H}\boldsymbol{N}_{i}{}^{T}\boldsymbol{t}\mathrm{d\Gamma}+2\int_{\mathcal{S}_{e}}\boldsymbol{N}_{i}{}^{T}\alpha\bar{\mbox{i}}\mathrm{d}\mathrm{\Gamma}\nonumber \\
 & - & \frac{1}{2}\int_{\mathcal{S}_{e}}\boldsymbol{B}_{i}^{T}\left[\boldsymbol{D}^{+}\right]^{T}\boldsymbol{n}\bar{\mbox{i}}\mathrm{d}\mathrm{\Gamma}+\frac{1}{2}\int_{\mathcal{S}_{e}}\boldsymbol{B}_{i}^{T}\left[\boldsymbol{D}^{-}\right]^{T}\boldsymbol{n}\bar{\mbox{i}}\mathrm{d}\mathrm{\Gamma}\label{eq:-27}
\end{eqnarray}
Rearranging (\ref{eq:-26}) and (\ref{eq:-27}), we have
\begin{eqnarray}
 &  & (\int_{\Omega_{e}}\boldsymbol{B}_{i}^{T}\boldsymbol{D}\boldsymbol{B}_{j}\mathrm{d\Omega})\boldsymbol{u}_{j}\nonumber \\
 &  & +(\int_{\Omega_{e}}\boldsymbol{B}_{i}^{T}\boldsymbol{D}\boldsymbol{B}_{j}\boldsymbol{H}\mathrm{d\Omega}-\int_{\Omega_{e}}\boldsymbol{B}_{i}^{T}\left[\boldsymbol{D}^{+}\right]^{T}\boldsymbol{n}\boldsymbol{N}_{j}\mathrm{d\Gamma}-\int_{\Omega_{e}}\boldsymbol{B}_{i}^{T}\left[\boldsymbol{D}^{-}\right]^{T}\boldsymbol{n}\boldsymbol{N}_{j}\mathrm{d\Gamma})\boldsymbol{a}_{j}\nonumber \\
 & = & \int_{\Omega_{e}}\boldsymbol{N}_{i}^{T}\boldsymbol{f}\mathrm{d\Omega}+\int_{\Gamma_{\mbox{n}_{e}}}\boldsymbol{N}_{i}^{T}\boldsymbol{t}\mathrm{d\Gamma}\nonumber \\
 &  & -\frac{1}{2}\int_{\mathcal{S}_{e}}\boldsymbol{B}_{i}^{T}\left[\boldsymbol{D}^{+}\right]^{T}\boldsymbol{n}\bar{\mbox{i}}\mathrm{d}\mathrm{\Gamma}-\frac{1}{2}\int_{\mathcal{S}_{e}}\boldsymbol{B}_{i}^{T}\left[\boldsymbol{D}^{-}\right]^{T}\boldsymbol{n}\bar{\mbox{i}}\mathrm{d}\mathrm{\Gamma}+\int_{\mathcal{S}_{e}}\boldsymbol{N}_{i}{}^{T}\bar{\mbox{j}}\mathrm{d}\mathrm{\Gamma}
\end{eqnarray}

\begin{eqnarray}
 &  & \left(\int_{\Omega_{e}}\boldsymbol{H}\boldsymbol{B}_{i}{}^{T}\boldsymbol{D}\boldsymbol{B}_{j}\mathrm{d\Omega}-\int_{\mathcal{S}_{e}}\boldsymbol{N}_{i}{}^{T}\boldsymbol{n}^{T}\boldsymbol{D}^{+}\boldsymbol{B}_{j}\mathrm{d\Gamma}-\int_{\mathcal{S}_{e}}\boldsymbol{N}_{i}{}^{T}\boldsymbol{n}^{T}\boldsymbol{D}^{-}\boldsymbol{B}_{j}\mathrm{d\Gamma}\right)\boldsymbol{u}_{j}\nonumber \\
 &  & +(\int_{\Omega_{e}}\boldsymbol{H}^{2}\boldsymbol{B}_{i}{}^{T}\boldsymbol{D}\boldsymbol{B}_{j}\mathrm{d\Omega}+\int_{\mathcal{S}_{e}}4\boldsymbol{N}_{i}{}^{T}\alpha\boldsymbol{N}_{j}\mathrm{d\Gamma}-\int_{\mathcal{S}_{e}}\boldsymbol{N}_{i}{}^{T}\boldsymbol{n}^{T}\boldsymbol{D}^{+}\boldsymbol{B}_{j}\mathrm{d\Gamma}+\int_{\mathcal{S}_{e}}\boldsymbol{N}_{i}\boldsymbol{n}^{T}\boldsymbol{D}^{-}\boldsymbol{B}_{j}\mathrm{d\Gamma}\nonumber \\
 &  & -\int_{\mathcal{S}_{e}}\boldsymbol{B}_{i}^{T}\left[\boldsymbol{D}^{+}\right]^{T}\boldsymbol{n}\boldsymbol{N}_{j}\mathrm{d\Gamma}+\int_{\mathcal{S}_{e}}\boldsymbol{B}_{i}^{T}\left[\boldsymbol{D}^{-}\right]^{T}\boldsymbol{n}\boldsymbol{N}_{j}\mathrm{d\Gamma})\boldsymbol{a}_{j}\nonumber \\
= &  & \int_{\Omega_{e}}\boldsymbol{H}\boldsymbol{N}_{i}^{T}\boldsymbol{f}\mathrm{d\Omega}-\frac{1}{2}\int_{\mathcal{S}_{e}}\boldsymbol{B}_{i}^{T}\left[\boldsymbol{D}^{+}\right]^{T}\boldsymbol{n}\bar{\mbox{i}}\mathrm{d}\mathrm{\Gamma}+\frac{1}{2}\int_{\mathcal{S}_{e}}\boldsymbol{B}_{i}^{T}\left[\boldsymbol{D}^{-}\right]^{T}\boldsymbol{n}\bar{\mbox{i}}\mathrm{d}\mathrm{\Gamma}\nonumber \\
 &  & +\int_{\Gamma_{\mbox{n}_{e}}}\boldsymbol{H}\boldsymbol{N}_{i}{}^{T}\boldsymbol{t}\mathrm{d\Gamma}+\int_{\mathcal{S}_{e}}2\boldsymbol{N}_{i}{}^{T}\alpha\bar{\mbox{i}}\mathrm{d}\mathrm{\Gamma}
\end{eqnarray}
\begin{eqnarray}
\left[\begin{array}{cc}
 & \sum_{e}\int_{\Omega_{e}}\boldsymbol{B}_{i}^{T}\boldsymbol{D}\boldsymbol{B}_{j}\boldsymbol{H}\mathrm{d\Omega}\\
\sum_{e}\int_{\Omega_{e}}\boldsymbol{B}_{i}^{T}\boldsymbol{D}\boldsymbol{B}_{j}\mathrm{d\Omega} & -\sum_{e}\int_{\Omega_{e}}\boldsymbol{B}_{i}^{T}\left[\boldsymbol{D}^{+}\right]^{T}\boldsymbol{n}\boldsymbol{N}_{j}\mathrm{d\Gamma}\\
 & -\sum_{e}\int_{\Omega_{e}}\boldsymbol{B}_{i}^{T}\left[\boldsymbol{D}^{-}\right]^{T}\boldsymbol{n}\boldsymbol{N}_{j}\mathrm{d\Gamma}\\
\\
 & \sum_{e}\int_{\Omega_{e}}\boldsymbol{B}_{i}{}^{T}\boldsymbol{D}\boldsymbol{B}_{j}\mathrm{d\Omega}\\
 & -\sum_{e}\int_{\mathcal{S}_{e}}\boldsymbol{N}_{i}{}^{T}\boldsymbol{n}^{T}\boldsymbol{D}^{+}\boldsymbol{B}_{j}\mathrm{d\Gamma}\\
\sum_{e}\int_{\Omega_{e}}\boldsymbol{H}\boldsymbol{B}_{i}{}^{T}\boldsymbol{D}\boldsymbol{B}_{j}\mathrm{d\Omega} & +\sum_{e}\int_{\mathcal{S}_{e}}\boldsymbol{N}_{i}\boldsymbol{n}^{T}\boldsymbol{D}^{-}\boldsymbol{B}_{j}\mathrm{d\Gamma}\\
-\sum_{e}\int_{\mathcal{S}_{e}}\boldsymbol{N}_{i}{}^{T}\boldsymbol{n}^{T}\boldsymbol{D}^{+}\boldsymbol{B}_{j}\mathrm{d\Gamma} & -\sum_{e}\int_{\mathcal{S}_{e}}\boldsymbol{B}_{i}^{T}\left[\boldsymbol{D}^{+}\right]^{T}\boldsymbol{n}\boldsymbol{N}_{j}\mathrm{d\Gamma}\\
-\sum_{e}\int_{\mathcal{S}_{e}}\boldsymbol{N}_{i}{}^{T}\boldsymbol{n}^{T}\boldsymbol{D}^{-}\boldsymbol{B}_{j}\mathrm{d\Gamma} & +\sum_{e}\int_{\mathcal{S}_{e}}\boldsymbol{B}_{i}^{T}\left[\boldsymbol{D}^{-}\right]^{T}\boldsymbol{n}\boldsymbol{N}_{j}\mathrm{d\Gamma}\\
 & +4\sum_{e}\int_{\mathcal{S}_{e}}\boldsymbol{N}_{i}{}^{T}\alpha\boldsymbol{N}_{j}\mathrm{d\Gamma}
\end{array}\right]\left\{ \begin{array}{c}
\boldsymbol{u}_{j}\\
\boldsymbol{a}_{j}
\end{array}\right\}  & = & \left\{ \begin{array}{c}
\sum_{e}\int_{\Omega_{e}}\boldsymbol{N}_{i}^{T}\boldsymbol{f}\mathrm{d\Omega}\\
+\sum_{e}\int_{\Gamma_{\mbox{n}_{e}}}\boldsymbol{N}_{i}^{T}\boldsymbol{t}\mathrm{d\Gamma}\\
-\frac{1}{2}\sum_{e}\int_{\mathcal{S}_{e}}\boldsymbol{B}_{i}^{T}\left[\boldsymbol{D}^{+}\right]^{T}\boldsymbol{n}\bar{\mbox{i}}\mathrm{d}\mathrm{\Gamma}\\
-\frac{1}{2}\sum_{e}\int_{\mathcal{S}_{e}}\boldsymbol{B}_{i}^{T}\left[\boldsymbol{D}^{-}\right]^{T}\boldsymbol{n}\bar{\mbox{i}}\mathrm{d}\mathrm{\Gamma}\\
+\sum_{e}\int_{\mathcal{S}_{e}}\boldsymbol{N}_{i}{}^{T}\bar{\mbox{j}}\mathrm{d}\mathrm{\Gamma}\\
\\
\sum_{e}\int_{\Omega_{e}}\boldsymbol{H}\boldsymbol{N}_{i}^{T}\boldsymbol{f}\mathrm{d\Omega}\\
+\sum_{e}\int_{\Gamma_{\mbox{n}_{e}}}\boldsymbol{H}\boldsymbol{N}_{i}{}^{T}\boldsymbol{t}\mathrm{d\Gamma}\\
-\frac{1}{2}\sum_{e}\int_{\mathcal{S}_{e}}\boldsymbol{B}_{i}^{T}\left[\boldsymbol{D}^{+}\right]^{T}\boldsymbol{n}\bar{\mbox{i}}\mathrm{d}\mathrm{\Gamma}\\
+\frac{1}{2}\sum_{e}\int_{\mathcal{S}_{e}}\boldsymbol{B}_{i}^{T}\left[\boldsymbol{D}^{-}\right]^{T}\boldsymbol{n}\bar{\mbox{i}}\mathrm{d}\mathrm{\Gamma}\\
+\sum_{e}\int_{\mathcal{S}_{e}}2\boldsymbol{N}_{i}{}^{T}\alpha\bar{\mbox{i}}\mathrm{d}\mathrm{\Gamma}
\end{array}\right\} \nonumber \\
\end{eqnarray}

\begin{equation}
\left[\begin{array}{cc}
\mathbf{K}_{b} & \mathbf{K}_{b}\boldsymbol{H}-\left[\mathbf{K}_{n}^{+}\right]^{T}-\left[\mathbf{K}_{n}^{-}\right]^{T}\\
\\
\boldsymbol{H}\mathbf{K}_{b}-\mathbf{K}_{n}^{+}-\mathbf{K}_{n}^{-} & \mathbf{K}_{b}-\mathbf{K}_{n}^{+}+\mathbf{K}_{n}^{-}\\
 & -\left[\mathbf{K}_{n}^{+}\right]^{T}+\left[\mathbf{K}_{n}^{-}\right]^{T}+4\mathbf{K}_{s}
\end{array}\right]\left\{ \begin{array}{c}
\boldsymbol{u}_{j}\\
\boldsymbol{a}_{j}
\end{array}\right\} =\left\{ \begin{array}{c}
\mathbf{f}_{b}+\mathbf{f}_{h}-\frac{1}{2}\mathbf{f}_{n}^{+}-\frac{1}{2}\mathbf{f}_{n}^{-}+\mathbf{f}_{\mbox{j}}\\
\boldsymbol{H}(\mathbf{f}_{b}+\mathbf{f}_{h})-\frac{1}{2}\mathbf{f}_{n}^{+}+\frac{1}{2}\mathbf{f}_{n}^{-}+2\mathbf{f}_{s}
\end{array}\right\} \label{eq:-19}
\end{equation}

We assemble all the elementary matrices to obtain the global system.
Here $\mathbf{K}_{b}$ is the bulk stiffness term, $\mathbf{K}_{n}^{\pm}$
is Nitsche's contribution to stiffness and $\mathbf{K}_{s}$ is the
stability term associated with the formulation. $\mathbf{f}_{b}$
is the bulk force term, $\mathbf{f}_{h}$ is Neumann's contribution
to the force term, $\mathbf{f}_{n}^{\pm}$ is  Nitsche's contribution
to the force term, $\mathbf{f}_{\mbox{j}}$ is the jump associated
with the flux and $\mathbf{f}_{s}$ is the stability parameter associated
to the force component. This can be compared to the tabular formulation
of the terms associated with the weak formulation. 

One can note here that when the interface is between two domains of
the same material, $\boldsymbol{D}^{+}=\boldsymbol{D}^{-}$, $\mathbf{K}_{n}^{+}=\mathbf{K}_{n}^{-}=\mathbf{K}_{n}$
and $\mathbf{f}_{n}^{+}=\mathbf{f}_{n}^{-}=\mathbf{f}_{n}$ resulting
in 
\begin{equation}
\left[\begin{array}{cc}
\mathbf{K}_{b} & \mathbf{K}_{b}\boldsymbol{H}-2\mathbf{K}_{n}^{T}\\
\\
\boldsymbol{H}\mathbf{K}_{b}-2\mathbf{K}_{n} & \mathbf{K}_{b}+4\mathbf{K}_{s}
\end{array}\right]\left\{ \begin{array}{c}
\boldsymbol{u}_{j}\\
\boldsymbol{a}_{j}
\end{array}\right\} =\left\{ \begin{array}{c}
\mathbf{f}_{b}+\mathbf{f}_{h}+\mathbf{f}_{\mbox{j}}-\mathbf{f}_{n}\\
\boldsymbol{H}(\mathbf{f}_{b}+\mathbf{f}_{h})+2\mathbf{f}_{s}
\end{array}\right\} \label{eq:-20}
\end{equation}

\subsubsection{Dirichlet condition\label{par:Dirichlet-condition}}

We can combine all the terms associated with the weak form of this
condition and try to obtain the system in a matrix form:

\begin{eqnarray}
 &  & \int_{\Omega_{e}}\varepsilon^{T}(\boldsymbol{v}_{h})\sigma(\boldsymbol{u}_{h})\mathrm{d\Omega}-\int_{\mathcal{S}_{e}}\boldsymbol{v}_{h}^{T+}\left(\sigma(\boldsymbol{u}_{h})^{+}\right).\mathbf{n}\mathrm{d}\mathrm{\Gamma}-\int_{\mathcal{S}_{e}}\left(\left(\sigma(\boldsymbol{v}_{h})^{+}\right).\mathbf{n}\right)^{T}\boldsymbol{u}_{h}^{+}\mathrm{d}\mathrm{\Gamma}\nonumber \\
 &  & +\int_{\mathcal{S}_{e}}\boldsymbol{v}_{h}^{T+}\alpha^{+}\boldsymbol{u}_{h}^{+}\mathrm{d}\mathrm{\Gamma}+\int_{\mathcal{S}_{e}}\boldsymbol{v}_{h}^{T-}\left(\sigma(\boldsymbol{u}_{h})^{-}\right).\mathbf{n}\mathrm{d}\mathrm{\Gamma}+\int_{\mathcal{S}_{e}}\left(\left(\sigma(\boldsymbol{v}_{h})^{-}\right).\mathbf{n}\right)^{T}\boldsymbol{u}_{h}^{-}\mathrm{d}\mathrm{\Gamma}+\int_{\mathcal{S}_{e}}\boldsymbol{v}_{h}^{T-}\alpha^{-}\boldsymbol{u}_{h}^{-}\mathrm{d}\mathrm{\Gamma}\nonumber \\
 & = & \int_{\Omega_{e}}\boldsymbol{v}_{h}^{T}\boldsymbol{f}\mathrm{d\Omega}+\int_{\Gamma_{\mbox{n}_{e}}}\boldsymbol{v}_{h}^{T}\boldsymbol{t}\mathrm{d\Gamma}+\int_{\mathcal{S}_{e}}\left(\left(\sigma(\boldsymbol{v}_{h})^{+}\right).\mathbf{n}\right)^{T}\boldsymbol{g}^{+}\mathrm{d}\mathrm{\Gamma}+\int_{\mathcal{S}_{e}}\boldsymbol{v}_{h}^{T+}\alpha^{+}\boldsymbol{g}^{+}\mathrm{d}\mathrm{\Gamma}\label{eq:-23}\\
 &  & -\int_{\mathcal{S}_{e}}\left(\left(\sigma(\boldsymbol{v}_{h})^{-}\right).\mathbf{n}\right)^{T}\boldsymbol{g}^{-}\mathrm{d}\mathrm{\Gamma}+\int_{\mathcal{S}_{e}}\boldsymbol{v}_{h}^{T-}\alpha^{-}\boldsymbol{g}^{-}\mathrm{d}\mathrm{\Gamma}\nonumber 
\end{eqnarray}
Similar to what was done for the jump condition formulation, we have:
\begin{eqnarray}
 &  & \int_{\Omega_{e}}\boldsymbol{v}_{i}^{T}\boldsymbol{B}_{i}^{T}\sigma(\boldsymbol{u}_{h})\mathrm{d\Omega}+\int_{\Omega_{e}}\boldsymbol{b}_{i}^{T}\boldsymbol{H}\boldsymbol{B}_{i}{}^{T}\sigma(\boldsymbol{u}_{h})\mathrm{d\Omega}+\int_{\mathcal{S}_{e}}\boldsymbol{v}_{i}^{T}\boldsymbol{N}_{i}{}^{T}\alpha^{-}\boldsymbol{u}_{h}^{-}\mathrm{d\Gamma}\nonumber \\
 &  & -\int_{\mathcal{S}_{e}}\boldsymbol{b}_{i}^{T}\boldsymbol{N}_{i}{}^{T}\alpha^{-}\boldsymbol{u}_{h}^{-}\mathrm{d\Gamma}+\int_{\mathcal{S}_{e}}\boldsymbol{v}_{i}^{T}\boldsymbol{N}_{i}{}^{T}\alpha^{+}\boldsymbol{u}_{h}^{+}\mathrm{d\Gamma}+\int_{\mathcal{S}_{e}}\boldsymbol{b}_{i}^{T}\boldsymbol{N}_{i}{}^{T}\alpha^{+}\boldsymbol{u}_{h}^{+}\mathrm{d\Gamma}\nonumber \\
 &  & -\int_{\mathcal{S}_{e}}\boldsymbol{v}_{i}^{T}\boldsymbol{N}_{i}{}^{T}(\sigma^{+}(\boldsymbol{u}_{h})).\mathbf{n}\mathrm{d\Gamma}-\int_{\mathcal{S}_{e}}\boldsymbol{b}_{i}^{T}\boldsymbol{N}_{i}{}^{T}(\sigma^{+}(\boldsymbol{u}_{h})).\mathbf{n}\mathrm{d\Gamma}-\int_{\mathcal{S}_{e}}\boldsymbol{v}_{i}^{T}\boldsymbol{B}_{i}{}^{T}\left[\boldsymbol{D}^{+}\right]^{T}\boldsymbol{n}\boldsymbol{u}_{h}^{+}\mathrm{d}\mathrm{\Gamma}\nonumber \\
 &  & -\int_{\mathcal{S}_{e}}\boldsymbol{b}_{i}^{T}\boldsymbol{B}_{i}{}^{T}\left[\boldsymbol{D}^{+}\right]^{T}\boldsymbol{n}\boldsymbol{u}_{h}^{+}\mathrm{d}\mathrm{\Gamma}+\int_{\mathcal{S}_{e}}\boldsymbol{v}_{i}^{T}\boldsymbol{B}_{i}{}^{T}\left[\boldsymbol{D}^{-}\right]^{T}\boldsymbol{n}\boldsymbol{u}_{h}^{-}\mathrm{d}\mathrm{\Gamma}\nonumber \\
 &  & +\int_{\mathcal{S}_{e}}\boldsymbol{b}_{i}^{T}\boldsymbol{B}_{i}{}^{T}\left[\boldsymbol{D}^{-}\right]^{T}\boldsymbol{n}\boldsymbol{u}_{h}^{-}\mathrm{d}\mathrm{\Gamma}+\int_{\mathcal{S}_{e}}\boldsymbol{v}_{i}^{T}\boldsymbol{N}_{i}{}^{T}(\sigma^{-}(\boldsymbol{u}_{h})).\mathbf{n}\mathrm{d\Gamma}-\int_{\mathcal{S}_{e}}\boldsymbol{b}_{i}^{T}\boldsymbol{N}_{i}{}^{T}(\sigma^{-}(\boldsymbol{u}_{h})).\mathbf{n}\mathrm{d\Gamma}\nonumber \\
 & = & \int_{\Omega_{e}}\boldsymbol{v}_{i}^{T}\boldsymbol{N}_{i}^{T}\boldsymbol{f}\mathrm{d\Omega}+\int_{\Omega_{e}}\boldsymbol{b}_{i}^{T}\boldsymbol{H}\boldsymbol{N}_{i}^{T}\boldsymbol{f}\mathrm{d\Omega}+\int_{\Gamma_{\mbox{n}_{e}}}\boldsymbol{v}_{i}^{T}\boldsymbol{N}_{i}^{T}\boldsymbol{t}\mathrm{d\Omega}+\int_{\Gamma_{\mbox{n}_{e}}}\boldsymbol{b}_{i}^{T}\boldsymbol{H}\boldsymbol{N}_{i}^{T}\boldsymbol{t}\mathrm{d\Omega}\nonumber \\
 &  & -\int_{\mathcal{S}_{e}}\boldsymbol{v}_{i}^{T}\boldsymbol{B}_{i}{}^{T}\left[\boldsymbol{D}^{+}\right]^{T}\boldsymbol{n}\boldsymbol{g}^{+}\mathrm{d}\mathrm{\Gamma}-\int_{\mathcal{S}_{e}}\boldsymbol{b}_{i}^{T}\boldsymbol{B}_{i}{}^{T}\left[\boldsymbol{D}^{+}\right]^{T}\boldsymbol{n}\boldsymbol{g}^{+}\mathrm{d}\mathrm{\Gamma}+\int_{\mathcal{S}_{e}}\boldsymbol{v}_{i}^{T}\boldsymbol{B}_{i}{}^{T}\left[\boldsymbol{D}^{-}\right]^{T}\boldsymbol{n}\boldsymbol{g}^{-}\mathrm{d}\mathrm{\Gamma}\nonumber \\
 &  & -\int_{\mathcal{S}_{e}}\boldsymbol{b}_{i}^{T}\boldsymbol{B}_{i}{}^{T}\left[\boldsymbol{D}^{-}\right]^{T}\boldsymbol{n}\boldsymbol{g}^{-}\mathrm{d}\mathrm{\Gamma}+\int_{\mathcal{S}_{e}}\boldsymbol{v}_{i}^{T}\boldsymbol{N}_{i}{}^{T}\alpha^{+}\boldsymbol{g}^{+}\mathrm{d\Gamma}+\int_{\mathcal{S}_{e}}\boldsymbol{b}_{i}^{T}\boldsymbol{N}_{i}{}^{T}\alpha^{+}\boldsymbol{g}^{+}\mathrm{d\Gamma}\nonumber \\
 &  & +\int_{\mathcal{S}_{e}}\boldsymbol{v}_{i}^{T}\boldsymbol{N}_{i}{}^{T}\alpha^{-}\boldsymbol{g}^{-}\mathrm{d\Gamma}-\int_{\mathcal{S}_{e}}\boldsymbol{b}_{i}^{T}\boldsymbol{N}_{i}{}^{T}\alpha^{-}\boldsymbol{g}^{-}\mathrm{d\Gamma}
\end{eqnarray}
Grouping $\boldsymbol{v}_{i}^{T}$and $\boldsymbol{b}_{i}^{T}$ and
knowing their arbitrariness, we can write:
\begin{eqnarray}
 &  & \int_{\Omega_{e}}\boldsymbol{B}_{i}^{T}\sigma(\boldsymbol{u}_{h})\mathrm{d\Omega}+\int_{\mathcal{S}_{e}}\boldsymbol{N}_{i}{}^{T}(\sigma^{+}(\boldsymbol{u}_{h})).\mathbf{n}\mathrm{d\Gamma}-\int_{\mathcal{S}_{e}}\boldsymbol{B}_{i}{}^{T}\left[\boldsymbol{D}^{+}\right]^{T}\boldsymbol{n}\boldsymbol{u}_{h}^{+}\mathrm{d}\mathrm{\Gamma}+\int_{\mathcal{S}_{e}}\boldsymbol{N}_{i}{}^{T}\alpha^{+}\boldsymbol{u}_{h}^{+}\mathrm{d\Gamma}\nonumber \\
 &  & +\int_{\mathcal{S}_{e}}\boldsymbol{N}_{i}{}^{T}(\sigma^{-}(\boldsymbol{u}_{h})).\mathbf{n}\mathrm{d\Gamma}+\int_{\mathcal{S}_{e}}\boldsymbol{B}_{i}{}^{T}\left[\boldsymbol{D}^{-}\right]^{T}\boldsymbol{n}\boldsymbol{u}_{h}^{-}\mathrm{d}\mathrm{\Gamma}+\int_{\mathcal{S}_{e}}\boldsymbol{N}_{i}{}^{T}\alpha^{-}\boldsymbol{u}_{h}^{-}\mathrm{d\Gamma}\nonumber \\
 & = & \int_{\Omega_{e}}\boldsymbol{N}_{i}^{T}\boldsymbol{f}\mathrm{d\Omega}+\int_{\Gamma_{\mbox{n}_{e}}}\boldsymbol{N}_{i}^{T}\boldsymbol{t}\mathrm{d\Omega}-\int_{\mathcal{S}_{e}}\boldsymbol{B}_{i}{}^{T}\left[\boldsymbol{D}^{+}\right]^{T}\boldsymbol{n}\boldsymbol{g}^{+}\mathrm{d}\mathrm{\Gamma}+\int_{\mathcal{S}_{e}}\boldsymbol{B}_{i}{}^{T}\left[\boldsymbol{D}^{-}\right]^{T}\boldsymbol{n}\boldsymbol{g}^{-}\mathrm{d}\mathrm{\Gamma}\nonumber \\
 &  & +\int_{\mathcal{S}_{e}}\boldsymbol{N}_{i}{}^{T}\alpha^{+}\boldsymbol{g}^{+}\mathrm{d\Gamma}+\int_{\mathcal{S}_{e}}\boldsymbol{N}_{i}{}^{T}\alpha^{-}\boldsymbol{g}^{-}\mathrm{d\Gamma}
\end{eqnarray}

\begin{eqnarray}
 &  & \int_{\Omega_{e}}\boldsymbol{H}\boldsymbol{B}_{i}{}^{T}\sigma(\boldsymbol{u}_{h})\mathrm{d\Omega}-\int_{\mathcal{S}_{e}}\boldsymbol{N}_{i}{}^{T}(\sigma^{+}(\boldsymbol{u}_{h})).\mathbf{n}\mathrm{d\Gamma}-\int_{\mathcal{S}_{e}}\boldsymbol{B}_{i}{}^{T}\left[\boldsymbol{D}^{+}\right]^{T}\boldsymbol{n}\boldsymbol{u}_{h}^{+}\mathrm{d}\mathrm{\Gamma}\nonumber \\
 &  & +\int_{\mathcal{S}_{e}}\boldsymbol{N}_{i}{}^{T}\alpha^{+}\boldsymbol{u}_{h}^{+}\mathrm{d\Gamma}-\int_{\mathcal{S}_{e}}\boldsymbol{N}_{i}{}^{T}(\sigma^{-}(\boldsymbol{u}_{h})).\mathbf{n}\mathrm{d\Gamma}\nonumber \\
 &  & -\int_{\mathcal{S}_{e}}\boldsymbol{B}_{i}{}^{T}\left[\boldsymbol{D}^{-}\right]^{T}\boldsymbol{n}\boldsymbol{u}_{h}^{-}\mathrm{d}\mathrm{\Gamma}-\int_{\mathcal{S}_{e}}\boldsymbol{N}_{i}{}^{T}\alpha^{-}\boldsymbol{u}_{h}^{-}\mathrm{d\Gamma}\nonumber \\
 & = & \int_{\Omega_{e}}\boldsymbol{H}\boldsymbol{N}_{i}^{T}\boldsymbol{f}\mathrm{d\Omega}+\int_{\Gamma_{\mbox{n}_{e}}}\boldsymbol{H}\boldsymbol{N}_{i}^{T}\boldsymbol{t}\mathrm{d\Omega}\nonumber \\
 &  & -\int_{\mathcal{S}_{e}}\boldsymbol{B}_{i}{}^{T}\left[\boldsymbol{D}^{+}\right]^{T}\boldsymbol{n}\boldsymbol{g}^{+}\mathrm{d}\mathrm{\Gamma}+\int_{\mathcal{S}_{e}}\boldsymbol{B}_{i}{}^{T}\left[\boldsymbol{D}^{-}\right]^{T}\boldsymbol{n}\boldsymbol{g}^{-}\mathrm{d}\mathrm{\Gamma}\nonumber \\
 &  & +\int_{\mathcal{S}_{e}}\boldsymbol{N}_{i}{}^{T}\alpha^{+}\boldsymbol{g}^{+}\mathrm{d\Gamma}-\int_{\mathcal{S}_{e}}\boldsymbol{N}_{i}{}^{T}\alpha^{-}\boldsymbol{g}^{-}\mathrm{d\Gamma}
\end{eqnarray}
Again applying the constitutive law and expanding gives us,
\begin{eqnarray}
 &  & \int_{\Omega_{e}}\boldsymbol{B}_{i}^{T}\boldsymbol{D}\boldsymbol{B}_{j}\boldsymbol{u}_{j}\mathrm{d\Omega}+\int_{\Omega_{e}}\boldsymbol{B}_{i}^{T}\boldsymbol{D}\boldsymbol{B}_{j}\boldsymbol{H}\boldsymbol{a}_{j}\mathrm{d\Omega}-\int_{\mathcal{S}_{e}}\boldsymbol{N}_{i}{}^{T}\boldsymbol{n}^{T}\boldsymbol{D}^{+}\boldsymbol{B}_{j}\boldsymbol{u}_{j}\mathrm{d\Gamma}-\int_{\mathcal{S}_{e}}\boldsymbol{N}_{i}{}^{T}\boldsymbol{n}^{T}\boldsymbol{D}^{+}\boldsymbol{B}_{j}\boldsymbol{a}_{j}\mathrm{d\Gamma}\nonumber \\
 &  & -\int_{\mathcal{S}_{e}}\boldsymbol{B}_{i}{}^{T}\left[\boldsymbol{D}^{+}\right]^{T}\boldsymbol{n}\boldsymbol{N}_{j}\boldsymbol{u}_{j}\mathrm{d}\mathrm{\Gamma}-\int_{\mathcal{S}_{e}}\boldsymbol{B}_{i}{}^{T}\left[\boldsymbol{D}^{+}\right]^{T}\boldsymbol{n}\boldsymbol{N}_{j}\boldsymbol{a}_{j}\mathrm{d}\mathrm{\Gamma}\nonumber \\
 &  & +\int_{\mathcal{S}_{e}}\boldsymbol{N}_{i}{}^{T}\alpha^{+}\boldsymbol{N}_{j}\boldsymbol{u}_{j}\mathrm{d\Gamma}+\int_{\mathcal{S}_{e}}\boldsymbol{N}_{i}{}^{T}\alpha^{+}\boldsymbol{N}_{j}\boldsymbol{a}_{j}\mathrm{d\Gamma}\nonumber \\
 &  & +\int_{\mathcal{S}_{e}}\boldsymbol{N}_{i}{}^{T}\boldsymbol{n}^{T}\boldsymbol{D}^{-}\boldsymbol{B}_{j}\boldsymbol{u}_{j}\mathrm{d\Gamma}-\int_{\mathcal{S}_{e}}\boldsymbol{N}_{i}{}^{T}\boldsymbol{n}^{T}\boldsymbol{D}^{-}\boldsymbol{B}_{j}\boldsymbol{a}_{j}\mathrm{d\Gamma}+\int_{\mathcal{S}_{e}}\boldsymbol{B}_{i}{}^{T}\left[\boldsymbol{D}^{-}\right]^{T}\boldsymbol{n}\boldsymbol{N}_{j}\boldsymbol{u}_{j}\mathrm{d}\mathrm{\Gamma}\nonumber \\
 &  & -\int_{\mathcal{S}_{e}}\boldsymbol{B}_{i}{}^{T}\left[\boldsymbol{D}^{-}\right]^{T}\boldsymbol{n}\boldsymbol{N}_{j}\boldsymbol{a}_{j}\mathrm{d}\mathrm{\Gamma}+\int_{\mathcal{S}_{e}}\boldsymbol{N}_{i}{}^{T}\alpha^{-}\boldsymbol{N}_{j}\boldsymbol{u}_{j}\mathrm{d\Gamma}-\int_{\mathcal{S}_{e}}\boldsymbol{N}_{i}{}^{T}\alpha^{-}\boldsymbol{N}_{j}\boldsymbol{a}_{j}\mathrm{d\Gamma}\nonumber \\
 & = & \int_{\Omega_{e}}\boldsymbol{N}_{i}^{T}\boldsymbol{f}\mathrm{d\Omega}+\int_{\Gamma_{\mbox{n}_{e}}}\boldsymbol{N}_{i}^{T}\boldsymbol{t}\mathrm{d\Omega}-\int_{\mathcal{S}_{e}}\boldsymbol{B}_{i}{}^{T}\left[\boldsymbol{D}^{+}\right]^{T}\boldsymbol{n}\boldsymbol{g}^{+}\mathrm{d}\mathrm{\Gamma}\nonumber \\
 &  & +\int_{\mathcal{S}_{e}}\boldsymbol{B}_{i}{}^{T}\left[\boldsymbol{D}^{-}\right]^{T}\boldsymbol{n}\boldsymbol{g}^{-}\mathrm{d}\mathrm{\Gamma}+\int_{\mathcal{S}_{e}}\boldsymbol{N}_{i}{}^{T}\alpha^{+}\boldsymbol{g}^{+}\mathrm{d\Gamma}+\int_{\mathcal{S}_{e}}\boldsymbol{N}_{i}{}^{T}\alpha^{-}\boldsymbol{g}^{-}\mathrm{d\Gamma}
\end{eqnarray}

\begin{eqnarray}
 &  & \int_{\Omega_{e}}\boldsymbol{H}\boldsymbol{B}_{i}{}^{T}\boldsymbol{D}\boldsymbol{B}_{j}\boldsymbol{u}_{j}\mathrm{d\Omega}+\int_{\Omega_{e}}\boldsymbol{B}_{i}{}^{T}\boldsymbol{D}\boldsymbol{B}_{j}\boldsymbol{a}_{j}\mathrm{d\Omega}-\int_{\mathcal{S}_{e}}\boldsymbol{N}_{i}{}^{T}\boldsymbol{n}^{T}\boldsymbol{D}^{+}\boldsymbol{B}_{j}\boldsymbol{u}_{j}\mathrm{d\Gamma}-\int_{\mathcal{S}_{e}}\boldsymbol{N}_{i}{}^{T}\boldsymbol{n}^{T}\boldsymbol{D}^{+}\boldsymbol{B}_{j}\boldsymbol{a}_{j}\mathrm{d\Gamma}\nonumber \\
 &  & -\int_{\mathcal{S}_{e}}\boldsymbol{B}_{i}{}^{T}\left[\boldsymbol{D}^{+}\right]^{T}\boldsymbol{n}\boldsymbol{N}_{j}\boldsymbol{u}_{j}\mathrm{d}\mathrm{\Gamma}-\int_{\mathcal{S}_{e}}\boldsymbol{B}_{i}{}^{T}\left[\boldsymbol{D}^{+}\right]^{T}\boldsymbol{n}\boldsymbol{N}_{j}\boldsymbol{a}_{j}\mathrm{d}\mathrm{\Gamma}\nonumber \\
 &  & +\int_{\mathcal{S}_{e}}\boldsymbol{N}_{i}{}^{T}\alpha^{+}\boldsymbol{N}_{j}\boldsymbol{u}_{j}\mathrm{d\Gamma}+\int_{\mathcal{S}_{e}}\boldsymbol{N}_{i}{}^{T}\alpha^{+}\boldsymbol{N}_{j}\boldsymbol{a}_{j}\mathrm{d\Gamma}\nonumber \\
 &  & -\int_{\mathcal{S}_{e}}\boldsymbol{N}_{i}{}^{T}\boldsymbol{n}^{T}\boldsymbol{D}^{-}\boldsymbol{B}_{j}\boldsymbol{u}_{j}\mathrm{d\Gamma}+\int_{\mathcal{S}_{e}}\boldsymbol{N}_{i}{}^{T}\boldsymbol{n}^{T}\boldsymbol{D}^{-}\boldsymbol{B}_{j}\boldsymbol{a}_{j}\mathrm{d\Gamma}-\int_{\mathcal{S}_{e}}\boldsymbol{B}_{i}{}^{T}\left[\boldsymbol{D}^{-}\right]^{T}\boldsymbol{n}\boldsymbol{N}_{j}\boldsymbol{u}_{j}\mathrm{d}\mathrm{\Gamma}\nonumber \\
 &  & +\int_{\mathcal{S}_{e}}\boldsymbol{B}_{i}{}^{T}\left[\boldsymbol{D}^{-}\right]^{T}\boldsymbol{n}\boldsymbol{N}_{j}\boldsymbol{a}_{j}\mathrm{d}\mathrm{\Gamma}-\int_{\mathcal{S}_{e}}\boldsymbol{N}_{i}{}^{T}\alpha^{-}\boldsymbol{N}_{j}\boldsymbol{u}_{j}\mathrm{d\Gamma}+\int_{\mathcal{S}_{e}}\boldsymbol{N}_{i}{}^{T}\alpha^{-}\boldsymbol{N}_{j}\boldsymbol{a}_{j}\mathrm{d\Gamma}\nonumber \\
 & = & \int_{\Omega_{e}}\boldsymbol{H}\boldsymbol{N}_{i}^{T}\boldsymbol{f}\mathrm{d\Omega}+\int_{\Gamma_{\mbox{n}_{e}}}\boldsymbol{H}\boldsymbol{N}_{i}^{T}\boldsymbol{t}\mathrm{d\Omega}-\int_{\mathcal{S}_{e}}\boldsymbol{B}_{i}{}^{T}\left[\boldsymbol{D}^{+}\right]^{T}\boldsymbol{n}\boldsymbol{g}^{+}\mathrm{d}\mathrm{\Gamma}\nonumber \\
 &  & -\int_{\mathcal{S}_{e}}\boldsymbol{B}_{i}{}^{T}\left[\boldsymbol{D}^{-}\right]^{T}\boldsymbol{n}\boldsymbol{g}^{-}\mathrm{d}\mathrm{\Gamma}+\int_{\mathcal{S}_{e}}\boldsymbol{N}_{i}{}^{T}\alpha^{+}\boldsymbol{g}^{+}\mathrm{d\Gamma}-\int_{\mathcal{S}_{e}}\boldsymbol{N}_{i}{}^{T}\alpha^{-}\boldsymbol{g}^{-}\mathrm{d\Gamma}
\end{eqnarray}
Rearranging and writing in a matrix form, similarly to (\ref{eq:-19}),
gives us
\begin{eqnarray}
 &  & (\int_{\Omega_{e}}\boldsymbol{B}_{i}^{T}\boldsymbol{D}\boldsymbol{B}_{j}\mathrm{d\Omega}-\int_{\mathcal{S}_{e}}\boldsymbol{N}_{i}{}^{T}\boldsymbol{n}^{T}\boldsymbol{D}^{+}\boldsymbol{B}_{j}\mathrm{d\Gamma}-\int_{\mathcal{S}_{e}}\boldsymbol{B}_{i}{}^{T}\left[\boldsymbol{D}^{+}\right]^{T}\boldsymbol{n}\boldsymbol{N}_{j}\mathrm{d}\mathrm{\Gamma}+\int_{\mathcal{S}_{e}}\boldsymbol{N}_{i}{}^{T}\alpha^{+}\boldsymbol{N}_{j}\mathrm{d\Gamma}\nonumber \\
 &  & +\int_{\mathcal{S}_{e}}\boldsymbol{N}_{i}{}^{T}\boldsymbol{n}^{T}\boldsymbol{D}^{-}\boldsymbol{B}_{j}\mathrm{d\Gamma}+\int_{\mathcal{S}_{e}}\boldsymbol{B}_{i}{}^{T}\left[\boldsymbol{D}^{-}\right]^{T}\boldsymbol{n}\boldsymbol{N}_{j}\mathrm{d}\mathrm{\Gamma}+\int_{\mathcal{S}_{e}}\boldsymbol{N}_{i}{}^{T}\alpha^{-}\boldsymbol{N}_{j}\mathrm{d\Gamma})\boldsymbol{u}_{j}\nonumber \\
 &  & +(\int_{\Omega_{e}}\boldsymbol{B}_{i}^{T}\boldsymbol{D}\boldsymbol{B}_{j}\boldsymbol{H}\mathrm{d\Omega}-\int_{\mathcal{S}_{e}}\boldsymbol{N}_{i}{}^{T}\boldsymbol{n}^{T}\boldsymbol{D}^{+}\boldsymbol{B}_{j}\mathrm{d\Gamma}-\int_{\mathcal{S}_{e}}\boldsymbol{B}_{i}{}^{T}\left[\boldsymbol{D}^{+}\right]^{T}\boldsymbol{n}\boldsymbol{N}_{j}\mathrm{d}\mathrm{\Gamma}+\int_{\mathcal{S}_{e}}\boldsymbol{N}_{i}{}^{T}\alpha^{+}\boldsymbol{N}_{j}\mathrm{d\Gamma}\nonumber \\
 &  & -\int_{\mathcal{S}_{e}}\boldsymbol{N}_{i}{}^{T}\boldsymbol{n}^{T}\boldsymbol{D}^{-}\boldsymbol{B}_{j}\mathrm{d\Gamma}-\int_{\mathcal{S}_{e}}\boldsymbol{B}_{i}{}^{T}\left[\boldsymbol{D}^{-}\right]^{T}\boldsymbol{n}\boldsymbol{N}_{j}\mathrm{d}\mathrm{\Gamma}-\int_{\mathcal{S}_{e}}\boldsymbol{N}_{i}{}^{T}\alpha^{-}\boldsymbol{N}_{j}\mathrm{d\Gamma})\boldsymbol{a}_{j}\nonumber \\
 & = & \int_{\Omega_{e}}\boldsymbol{N}_{i}^{T}\boldsymbol{f}\mathrm{d\Omega}+\int_{\Gamma_{\mbox{n}_{e}}}\boldsymbol{N}_{i}^{T}\boldsymbol{t}\mathrm{d\Omega}-\int_{\mathcal{S}_{e}}\boldsymbol{B}_{i}{}^{T}\left[\boldsymbol{D}^{+}\right]^{T}\boldsymbol{n}\boldsymbol{g}^{+}\mathrm{d}\mathrm{\Gamma}\\
 &  & +\int_{\mathcal{S}_{e}}\boldsymbol{B}_{i}{}^{T}\left[\boldsymbol{D}^{-}\right]^{T}\boldsymbol{n}\boldsymbol{g}^{-}\mathrm{d}\mathrm{\Gamma}+\int_{\mathcal{S}_{e}}\boldsymbol{N}_{i}{}^{T}\alpha^{+}\boldsymbol{g}^{+}\mathrm{d\Gamma}+\int_{\mathcal{S}_{e}}\boldsymbol{N}_{i}{}^{T}\alpha^{-}\boldsymbol{g}^{-}\mathrm{d\Gamma}\nonumber 
\end{eqnarray}
\begin{eqnarray}
 &  & (\int_{\Omega_{e}}\boldsymbol{H}\boldsymbol{B}_{i}{}^{T}\boldsymbol{D}\boldsymbol{B}_{j}\mathrm{d\Omega}-\int_{\mathcal{S}_{e}}\boldsymbol{N}_{i}{}^{T}\boldsymbol{n}^{T}\boldsymbol{D}^{+}\boldsymbol{B}_{j}\mathrm{d\Gamma}-\int_{\mathcal{S}_{e}}\boldsymbol{B}_{i}{}^{T}\left[\boldsymbol{D}^{+}\right]^{T}\boldsymbol{n}\boldsymbol{N}_{j}\mathrm{d}\mathrm{\Gamma}+\int_{\mathcal{S}_{e}}\boldsymbol{N}_{i}{}^{T}\alpha^{+}\boldsymbol{N}_{j}\mathrm{d\Gamma}\nonumber \\
 & - & \int_{\mathcal{S}_{e}}\boldsymbol{N}_{i}{}^{T}\boldsymbol{n}^{T}\boldsymbol{D}^{-}\boldsymbol{B}_{j}\mathrm{d\Gamma}-\int_{\mathcal{S}_{e}}\boldsymbol{B}_{i}{}^{T}\left[\boldsymbol{D}^{-}\right]^{T}\boldsymbol{n}\boldsymbol{N}_{j}\mathrm{d}\mathrm{\Gamma}-\int_{\mathcal{S}_{e}}\boldsymbol{N}_{i}{}^{T}\alpha^{-}\boldsymbol{N}_{j}\mathrm{d\Gamma})\boldsymbol{u}_{j}\nonumber \\
 & + & (\int_{\Omega_{e}}\boldsymbol{B}_{i}{}^{T}\boldsymbol{D}\boldsymbol{B}_{j}\mathrm{d\Omega}-\int_{\mathcal{S}_{e}}\boldsymbol{N}_{i}{}^{T}\boldsymbol{n}^{T}\boldsymbol{D}^{+}\boldsymbol{B}_{j}\mathrm{d\Gamma}-\int_{\mathcal{S}_{e}}\boldsymbol{B}_{i}{}^{T}\left[\boldsymbol{D}^{+}\right]^{T}\boldsymbol{n}\boldsymbol{N}_{j}\mathrm{d}\mathrm{\Gamma}+\int_{\mathcal{S}_{e}}\boldsymbol{N}_{i}{}^{T}\alpha^{+}\boldsymbol{N}_{j}\mathrm{d\Gamma}\nonumber \\
 & + & \int_{\mathcal{S}_{e}}\boldsymbol{N}_{i}{}^{T}\boldsymbol{n}^{T}\boldsymbol{D}^{-}\boldsymbol{B}_{j}\mathrm{d\Gamma}+\int_{\mathcal{S}_{e}}\boldsymbol{B}_{i}{}^{T}\left[\boldsymbol{D}^{-}\right]^{T}\boldsymbol{n}\boldsymbol{N}_{j}\mathrm{d}\mathrm{\Gamma}+\int_{\mathcal{S}_{e}}\boldsymbol{N}_{i}{}^{T}\alpha^{-}\boldsymbol{N}_{j}\mathrm{d\Gamma})\boldsymbol{a}_{j}\nonumber \\
= &  & \int_{\Omega_{e}}\boldsymbol{H}\boldsymbol{N}_{i}^{T}\boldsymbol{f}\mathrm{d\Omega}+\int_{\Gamma_{\mbox{n}_{e}}}\boldsymbol{H}\boldsymbol{N}_{i}^{T}\boldsymbol{t}\mathrm{d\Omega}-\int_{\mathcal{S}_{e}}\boldsymbol{B}_{i}{}^{T}\left[\boldsymbol{D}^{+}\right]^{T}\boldsymbol{n}\boldsymbol{g}^{+}\mathrm{d}\mathrm{\Gamma}\nonumber \\
 & - & \int_{\mathcal{S}_{e}}\boldsymbol{B}_{i}{}^{T}\left[\boldsymbol{D}^{-}\right]^{T}\boldsymbol{n}\boldsymbol{g}^{-}\mathrm{d}\mathrm{\Gamma}+\int_{\mathcal{S}_{e}}\boldsymbol{N}_{i}{}^{T}\alpha^{+}\boldsymbol{g}^{+}\mathrm{d\Gamma}-\int_{\mathcal{S}_{e}}\boldsymbol{N}_{i}{}^{T}\alpha^{-}\boldsymbol{g}^{-}\mathrm{d\Gamma}
\end{eqnarray}
\begin{eqnarray}
\left[\begin{array}{cc}
\sum_{e}\int_{\Omega_{e}}\boldsymbol{B}_{i}^{T}\boldsymbol{D}\boldsymbol{B}_{j}\mathrm{d\Omega} & \sum_{e}\int_{\Omega_{e}}\boldsymbol{B}_{i}^{T}\boldsymbol{D}\boldsymbol{B}_{j}\boldsymbol{H}\mathrm{d\Omega}\\
-\sum_{e}\int_{\mathcal{S}_{e}}\boldsymbol{N}_{i}{}^{T}\boldsymbol{n}^{T}\boldsymbol{D}^{+}\boldsymbol{B}_{j}\mathrm{d\Gamma} & -\sum_{e}\int_{\mathcal{S}_{e}}\boldsymbol{N}_{i}{}^{T}\boldsymbol{n}^{T}\boldsymbol{D}^{+}\boldsymbol{B}_{j}\mathrm{d\Gamma}\\
-\sum_{e}\int_{\mathcal{S}_{e}}\boldsymbol{B}_{i}{}^{T}\left[\boldsymbol{D}^{+}\right]^{T}\boldsymbol{n}\boldsymbol{N}_{j}\mathrm{d}\mathrm{\Gamma} & -\sum_{e}\int_{\mathcal{S}_{e}}\boldsymbol{B}_{i}{}^{T}\left[\boldsymbol{D}^{+}\right]^{T}\boldsymbol{n}\boldsymbol{N}_{j}\mathrm{d}\mathrm{\Gamma}\\
+\sum_{e}\int_{\mathcal{S}_{e}}\boldsymbol{N}_{i}{}^{T}\alpha^{+}\boldsymbol{N}_{j}\mathrm{d\Gamma} & +\sum_{e}\int_{\mathcal{S}_{e}}\boldsymbol{N}_{i}{}^{T}\alpha^{+}\boldsymbol{N}_{j}\mathrm{d\Gamma}\\
+\sum_{e}\int_{\mathcal{S}_{e}}\boldsymbol{N}_{i}{}^{T}\boldsymbol{n}^{T}\boldsymbol{D}^{-}\boldsymbol{B}_{j}\mathrm{d\Gamma} & -\sum_{e}\int_{\mathcal{S}_{e}}\boldsymbol{N}_{i}{}^{T}\boldsymbol{n}^{T}\boldsymbol{D}^{-}\boldsymbol{B}_{j}\mathrm{d\Gamma}\\
+\sum_{e}\int_{\mathcal{S}_{e}}\boldsymbol{B}_{i}{}^{T}\left[\boldsymbol{D}^{-}\right]^{T}\boldsymbol{n}\boldsymbol{N}_{j}\mathrm{d}\mathrm{\Gamma} & -\sum_{e}\int_{\mathcal{S}_{e}}\boldsymbol{B}_{i}{}^{T}\left[\boldsymbol{D}^{-}\right]^{T}\boldsymbol{n}\boldsymbol{N}_{j}\mathrm{d}\mathrm{\Gamma}\\
+\sum_{e}\int_{\mathcal{S}_{e}}\boldsymbol{N}_{i}{}^{T}\alpha^{-}\boldsymbol{N}_{j}\mathrm{d\Gamma} & -\sum_{e}\int_{\mathcal{S}_{e}}\boldsymbol{N}_{i}{}^{T}\alpha^{-}\boldsymbol{N}_{j}\mathrm{d\Gamma}\\
\\
\sum_{e}\int_{\Omega_{e}}\boldsymbol{H}\boldsymbol{B}_{i}{}^{T}\boldsymbol{D}\boldsymbol{B}_{j}\mathrm{d\Omega} & \sum_{e}\int_{\Omega_{e}}\boldsymbol{B}_{i}{}^{T}\boldsymbol{D}\boldsymbol{B}_{j}\mathrm{d\Omega}\\
-\sum_{e}\int_{\mathcal{S}_{e}}\boldsymbol{N}_{i}{}^{T}\boldsymbol{n}^{T}\boldsymbol{D}^{+}\boldsymbol{B}_{j}\mathrm{d\Gamma} & -\sum_{e}\int_{\mathcal{S}_{e}}\boldsymbol{N}_{i}{}^{T}\boldsymbol{n}^{T}\boldsymbol{D}^{+}\boldsymbol{B}_{j}\mathrm{d\Gamma}\\
-\sum_{e}\int_{\mathcal{S}_{e}}\boldsymbol{B}_{i}{}^{T}\left[\boldsymbol{D}^{+}\right]^{T}\boldsymbol{n}\boldsymbol{N}_{j}\mathrm{d}\mathrm{\Gamma} & -\sum_{e}\int_{\mathcal{S}_{e}}\boldsymbol{B}_{i}{}^{T}\left[\boldsymbol{D}^{+}\right]^{T}\boldsymbol{n}\boldsymbol{N}_{j}\mathrm{d}\mathrm{\Gamma}\\
+\sum_{e}\int_{\mathcal{S}_{e}}\boldsymbol{N}_{i}{}^{T}\alpha^{+}\boldsymbol{N}_{j}\mathrm{d\Gamma} & +\sum_{e}\int_{\mathcal{S}_{e}}\boldsymbol{N}_{i}{}^{T}\alpha^{+}\boldsymbol{N}_{j}\mathrm{d\Gamma}\\
-\sum_{e}\int_{\mathcal{S}_{e}}\boldsymbol{N}_{i}{}^{T}\boldsymbol{n}^{T}\boldsymbol{D}^{-}\boldsymbol{B}_{j}\mathrm{d\Gamma} & +\sum_{e}\int_{\mathcal{S}_{e}}\boldsymbol{N}_{i}{}^{T}\boldsymbol{n}^{T}\boldsymbol{D}^{-}\boldsymbol{B}_{j}\mathrm{d\Gamma}\\
-\sum_{e}\int_{\mathcal{S}_{e}}\boldsymbol{B}_{i}{}^{T}\left[\boldsymbol{D}^{-}\right]^{T}\boldsymbol{n}\boldsymbol{N}_{j}\mathrm{d}\mathrm{\Gamma} & +\sum_{e}\int_{\mathcal{S}_{e}}\boldsymbol{B}_{i}{}^{T}\left[\boldsymbol{D}^{-}\right]^{T}\boldsymbol{n}\boldsymbol{N}_{j}\mathrm{d}\mathrm{\Gamma}\\
-\sum_{e}\int_{\mathcal{S}_{e}}\boldsymbol{N}_{i}{}^{T}\alpha^{-}\boldsymbol{N}_{j}\mathrm{d\Gamma} & +\sum_{e}\int_{\mathcal{S}_{e}}\boldsymbol{N}_{i}{}^{T}\alpha^{-}\boldsymbol{N}_{j}\mathrm{d\Gamma})\boldsymbol{a}_{j}
\end{array}\right]\left\{ \begin{array}{c}
\boldsymbol{u}_{j}\\
\boldsymbol{a}_{j}
\end{array}\right\}  & = & \left\{ \begin{array}{c}
\sum_{e}\int_{\Omega_{e}}\boldsymbol{N}_{i}^{T}\boldsymbol{f}\mathrm{d\Omega}\\
+\sum_{e}\int_{\Gamma_{\mbox{n}_{e}}}\boldsymbol{N}_{i}^{T}\boldsymbol{t}\mathrm{d\Omega}\\
-\sum_{e}\int_{\mathcal{S}_{e}}\boldsymbol{B}_{i}{}^{T}\left[\boldsymbol{D}^{+}\right]^{T}\boldsymbol{n}\boldsymbol{g}^{+}\mathrm{d}\mathrm{\Gamma}\\
+\sum_{e}\int_{\mathcal{S}_{e}}\boldsymbol{B}_{i}{}^{T}\left[\boldsymbol{D}^{-}\right]^{T}\boldsymbol{n}\boldsymbol{g}^{-}\mathrm{d}\mathrm{\Gamma}\\
+\sum_{e}\int_{\mathcal{S}_{e}}\boldsymbol{N}_{i}{}^{T}\alpha^{+}\boldsymbol{g}^{+}\mathrm{d\Gamma}\\
+\sum_{e}\int_{\mathcal{S}_{e}}\boldsymbol{N}_{i}{}^{T}\alpha^{-}\boldsymbol{g}^{-}\mathrm{d\Gamma}\\
\\
\sum_{e}\int_{\Omega_{e}}\boldsymbol{H}\boldsymbol{N}_{i}^{T}\boldsymbol{f}\mathrm{d\Omega}\\
+\sum_{e}\int_{\Gamma_{\mbox{n}_{e}}}\boldsymbol{H}\boldsymbol{N}_{i}^{T}\boldsymbol{t}\mathrm{d\Omega}\\
-\sum_{e}\int_{\mathcal{S}_{e}}\boldsymbol{B}_{i}{}^{T}\left[\boldsymbol{D}^{+}\right]^{T}\boldsymbol{n}\boldsymbol{g}^{+}\mathrm{d}\mathrm{\Gamma}\\
-\sum_{e}\int_{\mathcal{S}_{e}}\boldsymbol{B}_{i}{}^{T}\left[\boldsymbol{D}^{-}\right]^{T}\boldsymbol{n}\boldsymbol{g}^{-}\mathrm{d}\mathrm{\Gamma}\\
+\sum_{e}\int_{\mathcal{S}_{e}}\boldsymbol{N}_{i}{}^{T}\alpha^{+}\boldsymbol{g}^{+}\mathrm{d\Gamma}\\
-\sum_{e}\int_{\mathcal{S}_{e}}\boldsymbol{N}_{i}{}^{T}\alpha^{-}\boldsymbol{g}^{-}\mathrm{d\Gamma}
\end{array}\right\} \nonumber \\
\end{eqnarray}
\begin{equation}
\left[\begin{array}{cc}
\mathbf{K}_{b}-\mathbf{K}_{n}^{+}-\left[\mathbf{K}_{n}^{+}\right]^{T}+\mathbf{K}_{s}^{+} & \mathbf{K}_{b}\boldsymbol{H}-\mathbf{K}_{n}^{+}-\left[\mathbf{K}_{n}^{+}\right]^{T}+\mathbf{K}_{s}^{+}\\
+\mathbf{K}_{n}^{-}+\left[\mathbf{K}_{n}^{-}\right]^{T}+\mathbf{K}_{s}^{-} & -\mathbf{K}_{n}^{-}-\left[\mathbf{K}_{n}^{-}\right]^{T}-\mathbf{K}_{s}^{-}\\
\\
\boldsymbol{H}\mathbf{K}_{b}-\mathbf{K}_{n}^{+}-\left[\mathbf{K}_{n}^{+}\right]^{T}+\mathbf{K}_{s}^{+} & \mathbf{K}_{b}-\mathbf{K}_{n}^{+}-\left[\mathbf{K}_{n}^{+}\right]^{T}+\mathbf{K}_{s}^{+}\\
-\mathbf{K}_{n}^{-}-\left[\mathbf{K}_{n}^{-}\right]^{T}-\mathbf{K}_{s}^{-} & +\mathbf{K}_{n}^{-}+\left[\mathbf{K}_{n}^{-}\right]^{T}+\mathbf{K}_{s}^{-}
\end{array}\right]\left\{ \begin{array}{c}
\boldsymbol{u}_{j}\\
\boldsymbol{a}_{j}
\end{array}\right\} =\left\{ \begin{array}{c}
\mathbf{f}_{b}+\mathbf{f}_{h}+\mathbf{f}_{s}^{+}+\mathbf{f}_{s}^{-}\\
-\mathbf{f}_{n}^{+}+\mathbf{f}_{n}^{-}\\
\\
\boldsymbol{H}(\mathbf{f}_{b}+\mathbf{f}_{h})+\mathbf{f}_{s}^{+}-\mathbf{f}_{s}^{-}\\
-\mathbf{f}_{n}^{+}-\mathbf{f}_{n}^{-}
\end{array}\right\} \label{eq:-28}
\end{equation}

Here, like in (\ref{eq:-19}) after assembly, $\mathbf{K}_{b}$ is
the bulk stiffness term, $\mathbf{K}_{n}$ is Nitsche's contribution
to stiffness and $\mathbf{K}_{s}^{\pm}$ is the stability term associated
with the formulation. $\mathbf{f}_{b}$ is the bulk force term, $\mathbf{f}_{h}$
is Neumann's contribution to force term, $\mathbf{f}_{n}^{\pm}$ is
Nitsche's contribution to the force term and $\mathbf{f}_{s}^{\pm}$
is the stability parameter associated with the force component. This
can be compared to the tabular formulation of the terms associated
with the weak formulation. 

To evaluate flux, we use the domain integral. Considering a set of
nodes $D$, whose supports intersect the interface $\mathcal{S}_{e}$,
the approximation $\mbox{j}_{h}$ to the interfacial flux is written:
\[
\mbox{j}_{h}=\sum_{I\in D}N_{I}(\mathbf{x})\mbox{j}_{I}\mbox{,\,\,\,\,\textbf{x\ensuremath{\in\mathcal{L}_{e}}}}
\]
where $\mbox{j}_{I}$ are to be determined. This can be implemented
as

\begin{eqnarray}
\int_{\mathcal{S}_{e}}\boldsymbol{v}_{h}^{T}\bar{\mathrm{j}}_{h}\mbox{d}\Gamma & = & \int_{\mathcal{B}_{e}}\boldsymbol{v}_{h}^{T}\boldsymbol{f}\mbox{d}\Omega+\int_{\Gamma_{\mbox{n}_{e}}}\boldsymbol{v}_{h}^{T}\boldsymbol{h}\mbox{d}\Omega-\int_{\mathcal{B}_{e}}\varepsilon^{T}(\boldsymbol{v}_{h})\sigma(\boldsymbol{u}_{h})\mbox{d}\Omega
\end{eqnarray}
\begin{eqnarray}
\int_{\mathcal{S}_{e}}(\boldsymbol{v}_{i}^{T}+\boldsymbol{b}_{i}^{T}H)\boldsymbol{N}_{i}{}^{T}\bar{\mathrm{j}}_{h}\mbox{d}\Gamma & = & \int_{\mathcal{B}_{e}}(\boldsymbol{v}_{i}^{T}+\boldsymbol{b}_{i}^{T}H)\boldsymbol{N}_{i}{}^{T}\boldsymbol{f}\mbox{d}\Omega+\nonumber \\
 &  & \int_{\Gamma_{\mbox{n}_{e}}}(\boldsymbol{v}_{i}^{T}+\boldsymbol{b}_{i}^{T}H)\boldsymbol{N}_{i}{}^{T}\boldsymbol{h}\mbox{d}\Omega-\int_{\mathcal{B}_{e}}(\boldsymbol{v}_{i}^{T}+\boldsymbol{b}_{i}^{T}H)\boldsymbol{B}_{i}{}^{T}\sigma(\boldsymbol{u}_{h})\mbox{d}\Omega
\end{eqnarray}
Knowing the arbitrariness of $\boldsymbol{v}_{i}$ and $\boldsymbol{b}_{i}$,
we can write,
\begin{eqnarray}
\int_{\mathcal{S}_{e}}\boldsymbol{N}_{i}{}^{T}\bar{\mathrm{j}}_{h}\mbox{d}\Gamma & = & \int_{\mathcal{B}_{e}}\boldsymbol{N}_{i}{}^{T}\boldsymbol{f}\mbox{d}\Omega+\int_{\Gamma_{\mbox{n}_{e}}}\boldsymbol{N}_{i}{}^{T}\boldsymbol{h}\mbox{d}\Omega-\int_{\mathcal{B}_{e}}\boldsymbol{B}_{i}{}^{T}\sigma(\boldsymbol{u}_{h})\mbox{d}\Omega\label{eq:-30}
\end{eqnarray}
From the constitutive law we can write,
\begin{eqnarray}
\int_{\mathcal{S}_{e}}\boldsymbol{N}_{i}{}^{T}\boldsymbol{N}_{j}\bar{\mathrm{j}}_{j}\mbox{d}\Gamma & = & \int_{\mathcal{B}_{e}}\boldsymbol{N}_{i}{}^{T}\boldsymbol{f}\mbox{d}\Omega+\int_{\Gamma_{\mbox{n}_{e}}}\boldsymbol{N}_{i}{}^{T}\boldsymbol{h}\mbox{d}\Gamma-\int_{\mathcal{B}_{e}}\boldsymbol{B}_{i}{}^{T}\boldsymbol{D}\boldsymbol{B}_{j}(\boldsymbol{u}_{j}+H\boldsymbol{a}_{j})\mbox{d}\Omega
\end{eqnarray}
\begin{equation}
\mathbf{M}_{d}\bar{\mathrm{j}}_{j}=\mathbf{f}_{b}+\mathbf{f}_{h}-\mathbf{K}_{b}\boldsymbol{u}_{j}-\mathbf{K}_{b}\boldsymbol{H}\boldsymbol{a}_{j}
\end{equation}
where $\mathbf{M}_{d}$ is the mass matrix over the interface. Knowing
the values of $\boldsymbol{u}_{j}$ and $\boldsymbol{a}_{j}$ from
the previous formulation we can approximate the value of the jump
in flux. We can now combine the two, the displacement system of equations
(\ref{eq:-28}) and the jump in flux, to get the following system:
\begin{equation}
\left[\begin{array}{ccc}
\mathbf{K}_{b}-\mathbf{K}_{n}^{+}-\left[\mathbf{K}_{n}^{+}\right]^{T}+\mathbf{K}_{s}^{+} & \mathbf{K}_{b}\boldsymbol{H}-\mathbf{K}_{n}^{+}-\left[\mathbf{K}_{n}^{+}\right]^{T}+\mathbf{K}_{s}^{+} & \mathbf{0}\\
+\mathbf{K}_{n}^{-}+\left[\mathbf{K}_{n}^{-}\right]^{T}+\mathbf{K}_{s}^{-} & -\mathbf{K}_{n}^{-}-\left[\mathbf{K}_{n}^{-}\right]^{T}-\mathbf{K}_{s}^{-}\\
\\
\boldsymbol{H}\mathbf{K}_{b}-\mathbf{K}_{n}^{+}-\left[\mathbf{K}_{n}^{+}\right]^{T}+\mathbf{K}_{s}^{+} & \mathbf{K}_{b}-\mathbf{K}_{n}^{+}-\left[\mathbf{K}_{n}^{+}\right]^{T}+\mathbf{K}_{s}^{+} & \mathbf{0}\\
-\mathbf{K}_{n}^{-}-\left[\mathbf{K}_{n}^{-}\right]^{T}-\mathbf{K}_{s}^{-} & +\mathbf{K}_{n}^{-}+\left[\mathbf{K}_{n}^{-}\right]^{T}+\mathbf{K}_{s}^{-}\\
\\
\mathbf{K}_{b} & \mathbf{K}_{b}\boldsymbol{H} & \mathbf{M}_{d}
\end{array}\right]\left\{ \begin{array}{c}
\boldsymbol{u}_{j}\\
\boldsymbol{a}_{j}\\
\bar{\mathrm{j}}_{j}
\end{array}\right\} =\left\{ \begin{array}{c}
\mathbf{f}_{b}+\mathbf{f}_{h}+\mathbf{f}_{s}^{+}+\mathbf{f}_{s}^{-}\\
-\mathbf{f}_{n}^{+}+\mathbf{f}_{n}^{-}\\
\\
\boldsymbol{H}(\mathbf{f}_{b}+\mathbf{f}_{h})+\mathbf{f}_{s}^{+}-\mathbf{f}_{s}^{-}\\
-\mathbf{f}_{n}^{+}-\mathbf{f}_{n}^{-}\\
\\
\mathbf{f}_{b}+\mathbf{f}_{h}
\end{array}\right\} 
\end{equation}
It can be noted that since the two systems are not coupled, the displacement
and jump in flux can be independently solved. 

One can note that when the interface is between two domains of the
same material, $\boldsymbol{D}^{+}=\boldsymbol{D}^{-}$, $\mathbf{K}_{n}^{+}=\mathbf{K}_{n}^{-}=\mathbf{K}_{n}$
resulting in
\begin{equation}
\left[\begin{array}{cc}
\mathbf{K}_{b}+\mathbf{K}_{s}^{+}+\mathbf{K}_{s}^{-} & \mathbf{K}_{b}\boldsymbol{H}-2(\mathbf{K}_{n}+\mathbf{K}_{n}^{T})\\
 & +\mathbf{K}_{s}^{+}-\mathbf{K}_{s}^{-}\\
\\
\boldsymbol{H}\mathbf{K}_{b}-2(\mathbf{K}_{n}+\mathbf{K}_{n}^{T})\\
+\mathbf{K}_{s}^{+}-\mathbf{K}_{s}^{-} & \mathbf{K}_{b}+\mathbf{K}_{s}^{+}+\mathbf{K}_{s}^{-}
\end{array}\right]\left\{ \begin{array}{c}
\boldsymbol{u}_{j}\\
\boldsymbol{a}_{j}
\end{array}\right\} =\left\{ \begin{array}{c}
\mathbf{f}_{b}+\mathbf{f}_{h}+\mathbf{f}_{s}^{+}+\mathbf{f}_{s}^{-}\\
-\mathbf{f}_{n}^{+}+\mathbf{f}_{n}^{-}\\
\\
\boldsymbol{H}(\mathbf{f}_{b}+\mathbf{f}_{h})+\mathbf{f}_{s}^{+}-\mathbf{f}_{s}^{-}\\
-\mathbf{f}_{n}^{+}-\mathbf{f}_{n}^{-}
\end{array}\right\} 
\end{equation}

To make a stronger comparison between the two conditions, Dirichlet
and jump, we can also impose $\alpha^{+}=\alpha^{-}=\alpha$ as the
stabilization parameter. This gives us $\mathbf{K}_{s}^{+}=\mathbf{K}_{s}^{-}=\mathbf{K}_{s}$.
The Dirichlet formulation now becomes
\begin{equation}
\left[\begin{array}{cc}
\mathbf{K}_{b}+2\mathbf{K}_{s} & \mathbf{K}_{b}\boldsymbol{H}-2(\mathbf{K}_{n}+\mathbf{K}_{n}^{T})\\
\\
\boldsymbol{H}\mathbf{K}_{b}-2(\mathbf{K}_{n}+\mathbf{K}_{n}^{T}) & \mathbf{K}_{b}+2\mathbf{K}_{s}
\end{array}\right]\left\{ \begin{array}{c}
\boldsymbol{u}_{j}\\
\boldsymbol{a}_{j}
\end{array}\right\} =\left\{ \begin{array}{c}
\mathbf{f}_{b}+\mathbf{f}_{h}+\mathbf{f}_{s}^{+}+\mathbf{f}_{s}^{-}\\
-\mathbf{f}_{n}^{+}+\mathbf{f}_{n}^{-}\\
\\
\boldsymbol{H}(\mathbf{f}_{b}+\mathbf{f}_{h})+\mathbf{f}_{s}^{+}-\mathbf{f}_{s}^{-}\\
-\mathbf{f}_{n}^{+}-\mathbf{f}_{n}^{-}
\end{array}\right\} \label{eq:-21}
\end{equation}

\subsection{Discretization with a shifted basis enrichment}

Consider a modification of the enrichment type known as shifted basis.\cite{key-9,key-7}
We can write
\begin{eqnarray}
\boldsymbol{u}_{h}(\mathbf{x}) & = & \sum_{i\in I}\boldsymbol{u}_{i}N_{i}(\mathbf{x})+\sum_{i\in L}\boldsymbol{a}_{i}N_{i}(\mathbf{x})(\tilde{H}_{i}(\mathbf{x}))\,\,\,\,\,\,\,\,\mathbb{U}_{h}\subset\mathbb{U}\\
\boldsymbol{v}_{h}(\mathbf{x}) & = & \sum_{i\in I}\boldsymbol{v}_{i}N_{i}(\mathbf{x})+\sum_{i\in L}\boldsymbol{b}_{i}N_{i}(\mathbf{x})(\tilde{H}_{i}(\mathbf{x}))\,\,\,\,\,\,\,\,\mathbb{U}_{0h}\subset\mathbb{U}_{0}\label{eq:-12-1}
\end{eqnarray}
where $\tilde{H}(\mathbf{x})$ is the shifted basis enrichment function:
\[
\tilde{H}_{i}(\mathbf{x})=H(\mathbf{x})-H(\mathbf{x}_{i})=H(\mathbf{x})-H{}_{i}
\]
 with $H(\mathbf{x})$, the global enrichment function and $H(\mathbf{x}_{i})$,
the local enrichment function. With this modification, and keeping
all the necessary spaces and conditions the same, we can now say that:
\[
\boldsymbol{u}_{h}^{+}=\sum_{i\in I}\boldsymbol{u}_{i}N_{i}(\mathbf{x})+\sum_{i\in L}\boldsymbol{a}_{i}\tilde{H}_{i}^{+}N_{i}(\mathbf{x})
\]
\[
\boldsymbol{u}_{h}^{-}=\sum_{i\in I}\boldsymbol{u}_{i}N_{i}(\mathbf{x})+\sum_{i\in L}\boldsymbol{a}_{i}\tilde{H}_{i}^{-}N_{i}(\mathbf{x})
\]
with $H_{i}=H(\mathbf{x}_{i})$ 
\begin{eqnarray*}
\tilde{H}_{i}^{+} & = & 1-H_{i}\\
\tilde{H}_{i}^{-} & = & -1-H_{i}
\end{eqnarray*}
and thus:
\begin{eqnarray*}
<\boldsymbol{u}_{h}> & = & \frac{1}{2}(\boldsymbol{u}_{h}^{+}+\boldsymbol{u}_{h}^{-})\\
 & = & \frac{1}{2}(\boldsymbol{u}_{i}+\tilde{H}_{i}^{+}\boldsymbol{a}_{i}+\boldsymbol{u}_{i}+\tilde{H}_{i}^{-}\boldsymbol{a}_{i})\\
 & = & \boldsymbol{u}_{i}+\frac{1}{2}(1-H_{i}+-1-H_{i})\boldsymbol{a}_{i}\\
 & = & \boldsymbol{u}_{i}+\frac{1}{2}(\tilde{H}_{i}^{+}+\tilde{H}_{i}^{-})\boldsymbol{a}_{i}\\
 & = & \boldsymbol{u}_{i}-H_{i}\boldsymbol{a}_{i}
\end{eqnarray*}
\begin{eqnarray*}
<\sigma(\boldsymbol{u}_{h})>.\mathbf{n} & = & \frac{1}{2}(\sigma(\boldsymbol{u}_{h}^{+}).\mathbf{n}+\sigma(\boldsymbol{u}_{h}^{-}).\mathbf{n})\\
 & = & \frac{1}{2}(\sigma(\boldsymbol{u}_{h}^{+}).\mathbf{n}+\sigma(\boldsymbol{u}_{h}^{-}).\mathbf{n})\\
 & = & \frac{1}{2}(\boldsymbol{n}^{T}\boldsymbol{D}^{+}\boldsymbol{B}_{i}\boldsymbol{u}_{i}+\boldsymbol{n}^{T}\boldsymbol{D}^{+}\boldsymbol{B}_{i}\tilde{\boldsymbol{H}}_{i}^{+}\boldsymbol{a}_{i}+\boldsymbol{n}^{T}\boldsymbol{D}^{-}\boldsymbol{B}_{i}\boldsymbol{u}_{i}+\boldsymbol{n}^{T}\boldsymbol{D}^{-}\boldsymbol{B}_{i}\tilde{\boldsymbol{H}}_{i}^{-}\boldsymbol{a}_{i})\\
 & = & \frac{1}{2}(\boldsymbol{n}^{T}\boldsymbol{D}^{+}\boldsymbol{B}_{i}+\boldsymbol{n}^{T}\boldsymbol{D}^{-}\boldsymbol{B}_{i})\boldsymbol{u}_{i}+\frac{1}{2}(\boldsymbol{n}^{T}\boldsymbol{D}^{+}\boldsymbol{B}_{i}\tilde{\boldsymbol{H}}_{i}^{+}+\boldsymbol{n}^{T}\boldsymbol{D}^{-}\boldsymbol{B}_{i}\tilde{\boldsymbol{H}}_{i}^{-})\boldsymbol{a}_{i}
\end{eqnarray*}
and
\begin{eqnarray*}
[[\boldsymbol{u}_{h}]] & = & \boldsymbol{u}_{h}^{+}-\boldsymbol{u}_{h}^{-}\\
 & = & \cancel{\boldsymbol{u}_{i}}+\tilde{\boldsymbol{H}}_{i}^{+}\boldsymbol{a}_{i}-\cancel{\boldsymbol{u}_{i}}-\tilde{\boldsymbol{H}}_{i}^{-}\boldsymbol{a}_{i}\\
 & = & (1-\cancel{H}_{i})\boldsymbol{a}_{i}+(1+\cancel{H}_{i})\boldsymbol{a}_{i}\\
 & = & 2\boldsymbol{a}_{i}
\end{eqnarray*}
We will continue our discussions using the notation of $\tilde{\boldsymbol{H}}$.

\subsubsection{Jump Condition}

Considering the jump in displacement condition from (\ref{eq:-24})
and implementing the above gives us 
\begin{eqnarray}
 &  & \int_{\Omega_{e}}\left(\boldsymbol{v}_{i}^{T}\boldsymbol{B}_{i}^{T}+\tilde{\boldsymbol{H}}_{i}^{T}\boldsymbol{b}_{i}^{T}\boldsymbol{B}_{i}{}^{T}\right)\sigma(\boldsymbol{u}_{h})\mathrm{d\Omega}-2\int_{\mathcal{S}_{e}}\boldsymbol{b}_{i}^{T}\boldsymbol{N}_{i}{}^{T}<\sigma(\boldsymbol{u}_{h})>.\mathbf{n}\mathrm{d\Gamma}\nonumber \\
 &  & -\frac{1}{2}\int_{\mathcal{S}_{e}}\boldsymbol{v}_{i}^{T}\left(\boldsymbol{B}_{i}^{T}\left[\boldsymbol{D}^{+}\right]^{T}\boldsymbol{n}+\boldsymbol{B}_{i}^{T}\left[\boldsymbol{D}^{-}\right]^{T}\boldsymbol{n}\right)[[\boldsymbol{u}_{h}]]\mathrm{d\Gamma}\nonumber \\
 &  & -\frac{1}{2}\int_{\mathcal{S}_{e}}\boldsymbol{b}_{i}^{T}\left(\boldsymbol{B}_{i}^{T}\left[\boldsymbol{D}^{+}\right]^{T}\boldsymbol{n}\tilde{\boldsymbol{H}}_{i}^{+}+\boldsymbol{B}_{i}^{T}\left[\boldsymbol{D}^{-}\right]^{T}\boldsymbol{n}\tilde{\boldsymbol{H}}_{i}^{-}\right)[[\boldsymbol{u}_{h}]]\mathrm{d\Gamma}+2\int_{\mathcal{S}_{e}}\boldsymbol{b}_{i}^{T}\boldsymbol{N}_{i}{}^{T}\alpha[[\boldsymbol{u}_{h}]]\mathrm{d\Gamma}\nonumber \\
 & = & \int_{\Omega_{e}}\left(\boldsymbol{v}_{i}^{T}\boldsymbol{N}_{i}^{T}+\tilde{\boldsymbol{H}}_{i}^{T}\boldsymbol{b}_{i}^{T}\boldsymbol{N}_{i}^{T}\right)\boldsymbol{f}\,\mathrm{d\Omega}+\int_{\Gamma_{\mbox{n}_{e}}}\left(\boldsymbol{v}_{i}^{T}\boldsymbol{N}_{i}^{T}+\tilde{\boldsymbol{H}}_{i}^{T}\boldsymbol{b}_{i}^{T}\boldsymbol{N}_{i}{}^{T}\right)\boldsymbol{t}\mathrm{d\Gamma}+2\int_{\mathcal{S}_{e}}\boldsymbol{b}_{i}^{T}\boldsymbol{N}_{i}{}^{T}\alpha\,\bar{\mbox{i}}\mathrm{\,d}\mathrm{\Gamma}\nonumber \\
 &  & -\frac{1}{2}\int_{\mathcal{S}_{e}}\boldsymbol{v}_{i}^{T}\left(\boldsymbol{B}_{i}^{T}\left[\boldsymbol{D}^{+}\right]^{T}\boldsymbol{n}+\boldsymbol{B}_{i}^{T}\left[\boldsymbol{D}^{-}\right]^{T}\boldsymbol{n}\right)\,\bar{\mbox{i}}\,\mathrm{d}\mathrm{\Gamma}-\int_{\mathcal{S}_{e}}\boldsymbol{H}_{i}^{T}\boldsymbol{b}_{i}^{T}\boldsymbol{N}_{i}{}^{T}\,\bar{\mbox{j}}\,\mathrm{d}\mathrm{\Gamma}\nonumber \\
 &  & -\frac{1}{2}\int_{\mathcal{S}_{e}}\boldsymbol{b}_{i}^{T}\left(\left[\tilde{\boldsymbol{H}}_{i}^{+}\right]^{T}\boldsymbol{B}_{i}^{T}\left[\boldsymbol{D}^{+}\right]^{T}\boldsymbol{n}+\left[\tilde{\boldsymbol{H}}_{i}^{-}\right]^{T}\boldsymbol{B}_{i}^{T}\left[\boldsymbol{D}^{-}\right]^{T}\boldsymbol{n}\right)\,\bar{\mbox{i}}\,\mathrm{d}\mathrm{\Gamma}+\int_{\mathcal{S}_{e}}\boldsymbol{v}_{i}^{T}\boldsymbol{N}_{i}{}^{T}\,\bar{\mbox{j}}\,\mathrm{d}\mathrm{\Gamma}\label{eq:-29}
\end{eqnarray}
Knowing the arbitrariness of $\boldsymbol{v}_{i}^{T}\mbox{ and }\boldsymbol{b}_{i}^{T}$
we can write (\ref{eq:-29}) in a two equation form, 
\begin{eqnarray}
 &  & \int_{\Omega_{e}}\boldsymbol{B}_{i}^{T}\sigma(\boldsymbol{u}_{h})\mathrm{d\Omega}-\frac{1}{2}\int_{\mathcal{S}_{e}}\left(\boldsymbol{B}_{i}^{T}\left[\boldsymbol{D}^{+}\right]^{T}\boldsymbol{n}+\boldsymbol{B}_{i}^{T}\left[\boldsymbol{D}^{-}\right]^{T}\boldsymbol{n}\right)[[\boldsymbol{u}_{h}]]\mathrm{d\Gamma}\nonumber \\
 & = & \int_{\Omega_{e}}\boldsymbol{N}_{i}^{T}\boldsymbol{f}\mathrm{d\Omega}+\int_{\Gamma_{\mbox{n}_{e}}}\boldsymbol{N}_{i}^{T}\boldsymbol{t}\mathrm{d\Gamma}-\frac{1}{2}\int_{\mathcal{S}_{e}}\left(\boldsymbol{B}_{i}^{T}\left[\boldsymbol{D}^{+}\right]^{T}\boldsymbol{n}+\boldsymbol{B}_{i}^{T}\left[\boldsymbol{D}^{-}\right]^{T}\boldsymbol{n}\right)\bar{\mbox{i}}\mathrm{d}\mathrm{\Gamma}+\int_{\mathcal{S}_{e}}\boldsymbol{N}_{i}{}^{T}\bar{\mbox{j}}\mathrm{d}\mathrm{\Gamma}
\end{eqnarray}
\begin{eqnarray}
 &  & \int_{\Omega_{e}}\tilde{\boldsymbol{H}}_{i}^{T}\boldsymbol{B}_{i}{}^{T}\sigma(\boldsymbol{u}_{h})\mathrm{d\Omega}-2\int_{\mathcal{S}_{e}}\boldsymbol{N}_{i}{}^{T}<\sigma(\boldsymbol{u}_{h}).n>\mathrm{d\Gamma}\nonumber \\
 &  & -\frac{1}{2}\int_{\mathcal{S}_{e}}\left(\left[\tilde{\boldsymbol{H}}_{i}^{+}\right]^{T}\boldsymbol{B}_{i}^{T}\left[\boldsymbol{D}^{+}\right]^{T}\boldsymbol{n}+\left[\tilde{\boldsymbol{H}}_{i}^{-}\right]^{T}\boldsymbol{B}_{i}^{T}\left[\boldsymbol{D}^{-}\right]^{T}\boldsymbol{n}\right)[[\boldsymbol{u}_{h}]]\mathrm{d\Gamma}+2\int_{\mathcal{S}_{e}}\boldsymbol{N}_{i}{}^{T}\alpha[[\boldsymbol{u}_{h}]]\mathrm{d\Gamma}\nonumber \\
 & = & \int_{\Omega_{e}}\tilde{\boldsymbol{H}}_{i}^{T}\boldsymbol{N}_{i}^{T}\boldsymbol{f}\mathrm{d\Omega}+\int_{\Gamma_{\mbox{n}_{e}}}\tilde{\boldsymbol{H}}_{i}^{T}\boldsymbol{N}_{i}{}^{T}\boldsymbol{t}\mathrm{d\Gamma}\nonumber \\
 &  & -\frac{1}{2}\int_{\mathcal{S}_{e}}\left(\left[\tilde{\boldsymbol{H}}_{i}^{+}\right]^{T}\boldsymbol{B}_{i}^{T}\left[\boldsymbol{D}^{+}\right]^{T}\boldsymbol{n}+\left[\tilde{\boldsymbol{H}}_{i}^{-}\right]^{T}\boldsymbol{B}_{i}^{T}\left[\boldsymbol{D}^{-}\right]^{T}\boldsymbol{n}\right)\,\bar{\mbox{i}}\,\mathrm{d}\mathrm{\Gamma}+2\int_{\mathcal{S}_{e}}\boldsymbol{N}_{i}{}^{T}\alpha\,\bar{\mbox{i}}\,\mathrm{d}\mathrm{\Gamma}\nonumber \\
 &  & +\frac{1}{2}\int_{\mathcal{S}_{e}}\boldsymbol{H}_{i}^{T}\boldsymbol{N}_{i}{}^{T}\,\bar{\mbox{j}}\,\mathrm{d}\mathrm{\Gamma}
\end{eqnarray}
Applying the constitutive law,
\begin{eqnarray}
 &  & \int_{\Omega_{e}}\boldsymbol{B}_{i}^{T}\boldsymbol{D}\boldsymbol{B}_{j}\left(\boldsymbol{u}_{j}+\tilde{\boldsymbol{H}}_{j}\boldsymbol{a}_{j}\right)\mathrm{d\Omega}-\int_{\mathcal{S}_{e}}\left(\boldsymbol{B}_{i}^{T}\left[\boldsymbol{D}^{+}\right]^{T}\boldsymbol{n}+\boldsymbol{B}_{i}^{T}\left[\boldsymbol{D}^{-}\right]^{T}\boldsymbol{n}\right)\boldsymbol{N}_{j}\boldsymbol{a}_{j}\mathrm{d\Gamma}\nonumber \\
 & = & \int_{\Omega_{e}}\boldsymbol{N}_{i}^{T}\boldsymbol{f}\mathrm{d\Omega}+\int_{\Gamma_{\mbox{n}_{e}}}\boldsymbol{N}_{i}^{T}\boldsymbol{t}\mathrm{d\Gamma}-\frac{1}{2}\int_{\mathcal{S}_{e}}\left(\boldsymbol{B}_{i}^{T}\left[\boldsymbol{D}^{+}\right]^{T}\boldsymbol{n}+\boldsymbol{B}_{i}^{T}\left[\boldsymbol{D}^{-}\right]^{T}\boldsymbol{n}\right)\,\bar{\mbox{i}}\,\mathrm{d}\mathrm{\Gamma}\nonumber \\
 &  & +\int_{\mathcal{S}_{e}}\boldsymbol{N}_{i}{}^{T}\,\bar{\mbox{j}}\,\mathrm{d}\mathrm{\Gamma}
\end{eqnarray}
\begin{eqnarray}
 &  & \int_{\Omega_{e}}\tilde{\boldsymbol{H}}_{i}^{T}\boldsymbol{B}_{i}{}^{T}\boldsymbol{D}\boldsymbol{B}_{j}\left(\boldsymbol{u}_{j}+\tilde{\boldsymbol{H}}_{j}\boldsymbol{a}_{j}\right)\mathrm{d\Omega}-\int_{\mathcal{S}_{e}}\boldsymbol{N}_{i}{}^{T}\left(\boldsymbol{n}^{T}\boldsymbol{D}^{+}\boldsymbol{B}_{j}+\boldsymbol{n}^{T}\boldsymbol{D}^{-}\boldsymbol{B}_{j}\right)\boldsymbol{u}_{j}\mathrm{d\Gamma}\nonumber \\
 &  & -\int_{\mathcal{S}_{e}}\boldsymbol{N}_{i}{}^{T}\left(\boldsymbol{n}^{T}\boldsymbol{D}^{+}\boldsymbol{B}_{j}\tilde{\boldsymbol{H}}_{j}^{+}+\boldsymbol{n}^{T}\boldsymbol{D}^{-}\boldsymbol{B}_{j}\tilde{\boldsymbol{H}}_{j}^{-}\right)\boldsymbol{a}_{j}\mathrm{d\Gamma}\nonumber \\
 &  & -\int_{\mathcal{S}_{e}}\left(\left[\tilde{\boldsymbol{H}}_{i}^{+}\right]^{T}\boldsymbol{B}_{i}^{T}\left[\boldsymbol{D}^{+}\right]^{T}\boldsymbol{n}+\left[\tilde{\boldsymbol{H}}_{i}^{-}\right]^{T}\boldsymbol{B}_{i}^{T}\left[\boldsymbol{D}^{-}\right]^{T}\boldsymbol{n}\right)\boldsymbol{N}_{j}\boldsymbol{a}_{j}\mathrm{d\Gamma}+4\int_{\mathcal{S}_{e}}\boldsymbol{N}_{i}{}^{T}\alpha\boldsymbol{N}_{j}\boldsymbol{a}_{j}\mathrm{d\Gamma}\nonumber \\
 & = & \int_{\Omega_{e}}\tilde{\boldsymbol{H}}_{i}^{T}\boldsymbol{N}_{i}^{T}\boldsymbol{f}\mathrm{d\Omega}+\int_{\Gamma_{\mbox{n}_{e}}}\tilde{\boldsymbol{H}}_{i}^{T}\boldsymbol{N}_{i}{}^{T}\boldsymbol{t}\mathrm{d\Gamma}+2\int_{\mathcal{S}_{e}}\boldsymbol{N}_{i}{}^{T}\alpha\,\bar{\mbox{i}}\,\mathrm{d}\mathrm{\Gamma}-\int_{\mathcal{S}_{e}}\boldsymbol{H}_{i}^{T}\boldsymbol{N}_{i}{}^{T}\,\bar{\mbox{j}}\,\mathrm{d}\mathrm{\Gamma}\nonumber \\
 &  & -\frac{1}{2}\int_{\mathcal{S}_{e}}\left(\left[\tilde{\boldsymbol{H}}_{i}^{+}\right]^{T}\boldsymbol{B}_{i}^{T}\left[\boldsymbol{D}^{+}\right]^{T}\boldsymbol{n}+\left[\tilde{\boldsymbol{H}}_{i}^{-}\right]^{T}\boldsymbol{B}_{i}^{T}\left[\boldsymbol{D}^{-}\right]^{T}\boldsymbol{n}\right)\,\bar{\mbox{i}}\,\mathrm{d}\mathrm{\Gamma}
\end{eqnarray}
Separating the variables, and writing in a matrix form as before gives
us the following system:\label{Considering-the-jump}{\footnotesize{}
\begin{eqnarray}
\left[\begin{array}{cc}
\sum_{e}\int_{\Omega_{e}}\boldsymbol{B}_{i}^{T}\boldsymbol{D}\boldsymbol{B}_{j}\mathrm{d\Omega} & \sum_{e}\int_{\Omega_{e}}\boldsymbol{B}_{i}^{T}\boldsymbol{D}\boldsymbol{B}_{j}\tilde{\boldsymbol{H}_{j}}\mathrm{d\Omega}\\
 & -\sum_{e}\int_{\mathcal{S}_{e}}\boldsymbol{B}_{i}^{T}\left[\boldsymbol{D}^{+}\right]^{T}\boldsymbol{n}\boldsymbol{N}_{j}\mathrm{d\Gamma}\\
 & -\sum_{e}\int_{\mathcal{S}_{e}}\boldsymbol{B}_{i}^{T}\left[\boldsymbol{D}^{-}\right]^{T}\boldsymbol{n}\boldsymbol{N}_{j}\mathrm{d\Gamma}\\
\\
\sum_{e}\int_{\Omega_{e}}\tilde{\boldsymbol{H}}_{i}^{T}\boldsymbol{B}_{i}{}^{T}\boldsymbol{D}\boldsymbol{B}_{j}\mathrm{d\Gamma} & \sum_{e}\int_{\Omega_{e}}\tilde{\boldsymbol{H}}_{i}^{T}\boldsymbol{B}_{i}{}^{T}\boldsymbol{D}\boldsymbol{B}_{j}\tilde{\boldsymbol{H}}_{j}\mathrm{d\Omega}\\
-\sum_{e}\int_{\mathcal{S}_{e}}\boldsymbol{N}_{i}{}^{T}\boldsymbol{n}^{T}\boldsymbol{D}^{+}\boldsymbol{B}_{j}\mathrm{d\Gamma} & -\sum_{e}\int_{\mathcal{S}_{e}}\boldsymbol{N}_{i}{}^{T}\boldsymbol{n}^{T}\boldsymbol{D}^{+}\boldsymbol{B}_{j}\tilde{\boldsymbol{H}}_{j}^{+}\mathrm{d\Gamma}\\
-\sum_{e}\int_{\mathcal{S}_{e}}\boldsymbol{N}_{i}{}^{T}\boldsymbol{n}^{T}\boldsymbol{D}^{-}\boldsymbol{B}_{j}\mathrm{d\Gamma} & -\sum_{e}\int_{\mathcal{S}_{e}}\boldsymbol{N}_{i}{}^{T}\boldsymbol{n}^{T}\boldsymbol{D}^{-}\boldsymbol{B}_{j}\tilde{\boldsymbol{H}}_{j}^{-}\mathrm{d\Gamma}\\
 & -\sum_{e}\int_{\mathcal{S}_{e}}\left[\tilde{\boldsymbol{H}}_{i}^{+}\right]^{T}\boldsymbol{B}_{i}^{T}\left[\boldsymbol{D}^{+}\right]^{T}\boldsymbol{n}\boldsymbol{N}_{j}\mathrm{d\Gamma}\\
 & -\sum_{e}\int_{\mathcal{S}_{e}}\left[\tilde{\boldsymbol{H}}_{i}^{-}\right]^{T}\boldsymbol{B}_{i}^{T}\left[\boldsymbol{D}^{-}\right]^{T}\boldsymbol{n}\boldsymbol{N}_{j}\mathrm{d\Gamma}\\
 & +4\sum_{e}\int_{\mathcal{S}_{e}}\boldsymbol{N}_{i}{}^{T}\alpha\boldsymbol{N}_{j}\mathrm{d\Gamma}
\end{array}\right]\left\{ \begin{array}{c}
\boldsymbol{u}_{j}\\
\boldsymbol{a}_{j}
\end{array}\right\}  & = & \left\{ \begin{array}{c}
\sum_{e}\int_{\Omega_{e}}\boldsymbol{N}_{i}^{T}\boldsymbol{f}\mathrm{d\Omega}\\
+\sum_{e}\int_{\Gamma_{\mbox{n}_{e}}}\boldsymbol{N}_{i}^{T}\boldsymbol{h}\mathrm{d\Gamma}\\
-\frac{1}{2}\sum_{e}\int_{\mathcal{S}_{e}}\boldsymbol{B}_{i}^{T}\left[\boldsymbol{D}^{+}\right]^{T}\boldsymbol{n}\,\bar{\mbox{i}}\,\mathrm{d}\mathrm{\Gamma}\\
-\frac{1}{2}\sum_{e}\int_{\mathcal{S}_{e}}\boldsymbol{B}_{i}^{T}\left[\boldsymbol{D}^{-}\right]^{T}\boldsymbol{n}\,\bar{\mbox{i}}\,\mathrm{d}\mathrm{\Gamma}\\
+\sum_{e}\int_{\mathcal{S}_{e}}\boldsymbol{N}_{i}{}^{T}\,\bar{\mbox{j}}\,\mathrm{d}\mathrm{\Gamma}\\
\\
\sum_{e}\int_{\Omega_{e}}\tilde{\boldsymbol{H}}_{i}^{T}\boldsymbol{N}_{i}^{T}\boldsymbol{f}\mathrm{d\Omega}\\
+\sum_{e}\int_{\Gamma_{\mbox{n}_{e}}}\tilde{\boldsymbol{H}}_{i}^{T}\boldsymbol{N}_{i}{}^{T}\boldsymbol{t}\mathrm{d\Gamma}\\
+2\sum_{e}\int_{\mathcal{S}_{e}}\boldsymbol{N}_{i}{}^{T}\alpha\,\bar{\mbox{i}}\,\mathrm{d}\mathrm{\Gamma}\\
-\sum_{e}\int_{\mathcal{S}_{e}}\boldsymbol{H}_{i}^{T}\boldsymbol{N}_{i}{}^{T}\,\bar{\mbox{j}}\,\mathrm{d}\mathrm{\Gamma}\\
-\frac{1}{2}\sum_{e}\int_{\mathcal{S}_{e}}\left[\tilde{\boldsymbol{H}}_{i}^{+}\right]^{T}\boldsymbol{B}_{i}^{T}\left[\boldsymbol{D}^{+}\right]^{T}\boldsymbol{n}\,\bar{\mbox{i}}\,\mathrm{d}\mathrm{\Gamma}\\
-\frac{1}{2}\sum_{e}\int_{\mathcal{S}_{e}}\left[\tilde{\boldsymbol{H}}_{i}^{-}\right]^{T}\boldsymbol{B}_{i}^{T}\left[\boldsymbol{D}^{-}\right]^{T}\boldsymbol{n}\,\bar{\mbox{i}}\,\mathrm{d}\mathrm{\Gamma}
\end{array}\right\} \nonumber \\
\end{eqnarray}
}
\begin{equation}
\left[\begin{array}{cc}
\mathbf{K}_{b} & \mathbf{K}_{b}\tilde{\boldsymbol{H}}-\left[\mathbf{K}_{n}^{+}\right]^{T}-\left[\mathbf{K}_{n}^{-}\right]^{T}\\
\\
 & \tilde{\boldsymbol{H}}^{T}\mathbf{K}_{b}\tilde{\boldsymbol{H}}+4\mathbf{K}_{s}\\
\tilde{\boldsymbol{H}}^{T}\mathbf{K}_{b}-\mathbf{K}_{n}^{+}-\mathbf{K}_{n}^{-} & -(\mathbf{K}_{n}^{+}\tilde{\boldsymbol{H}}^{+}+\mathbf{K}_{n}^{-}\tilde{\boldsymbol{H}}^{-})\\
 & -\left(\left[\mathbf{K}_{n}^{+}\tilde{\boldsymbol{H}}^{+}\right]^{T}+\left[\mathbf{K}_{n}^{-}\tilde{\boldsymbol{H}}^{-}\right]^{T}\right)
\end{array}\right]\left\{ \begin{array}{c}
\boldsymbol{u}_{j}\\
\boldsymbol{a}_{j}
\end{array}\right\} =\left\{ \begin{array}{c}
\mathbf{f}_{b}+\mathbf{f}_{h}-\frac{1}{2}\mathbf{f}_{n}^{+}-\frac{1}{2}\mathbf{f}_{n}^{-}+\mathbf{f}_{\mbox{j}}\\
\\
\tilde{\boldsymbol{H}}^{T}(\mathbf{f}_{b}+\mathbf{f}_{h})\\
-\frac{1}{2}\left(\left[\tilde{\boldsymbol{H}}^{+}\right]^{T}\mathbf{f}_{n}^{+}+\left[\tilde{\boldsymbol{H}}^{-}\right]^{T}\mathbf{f}_{n}^{-}\right)\\
+2\mathbf{f}_{s}-\boldsymbol{H}^{T}\mathbf{f}_{\mbox{j}}
\end{array}\right\} \label{eq:-19-1}
\end{equation}
where:
\[
\mathbf{K}_{b}\tilde{\boldsymbol{H}}=\sum_{e}\int_{\Omega_{e}}\boldsymbol{B}_{i}^{T}\boldsymbol{D}\boldsymbol{B}_{j}\tilde{\boldsymbol{H}}_{j}\mathrm{d\Omega}
\]
\[
\tilde{\boldsymbol{H}}^{T}\mathbf{K}_{b}=\sum_{e}\int_{\Omega_{e}}\tilde{\boldsymbol{H}}_{i}^{T}\boldsymbol{B}_{i}^{T}\boldsymbol{D}\boldsymbol{B}_{j}\mathrm{d\Omega}
\]
and:
\[
\tilde{\boldsymbol{H}}^{T}\mathbf{K}_{b}\tilde{\boldsymbol{H}}=\sum_{e}\int_{\Omega_{e}}\tilde{\boldsymbol{H}}_{i}^{T}\boldsymbol{B}_{i}^{T}\boldsymbol{D}\boldsymbol{B}_{j}\tilde{\boldsymbol{H}}_{j}\mathrm{d\Omega}
\]
It can be seen that $\tilde{\boldsymbol{H}}$ and $\boldsymbol{H}$
are similar but does not imply the same type of enrichment.

\subsubsection{Dirichlet Condition}

Consider the Dirichlet condition case from (\ref{eq:-23}) and implementing
the shifted enrichment gives:

\begin{eqnarray}
 &  & \int_{\Omega_{e}}\left[\boldsymbol{v}_{i}+\tilde{\boldsymbol{H}}_{i}\boldsymbol{b}_{i}\right]{}^{T}\boldsymbol{B}_{i}^{T}\sigma(\boldsymbol{u}_{h})\mathrm{d\Omega}-\int_{\mathcal{S}_{e}}\left[\boldsymbol{v}_{i}+\tilde{\boldsymbol{H}}_{i}^{+}\boldsymbol{b}_{i}\right]{}^{T}\boldsymbol{N}_{i}^{T}\left(\sigma(\boldsymbol{u}_{h})^{+}.\mathbf{n}\right)\mathrm{d}\mathrm{\Gamma}\nonumber \\
 &  & -\int_{\mathcal{S}_{e}}\left[\boldsymbol{v}_{i}+\tilde{\boldsymbol{H}}_{i}^{+}\boldsymbol{b}_{i}\right]{}^{T}\boldsymbol{B}_{i}^{T}\left[\boldsymbol{D}^{+}\right]^{T}\boldsymbol{n}\boldsymbol{u}_{h}^{+}\mathrm{d}\mathrm{\Gamma}+\int_{\mathcal{S}_{e}}\left[\boldsymbol{v}_{i}+\tilde{\boldsymbol{H}}_{i}^{+}\boldsymbol{b}_{i}\right]{}^{T}\boldsymbol{N}_{i}^{T}\alpha^{+}\boldsymbol{u}_{h}^{+}\mathrm{d}\mathrm{\Gamma}\nonumber \\
 &  & +\int_{\mathcal{S}_{e}}\left[\boldsymbol{v}_{i}+\tilde{\boldsymbol{H}}_{i}^{-}\boldsymbol{b}_{i}\right]{}^{T}\boldsymbol{N}_{i}^{T}\left(\sigma(\boldsymbol{u}_{h})^{-}.\mathbf{n}\right)\mathrm{d}\mathrm{\Gamma}+\int_{\mathcal{S}_{e}}\left[\boldsymbol{v}_{i}+\tilde{\boldsymbol{H}}_{i}^{-}\boldsymbol{b}_{i}\right]{}^{T}\boldsymbol{B}_{i}^{T}\left[\boldsymbol{D}^{-}\right]^{T}\boldsymbol{n}\boldsymbol{u}_{h}^{-}\mathrm{d}\mathrm{\Gamma}\nonumber \\
 &  & +\int_{\mathcal{S}_{e}}\left[\boldsymbol{v}_{i}+\tilde{\boldsymbol{H}}_{i}^{-}\boldsymbol{b}_{i}\right]{}^{T}\boldsymbol{N}_{i}^{T}\alpha^{-}\boldsymbol{u}_{h}^{-}\mathrm{d}\mathrm{\Gamma}\nonumber \\
 & = & \int_{\Omega_{e}}\left[\boldsymbol{v}_{i}+\tilde{\boldsymbol{H}}_{i}\boldsymbol{b}_{i}\right]{}^{T}\boldsymbol{N}_{i}^{T}\boldsymbol{f}\mathrm{d\Omega}+\int_{\Gamma_{\mbox{n}_{e}}}\left[\boldsymbol{v}_{i}+\tilde{\boldsymbol{H}}_{i}\boldsymbol{b}_{i}\right]{}^{T}\boldsymbol{N}_{i}^{T}\boldsymbol{t}\mathrm{d\Gamma}\nonumber \\
 &  & -\int_{\mathcal{S}_{e}}\left[\boldsymbol{v}_{i}+\tilde{\boldsymbol{H}}_{i}^{+}\boldsymbol{b}_{i}\right]{}^{T}\boldsymbol{B}_{i}^{T}\left[\boldsymbol{D}^{+}\right]^{T}\boldsymbol{n}\boldsymbol{g}^{+}\mathrm{d}\mathrm{\Gamma}+\int_{\mathcal{S}_{e}}\left[\boldsymbol{v}_{i}+\tilde{\boldsymbol{H}}_{i}^{+}\boldsymbol{b}_{i}\right]{}^{T}\boldsymbol{N}_{i}^{T}\alpha^{+}\boldsymbol{g}^{+}\mathrm{d}\mathrm{\Gamma}\nonumber \\
 &  & +\int_{\mathcal{S}_{e}}\left[\boldsymbol{v}_{i}+\tilde{\boldsymbol{H}}_{i}^{-}\boldsymbol{b}_{i}\right]{}^{T}\boldsymbol{B}_{i}^{T}\left[\boldsymbol{D}^{-}\right]^{T}\boldsymbol{n}\boldsymbol{g}^{-}\mathrm{d}\mathrm{\Gamma}+\int_{\mathcal{S}_{e}}\left[\boldsymbol{v}_{i}+\tilde{\boldsymbol{H}}_{i}^{-}\boldsymbol{b}_{i}\right]{}^{T}\boldsymbol{N}_{i}^{T}\alpha^{-}\boldsymbol{g}^{-}\mathrm{d}\mathrm{\Gamma}
\end{eqnarray}
Knowing the arbitrariness of $\boldsymbol{v}_{i}\mbox{ and }\boldsymbol{b}_{i}$,
we can separate the system into two equations: 
\begin{eqnarray}
 &  & \int_{\Omega_{e}}\boldsymbol{B}_{i}^{T}\boldsymbol{D}\boldsymbol{B}_{j}\left(\boldsymbol{u}_{j}+\tilde{\boldsymbol{H}}_{j}\boldsymbol{a}_{j}\right)\mathrm{d\Omega}-\int_{\mathcal{S}_{e}}\boldsymbol{N}_{i}^{T}\boldsymbol{n}^{T}\boldsymbol{D}^{+}\boldsymbol{B}_{j}\left(\boldsymbol{u}_{j}+\tilde{\boldsymbol{H}}_{j}^{+}\boldsymbol{a}_{j}\right)\mathrm{d}\mathrm{\Gamma}\nonumber \\
 &  & -\int_{\mathcal{S}_{e}}\boldsymbol{B}_{i}^{T}\left[\boldsymbol{D}^{+}\right]^{T}\boldsymbol{n}\boldsymbol{N}_{j}\left(\boldsymbol{u}_{j}+\tilde{\boldsymbol{H}}_{j}^{+}\boldsymbol{a}_{j}\right)\mathrm{d}\mathrm{\Gamma}+\int_{\mathcal{S}_{e}}\boldsymbol{N}_{i}^{T}\alpha^{+}\boldsymbol{N}_{j}\left(\boldsymbol{u}_{j}+\tilde{\boldsymbol{H}}_{j}^{+}\boldsymbol{a}_{j}\right)\mathrm{d}\mathrm{\Gamma}\nonumber \\
 &  & +\int_{\mathcal{S}_{e}}\boldsymbol{N}_{i}^{T}\boldsymbol{n}^{T}\boldsymbol{D}^{-}\boldsymbol{B}_{j}\left(\boldsymbol{u}_{j}+\tilde{\boldsymbol{H}}_{j}^{-}\boldsymbol{a}_{j}\right)\mathrm{d}\mathrm{\Gamma}+\int_{\mathcal{S}_{e}}\boldsymbol{B}_{i}^{T}\left[\boldsymbol{D}^{-}\right]^{T}\boldsymbol{n}\boldsymbol{N}_{j}\left(\boldsymbol{u}_{j}+\tilde{\boldsymbol{H}}_{j}^{-}\boldsymbol{a}_{j}\right)\mathrm{d}\mathrm{\Gamma}\nonumber \\
 &  & +\int_{\mathcal{S}_{e}}\boldsymbol{N}_{i}^{T}\alpha^{-}\boldsymbol{N}_{j}\left(\boldsymbol{u}_{j}+\tilde{\boldsymbol{H}}_{j}^{-}\boldsymbol{a}_{j}\right)\mathrm{d}\mathrm{\Gamma}\nonumber \\
 & = & \int_{\Omega_{e}}\boldsymbol{N}_{i}^{T}\boldsymbol{f}\mathrm{d\Omega}+\int_{\Gamma_{\mbox{n}_{e}}}\boldsymbol{N}_{i}^{T}\boldsymbol{t}\mathrm{d\Gamma}-\int_{\mathcal{S}_{e}}\boldsymbol{B}_{i}^{T}\left[\boldsymbol{D}^{+}\right]^{T}\boldsymbol{n}\boldsymbol{g}^{+}\mathrm{d}\mathrm{\Gamma}\nonumber \\
 &  & +\int_{\mathcal{S}_{e}}\boldsymbol{N}_{i}^{T}\alpha^{+}\boldsymbol{g}^{+}\mathrm{d}\mathrm{\Gamma}+\int_{\mathcal{S}_{e}}\boldsymbol{B}_{i}^{T}\left[\boldsymbol{D}^{-}\right]^{T}\boldsymbol{n}\boldsymbol{g}^{-}\mathrm{d}\mathrm{\Gamma}+\int_{\mathcal{S}_{e}}\boldsymbol{N}_{i}^{T}\alpha^{-}\boldsymbol{g}^{-}\mathrm{d}\mathrm{\Gamma}
\end{eqnarray}
\begin{eqnarray}
 &  & \int_{\Omega_{e}}\left[\tilde{\boldsymbol{H}}_{i}^{-}\right]^{T}\boldsymbol{B}_{i}^{T}\boldsymbol{D}\boldsymbol{B}_{j}\left(\boldsymbol{u}_{j}+\tilde{\boldsymbol{H}}_{j}\boldsymbol{a}_{j}\right)\mathrm{d\Omega}-\int_{\mathcal{S}_{e}}\left[\tilde{\boldsymbol{H}}_{i}^{-}\right]^{T}{}^{+}\boldsymbol{N}_{i}^{T}\boldsymbol{n}^{T}\boldsymbol{D}^{+}\boldsymbol{B}_{j}\left(\boldsymbol{u}_{j}+\tilde{\boldsymbol{H}}_{j}^{+}\boldsymbol{a}_{j}\right)\mathrm{d}\mathrm{\Gamma}\nonumber \\
 &  & -\int_{\mathcal{S}_{e}}\left[\tilde{\boldsymbol{H}}_{i}^{-}\right]^{T}{}^{+}\boldsymbol{B}_{i}^{T}\left[\boldsymbol{D}^{+}\right]^{T}\boldsymbol{n}\boldsymbol{N}_{j}\left(\boldsymbol{u}_{j}+\tilde{\boldsymbol{H}}_{j}^{+}\boldsymbol{a}_{j}\right)\mathrm{d}\mathrm{\Gamma}+\int_{\mathcal{S}_{e}}\left[\tilde{\boldsymbol{H}}_{i}^{-}\right]^{T}{}^{+}\boldsymbol{N}_{i}^{T}\alpha^{+}\boldsymbol{N}_{j}\left(\boldsymbol{u}_{j}+\tilde{\boldsymbol{H}}_{j}^{+}\boldsymbol{a}_{j}\right)\mathrm{d}\mathrm{\Gamma}\nonumber \\
 &  & +\int_{\mathcal{S}_{e}}\left[\tilde{\boldsymbol{H}}_{i}^{-}\right]^{T}{}^{-}\boldsymbol{N}_{i}^{T}\boldsymbol{n}^{T}\boldsymbol{D}^{-}\boldsymbol{B}_{j}\left(\boldsymbol{u}_{j}+\tilde{\boldsymbol{H}}_{j}^{-}\boldsymbol{a}_{j}\right)\mathrm{d}\mathrm{\Gamma}+\int_{\mathcal{S}_{e}}\left[\tilde{\boldsymbol{H}}_{i}^{-}\right]^{T}{}^{-}\boldsymbol{B}_{i}^{T}\left[\boldsymbol{D}^{-}\right]^{T}\boldsymbol{n}\boldsymbol{N}_{j}\left(\boldsymbol{u}_{j}+\tilde{\boldsymbol{H}}_{j}^{-}\boldsymbol{a}_{j}\right)\mathrm{d}\mathrm{\Gamma}\nonumber \\
 &  & +\int_{\mathcal{S}_{e}}\left[\tilde{\boldsymbol{H}}_{i}^{-}\right]^{T}{}^{-}\boldsymbol{N}_{i}^{T}\alpha^{-}\boldsymbol{N}_{j}\left(\boldsymbol{u}_{j}+\tilde{\boldsymbol{H}}_{j}^{-}\boldsymbol{a}_{j}\right)\mathrm{d}\mathrm{\Gamma}\nonumber \\
 & = & \int_{\Omega_{e}}\left[\tilde{\boldsymbol{H}}_{i}^{-}\right]^{T}\boldsymbol{N}_{i}^{T}\boldsymbol{f}\mathrm{d\Omega}+\int_{\Gamma_{\mbox{n}_{e}}}\left[\tilde{\boldsymbol{H}}_{i}^{-}\right]^{T}\boldsymbol{N}_{i}^{T}\boldsymbol{t}\mathrm{d\Gamma}-\int_{\mathcal{S}_{e}}\left[\tilde{\boldsymbol{H}}_{i}^{-}\right]^{T}{}^{+}\boldsymbol{B}_{i}^{T}\left[\boldsymbol{D}^{+}\right]^{T}\boldsymbol{n}\boldsymbol{g}^{+}\mathrm{d}\mathrm{\Gamma}\nonumber \\
 &  & +\int_{\mathcal{S}_{e}}\left[\tilde{\boldsymbol{H}}_{i}^{-}\right]^{T}{}^{+}\boldsymbol{N}_{i}^{T}\alpha^{+}\boldsymbol{g}^{+}\mathrm{d}\mathrm{\Gamma}+\int_{\mathcal{S}_{e}}\left[\tilde{\boldsymbol{H}}_{i}^{-}\right]^{T}{}^{-}\boldsymbol{B_{i}}^{T}\left[\boldsymbol{D}^{-}\right]^{T}\boldsymbol{n}\boldsymbol{g}^{-}\mathrm{d}\mathrm{\Gamma}+\int_{\mathcal{S}_{e}}\left[\tilde{\boldsymbol{H}}_{i}^{-}\right]^{T}{}^{-}\boldsymbol{N}_{i}^{T}\alpha^{-}\boldsymbol{g}^{-}\mathrm{d}\mathrm{\Gamma}\nonumber \\
\end{eqnarray}
Writing the system of equations in the form of a matrix gives us the
following:{\footnotesize{}
\begin{eqnarray}
\left[\begin{array}{cc}
\sum_{e}\int_{\Omega_{e}}\boldsymbol{B}_{i}^{T}\boldsymbol{D}\boldsymbol{B}_{j}\mathrm{d\Omega} & \sum_{e}\int_{\Omega_{e}}\boldsymbol{B}_{i}^{T}\boldsymbol{D}\boldsymbol{B}_{j}\tilde{\boldsymbol{H}}_{j}\mathrm{d\Omega}\\
-\sum_{e}\int_{\mathcal{S}_{e}}\boldsymbol{N}_{i}^{T}\boldsymbol{n}^{T}\boldsymbol{D}^{+}\boldsymbol{B}_{j}\mathrm{d}\mathrm{\Gamma} & -\sum_{e}\int_{\mathcal{S}_{e}}\boldsymbol{N}_{i}^{T}\boldsymbol{n}^{T}\boldsymbol{D}^{+}\boldsymbol{B}_{j}\tilde{\boldsymbol{H}}_{j}^{+}\mathrm{d}\mathrm{\Gamma}\\
-\sum_{e}\int_{\mathcal{S}_{e}}\boldsymbol{B}_{i}^{T}\left[\boldsymbol{D}^{+}\right]^{T}\boldsymbol{n}\boldsymbol{N}_{j}\mathrm{d}\mathrm{\Gamma} & -\sum_{e}\int_{\mathcal{S}_{e}}\boldsymbol{B}_{i}^{T}\left[\boldsymbol{D}^{+}\right]^{T}\boldsymbol{n}\boldsymbol{N}_{j}\tilde{\boldsymbol{H}}_{j}^{+}\mathrm{d}\mathrm{\Gamma}\\
+\sum_{e}\int_{\mathcal{S}_{e}}\boldsymbol{N}_{i}^{T}\alpha^{+}\boldsymbol{N}_{j}\mathrm{d}\mathrm{\Gamma} & +\sum_{e}\int_{\mathcal{S}_{e}}\boldsymbol{N}_{i}^{T}\alpha^{+}\boldsymbol{N}_{j}\tilde{\boldsymbol{H}}_{j}^{+}\mathrm{d}\mathrm{\Gamma}\\
+\sum_{e}\int_{\mathcal{S}_{e}}\boldsymbol{N}_{i}^{T}\boldsymbol{n}^{T}\boldsymbol{D}^{-}\boldsymbol{B}_{j}\mathrm{d}\mathrm{\Gamma} & +\sum_{e}\int_{\mathcal{S}_{e}}\boldsymbol{N}_{i}^{T}\boldsymbol{n}^{T}\boldsymbol{D}^{-}\boldsymbol{B}_{j}\tilde{\boldsymbol{H}}_{j}^{-}\mathrm{d}\mathrm{\Gamma}\\
+\sum_{e}\int_{\mathcal{S}_{e}}\boldsymbol{B}_{i}^{T}\left[\boldsymbol{D}^{-}\right]^{T}\boldsymbol{n}\boldsymbol{N}_{j}\mathrm{d}\mathrm{\Gamma} & +\sum_{e}\int_{\mathcal{S}_{e}}\boldsymbol{B}_{i}^{T}\left[\boldsymbol{D}^{-}\right]^{T}\boldsymbol{n}\boldsymbol{N}_{j}\tilde{\boldsymbol{H}}_{j}^{-}\mathrm{d}\mathrm{\Gamma}\\
+\sum_{e}\int_{\mathcal{S}_{e}}\boldsymbol{N}_{i}^{T}\alpha^{-}\boldsymbol{N}_{j}\mathrm{d}\mathrm{\Gamma} & +\sum_{e}\int_{\mathcal{S}_{e}}\boldsymbol{N}_{i}^{T}\alpha^{-}\boldsymbol{N}_{j}\tilde{\boldsymbol{H}}_{j}^{-}\mathrm{d}\mathrm{\Gamma}\\
\\
\sum_{e}\int_{\Omega_{e}}\left[\tilde{\boldsymbol{H}}_{i}^{-}\right]^{T}\boldsymbol{B}_{i}^{T}\boldsymbol{D}\boldsymbol{B}_{j}\mathrm{d\Omega} & \sum_{e}\int_{\Omega_{e}}\left[\tilde{\boldsymbol{H}}_{i}^{-}\right]^{T}\boldsymbol{B}_{i}^{T}\boldsymbol{D}\boldsymbol{B}_{j}\tilde{\boldsymbol{H}}_{j}\mathrm{d\Omega}\\
-\sum_{e}\int_{\mathcal{S}_{e}}\left[\tilde{\boldsymbol{H}}_{i}^{-}\right]^{T}{}^{+}\boldsymbol{N}_{i}^{T}\boldsymbol{n}^{T}\boldsymbol{D}^{+}\boldsymbol{B}_{j}\mathrm{d}\mathrm{\Gamma} & -\sum_{e}\int_{\mathcal{S}_{e}}\left[\tilde{\boldsymbol{H}}_{i}^{-}\right]^{T}{}^{+}\boldsymbol{N}_{i}^{T}\boldsymbol{n}^{T}\boldsymbol{D}^{+}\boldsymbol{B}_{j}\tilde{\boldsymbol{H}}_{j}^{+}\mathrm{d}\mathrm{\Gamma}\\
-\sum_{e}\int_{\mathcal{S}_{e}}\left[\tilde{\boldsymbol{H}}_{i}^{-}\right]^{T}{}^{+}\boldsymbol{B}_{i}^{T}\left[\boldsymbol{D}^{+}\right]^{T}\boldsymbol{n}\boldsymbol{N}_{j}\mathrm{d}\mathrm{\Gamma} & -\sum_{e}\int_{\mathcal{S}_{e}}\left[\tilde{\boldsymbol{H}}_{i}^{-}\right]^{T}{}^{+}\boldsymbol{B}_{i}^{T}\left[\boldsymbol{D}^{+}\right]^{T}\boldsymbol{n}\boldsymbol{N}_{j}\tilde{\boldsymbol{H}}_{j}^{+}\mathrm{d}\mathrm{\Gamma}\\
+\sum_{e}\int_{\mathcal{S}_{e}}\left[\tilde{\boldsymbol{H}}_{i}^{-}\right]^{T}{}^{+}\boldsymbol{N}_{i}^{T}\alpha^{+}\boldsymbol{N}_{j}\mathrm{d}\mathrm{\Gamma} & +\sum_{e}\int_{\mathcal{S}_{e}}\left[\tilde{\boldsymbol{H}}_{i}^{-}\right]^{T}{}^{+}\boldsymbol{N}_{i}^{T}\alpha^{+}\boldsymbol{N}_{j}\tilde{\boldsymbol{H}}_{j}^{+}\mathrm{d}\mathrm{\Gamma}\\
+\sum_{e}\int_{\mathcal{S}_{e}}\left[\tilde{\boldsymbol{H}}_{i}^{-}\right]^{T}{}^{-}\boldsymbol{N}_{i}^{T}\boldsymbol{n}^{T}\boldsymbol{D}^{-}\boldsymbol{B}_{j}\mathrm{d}\mathrm{\Gamma} & +\sum_{e}\int_{\mathcal{S}_{e}}\left[\tilde{\boldsymbol{H}}_{i}^{-}\right]^{T}{}^{-}\boldsymbol{N}_{i}^{T}\boldsymbol{n}^{T}\boldsymbol{D}^{-}\boldsymbol{B}_{j}\tilde{\boldsymbol{H}}_{j}^{-}\mathrm{d}\mathrm{\Gamma}\\
+\sum_{e}\int_{\mathcal{S}_{e}}\left[\tilde{\boldsymbol{H}}_{i}^{-}\right]^{T}{}^{-}\boldsymbol{B}_{i}^{T}\left[\boldsymbol{D}^{-}\right]^{T}\boldsymbol{n}\boldsymbol{N}_{j}\mathrm{d}\mathrm{\Gamma} & +\sum_{e}\int_{\mathcal{S}_{e}}\left[\tilde{\boldsymbol{H}}_{i}^{-}\right]^{T}{}^{-}\boldsymbol{B}_{i}^{T}\left[\boldsymbol{D}^{-}\right]^{T}\boldsymbol{n}\boldsymbol{N}_{j}\tilde{\boldsymbol{H}}_{j}^{-}\mathrm{d}\mathrm{\Gamma}\\
+\sum_{e}\int_{\mathcal{S}_{e}}\left[\tilde{\boldsymbol{H}}_{i}^{-}\right]^{T}{}^{-}\boldsymbol{N}_{i}^{T}\alpha^{-}\boldsymbol{N}_{j}\mathrm{d}\mathrm{\Gamma} & +\sum_{e}\int_{\mathcal{S}_{e}}\left[\tilde{\boldsymbol{H}}_{i}^{-}\right]^{T}{}^{-}\boldsymbol{N}_{i}^{T}\alpha^{-}\boldsymbol{N}_{j}\tilde{\boldsymbol{H}}_{j}^{-}\mathrm{d}\mathrm{\Gamma}
\end{array}\right]\left\{ \begin{array}{c}
\boldsymbol{u}_{j}\\
\boldsymbol{a}_{j}
\end{array}\right\} \nonumber \\
=\left\{ \begin{array}{c}
\sum_{e}\int_{\Omega_{e}}\boldsymbol{N}_{i}^{T}\boldsymbol{f}\mathrm{d\Omega}\\
+\sum_{e}\int_{\Gamma_{\mbox{n}_{e}}}\boldsymbol{N}_{i}^{T}\boldsymbol{t}\mathrm{d\Gamma}\\
-\sum_{e}\int_{\mathcal{S}_{e}}\boldsymbol{B}_{i}^{T}\left[\boldsymbol{D}^{+}\right]^{T}\boldsymbol{n}\boldsymbol{g}^{+}\mathrm{d}\mathrm{\Gamma}\\
+\sum_{e}\int_{\mathcal{S}_{e}}\boldsymbol{N}_{i}^{T}\alpha^{+}\boldsymbol{g}^{+}\mathrm{d}\mathrm{\Gamma}\\
+\sum_{e}\int_{\mathcal{S}_{e}}\boldsymbol{B}_{i}^{T}\left[\boldsymbol{D}^{-}\right]^{T}\boldsymbol{n}\boldsymbol{g}^{-}\mathrm{d}\mathrm{\Gamma}\\
+\sum_{e}\int_{\mathcal{S}_{e}}\boldsymbol{N}_{i}^{T}\alpha^{-}\boldsymbol{g}^{-}\mathrm{d}\mathrm{\Gamma}\\
\\
\sum_{e}\int_{\Omega_{e}}\left[\tilde{\boldsymbol{H}}_{i}^{-}\right]^{T}\boldsymbol{N}_{i}^{T}\boldsymbol{f}\mathrm{d\Omega}\\
+\sum_{e}\int_{\Gamma_{\mbox{n}_{e}}}\left[\tilde{\boldsymbol{H}}_{i}^{-}\right]^{T}\boldsymbol{N}_{i}^{T}\boldsymbol{t}\mathrm{d\Gamma}\\
-\sum_{e}\int_{\mathcal{S}_{e}}\left[\tilde{\boldsymbol{H}}_{i}^{-}\right]^{T}{}^{+}\boldsymbol{B}_{i}^{T}\left[\boldsymbol{D}^{+}\right]^{T}\boldsymbol{n}\boldsymbol{g}^{+}\mathrm{d}\mathrm{\Gamma}\\
+\sum_{e}\int_{\mathcal{S}_{e}}\left[\tilde{\boldsymbol{H}}_{i}^{-}\right]^{T}{}^{+}\boldsymbol{N}_{i}^{T}\alpha^{+}\boldsymbol{g}^{+}\mathrm{d}\mathrm{\Gamma}\\
+\sum_{e}\int_{\mathcal{S}_{e}}\left[\tilde{\boldsymbol{H}}_{i}^{-}\right]^{T}{}^{-}\boldsymbol{B_{i}}^{T}\left[\boldsymbol{D}^{-}\right]^{T}\boldsymbol{n}\boldsymbol{g}^{-}\mathrm{d}\mathrm{\Gamma}\\
+\sum_{e}\int_{\mathcal{S}_{e}}\left[\tilde{\boldsymbol{H}}_{i}^{-}\right]^{T}{}^{-}\boldsymbol{N}_{i}^{T}\alpha^{-}\boldsymbol{g}^{-}\mathrm{d}\mathrm{\Gamma}
\end{array}\right\} 
\end{eqnarray}
}
\begin{equation}
\left[\begin{array}{cc}
\mathbf{K}_{b}^{+}+\mathbf{K}_{n}^{+} & \left[\mathbf{K}_{b}^{+}+\mathbf{K}_{n}^{+}+\mathbf{K}_{s}^{+}\right]\tilde{\boldsymbol{H}}^{+}\\
+\left[\mathbf{K}_{n}^{+}\right]^{T}+\mathbf{K}_{s}^{+} & +\left[\mathbf{K}_{n}^{+}\right]^{T}\tilde{\boldsymbol{H}}^{+}\\
+\mathbf{K}_{b}^{-}-\mathbf{K}_{n}^{-} & +\left[\mathbf{K}_{b}^{-}-\mathbf{K}_{n}^{-}+\mathbf{K}_{s}^{-}\right]\tilde{\boldsymbol{H}}^{-}\\
-\left[\mathbf{K}_{n}^{-}\right]^{T}+\mathbf{K}_{s}^{-} & -\left[\mathbf{K}_{n}^{-}\right]^{T}\tilde{\boldsymbol{H}}^{-}\\
\\
\tilde{\boldsymbol{H}}^{+}\left[\mathbf{K}_{b}^{+}+\mathbf{K}_{n}^{+}+\mathbf{K}_{s}^{+}\right] & \tilde{\boldsymbol{H}}^{+}\left[\mathbf{K}_{b}^{+}+\mathbf{K}_{n}^{+}+\mathbf{K}_{s}^{+}\right]\tilde{\boldsymbol{H}}^{+}\\
+\tilde{\boldsymbol{H}}^{+}\left[\mathbf{K}_{n}^{+}\right]^{T} & +\tilde{\boldsymbol{H}}^{+}\left[\mathbf{K}_{n}^{+}\right]^{T}\tilde{\boldsymbol{H}}^{+}\\
+\tilde{\boldsymbol{H}}^{-}\left[\mathbf{K}_{b}^{-}-\mathbf{K}_{n}^{-}+\mathbf{K}_{s}^{-}\right] & +\tilde{\boldsymbol{H}}^{-}\left[\mathbf{K}_{b}^{-}-\mathbf{K}_{n}^{-}+\mathbf{K}_{s}^{-}\right]\tilde{\boldsymbol{H}}^{-}\\
-\tilde{\boldsymbol{H}}^{-}\left[\mathbf{K}_{n}^{-}\right]^{T} & -\tilde{\boldsymbol{H}}^{-}\left[\mathbf{K}_{n}^{-}\right]^{T}\tilde{\boldsymbol{H}}^{-}
\end{array}\right]\left\{ \begin{array}{c}
\boldsymbol{u}_{j}\\
\boldsymbol{a}_{j}
\end{array}\right\} =\left\{ \begin{array}{c}
\mathbf{f}_{b}^{+}+\mathbf{f}_{h}^{+}+\mathbf{f}_{s}^{+}+\mathbf{f}_{s}^{-}\\
+\mathbf{f}_{b}^{-}+\mathbf{f}_{h}^{-}+\mathbf{f}_{n}^{+}-\mathbf{f}_{n}^{-}\\
\\
\tilde{\boldsymbol{H}}^{+}\left(\mathbf{f}_{b}^{+}+\mathbf{f}_{h}^{+}+\mathbf{f}_{s}^{+}+\mathbf{f}_{n}^{+}\right)\\
+\tilde{\boldsymbol{H}}^{-}\left(\mathbf{f}_{b}^{-}+\mathbf{f}_{h}^{-}+\mathbf{f}_{s}^{-}-\mathbf{f}_{n}^{-}\right)
\end{array}\right\} 
\end{equation}
where $\tilde{\boldsymbol{H}}$ is a diagonal matrix of shifted basis
enrichments for the corresponding nodes. 

Knowing the arbitrariness of $\boldsymbol{v}_{i}\mbox{ and }\boldsymbol{b}_{i}$,
we can write the jump in flux from (\ref{eq:-30}):
\begin{equation}
\int_{\mathcal{S}_{e}}\boldsymbol{N}_{i}{}^{T}\bar{\mathrm{j}}_{h}\mbox{d}\Gamma=\int_{\mathcal{B}_{e}}\boldsymbol{N}_{i}{}^{T}\boldsymbol{f}\mbox{d}\Omega+\int_{\Gamma_{\mbox{n}_{e}}}\boldsymbol{N}_{i}{}^{T}\boldsymbol{t}\mbox{d}\Omega-\int_{\mathcal{B}_{e}}\boldsymbol{B}_{i}{}^{T}\sigma(\boldsymbol{u}_{h})\mbox{d}\Omega
\end{equation}
Similarly, as we have done in section (\ref{par:Dirichlet-condition}),
we can discretize the jump in flux:
\begin{eqnarray}
\int_{\mathcal{S}_{e}}\boldsymbol{N}_{i}{}^{T}\boldsymbol{N}_{j}\bar{\mathrm{j}}_{j}\mbox{d}\Gamma & = & \int_{\mathcal{B}_{e}}\boldsymbol{N}_{i}{}^{T}\boldsymbol{f}\mbox{d}\Omega+\int_{\Gamma_{\mbox{n}_{e}}}\boldsymbol{N}_{i}{}^{T}\boldsymbol{t}\mbox{d}\Gamma-\int_{\mathcal{B}_{e}}\boldsymbol{B}_{i}{}^{T}\boldsymbol{D}\boldsymbol{B}_{j}\left(\boldsymbol{u}_{j}+\tilde{\boldsymbol{H}}_{j}\boldsymbol{a}_{j}\right)\mbox{d}\Omega
\end{eqnarray}
\begin{equation}
\mathbf{M}_{d}\bar{\mathrm{j}}_{j}=\mathbf{f}_{b}+\mathbf{f}_{h}-\mathbf{K}_{b}\boldsymbol{u}_{j}-\mathbf{K}_{b}\tilde{\boldsymbol{H}}\boldsymbol{a}_{j}
\end{equation}
 where, as before, $\mathbf{M}_{d}$ is the mass matrix.

\subsection{Non-Linear Iteration}

Consider the residual of the system to be solved:
\begin{equation}
R\left(\left\{ u\right\} \right)=\left\{ f_{int}\left(\left\{ u\right\} \right)\right\} -\left\{ f_{ext}\right\} 
\end{equation}

\subsubsection{Jump conditions}

We introduce the penalization term and the term for balancing variational
consistency in the system with jump in displacement conditions (\ref{eq:-24})
and (\ref{eq:-29}).
\begin{eqnarray}
R\left(\left\{ u\right\} \right) & = & \int_{\Omega}\sigma\left(\left\{ u\right\} \right):\varepsilon\left(\left\{ u\right\} \right)\mbox{d}V+\int_{\Gamma^{*}}[[v]]^{T}\alpha\left(\left\{ \left[\left[u\right]\right]\right\} -\bar{\mbox{i}}\right)\mathrm{d\Gamma}-\int_{\Gamma^{*}}<v>^{T}\bar{\mbox{j}}\,\mathrm{d}\mathrm{\Gamma}\nonumber \\
 &  & -\int_{\Gamma^{*}}[[v]]^{T}<\sigma(\left\{ u\right\} )>.\mathbf{n}\mathrm{d\Gamma}-\int_{\Omega}v^{T}\boldsymbol{f}\,\mathrm{d\Omega}-\int_{\Gamma}v^{T}\boldsymbol{t}\,\mbox{d}\Gamma
\end{eqnarray}
Discretizing the above system gives us:
\begin{eqnarray}
R\left(\left\{ u\right\} \right) & = & \int_{\Omega}\left[\begin{array}{c}
\mathbf{B}^{T}\\
\tilde{\boldsymbol{H}}\mathbf{B}^{T}
\end{array}\right]\sigma\left(\left\{ u\right\} \right)\mbox{d}V+\int_{\Gamma^{*}}\left[\begin{array}{c}
\mathbf{0}\\
2\mathbf{N}^{T}
\end{array}\right]\alpha\left(\left[\begin{array}{cc}
\mathbf{0} & 2\mathbf{N}\end{array}\right]\left\{ u\right\} -\bar{\mbox{i}}\right)\mathrm{d\Gamma}\nonumber \\
 &  & -\int_{\Gamma^{*}}\left[\begin{array}{c}
\mathbf{N}^{T}\\
\boldsymbol{H}\mathbf{N}^{T}
\end{array}\right]\bar{\mbox{j}}\,\mathrm{d}\mathrm{\Gamma}-\int_{\Gamma^{*}}\left[\begin{array}{c}
\mathbf{0}\\
2\mathbf{N}^{T}
\end{array}\right]<\sigma(\left\{ u\right\} )>.\mathbf{n}\mathrm{d\Gamma}\nonumber \\
 &  & -\int_{\Omega}\left[\begin{array}{c}
\mathbf{N}^{T}\\
\tilde{\boldsymbol{H}}\mathbf{N}^{T}
\end{array}\right]\boldsymbol{f}\,\mathrm{d\Omega}-\int_{\Gamma}\left[\begin{array}{c}
\mathbf{N}^{T}\\
\tilde{\boldsymbol{H}}\mathbf{N}^{T}
\end{array}\right]\boldsymbol{t}\,\mbox{d}\Gamma
\end{eqnarray}
The Jacobian of the system is obtained by differentiating the system
with respect to the discrete nodal displacements.
\begin{eqnarray}
\frac{\partial R}{\partial\left\{ u\right\} }\left\{ u_{k}\right\}  & = & \int_{\Omega}\left[\begin{array}{c}
\mathbf{B}^{T}\\
\tilde{\boldsymbol{H}}\mathbf{B}^{T}
\end{array}\right]\frac{\partial\sigma}{\partial\varepsilon}\frac{\partial\varepsilon\left(\left\{ u\right\} \right)}{\partial\left\{ u\right\} }\mbox{d}V+\int_{\Gamma^{*}}\left[\begin{array}{c}
\mathbf{0}\\
2\mathbf{N}^{T}
\end{array}\right]\alpha\left[\begin{array}{cc}
\mathbf{0} & 2\mathbf{N}\end{array}\right]\mathrm{d\Gamma}\nonumber \\
 &  & -\int_{\Gamma^{*}}\left[\begin{array}{c}
\mathbf{0}\\
2\mathbf{N}^{T}
\end{array}\right]\boldsymbol{n}<\frac{\partial\sigma}{\partial\varepsilon}\frac{\partial\varepsilon\left(\left\{ u\right\} \right)}{\partial\left\{ u\right\} }>\mathrm{d\Gamma}\label{eq:-32}
\end{eqnarray}
\begin{eqnarray}
\frac{\partial R}{\partial\left\{ u\right\} }\left\{ u_{k}\right\}  & = & \int_{\Omega}\left[\begin{array}{c}
\mathbf{B}^{T}\\
\tilde{\boldsymbol{H}}\mathbf{B}^{T}
\end{array}\right]\mathbf{D_{T}}\left[\begin{array}{cc}
\mathbf{B} & \mathbf{B}\tilde{\boldsymbol{H}}\end{array}\right]\mbox{d}V+\int_{\Gamma^{*}}\left[\begin{array}{c}
\mathbf{0}\\
2\mathbf{N}^{T}
\end{array}\right]\alpha\left[\begin{array}{cc}
\mathbf{0} & 2\mathbf{N}\end{array}\right]\mathrm{d\Gamma}\nonumber \\
 &  & -\int_{\Gamma^{*}}\left[\begin{array}{c}
\mathbf{0}\\
2\mathbf{N}^{T}
\end{array}\right]\boldsymbol{n}\mathbf{D_{T}}\left[\begin{array}{cc}
\mathbf{B} & \frac{1}{2}\mathbf{B}\left(\tilde{\boldsymbol{H}}^{+}+\tilde{\boldsymbol{H}}^{-}\right)\end{array}\right]\mathrm{d\Gamma}
\end{eqnarray}
where $\mathbf{D_{T}}=\left(\partial\sigma/\partial\varepsilon\right)$
and is equal to the Hook Tensor in the case of linear elasticity.
As we have learnt before regarding penalization method, the term representing
variational consistency results in an asymmetric matrix. We seek a
solution $\left\{ \delta u\right\} $ of the system:
\begin{equation}
\left(\frac{\partial R}{\partial\left\{ u\right\} }\left\{ u_{k}\right\} \right)\left\{ \delta u\right\} =-R\left(\left\{ u_{k}\right\} \right)
\end{equation}
We can now introduce the symmetric part of the matrix in the system
to obtain Nitsche's system.
\begin{eqnarray}
 &  & \left(\frac{\partial R}{\partial\left\{ u\right\} }\left\{ u_{k}\right\} \right)\left\{ \delta u\right\} -\left(\int_{\Gamma^{*}}\left[\begin{array}{c}
\mathbf{B}^{T}\\
\frac{1}{2}\left[\tilde{\boldsymbol{H}}^{+}+\tilde{\boldsymbol{H}}^{-}\right]^{T}\mathbf{B}^{T}
\end{array}\right]\mathbf{D_{T}}\boldsymbol{n}^{T}\left[\begin{array}{cc}
\mathbf{0} & 2\mathbf{N}\end{array}\right]\mathrm{d\Gamma}\right)\left\{ u_{k+1}\right\} \nonumber \\
 & = & -R\left(\left\{ u_{k}\right\} \right)-\int_{\Gamma^{*}}\left[\begin{array}{c}
\mathbf{B}^{T}\\
\frac{1}{2}\left[\tilde{\boldsymbol{H}}^{+}+\tilde{\boldsymbol{H}}^{-}\right]^{T}\mathbf{B}^{T}
\end{array}\right]\mathbf{D_{T}}\boldsymbol{n}^{T}\bar{\mbox{i}}\,\mathrm{d\Gamma}
\end{eqnarray}
With $u_{k+1}=u_{k}+\delta u$ we can rearrange the terms to obtain:
\begin{eqnarray}
 &  & \left(\frac{\partial R}{\partial\left\{ u\right\} }\left\{ u_{k}\right\} \right)\left\{ \delta u\right\} -\left(\int_{\Gamma^{*}}\left[\begin{array}{c}
\mathbf{B}^{T}\\
\frac{1}{2}\left[\tilde{\boldsymbol{H}}^{+}+\tilde{\boldsymbol{H}}^{-}\right]^{T}\mathbf{B}^{T}
\end{array}\right]\mathbf{D_{T}}\boldsymbol{n}^{T}\left[\begin{array}{cc}
\mathbf{0} & 2\mathbf{N}\end{array}\right]\mathrm{d\Gamma}\right)\left\{ \delta u\right\} \nonumber \\
 & = & -R\left(\left\{ u_{k}\right\} \right)+\int_{\Gamma^{*}}\left[\begin{array}{c}
\mathbf{B}^{T}\\
\frac{1}{2}\left[\tilde{\boldsymbol{H}}^{+}+\tilde{\boldsymbol{H}}^{-}\right]^{T}\mathbf{B}^{T}
\end{array}\right]\mathbf{D_{T}}\boldsymbol{n}^{T}\left(\left[\begin{array}{cc}
\mathbf{0} & 2\mathbf{N}\end{array}\right]\left\{ u_{k}\right\} -\bar{\mbox{i}}\right)\mathrm{d\Gamma}
\end{eqnarray}
since $\left[\left[u\right]\right]=\bar{\mbox{i}}\mbox{ on }\Gamma^{*}$.
Finally we can write:
\begin{equation}
\left[\mathbf{K_{T}}\right]\left\{ \delta u\right\} =-R'\left(\left\{ u_{k}\right\} \right)
\end{equation}
where,
\begin{eqnarray}
\left[\mathbf{K_{T}}\right] & = & \int_{\Omega}\left[\begin{array}{c}
\mathbf{B}^{T}\\
\tilde{\boldsymbol{H}}\mathbf{B}^{T}
\end{array}\right]\mathbf{D_{T}}\left[\begin{array}{cc}
\mathbf{B} & \mathbf{B}\tilde{\boldsymbol{H}}\end{array}\right]\mbox{d}V+\int_{\Gamma^{*}}\left[\begin{array}{c}
\mathbf{0}\\
2\mathbf{N}^{T}
\end{array}\right]\alpha\left[\begin{array}{cc}
\mathbf{0} & 2\mathbf{N}\end{array}\right]\mathrm{d\Gamma}\nonumber \\
 &  & -\int_{\Gamma^{*}}\left[\begin{array}{c}
\mathbf{0}\\
2\mathbf{N}^{T}
\end{array}\right]\boldsymbol{n}\mathbf{D_{T}}\left[\begin{array}{cc}
\mathbf{B} & \frac{1}{2}\mathbf{B}\left(\tilde{\boldsymbol{H}}^{+}+\tilde{\boldsymbol{H}}^{-}\right)\end{array}\right]\mathrm{d\Gamma}\nonumber \\
 &  & -\int_{\Gamma^{*}}\left[\begin{array}{c}
\mathbf{B}^{T}\\
\frac{1}{2}\left[\tilde{\boldsymbol{H}}^{+}+\tilde{\boldsymbol{H}}^{-}\right]^{T}\mathbf{B}^{T}
\end{array}\right]\mathbf{D_{T}}\boldsymbol{n}^{T}\left[\begin{array}{cc}
\mathbf{0} & 2\mathbf{N}\end{array}\right]\mathrm{d\Gamma}
\end{eqnarray}
and
\begin{eqnarray}
-R'\left(\left\{ u_{k}\right\} \right) & = & -\int_{\Omega}\left[\begin{array}{c}
\mathbf{B}^{T}\\
\tilde{\boldsymbol{H}}\mathbf{B}^{T}
\end{array}\right]\sigma\left(\left\{ u_{k}\right\} \right)\mbox{d}V+\int_{\Gamma^{*}}\left[\begin{array}{c}
\mathbf{N}^{T}\\
\frac{1}{2}\boldsymbol{H}\mathbf{N}^{T}
\end{array}\right]\bar{\mbox{j}}\mathrm{\,d}\mathrm{\Gamma}\nonumber \\
 &  & +\int_{\Gamma^{*}}\left[\begin{array}{c}
\mathbf{0}\\
2\mathbf{N}^{T}
\end{array}\right]<\sigma(\left\{ u_{k}\right\} )>.\mathbf{n}\mathrm{d\Gamma}+\int_{\Omega}\left[\begin{array}{c}
\mathbf{N}^{T}\\
\tilde{\boldsymbol{H}}\mathbf{N}^{T}
\end{array}\right]\boldsymbol{f}\,\mathrm{d\Omega}+\int_{\Gamma}\left[\begin{array}{c}
\mathbf{N}^{T}\\
\tilde{\boldsymbol{H}}\mathbf{N}^{T}
\end{array}\right]\boldsymbol{t}\,\mbox{d}\Gamma\nonumber \\
 &  & -\int_{\Gamma^{*}}\left[\begin{array}{c}
\mathbf{0}\\
2\mathbf{N}^{T}
\end{array}\right]\alpha\left(\left[\begin{array}{cc}
\mathbf{0} & 2\mathbf{N}\end{array}\right]\left\{ u_{k}\right\} -\bar{\mbox{i}}\right)\mathrm{d\Gamma}\nonumber \\
 &  & -\int_{\Gamma^{*}}\left[\begin{array}{c}
\mathbf{B}^{T}\\
\frac{1}{2}\left[\tilde{\boldsymbol{H}}^{+}+\tilde{\boldsymbol{H}}^{-}\right]^{T}\mathbf{B}^{T}
\end{array}\right]\mathbf{D_{T}}\boldsymbol{n}^{T}\left(\left[\begin{array}{cc}
\mathbf{0} & 2\mathbf{N}\end{array}\right]\left\{ u_{k}\right\} -\bar{\mbox{i}}\right)\mathrm{d\Gamma}
\end{eqnarray}

\subsubsection{Dirichlet Conditions}

Now, considering Dirichlet conditions (\ref{eq:-23}):
\begin{eqnarray}
R\left(\left\{ u\right\} \right) & = & \int_{\Omega}\sigma\left(\left\{ u\right\} \right):\varepsilon\left(v\right)\mbox{d}\Omega-\int_{\Omega}v^{T}\boldsymbol{f}\,\mbox{d}\Omega-\int_{\Gamma}v^{T}\boldsymbol{t}\,\mbox{d}\Gamma+\int_{\Gamma^{*}}\left(v^{+}\right)^{T}\alpha^{+}\left(\left\{ u^{+}\right\} -g^{+}\right)\\
 &  & +\int_{\Gamma^{*}}\left(v^{-}\right)^{T}\beta\alpha^{-}\left(\left\{ u^{-}\right\} -g^{-}\right)-\int_{\Gamma^{*}}\left(v^{+}\right)^{T}\sigma^{+}(\left\{ u\right\} ).\mathbf{n}\mathrm{d\Gamma}+\int_{\Gamma^{*}}\left(v^{-}\right)^{T}\sigma^{-}(\left\{ u\right\} ).\mathbf{n}\mathrm{d\Gamma}\nonumber 
\end{eqnarray}
We discretize the system to obtain
\begin{eqnarray}
R\left(\left\{ u\right\} \right) & = & \int_{\Omega}\left[\begin{array}{c}
\mathbf{B}^{T}\\
\tilde{\boldsymbol{H}}^{T}\mathbf{B}^{T}
\end{array}\right]\sigma\left(\left\{ u\right\} \right)\mbox{d}\Omega-\int_{\Omega}\left[\begin{array}{c}
\mathbf{N}^{T}\\
\tilde{\boldsymbol{H}}^{T}\mathbf{N}^{T}
\end{array}\right]\boldsymbol{f}\,\mbox{d}\Omega-\int_{\Gamma}\left[\begin{array}{c}
\mathbf{N}^{T}\\
\tilde{\boldsymbol{H}}^{T}\mathbf{N}^{T}
\end{array}\right]\boldsymbol{t}\,\mbox{d}\Gamma\nonumber \\
 &  & +\int_{\Gamma^{*}}\left[\begin{array}{c}
\mathbf{N}^{T}\\
\left[\tilde{\boldsymbol{H}}^{+}\right]^{T}\mathbf{N}^{T}
\end{array}\right]\alpha^{+}\left(\left[\begin{array}{cc}
\mathbf{N} & \mathbf{N}\tilde{\boldsymbol{H}}^{+}\end{array}\right]\left\{ u\right\} -g^{+}\right)\mathrm{d\Gamma}\nonumber \\
 &  & +\int_{\Gamma^{*}}\left[\begin{array}{c}
\mathbf{N}^{T}\\
\left[\tilde{\boldsymbol{H}}^{-}\right]^{T}\mathbf{N}^{T}
\end{array}\right]\alpha^{-}\left(\left[\begin{array}{cc}
\mathbf{N} & \mathbf{N}\tilde{\boldsymbol{H}}^{-}\end{array}\right]\left\{ u\right\} -g^{-}\right)\mathrm{d\Gamma}\nonumber \\
 &  & -\int_{\Gamma^{*}}\left[\begin{array}{c}
\mathbf{N}^{T}\\
\left[\tilde{\boldsymbol{H}}^{+}\right]^{T}\mathbf{N}^{T}
\end{array}\right]\sigma^{+}(\left\{ u\right\} ).\mathbf{n}\mathrm{d\Gamma}+\int_{\Gamma^{*}}\left[\begin{array}{c}
\mathbf{N}^{T}\\
\left[\tilde{\boldsymbol{H}}^{-}\right]^{T}\mathbf{N}^{T}
\end{array}\right]\sigma^{-}(\left\{ u\right\} ).\mathbf{n}\mathrm{d\Gamma}
\end{eqnarray}
Similar to what was done in the case of jump in displacement (\ref{eq:-32}),
we differentiate with respect to the nodal displacement.
\begin{eqnarray}
\frac{\partial R}{\partial\left\{ u\right\} }\left\{ u_{k}\right\}  & = & \int_{\Omega}\left[\begin{array}{c}
\mathbf{B}^{T}\\
\tilde{\boldsymbol{H}^{T}}\mathbf{B}^{T}
\end{array}\right]\frac{\partial\sigma}{\partial\varepsilon}\frac{\partial\varepsilon\left(\left\{ u\right\} \right)}{\partial\left\{ u\right\} }\mbox{d}\Omega+\int_{\Gamma^{*}}\left[\begin{array}{c}
\mathbf{N}^{T}\\
\left[\tilde{\boldsymbol{H}}^{+}\right]^{T}\mathbf{N}^{T}
\end{array}\right]\alpha^{+}\left[\begin{array}{cc}
\mathbf{N} & \mathbf{N}\tilde{\boldsymbol{H}}^{+}\end{array}\right]\mathrm{d\Gamma}\nonumber \\
 &  & +\int_{\Gamma^{*}}\left[\begin{array}{c}
\mathbf{N}^{T}\\
\left[\tilde{\boldsymbol{H}}^{-}\right]^{T}\mathbf{N}^{T}
\end{array}\right]\alpha^{-}\left[\begin{array}{cc}
\mathbf{N} & \mathbf{N}\tilde{\boldsymbol{H}}^{-}\end{array}\right]\mathrm{d\Gamma}-\int_{\Gamma^{*}}\left[\begin{array}{c}
\mathbf{N}^{T}\\
\left[\tilde{\boldsymbol{H}}^{+}\right]^{T}\mathbf{N}^{T}
\end{array}\right]\boldsymbol{n}\frac{\partial\sigma}{\partial\varepsilon}\frac{\partial\varepsilon^{+}\left(\left\{ u\right\} \right)}{\partial\left\{ u\right\} }\mathrm{d\Gamma}\nonumber \\
 &  & +\int_{\Gamma^{*}}\left[\begin{array}{c}
\mathbf{N}^{T}\\
\left[\tilde{\boldsymbol{H}}^{-}\right]^{T}\mathbf{N}^{T}
\end{array}\right]\boldsymbol{n}\frac{\partial\sigma}{\partial\varepsilon}\frac{\partial\varepsilon^{-}\left(\left\{ u\right\} \right)}{\partial\left\{ u\right\} }\mathrm{d\Gamma}
\end{eqnarray}
\begin{eqnarray}
\frac{\partial R}{\partial\left\{ u\right\} }\left\{ u_{k}\right\}  & = & \int_{\Omega}\left[\begin{array}{c}
\mathbf{B}^{T}\\
\tilde{\boldsymbol{H}}^{T}\mathbf{B}^{T}
\end{array}\right]\mathbf{D_{T}}\left[\begin{array}{cc}
\mathbf{B} & \mathbf{B}\tilde{\boldsymbol{H}}\end{array}\right]\mbox{d}\Omega+\int_{\Gamma^{*}}\left[\begin{array}{c}
\mathbf{N}^{T}\\
\left[\tilde{\boldsymbol{H}}^{+}\right]^{T}\mathbf{N}^{T}
\end{array}\right]\alpha^{+}\left[\begin{array}{cc}
\mathbf{N} & \mathbf{N}\tilde{\boldsymbol{H}}^{+}\end{array}\right]\mathrm{d\Gamma}\nonumber \\
 &  & +\int_{\Gamma^{*}}\left[\begin{array}{c}
\mathbf{N}^{T}\\
\left[\tilde{\boldsymbol{H}}^{-}\right]^{T}\mathbf{N}^{T}
\end{array}\right]\alpha^{-}\left[\begin{array}{cc}
\mathbf{N} & \mathbf{N}\tilde{\boldsymbol{H}}^{-}\end{array}\right]\mathrm{d\Gamma}-\int_{\Gamma^{*}}\left[\begin{array}{c}
\mathbf{N}^{T}\\
\left[\tilde{\boldsymbol{H}}^{+}\right]^{T}\mathbf{N}^{T}
\end{array}\right]\boldsymbol{n}\mathbf{D_{T}}\left[\begin{array}{cc}
\mathbf{B} & \mathbf{B}\tilde{\boldsymbol{H}}^{+}\end{array}\right]\mathrm{d\Gamma}\nonumber \\
 &  & +\int_{\Gamma^{*}}\left[\begin{array}{c}
\mathbf{N}^{T}\\
\left[\tilde{\boldsymbol{H}}^{-}\right]^{T}\mathbf{N}^{T}
\end{array}\right]\boldsymbol{n}\mathbf{D_{T}\frac{\partial\sigma}{\partial\varepsilon}}\left[\begin{array}{cc}
\mathbf{B} & \mathbf{B}\tilde{\boldsymbol{H}}^{-}\end{array}\right]\mathrm{d\Gamma}
\end{eqnarray}
To obtain a symmetric form of the matrix in the system of equations,
we introduce Nitsche's terms in the system,
\begin{eqnarray}
 &  & \left(\frac{\partial R}{\partial\left\{ u\right\} }\left\{ u_{k}\right\} \right)\left\{ \delta u\right\} \nonumber \\
 &  & -\left(\int_{\Gamma^{*}}\left[\begin{array}{c}
\mathbf{B}^{T}\\
\left[\tilde{\boldsymbol{H}}^{+}\right]^{T}\mathbf{B}^{T}
\end{array}\right]\mathbf{D_{T}}\boldsymbol{n}^{T}\left[\begin{array}{cc}
\mathbf{N} & \mathbf{N}\tilde{\boldsymbol{H}}^{+}\end{array}\right]\mathrm{d\Gamma}-\int_{\Gamma^{*}}\left[\begin{array}{c}
\mathbf{B}^{T}\\
\left[\tilde{\boldsymbol{H}}^{-}\right]^{T}\mathbf{B}^{T}
\end{array}\right]\mathbf{D_{T}}\boldsymbol{n}^{T}\left[\begin{array}{cc}
\mathbf{N} & \mathbf{N}\tilde{\boldsymbol{H}}^{-}\end{array}\right]\mathrm{d\Gamma}\right)\left\{ u_{k+1}\right\} \nonumber \\
 & = & -R\left(\left\{ u_{k}\right\} \right)-\int_{\Gamma^{*}}\left[\begin{array}{c}
\mathbf{B}^{T}\\
\left[\tilde{\boldsymbol{H}}^{+}\right]^{T}\mathbf{B}^{T}
\end{array}\right]\mathbf{D_{T}}\boldsymbol{n}^{T}g^{+}\mathrm{d\Gamma}+\int_{\Gamma^{*}}\left[\begin{array}{c}
\mathbf{B}^{T}\\
\left[\tilde{\boldsymbol{H}}^{-}\right]^{T}\mathbf{B}^{T}
\end{array}\right]\mathbf{D_{T}}\boldsymbol{n}^{T}g^{-}\mathrm{d\Gamma}
\end{eqnarray}
where $u_{k+1}=u_{k}+\delta u$ and since $u^{+}=g^{+},\,u^{-}=g^{-}\mbox{ on }\Gamma^{*}$,
we obtain 
\begin{eqnarray}
 &  & \left(\frac{\partial R}{\partial\left\{ u\right\} }\left\{ u_{k}\right\} \right)\left\{ \delta u\right\} \nonumber \\
 &  & -\left(\int_{\Gamma^{*}}\left[\begin{array}{c}
\mathbf{B}^{T}\\
\left[\tilde{\boldsymbol{H}}^{+}\right]^{T}\mathbf{B}^{T}
\end{array}\right]\mathbf{D_{T}}\boldsymbol{n}^{T}\left[\begin{array}{cc}
\mathbf{N} & \mathbf{N}\tilde{\boldsymbol{H}}^{+}\end{array}\right]\mathrm{d\Gamma}-\int_{\Gamma^{*}}\left[\begin{array}{c}
\mathbf{B}^{T}\\
\left[\tilde{\boldsymbol{H}}^{-}\right]^{T}\mathbf{B}^{T}
\end{array}\right]\mathbf{D_{T}}\boldsymbol{n}^{T}\left[\begin{array}{cc}
\mathbf{N} & \mathbf{N}\tilde{\boldsymbol{H}}^{-}\end{array}\right]\mathrm{d\Gamma}\right)\left\{ \delta u\right\} \nonumber \\
 & = & -R\left(\left\{ u_{k}\right\} \right)-\int_{\Gamma^{*}}\left[\begin{array}{c}
\mathbf{B}^{T}\\
\left[\tilde{\boldsymbol{H}}^{+}\right]^{T}\mathbf{B}^{T}
\end{array}\right]\mathbf{D_{T}}\boldsymbol{n}^{T}\left(\left[\begin{array}{cc}
\mathbf{N} & \mathbf{N}\tilde{\boldsymbol{H}}^{+}\end{array}\right]\left\{ u_{k}\right\} -g^{+}\right)\mathrm{d\Gamma}\nonumber \\
 &  & +\int_{\Gamma^{*}}\left[\begin{array}{c}
\mathbf{B}^{T}\\
\left[\tilde{\boldsymbol{H}}^{-}\right]^{T}\mathbf{B}^{T}
\end{array}\right]\mathbf{D_{T}}\boldsymbol{n}^{T}\left(\left[\begin{array}{cc}
\mathbf{N} & \mathbf{N}\tilde{\boldsymbol{H}}^{-}\end{array}\right]\left\{ u_{k}\right\} -g^{-}\right)\mathrm{d\Gamma}
\end{eqnarray}
which can be written as
\[
\left[\mathbf{K_{T}}\right]\left\{ \delta u\right\} =-R\left(\left\{ u_{k}\right\} \right)
\]
with
\begin{eqnarray}
\left[\mathbf{K_{T}}\right] & = & \int_{\Omega}\left[\begin{array}{c}
\mathbf{B}^{T}\\
\tilde{\boldsymbol{H}}^{T}\mathbf{B}^{T}
\end{array}\right]\mathbf{D_{T}}\left[\begin{array}{cc}
\mathbf{B} & \mathbf{B}\tilde{\boldsymbol{H}}\end{array}\right]\mbox{d}\Omega+\int_{\Gamma^{*}}\left[\begin{array}{c}
\mathbf{N}^{T}\\
\left[\tilde{\boldsymbol{H}}^{+}\right]^{T}\mathbf{N}^{T}
\end{array}\right]\alpha^{+}\left[\begin{array}{cc}
\mathbf{N} & \mathbf{N}\tilde{\boldsymbol{H}}^{+}\end{array}\right]\mathrm{d\Gamma}\nonumber \\
 &  & -\int_{\Gamma^{*}}\left[\begin{array}{c}
\mathbf{B}^{T}\\
\left[\tilde{\boldsymbol{H}}^{+}\right]^{T}\mathbf{B}^{T}
\end{array}\right]\mathbf{D_{T}}\boldsymbol{n}^{T}\left[\begin{array}{cc}
\mathbf{N} & \mathbf{N}\tilde{\boldsymbol{H}}^{+}\end{array}\right]\mathrm{d\Gamma}+\int_{\Gamma^{*}}\left[\begin{array}{c}
\mathbf{B}^{T}\\
\left[\tilde{\boldsymbol{H}}^{-}\right]^{T}\mathbf{B}^{T}
\end{array}\right]\mathbf{D_{T}}\boldsymbol{n}^{T}\left[\begin{array}{cc}
\mathbf{N} & \mathbf{N}\tilde{\boldsymbol{H}}^{-}\end{array}\right]\mathrm{d\Gamma}\nonumber \\
 &  & +\int_{\Gamma^{*}}\left[\begin{array}{c}
\mathbf{N}^{T}\\
\left[\tilde{\boldsymbol{H}}^{-}\right]^{T}\mathbf{N}^{T}
\end{array}\right]\alpha^{-}\left[\begin{array}{cc}
\mathbf{N} & \mathbf{N}\tilde{\boldsymbol{H}}^{-}\end{array}\right]\mathrm{d\Gamma}-\int_{\Gamma^{*}}\left[\begin{array}{c}
\mathbf{N}^{T}\\
\left[\tilde{\boldsymbol{H}}^{+}\right]^{T}\mathbf{N}^{T}
\end{array}\right]\boldsymbol{n}\mathbf{D_{T}}\left[\begin{array}{cc}
\mathbf{B} & \mathbf{B}\tilde{\boldsymbol{H}}^{+}\end{array}\right]\mathrm{d\Gamma}\nonumber \\
 &  & +\int_{\Gamma^{*}}\left[\begin{array}{c}
\mathbf{N}^{T}\\
\left[\tilde{\boldsymbol{H}}^{-}\right]^{T}\mathbf{N}^{T}
\end{array}\right]\boldsymbol{n}\mathbf{D_{T}\frac{\partial\sigma}{\partial\varepsilon}}\left[\begin{array}{cc}
\mathbf{B} & \mathbf{B}\tilde{\boldsymbol{H}}^{-}\end{array}\right]\mathrm{d\Gamma}
\end{eqnarray}
\begin{eqnarray}
R\left(\left\{ u_{k}\right\} \right) & = & \int_{\Omega}\left[\begin{array}{c}
\mathbf{B}^{T}\\
\tilde{\boldsymbol{H}}^{T}\mathbf{B}^{T}
\end{array}\right]\sigma\left(\left\{ u_{k}\right\} \right)\mbox{d}\Omega-\int_{\Omega}\left[\begin{array}{c}
\mathbf{N}^{T}\\
\tilde{\boldsymbol{H}}^{T}\mathbf{N}^{T}
\end{array}\right]\boldsymbol{f}\,\mbox{d}\Omega-\int_{\Gamma}\left[\begin{array}{c}
\mathbf{N}^{T}\\
\tilde{\boldsymbol{H}}^{T}\mathbf{N}^{T}
\end{array}\right]\boldsymbol{t}\,\mbox{d}\Gamma\nonumber \\
 &  & +\int_{\Gamma^{*}}\left[\begin{array}{c}
\mathbf{N}^{T}\\
\left[\tilde{\boldsymbol{H}}^{+}\right]^{T}\mathbf{N}^{T}
\end{array}\right]\alpha^{+}\left(\left[\begin{array}{cc}
\mathbf{N} & \mathbf{N}\tilde{\boldsymbol{H}}^{+}\end{array}\right]\left\{ u_{k}\right\} -g^{+}\right)\mathrm{d\Gamma}\nonumber \\
 &  & +\int_{\Gamma^{*}}\left[\begin{array}{c}
\mathbf{N}^{T}\\
\left[\tilde{\boldsymbol{H}}^{-}\right]^{T}\mathbf{N}^{T}
\end{array}\right]\alpha^{-}\left(\left[\begin{array}{cc}
\mathbf{N} & \mathbf{N}\tilde{\boldsymbol{H}}^{-}\end{array}\right]\left\{ u_{k}\right\} -g^{-}\right)\mathrm{d\Gamma}\nonumber \\
 &  & -\int_{\Gamma^{*}}\left[\begin{array}{c}
\mathbf{N}^{T}\\
\left[\tilde{\boldsymbol{H}}^{+}\right]^{T}\mathbf{N}^{T}
\end{array}\right]\sigma^{+}(\left\{ u_{k}\right\} ).\mathbf{n}\mathrm{d\Gamma}+\int_{\Gamma^{*}}\left[\begin{array}{c}
\mathbf{N}^{T}\\
\left[\tilde{\boldsymbol{H}}^{-}\right]^{T}\mathbf{N}^{T}
\end{array}\right]\sigma^{-}(\left\{ u_{k}\right\} ).\mathbf{n}\mathrm{d\Gamma}\nonumber \\
 &  & -\int_{\Gamma^{*}}\left[\begin{array}{c}
\mathbf{B}^{T}\\
\left[\tilde{\boldsymbol{H}}^{+}\right]^{T}\mathbf{B}^{T}
\end{array}\right]\mathbf{D_{T}}\boldsymbol{n}^{T}\left(\left[\begin{array}{cc}
\mathbf{N} & \mathbf{N}\tilde{\boldsymbol{H}}^{+}\end{array}\right]\left\{ u_{k}\right\} -g^{+}\right)\mathrm{d\Gamma}\nonumber \\
 &  & +\int_{\Gamma^{*}}\left[\begin{array}{c}
\mathbf{B}^{T}\\
\left[\tilde{\boldsymbol{H}}^{-}\right]^{T}\mathbf{B}^{T}
\end{array}\right]\mathbf{D_{T}}\boldsymbol{n}^{T}\left(\left[\begin{array}{cc}
\mathbf{N} & \mathbf{N}\tilde{\boldsymbol{H}}^{-}\end{array}\right]\left\{ u_{k}\right\} -g^{-}\right)\mathrm{d\Gamma}
\end{eqnarray}

\subsection{\label{subsec:Weighted-Discretization}Weighted Discretization}

We continue the XFEM discretization with a novel weighting for the
interfacial consistency terms arising in Nitsche's variational form.
We recollect the part from section (\ref{par:Jump-Condition}) regarding
jump condition.\cite{key-4} We now use a weighted approach with a
weight $0\leq\gamma\leq1$.
\[
\mbox{div}\sigma^{+}=f^{+}\mbox{ in }\Omega^{+}
\]
\[
\mbox{div}\sigma^{-}=f^{-}\mbox{ in }\Omega^{-}
\]
with:
\[
u=u_{0}\,\,\,\,\mathrm{on}\,\Gamma=\partial\Omega
\]
\begin{eqnarray*}
\left[\left[u\right]\right] & = & \bar{\mathrm{i}}\,\,\,\,\mathrm{on}\,\mathcal{S}\\
u^{+}-u^{-} & = & \bar{\mathrm{i}}
\end{eqnarray*}
and:
\begin{eqnarray*}
\left[\left[\sigma\right]\right]\centerdot\mathbf{n} & = & \bar{\mathrm{j}}\,\,\,\,\mathrm{on}\,\mathcal{S}\\
(\sigma^{+}-\sigma^{-}).\mathbf{n} & = & \bar{\mathrm{j}}
\end{eqnarray*}
where $\mathbf{n}$ is the normal vector pointing outwards of $\Omega^{-}$.
We calculate the averages of the flux and the displacement as:
\begin{equation}
\left\langle \sigma\right\rangle _{\gamma}.\mathbf{n}=\gamma\sigma^{+}.\mathbf{n}+(1-\gamma)\sigma^{-}.\mathbf{n}
\end{equation}
and: 
\[
\left\langle v\right\rangle =\gamma v^{+}+(1-\gamma)v^{-}
\]

From the weak Galerkin formulation of our problem statement we can
write: find $u\in\mathbb{U}$ such that:
\begin{equation}
\int_{\Omega}\varepsilon(v)\sigma\mbox{d}\Omega-\int_{\mathcal{S}}\left[v^{+}(\sigma^{+}.\mathbf{n})-v^{-}(\sigma^{-}.\mathbf{n})\right]\mbox{d}\Gamma=\int_{\Omega}vf\mbox{d}\Omega\mbox{ }\forall v\in\mathbb{U}_{0}\label{eq:-31}
\end{equation}
We can write:
\[
\gamma\sigma^{+}.\mathbf{n}-\gamma\sigma^{-}.\mathbf{n}=\gamma\bar{\mathrm{j}}
\]
which gives us:
\[
\left\langle \sigma\right\rangle _{\gamma}.\mathbf{n}-\gamma\bar{\mathrm{j}}=\cancel{\gamma\sigma^{+}.\mathbf{n}}+(1-\cancel{\gamma})\sigma^{-}.\mathbf{n}-\cancel{\gamma\sigma^{+}.\mathbf{n}}+\cancel{\gamma\sigma^{-}.\mathbf{n}}
\]
\[
\sigma^{-}.\mathbf{n}=\left\langle \sigma\right\rangle _{\gamma}.\mathbf{n}-\gamma\bar{\mathrm{j}}
\]
Similarly:
\[
\sigma^{+}.\mathbf{n}=\left\langle \sigma\right\rangle _{\gamma}.\mathbf{n}+(1-\gamma)\bar{\mathrm{j}}
\]
We can put these two results in equation (\ref{eq:-31}) to obtain:
\[
\int_{\Omega}\varepsilon(v)\sigma\mbox{d}\Omega-\int_{\mathcal{S}}\left[v^{+}(\left\langle \sigma\right\rangle _{\gamma}.\mathbf{n}+(1-\gamma)\bar{\mathrm{j}})-v^{-}(\left\langle \sigma\right\rangle _{\gamma}.\mathbf{n}-\gamma\bar{\mathrm{j}})\right]\mbox{d}\Gamma=\int_{\Omega}vf\,\mbox{d}\Omega
\]
\[
\int_{\Omega}\varepsilon(v)\sigma\mbox{d}\Omega-\int_{\mathcal{S}}\left[v^{+}\left\langle \sigma\right\rangle _{\gamma}.\mathbf{n}+(1-\gamma)v^{+}\bar{\mathrm{j}}-v^{-}\left\langle \sigma\right\rangle _{\gamma}.\mathbf{n}+\gamma v^{-}\bar{\mathrm{j}}\right]\mbox{d}\Gamma=\int_{\Omega}vf\mbox{\,d}\Omega
\]
\[
\int_{\Omega}\varepsilon(v)\sigma\mbox{d}\Omega-\int_{\mathcal{S}}\left[(v^{+}-v^{-})\left\langle \sigma\right\rangle _{\gamma}.\mathbf{n}+((1-\gamma)v^{+}+\gamma v^{-})\bar{\mathrm{j}}\right]\mbox{d}\Gamma=\int_{\Omega}vf\,\mbox{d}\Omega
\]
\[
\int_{\Omega}\varepsilon(v)\sigma\mbox{d}\Omega-\int_{\mathcal{S}}\left[[[v]]\left\langle \sigma\right\rangle _{\gamma}.\mathbf{n}+\left\langle v\right\rangle _{1-\gamma}\bar{\mathrm{j}}\right]\mbox{d}\Gamma=\int_{\Omega}vf\,\mbox{d}\Omega
\]
 Adding Nitsche's terms, stabilization terms and terms for variational
consistency and symmetry, we have the weighted Nitsche's formulation
for jump conditions:
\begin{eqnarray}
\int_{\Omega}\varepsilon(v)\sigma\mbox{d}\Omega-\int_{\mathcal{S}}[[v]]\left\langle \sigma\right\rangle _{\gamma}.\mathbf{n}\mbox{d}\Gamma-\int_{\mathcal{S}}\left\langle \sigma(v)\right\rangle _{\gamma}.\mathbf{n}[[u]]\mbox{d}\Gamma+\int_{\mathcal{S}}[[v]]\epsilon[[u]]\mbox{d}\Gamma\nonumber \\
=\int_{\Omega}vf\mbox{d}\Omega+\int_{\mathcal{S}}\left\langle v\right\rangle _{1-\gamma}\bar{\mathrm{j}}\mbox{d}\Gamma-\int_{\mathcal{S}}\left\langle \sigma\right\rangle _{\gamma}.\mathbf{n}\bar{\mathrm{i}}\mbox{d}\Gamma+\int_{\mathcal{S}}[[v]]\epsilon\bar{\mathrm{i}}\mbox{d}\Gamma
\end{eqnarray}
As we can see here, the weighting terms influence Nitsche's terms
of the formulation as well as the term corresponding to the jump in
the flux. This can let us safely conclude that the discretized system
with shifted basis enrichment will be as follows.{\footnotesize{}
\begin{equation}
\left[\begin{array}{cc}
\mathbf{K}_{b} & \mathbf{K}_{b}\tilde{\boldsymbol{H}}-2(\mathbf{G}\mathbf{K}_{n}^{T+}+\tilde{\mathbf{G}}\mathbf{K}_{n}^{T-})\\
\\
 & \tilde{\boldsymbol{H}}^{T}\mathbf{K}_{b}\tilde{\boldsymbol{H}}_{j}+4\mathbf{K}_{s}\\
\tilde{\boldsymbol{H}}^{T}\mathbf{K}_{b} & -2(\mathbf{G}\mathbf{K}_{n}^{+}\tilde{\boldsymbol{H}}^{+}+\tilde{\mathbf{G}}\mathbf{K}_{n}^{-}\tilde{\boldsymbol{H}}^{-})\\
-2(\mathbf{G}\mathbf{K}_{n}^{+}+\tilde{\mathbf{G}}\mathbf{K}_{n}^{-}) & -2\left(\mathbf{G}\left[\tilde{\boldsymbol{H}}^{+}\right]^{T}\mathbf{K}_{n}^{T+}-\tilde{\mathbf{G}}\left[\tilde{\boldsymbol{H}}^{-}\right]^{T}\mathbf{K}_{n}^{T-}\right)
\end{array}\right]\left\{ \begin{array}{c}
\boldsymbol{u}_{j}\\
\boldsymbol{a}_{j}
\end{array}\right\} =\left\{ \begin{array}{c}
\mathbf{f}_{b}+\mathbf{f}_{h}-(\mathbf{G}\mathbf{f}_{n}^{+}+\tilde{\mathbf{G}}\mathbf{f}_{n}^{-})+\mathbf{f}_{\mbox{j}}\\
\\
\tilde{\boldsymbol{H}}(\mathbf{f}_{b}+\mathbf{f}_{h})+2\mathbf{f}_{s}\\
-\left(\mathbf{G}\left[\tilde{\boldsymbol{H}}^{+}\right]^{T}\mathbf{f}_{n}^{+}+\tilde{\mathbf{G}}\left[\tilde{\boldsymbol{H}}^{-}\right]^{T}\mathbf{f}_{n}^{-}\right)\\
+\left(\mathbf{G}\left[\tilde{\boldsymbol{H}}_{i}^{+}\right]^{T}+\tilde{\mathbf{G}}\left[\tilde{\boldsymbol{H}}^{-}\right]^{T}\right)\mathbf{f}_{\mbox{j}}
\end{array}\right\} 
\end{equation}
}where $\mathbf{G}$ and $\tilde{\mathbf{G}}$ are the diagonal matrices
with the corresponding values of the weights ($\gamma$ and $1-\gamma$)
for that element. We can see here that we will recover the system
(\ref{eq:-19-1}) if we chose $\gamma=\frac{1}{2}$ and this is indeed
the classical Nitsche's algorithm. Following the works of Annavarapu
et al., we implement the weights as follows for an element $e$
\begin{equation}
\gamma_{e}=\frac{\mbox{meas}(\Omega_{e}^{+})/\left|\mathbf{D}^{+}\right|}{\mbox{meas}(\Omega_{e}^{+})/\left|\mathbf{D}^{+}\right|+\mbox{meas}(\Omega_{e}^{-})/\left|\mathbf{D}^{-}\right|}\label{eq:-40}
\end{equation}
\[
1-\gamma_{e}=\frac{\mbox{meas}(\Omega_{e}^{-})/\left|\mathbf{D}^{-}\right|}{\mbox{meas}(\Omega_{e}^{+})/\left|\mathbf{D}^{+}\right|+\mbox{meas}(\Omega_{e}^{-})/\left|\mathbf{D}^{-}\right|}
\]

\newpage{}

\section{Analytical Solution of a simple problem}

\begin{figure}
\begin{centering}
\includegraphics[bb=350bp 500bp 575bp 710bp,clip,height=6cm]{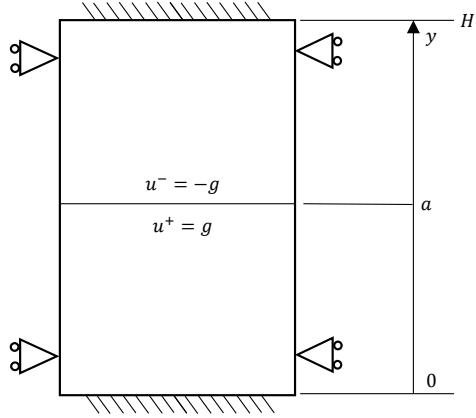}
\par\end{centering}
\caption{\label{fig:Discretization-of-the}simple problem}
\end{figure}

Consider the problem given in figure (\ref{fig:Discretization-of-the})

If we consider there is a jump in the displacement at the interface,
we can define the given problem by:
\[
\boldsymbol{u}(y=0)=\boldsymbol{0}
\]
\[
\boldsymbol{u}(y=H)=\boldsymbol{0}
\]
\[
\mbox{i}(y=a)=[[u^{+}-u^{-}]]=2g\boldsymbol{e_{y}}\,\,\,\,\forall g\in\mathbb{R}
\]
\[
\mbox{j}(y=a)=[[\underline{\underline{\sigma}}(y=a)]].n_{y}=\boldsymbol{0}
\]

\subsection{Exact numerical solution}

Here the bar is divided by an interface at $y=a$. All horizontal
movements are restricted, thus letting us simplify the given problem
into a 1D bar problem. The entire bar is made of a single material
with Young's modulus $E$. We assume $\nu$, Poisson's ratio, to be
zero. From equilibrium equations, we have, in the domain $\Omega^{+}$,
\[
\underbar{\mbox{div}}\,\underline{\underline{\sigma}}=\boldsymbol{0}
\]
Integrating the equilibrium equation gives:
\[
\underline{\underline{\sigma}}=A\boldsymbol{e_{y}}\otimes\boldsymbol{e_{y}}
\]
or:
\begin{eqnarray*}
\sigma & = & A\\
E\varepsilon & = & A
\end{eqnarray*}
Integrating again gives:
\[
u=\frac{1}{E}\left(Ay+B\right)
\]
From boundary conditions we have:
\begin{eqnarray*}
u(y=H) & = & 0\\
0 & = & AH+B
\end{eqnarray*}
and also:
\[
u(y=a^{+})=\frac{1}{E}\left(Aa+B\right)
\]
In the domain $\Omega^{-}$, we have:
\[
\underbar{\mbox{div}}\,\underline{\underline{\sigma}}=\boldsymbol{0}
\]
\[
\underline{\underline{\sigma}}=C
\]
\[
u=\frac{1}{E}\left(Cy+D\right)
\]
Applying the boundary conditions gives:
\begin{eqnarray*}
u\left(y=0\right) & = & 0\\
0 & = & D
\end{eqnarray*}
and also:
\[
u(y=a^{-})=\frac{1}{E}\left(Ca\right)
\]
From the condition of jump in the stress at the interface, we have:
\[
\mbox{j}(y=a)=A-C=0
\]
and jump in displacement:
\begin{eqnarray*}
\mbox{i}(y=a) & = & \frac{1}{E}\left(Aa+B\right)-\frac{1}{E}\left(Ca\right)\\
2gE & = & Aa+B-Ca
\end{eqnarray*}
which gives:
\[
B=2Eg
\]
\[
A=-\frac{2E}{H}g
\]
\[
C=-\frac{2E}{H}g
\]
and: 
\[
D=0
\]
Thus, we have:
\[
\sigma=-\frac{2E}{H}g\,\,\,\,\mbox{in}\,\Omega^{+}
\]
\[
u=-\frac{2g}{H}y+2g\,\,\,\,\mbox{in}\,\Omega^{+}
\]
\[
\sigma=-\frac{2E}{H}g\,\,\,\,\mbox{in}\,\Omega^{-}
\]
\[
u=-\frac{2g}{H}y\,\,\,\,\mbox{in}\,\Omega^{-}
\]

\subsection{Discretized solution}

\begin{figure}
\begin{centering}
\includegraphics[bb=200bp 580bp 435bp 700bp,clip]{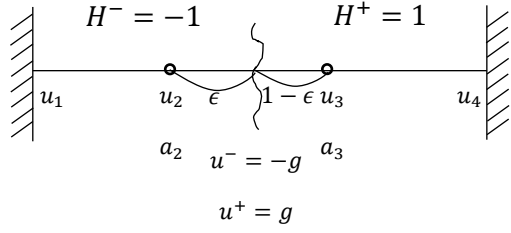}
\par\end{centering}
\caption{\label{fig:Discretization-of-the-1}simple problem discretized}
\end{figure}

Discretizing the above problem by XFEM gives us the system as in figure
(\ref{fig:Discretization-of-the-1}).We now divide the bar into three
elements of length $h$. The total length of the bar is $3h$. From
the analytic solution, we have: 
\[
u\left(y=h\right)=-\frac{2g}{3}
\]
\[
u\left(y=2h\right)=\frac{2g}{3}
\]
Since $a=h+\epsilon h$, we have:
\[
u\left(y=a^{-}\right)=-g
\]
\[
u\left(y=a^{+}\right)=g
\]
It is to note here that since the entire domain is of a single material
and the jump in stress is zero, the solution will be symmetric and
this is the same solution on applying a Dirichlet type of boundary
condition at the interface. 

On solving by the method of shifted enrichment X-FEM, we see that
the element 2, supported by nodes 2 and 3, includes the interface
and we consider the enrichment $H(\Omega^{-})=-1\mbox{ and }H(\Omega^{+})=+1$.
We consider the interface to be at a distance of $\epsilon h$ from
node 2 with $0\leqslant\epsilon\leqslant1$. We have:
\begin{eqnarray*}
u(x) & = & \cancel{u_{1}N_{1}(x)}+u_{2}N_{2}(x)+u_{3}N_{3}(x)+a_{2}\tilde{H}_{2}(x)N_{2}(x)+a_{3}\tilde{H}_{3}(x)N_{3}(x)+\cancel{u_{4}N_{4}(x)}\\
 & = & u_{2}N_{2}(x)+u_{3}N_{3}(x)+(H(x)+1)a_{2}N_{2}(x)+(H(x)-1)a_{3}N_{3}(x)\\
v(x) & = & v_{2}N_{2}(x)+v_{3}N_{3}(x)+(H(x)+1)b_{2}N_{2}(x)+(H(x)-1)b_{3}N_{3}(x)
\end{eqnarray*}
since $u_{1}=0$. Consider element 1 supported by nodes 1 and 2:
\[
E\int_{x_{1}}^{x_{2}}v_{,x}u_{,x}\mbox{d}x=E\int_{x_{1}}^{x_{2}}(v_{2}N_{2,x}+(-1+1)b_{2}N_{2,x})(u_{2}N_{2,x}+(-1+1)a_{2}N_{2,x})\mbox{d}x
\]
\begin{equation}
=E\int_{x_{1}}^{x_{2}}(v_{2}N_{2,x})(u_{2}N_{2,x})\mbox{d}x=\frac{E}{h}\left[1\right]\left\{ u_{2}\right\} \label{eq:-34}
\end{equation}
with the force contribution term equaling $0$. Similarly, for element
3 the contribution towards the element stiffness matrix is:
\[
\frac{E}{h}\left[1\right]\left\{ u_{3}\right\} 
\]
again with the force contribution term equaling $0$.

\subsubsection{Jump in Displacement type solution}

We apply Nitsche's method for jump in displacement boundary conditions
for element 2. We take:
\[
N_{2}=\frac{x_{3}-x}{x_{3}-x_{2}}
\]
\[
N_{3}=\frac{x-x_{2}}{x_{3}-x_{2}}
\]
\[
B_{2}=\frac{-1}{x_{3}-x_{2}}
\]
\[
B_{3}=\frac{1}{x_{3}-x_{2}}
\]
From the variational form for jump in displacement (\ref{eq:-24}),
\[
\int_{x_{2}}^{x_{3}}\varepsilon^{T}(\boldsymbol{v}_{h})\sigma(\boldsymbol{u}_{h})\mathrm{d\Omega}=E\left\{ \boldsymbol{v}\right\} ^{T}\left[\begin{array}{cc}
\mathbf{B} & \mathbf{B\tilde{H}}\end{array}\right]^{T}\left[\begin{array}{cc}
\mathbf{B} & \mathbf{B\tilde{H}}\end{array}\right]\left\{ \boldsymbol{u}\right\} 
\]
where $\mathbf{\tilde{H}}=\gamma\tilde{\boldsymbol{H}}^{+}+\left(1-\gamma\right)\tilde{\boldsymbol{H}}^{-}$,
$\gamma$ being a weighting parameter. The FEM element stiffness matrix
is given by:
\[
\boldsymbol{K}_{b}=\frac{E}{h}\left[\begin{array}{cc}
1 & -1\\
-1 & 1
\end{array}\right]
\]
The heavyside enrichment matrices for the element is given by:
\[
\boldsymbol{\tilde{H}}^{+}=\left[\begin{array}{cc}
(H(x)+1) & 0\\
0 & (H(x)-1)
\end{array}\right]=\left[\begin{array}{cc}
2 & 0\\
0 & 0
\end{array}\right]
\]
and: 
\[
\tilde{\boldsymbol{H}}^{-}=\left[\begin{array}{cc}
(H(x)+1) & 0\\
0 & (H(x)-1)
\end{array}\right]=\left[\begin{array}{cc}
0 & 0\\
0 & -2
\end{array}\right]
\]
By considering: 
\[
\boldsymbol{K}_{b}=\mathbf{B}^{T}E\mathbf{B}
\]
the XFEM element stiffness matrix is:
\begin{eqnarray*}
\mathbf{K}_{b} & = & \left[\begin{array}{cc}
\boldsymbol{K}_{b} & \epsilon\boldsymbol{K}_{b}\tilde{\boldsymbol{H}}^{-}+\left(1-\epsilon\right)\boldsymbol{K}_{b}\tilde{\boldsymbol{H}}^{+}\\
\left[\tilde{\boldsymbol{H}}^{-}\right]^{T}\epsilon\boldsymbol{K}_{b}+\left[\tilde{\boldsymbol{H}}^{+}\right]^{T}\left(1-\epsilon\right)\boldsymbol{K}_{b} & \left[\tilde{\boldsymbol{H}}^{-}\right]^{T}\epsilon\boldsymbol{K}_{b}\tilde{\boldsymbol{H}}^{-}+\left[\tilde{\boldsymbol{H}}^{+}\right]^{T}\left(1-\epsilon\right)\boldsymbol{K}_{b}\tilde{\boldsymbol{H}}^{+}
\end{array}\right]\\
 & = & \frac{E}{h}\left[\begin{array}{cccc}
1 & -1 & 2\left(1-\epsilon\right) & 2\epsilon\\
-1 & 1 & -2\left(1-\epsilon\right) & -2\epsilon\\
2\left(1-\epsilon\right) & -2\left(1-\epsilon\right) & 4\left(1-\epsilon\right) & 0\\
2\epsilon & -2\epsilon & 0 & 4\epsilon
\end{array}\right]
\end{eqnarray*}
where $\left(1-\epsilon\right)=\gamma$. For the stabilization part:
\begin{eqnarray*}
\left[\left[u\left(x=\left(1+\epsilon\right)h\right)\right]\right] & = & 2a_{2}N_{2}\left(x=\left(1+\epsilon\right)h\right)+2a_{3}N_{3}\left(x=\left(1+\epsilon\right)h\right)
\end{eqnarray*}
Thus:
\[
[[\boldsymbol{v}\left(x=\left(1+\epsilon\right)h\right)]]^{T}\alpha[[\boldsymbol{u}\left(x=\left(1+\epsilon\right)h\right)]]=\alpha\left\{ \boldsymbol{v}\right\} ^{T}\left[\begin{array}{c}
\mathbf{0}\\
2\mathbf{N}^{T}\left(x=\left(1+\epsilon\right)h\right)
\end{array}\right]\left[\begin{array}{cc}
\mathbf{0} & 2\mathbf{N}\left(x=\left(1+\epsilon\right)h\right)\end{array}\right]\left\{ \boldsymbol{u}\right\} 
\]
Let us look at:
\[
\mathbf{N}^{T}\mathbf{N}=\left[\begin{array}{cc}
\left(\frac{2h-h-h\epsilon}{h}\right)\left(\frac{2h-h-h\epsilon}{h}\right) & \left(\frac{2h-h-h\epsilon}{h}\right)\left(\frac{h-h+h\epsilon}{h}\right)\\
\left(\frac{h-h+h\epsilon}{h}\right)\left(\frac{2h-h-h\epsilon}{h}\right) & \left(\frac{h-h+h\epsilon}{h}\right)\left(\frac{h-h+h\epsilon}{h}\right)
\end{array}\right]
\]
If we consider:
\[
\boldsymbol{K}_{s}=\mathbf{N}^{T}\mathbf{N}=\left[\begin{array}{cc}
\left(1-\epsilon\right)^{2} & \left(1-\epsilon\right)\epsilon\\
\left(1-\epsilon\right)\epsilon & \epsilon^{2}
\end{array}\right]
\]
then the stabilization matrix is:
\begin{eqnarray*}
\mathbf{K}_{s} & = & \alpha\left[\begin{array}{cc}
\mathbf{0} & \mathbf{0}\\
\mathbf{0} & 4\boldsymbol{K}_{s}
\end{array}\right]\\
 & = & \alpha\left[\begin{array}{cccc}
0 & 0 & 0 & 0\\
0 & 0 & 0 & 0\\
0 & 0 & 4\left(1-\epsilon\right)^{2} & 4\left(1-\epsilon\right)\epsilon\\
0 & 0 & 4\left(1-\epsilon\right)\epsilon & 4\epsilon^{2}
\end{array}\right]
\end{eqnarray*}
We have two terms in Nitsche's matrix, the variational consistency
term:
\begin{eqnarray*}
[[\boldsymbol{v}\left(x=\left(1+\epsilon\right)h\right)]]^{T}<\sigma(\boldsymbol{u}\left(x=\left(1+\epsilon\right)h\right))>.\mathbf{n} & = & E\left\{ \boldsymbol{v}\right\} ^{T}\left[\begin{array}{c}
\mathbf{0}\\
2\mathbf{N}^{T}\left(x=\left(1+\epsilon\right)h\right)
\end{array}\right]\left[\begin{array}{cc}
\boldsymbol{n}^{T}\mathbf{B} & \boldsymbol{n}^{T}\mathbf{B\tilde{H}}\end{array}\right]\left\{ \boldsymbol{u}\right\} 
\end{eqnarray*}
and the symmetric term:
\[
(<\sigma(\boldsymbol{v}\left(x=\left(1+\epsilon\right)h\right))>.\mathbf{n})^{T}[[\boldsymbol{u}\left(x=\left(1+\epsilon\right)h\right)]]=E\left\{ \boldsymbol{v}\right\} ^{T}\left[\begin{array}{c}
\left[\boldsymbol{n}^{T}\mathbf{B}\right]^{T}\\
\left[\boldsymbol{n}^{T}\mathbf{B\tilde{H}}\right]^{T}
\end{array}\right]\left[\begin{array}{cc}
\mathbf{0} & 2\mathbf{N}\left(x=\left(1+\epsilon\right)h\right)\end{array}\right]\left\{ \boldsymbol{u}\right\} 
\]
Since we are considering a 1D example, $\mathbf{n}$ is a unit vector
in (-x)-direction. Thus $\boldsymbol{n}^{T}$ is a matrix of 1x1 dimension
with a unit value. 
\[
E\left[\begin{array}{c}
\mathbf{0}\\
2\mathbf{N}^{T}\left(x=\left(1+\epsilon\right)h\right)
\end{array}\right]\left[\begin{array}{cc}
\mathbf{B} & \mathbf{B}\left(\left(1-\epsilon\right)\boldsymbol{\tilde{H}}^{+}+\epsilon\boldsymbol{\tilde{H}}^{-}\right)\end{array}\right]=E\left[\begin{array}{cc}
\mathbf{0} & \mathbf{0}\\
2\mathbf{N}^{T}\mathbf{B} & 2\mathbf{N}^{T}\mathbf{B}\left(\left(1-\epsilon\right)\tilde{\boldsymbol{H}}^{+}+\epsilon\tilde{\boldsymbol{H}}^{-}\right)
\end{array}\right]
\]
If we look at:
\[
\mathbf{N}^{T}\mathbf{B}=\left[\begin{array}{cc}
\frac{2h-\left(1+\epsilon\right)h}{h}\frac{-1}{h} & \frac{2h-\left(1+\epsilon\right)h}{h}\frac{1}{h}\\
\frac{\left(1+\epsilon\right)h-h}{h}\frac{-1}{h} & \frac{\left(1+\epsilon\right)h-h}{h}\frac{1}{h}
\end{array}\right]=\frac{1}{h}\left[\begin{array}{cc}
-\left(1-\epsilon\right) & \epsilon\\
-\left(1-\epsilon\right) & \epsilon
\end{array}\right]
\]
Let us consider:
\[
\boldsymbol{K}_{n}=E\mathbf{N}^{T}\mathbf{B}=\frac{E}{h}\left[\begin{array}{cc}
-\left(1-\epsilon\right) & \epsilon\\
-\left(1-\epsilon\right) & \epsilon
\end{array}\right]
\]
and thus Nitsche's term of the matrix is:
\begin{eqnarray*}
\mathbf{K}_{n} & = & \left[\begin{array}{cc}
\mathbf{0} & \mathbf{0}\\
2\boldsymbol{K}_{n} & 2\left(\left(1-\epsilon\right)\boldsymbol{K}_{n}\tilde{\boldsymbol{H}}^{+}+\epsilon\boldsymbol{K}_{n}\tilde{\boldsymbol{H}}^{-}\right)
\end{array}\right]+\left[\begin{array}{cc}
\mathbf{0} & 2\boldsymbol{K}_{n}^{T}\\
\mathbf{0} & 2\left(\left(1-\epsilon\right)\tilde{\boldsymbol{H}}^{+}\boldsymbol{K}_{n}^{T}+\epsilon\tilde{\boldsymbol{H}}^{-}\boldsymbol{K}_{n}^{T}\right)
\end{array}\right]\\
 & = & +\frac{E}{h}\left[\begin{array}{cccc}
0 & 0 & 0 & 0\\
0 & 0 & 0 & 0\\
-2\left(1-\epsilon\right) & 2\epsilon & -4\left(1-\epsilon\right)^{2} & -4\epsilon^{2}\\
-2\left(1-\epsilon\right) & 2\epsilon & -4\left(1-\epsilon\right)^{2} & -4\epsilon^{2}
\end{array}\right]+\frac{E}{h}\left[\begin{array}{cccc}
0 & 0 & -2\left(1-\epsilon\right) & -2\left(1-\epsilon\right)\\
0 & 0 & 2\epsilon & 2\epsilon\\
0 & 0 & -4\left(1-\epsilon\right)^{2} & -4\left(1-\epsilon\right)^{2}\\
0 & 0 & -4\epsilon^{2} & -4\epsilon^{2}
\end{array}\right]
\end{eqnarray*}
In the right hand side, we have the bulk force term,
\[
\int_{x_{2}}^{x_{3}}\boldsymbol{v}_{h}^{T}\boldsymbol{f}\mathrm{d\Omega}=0
\]
the stabilization term:
\begin{eqnarray*}
[[\boldsymbol{v}\left(x=\left(1+\epsilon\right)h\right)]]^{T}\alpha\,\bar{\mbox{i}} & = & 2g\alpha\left\{ v\right\} ^{T}\left[\begin{array}{c}
\mathbf{0}\\
2\mathbf{N}^{T}\left(x=\left(1+\epsilon\right)h\right)
\end{array}\right]\\
 & = & 4g\alpha\left\{ \begin{array}{c}
0\\
0\\
\left(1-\epsilon\right)\\
\left(\epsilon\right)
\end{array}\right\} 
\end{eqnarray*}
and Nitsche's term:
\begin{eqnarray*}
\left(<\sigma\left(\boldsymbol{v}\left(x=\left(1+\epsilon\right)h\right)\right)>.\mathbf{n}\right){}^{T}\,\bar{\mbox{i}} & = & 2gE\left\{ v\right\} ^{T}\left[\begin{array}{c}
\left[\boldsymbol{n}^{T}\mathbf{B}\right]^{T}\\
\left[\boldsymbol{n}^{T}\mathbf{B\tilde{H}}\right]^{T}
\end{array}\right]\\
 & = & 2gE\left\{ v\right\} ^{T}\left[\begin{array}{c}
\begin{array}{c}
-1\\
1
\end{array}\\
\left(\left(1-\epsilon\right)\tilde{\boldsymbol{H}}^{+}+\epsilon\tilde{\boldsymbol{H}}^{-}\right)^{T}\left\{ \begin{array}{c}
-1\\
1
\end{array}\right\} 
\end{array}\right]\\
 & = & 2g\frac{E}{h}\left\{ \begin{array}{c}
-1\\
1\\
-2\left(1-\epsilon\right)\\
-2\epsilon
\end{array}\right\} 
\end{eqnarray*}
Assembling all the terms of all element 2 gives us the combined stiffness
matrix:
\begin{eqnarray*}
\mathbf{A}_{2} & = & \frac{E}{h}\left[\begin{array}{cccc}
1 & -1 & 2\left(1-\epsilon\right) & 2\epsilon\\
-1 & 1 & -2\left(1-\epsilon\right) & -2\epsilon\\
2\left(1-\epsilon\right) & -2\left(1-\epsilon\right) & 4\left(1-\epsilon\right) & 0\\
2\epsilon & -2\epsilon & 0 & 4\epsilon
\end{array}\right]+\alpha\left[\begin{array}{cccc}
0 & 0 & 0 & 0\\
0 & 0 & 0 & 0\\
0 & 0 & 4\left(1-\epsilon\right)^{2} & 4\left(1-\epsilon\right)\epsilon\\
0 & 0 & 4\left(1-\epsilon\right)\epsilon & 4\epsilon^{2}
\end{array}\right]\\
 &  & +\frac{E}{h}\left[\begin{array}{cccc}
0 & 0 & 0 & 0\\
0 & 0 & 0 & 0\\
-2\left(1-\epsilon\right) & 2\epsilon & -4\left(1-\epsilon\right)^{2} & -4\epsilon^{2}\\
-2\left(1-\epsilon\right) & 2\epsilon & -4\left(1-\epsilon\right)^{2} & -4\epsilon^{2}
\end{array}\right]+\frac{E}{h}\left[\begin{array}{cccc}
0 & 0 & -2\left(1-\epsilon\right) & -2\left(1-\epsilon\right)\\
0 & 0 & 2\epsilon & 2\epsilon\\
0 & 0 & -4\left(1-\epsilon\right)^{2} & -4\left(1-\epsilon\right)^{2}\\
0 & 0 & -4\epsilon^{2} & -4\epsilon^{2}
\end{array}\right]
\end{eqnarray*}
with the right hand side contribution
\[
\mathbf{f}=4g\alpha\left\{ \begin{array}{c}
0\\
0\\
\left(1-\epsilon\right)\\
\left(\epsilon\right)
\end{array}\right\} +2g\frac{E}{h}\left\{ \begin{array}{c}
-1\\
1\\
-2\left(1-\epsilon\right)\\
-2\epsilon
\end{array}\right\} 
\]
If:
\[
\mathbf{A}=\mathbf{A}_{2}+\left[\begin{array}{cccc}
1 & 0 & 0 & 0\\
0 & 1 & 0 & 0\\
0 & 0 & 0 & 0\\
0 & 0 & 0 & 0
\end{array}\right]
\]
the global stiffness matrix for the 1D problem, then:
\[
\mathbf{A}u=\mathbf{f}
\]
We see that the sign convention used here is opposite to that of what
we used in the initial derivation. This is because of our choice of
$\mathbf{n}$. Let us try to analyses the matrix \textbf{$\mathbf{A}_{2}$.}
\[
\mathbf{A}_{2}=\left[\begin{array}{cccc}
\frac{E}{h} & \frac{-E}{h} & 0 & \frac{\left(4\epsilon-2\right)E}{h}\\
\frac{-E}{h} & \frac{E}{h} & \frac{\left(4\epsilon-2\right)E}{h} & 0\\
0 & \frac{\left(4\epsilon-2\right)E}{h} & \frac{4E}{h}\left(1-\epsilon\right)\left(2\epsilon-1\right)+4\alpha\left(1-\epsilon\right)^{2} & \frac{-4\left(\epsilon^{2}+\left(1-\epsilon\right)^{2}\right)E}{h}+4\alpha\left(1-\epsilon\right)\epsilon\\
\frac{\left(4\epsilon-2\right)E}{h} & 0 & \frac{-4\left(\epsilon^{2}+\left(1-\epsilon\right)^{2}\right)E}{h}+4\alpha\left(1-\epsilon\right)\epsilon & \frac{4E}{h}\epsilon\left(1-2\epsilon\right)+4\alpha\epsilon^{2}
\end{array}\right]
\]
If $\lambda_{min}\left(\mathbf{A}_{2}\right)$ is the minimum eigen
value of $\mathbf{A}_{2}$, then to maintain coercivity, we need:
\[
\mathbf{x}^{T}\mathbf{A}_{2}\mathbf{x}\ge\lambda_{min}\left(\mathbf{A}_{2}\right)\mathbf{x}^{T}\mathbf{x}
\]
for any nonzero \textbf{$\mathbf{x}\in\mathbb{R}^{n}$}.\textbf{ }Thus
we try to find an optimal $\alpha$ that gives a coercive behavior
for $\mathbf{A}_{2}$. Since a positive definite matrix has $\lambda\left(\mathbf{A}_{2}\right)>0$,
we consider: 
\[
\left|\mathbf{A}_{2}\right|>0
\]
\[
\left|\begin{array}{cccc}
\frac{E}{h} & \frac{-E}{h} & \frac{\left(4\epsilon-2\right)E}{h} & 0\\
\frac{-E}{h} & \frac{E}{h} & 0 & \frac{\left(4\epsilon-2\right)E}{h}\\
\frac{\left(4\epsilon-2\right)E}{h} & 0 & \frac{\left(12\epsilon-8\right)E}{h}+4\alpha\left(1-\epsilon\right)^{2} & \frac{-4E}{h}+4\alpha\left(1-\epsilon\right)\epsilon\\
0 & \frac{\left(4\epsilon-2\right)E}{h} & \frac{-4E}{h}+4\alpha\left(1-\epsilon\right)\epsilon & \frac{\left(4-12\epsilon\right)E}{h}+4\alpha\epsilon^{2}
\end{array}\right|>0
\]
\[
\alpha>\frac{E}{h}
\]

To compare the XFEM solution with the analytical solution, let us
introduce the values of $h=1$ and $\epsilon=\frac{1}{2}$. This gives
$\alpha>\frac{E}{h}$. Solving the system with these values gives
us 
\[
u=\left\{ \begin{array}{c}
-\frac{2g}{3}\\
\frac{2g}{3}\\
g\\
g
\end{array}\right\} 
\]
which is the same as the analytical solution. We base our calculations
of $\alpha$ as devised by Dolbow. 

\subsubsection{Dirichlet type solution}

By considering the above problem with Dirichlet conditions, we obtain
the same solution if we consider:
\begin{eqnarray*}
u^{+}(y=a^{+}) & = & g\\
u^{-}(y=a^{-}) & = & -g
\end{eqnarray*}
We now discretize this problem for element 2 by the help of the variational
form from (\ref{eq:-23}). The XFEM element stiffness matrix is given
by:
\begin{eqnarray*}
\mathbf{K}_{b} & = & \epsilon\left[\begin{array}{cc}
\boldsymbol{K}_{b} & \boldsymbol{K}_{b}\tilde{\boldsymbol{H}}^{-}\\
\left[\tilde{\boldsymbol{H}}^{-}\right]^{T}\boldsymbol{K}_{b} & \left[\tilde{\boldsymbol{H}}^{-}\right]^{T}\boldsymbol{K}_{b}\tilde{\boldsymbol{H}}^{-}
\end{array}\right]+\left(1-\epsilon\right)\left[\begin{array}{cc}
\boldsymbol{K}_{b} & \boldsymbol{K}_{b}\tilde{\boldsymbol{H}}^{+}\\
\left[\tilde{\boldsymbol{H}}^{+}\right]^{T}\boldsymbol{K}_{b} & \left[\tilde{\boldsymbol{H}}^{+}\right]^{T}\boldsymbol{K}_{b}\tilde{\boldsymbol{H}}^{+}
\end{array}\right]\\
 & = & \frac{\epsilon E}{h}\left[\begin{array}{cccc}
1 & -1 & 0 & 2\\
-1 & 1 & 0 & -2\\
0 & 0 & 0 & 0\\
2 & -2 & 0 & 4
\end{array}\right]+\frac{\left(1-\epsilon\right)E}{h}\left[\begin{array}{cccc}
1 & -1 & 2 & 0\\
-1 & 1 & -2 & 0\\
2 & -2 & 4 & 0\\
0 & 0 & 0 & 0
\end{array}\right]
\end{eqnarray*}
The stabilization matrix is given by:
\begin{eqnarray*}
 &  & \boldsymbol{v}\left(x^{+}=\left(1+\epsilon\right)h\right)_{h}\alpha^{+}\boldsymbol{u}\left(x^{+}=\left(1+\epsilon\right)h\right)_{h}+\boldsymbol{v}\left(x^{\text{-}}=\left(1+\epsilon\right)h\right)_{h}\alpha^{-}\boldsymbol{u}\left(x^{-}=\left(1+\epsilon\right)h\right)_{h}\\
 & = & \alpha^{+}\left\{ v\right\} ^{T}\left[\begin{array}{c}
\mathbf{N}^{T}\\
\left[\mathbf{N}\tilde{\boldsymbol{H}}^{+}\right]^{T}
\end{array}\right]\left[\begin{array}{cc}
\mathbf{N} & \mathbf{N}\tilde{\boldsymbol{H}}^{+}\end{array}\right]\left\{ u\right\} +\alpha^{-}\left\{ v\right\} ^{T}\left[\begin{array}{c}
\mathbf{N}^{T}\\
\left[\mathbf{N}\tilde{\boldsymbol{H}}^{-}\right]^{T}
\end{array}\right]\left[\begin{array}{cc}
\mathbf{N} & \mathbf{N}\tilde{\boldsymbol{H}}^{-}\end{array}\right]\left\{ u\right\} \\
\mathbf{K}_{s} & = & \alpha^{+}\left[\begin{array}{cc}
\boldsymbol{K_{s}} & \boldsymbol{K}_{s}\boldsymbol{\tilde{H}}^{+}\\
\left[\tilde{\boldsymbol{H}}^{+}\right]^{T}\boldsymbol{K}_{s} & \left[\tilde{\boldsymbol{H}}^{+}\right]^{T}\boldsymbol{K}_{s}\tilde{\boldsymbol{H}}^{+}
\end{array}\right]+\alpha^{-}\left[\begin{array}{cc}
\boldsymbol{K}_{s} & \boldsymbol{K}_{s}\tilde{\boldsymbol{H}}^{-}\\
\left[\tilde{\boldsymbol{H}}^{-}\right]^{T}\boldsymbol{K}_{s} & \left[\tilde{\boldsymbol{H}}^{-}\right]^{T}\boldsymbol{K}_{s}\tilde{\boldsymbol{H}}^{-}
\end{array}\right]\\
 & = & \alpha^{+}\left[\begin{array}{cccc}
\left(1-\epsilon\right)^{2} & \left(1-\epsilon\right)\epsilon & 2\left(1-\epsilon\right)^{2} & 0\\
\left(1-\epsilon\right)\epsilon & \epsilon^{2} & 2\left(1-\epsilon\right)\epsilon & 0\\
2\left(1-\epsilon\right)^{2} & 2\left(1-\epsilon\right)\epsilon & 4\left(1-\epsilon\right)^{2} & 0\\
0 & 0 & 0 & 0
\end{array}\right]+\alpha^{-}\left[\begin{array}{cccc}
\left(1-\epsilon\right)^{2} & \left(1-\epsilon\right)\epsilon & 0 & -2\left(1-\epsilon\right)\epsilon\\
\left(1-\epsilon\right)\epsilon & \epsilon^{2} & 0 & -2\epsilon^{2}\\
0 & 0 & 0 & 0\\
-2\left(1-\epsilon\right)\epsilon & -2\epsilon^{2} & 0 & 4\epsilon^{2}
\end{array}\right]
\end{eqnarray*}
Nitsche's term of the matrix is obtained as:
\begin{eqnarray*}
 &  & \left[\boldsymbol{v}_{h}^{T+}\left(\sigma(\boldsymbol{u}_{h})^{+}\right).\mathbf{n}+\left(\left(\sigma(\boldsymbol{v}_{h})^{+}\right).\mathbf{n}\right)^{T}\boldsymbol{u}_{h}^{+}-\boldsymbol{v}_{h}^{T-}\left(\sigma(\boldsymbol{u}_{h})^{-}\right).\mathbf{n}-\left(\left(\sigma(\boldsymbol{v}_{h})^{-}\right).\mathbf{n}\right)^{T}\boldsymbol{u}_{h}^{-}\right]_{x=\left(1+\epsilon\right)h}\\
 & = & E\left\{ v\right\} ^{T}\left[\begin{array}{cc}
\mathbf{N} & \mathbf{N}\tilde{\boldsymbol{H}}^{+}\end{array}\right]^{T}\left[\begin{array}{cc}
\mathbf{B} & \mathbf{B}\tilde{\boldsymbol{H}}^{+}\end{array}\right]\left\{ u\right\} +E\left\{ v\right\} ^{T}\left[\begin{array}{cc}
\mathbf{B} & \mathbf{B}\tilde{\boldsymbol{H}}^{+}\end{array}\right]^{T}\left[\begin{array}{cc}
\mathbf{N} & \mathbf{N}\tilde{\boldsymbol{H}}^{+}\end{array}\right]\left\{ u\right\} \\
 & - & E\left\{ v\right\} ^{T}\left[\begin{array}{cc}
\mathbf{N} & \mathbf{N}\tilde{\boldsymbol{H}}^{-}\end{array}\right]^{T}\left[\begin{array}{cc}
\mathbf{B} & \mathbf{B}\tilde{\boldsymbol{H}}^{-}\end{array}\right]\left\{ u\right\} -E\left\{ v\right\} ^{T}\left[\begin{array}{cc}
\mathbf{B} & \mathbf{B}\tilde{\boldsymbol{H}}^{-}\end{array}\right]^{T}\left[\begin{array}{cc}
\mathbf{N} & \mathbf{N}\tilde{\boldsymbol{H}}^{-}\end{array}\right]\left\{ u\right\} \\
 & = & \left[\begin{array}{cc}
\boldsymbol{K}_{n} & \boldsymbol{K}_{n}\tilde{\boldsymbol{H}}^{+}\\
\left[\tilde{\boldsymbol{H}}^{+}\right]^{T}\boldsymbol{K}_{n} & \left[\tilde{\boldsymbol{H}}^{+}\right]^{T}\boldsymbol{K}_{n}\tilde{\boldsymbol{H}}^{+}
\end{array}\right]+\left[\begin{array}{cc}
\boldsymbol{K}_{n}^{T} & \left[\tilde{\boldsymbol{H}}^{+}\right]^{T}\boldsymbol{K}_{n}^{T}\\
\boldsymbol{K}_{n}^{T}\tilde{\boldsymbol{H}}^{+} & \left[\tilde{\boldsymbol{H}}^{+}\right]^{T}\boldsymbol{K}_{n}^{T}\tilde{\boldsymbol{H}}^{+}
\end{array}\right]\\
 & - & \left[\begin{array}{cc}
\boldsymbol{K}_{n} & \boldsymbol{K}_{n}\tilde{\boldsymbol{H}}^{-}\\
\left[\tilde{\boldsymbol{H}}^{-}\right]^{T}\boldsymbol{K}_{n} & \left[\tilde{\boldsymbol{H}}^{-}\right]^{T}\boldsymbol{K}_{n}\tilde{\boldsymbol{H}}^{-}
\end{array}\right]-\left[\begin{array}{cc}
\boldsymbol{K}_{n}^{T} & \tilde{\boldsymbol{H}}^{-}\boldsymbol{K}_{n}^{T}\\
\boldsymbol{K}_{n}^{T}\tilde{\boldsymbol{H}}^{-} & \tilde{\boldsymbol{H}}^{-}\boldsymbol{K}_{n}^{T}\tilde{\boldsymbol{H}}^{-}
\end{array}\right]\\
 & = & \frac{E}{h}\left[\begin{array}{cccc}
-\left(1-\epsilon\right) & \epsilon & -2\left(1-\epsilon\right) & 0\\
-\left(1-\epsilon\right) & \epsilon & -2\left(1-\epsilon\right) & 0\\
-2\left(1-\epsilon\right) & 2\epsilon & -4\left(1-\epsilon\right) & 0\\
0 & 0 & 0 & 0
\end{array}\right]+\frac{E}{h}\left[\begin{array}{cccc}
-\left(1-\epsilon\right) & -\left(1-\epsilon\right) & -2\left(1-\epsilon\right) & 0\\
\epsilon & \epsilon & 2\epsilon & 0\\
-2\left(1-\epsilon\right) & -2\left(1-\epsilon\right) & -4\left(1-\epsilon\right) & 0\\
0 & 0 & 0 & 0
\end{array}\right]\\
 & - & \frac{E}{h}\left[\begin{array}{cccc}
-\left(1-\epsilon\right) & \epsilon & 0 & -2\epsilon\\
-\left(1-\epsilon\right) & \epsilon & 0 & -2\epsilon\\
0 & 0 & 0 & 0\\
2\left(1-\epsilon\right) & -2\epsilon & 0 & 4\epsilon
\end{array}\right]-\frac{E}{h}\left[\begin{array}{cccc}
-\left(1-\epsilon\right) & -\left(1-\epsilon\right) & 0 & 2\left(1-\epsilon\right)\\
\epsilon & \epsilon & 0 & -2\epsilon\\
0 & 0 & 0 & 0\\
-2\epsilon & -2\epsilon & 0 & 4\epsilon
\end{array}\right]
\end{eqnarray*}
Nitsche's part of the right hand terms are given by:
\begin{eqnarray*}
 &  & \left[\left(\left(\sigma(\boldsymbol{v}_{h})^{+}\right).\mathbf{n}\right)^{T}\boldsymbol{g}^{+}-\left(\left(\sigma(\boldsymbol{v}_{h})^{-}\right).\mathbf{n}\right)^{T}\boldsymbol{g}^{-}\right]_{x=\left(1+\epsilon\right)h}\\
 & = & gE\left\{ v\right\} ^{T}\left[\begin{array}{cc}
\mathbf{B} & \mathbf{B}\tilde{\boldsymbol{H}}^{+}\end{array}\right]^{T}+gE\left\{ v\right\} ^{T}\left[\begin{array}{cc}
\mathbf{B} & \mathbf{B}\tilde{\boldsymbol{H}}^{-}\end{array}\right]^{T}\\
 & = & g\alpha^{+}\left\{ \begin{array}{c}
1-\epsilon\\
\epsilon\\
2\left(1-\epsilon\right)\\
0
\end{array}\right\} -g\alpha^{-}\left\{ \begin{array}{c}
1-\epsilon\\
\epsilon\\
0\\
-2\epsilon
\end{array}\right\} 
\end{eqnarray*}
and the stabilization part of the right hand side are given by:
\begin{eqnarray*}
 &  & \left[\boldsymbol{v}_{h}^{T+}\alpha^{+}\boldsymbol{g}^{+}+\boldsymbol{v}_{h}^{T-}\alpha^{-}\boldsymbol{g}^{-}\right]_{x=\left(1+\epsilon\right)h}\\
 & = & g\alpha^{+}\left\{ v\right\} ^{T}\left[\begin{array}{cc}
\mathbf{N} & \mathbf{N}\tilde{\boldsymbol{H}}^{+}\end{array}\right]^{T}-g\alpha^{-}\left\{ v\right\} ^{T}\left[\begin{array}{cc}
\mathbf{N} & \mathbf{N}\tilde{\boldsymbol{H}}^{-}\end{array}\right]^{T}\\
 & = & g\frac{E}{h}\left\{ \begin{array}{c}
-1\\
1\\
-2\\
0
\end{array}\right\} +g\frac{E}{h}\left\{ \begin{array}{c}
-1\\
1\\
0\\
-2
\end{array}\right\} 
\end{eqnarray*}
Assembling all the terms of all the elements gives us the system:
\begin{eqnarray*}
\mathbf{K} & = & \frac{E}{h}\left[\begin{array}{cccc}
1 & 0 & 0 & 0\\
0 & 1 & 0 & 0\\
0 & 0 & 0 & 0\\
0 & 0 & 0 & 0
\end{array}\right]\\
 & + & \frac{\epsilon E}{h}\left[\begin{array}{cccc}
1 & -1 & 0 & 2\\
-1 & 1 & 0 & -2\\
0 & 0 & 0 & 0\\
2 & -2 & 0 & 4
\end{array}\right]+\frac{\left(1-\epsilon\right)E}{h}\left[\begin{array}{cccc}
1 & -1 & 2 & 0\\
-1 & 1 & -2 & 0\\
2 & -2 & 4 & 0\\
0 & 0 & 0 & 0
\end{array}\right]\\
 & + & \alpha^{+}\left[\begin{array}{cccc}
\left(1-\epsilon\right)^{2} & \left(1-\epsilon\right)\epsilon & 2\left(1-\epsilon\right)^{2} & 0\\
\left(1-\epsilon\right)\epsilon & \epsilon^{2} & 2\left(1-\epsilon\right)\epsilon & 0\\
2\left(1-\epsilon\right)^{2} & 2\left(1-\epsilon\right)\epsilon & 4\left(1-\epsilon\right)^{2} & 0\\
0 & 0 & 0 & 0
\end{array}\right]+\alpha^{-}\left[\begin{array}{cccc}
\left(1-\epsilon\right)^{2} & \left(1-\epsilon\right)\epsilon & 0 & -2\left(1-\epsilon\right)\epsilon\\
\left(1-\epsilon\right)\epsilon & \epsilon^{2} & 0 & -2\epsilon^{2}\\
0 & 0 & 0 & 0\\
-2\left(1-\epsilon\right)\epsilon & -2\epsilon^{2} & 0 & 4\epsilon^{2}
\end{array}\right]\\
 & + & \frac{E}{h}\left[\begin{array}{cccc}
-\left(1-\epsilon\right) & \epsilon & -2\left(1-\epsilon\right) & 0\\
-\left(1-\epsilon\right) & \epsilon & -2\left(1-\epsilon\right) & 0\\
-2\left(1-\epsilon\right) & 2\epsilon & -4\left(1-\epsilon\right) & 0\\
0 & 0 & 0 & 0
\end{array}\right]+\frac{E}{h}\left[\begin{array}{cccc}
-\left(1-\epsilon\right) & -\left(1-\epsilon\right) & -2\left(1-\epsilon\right) & 0\\
\epsilon & \epsilon & 2\epsilon & 0\\
-2\left(1-\epsilon\right) & -2\left(1-\epsilon\right) & -4\left(1-\epsilon\right) & 0\\
0 & 0 & 0 & 0
\end{array}\right]\\
 & - & \frac{E}{h}\left[\begin{array}{cccc}
-\left(1-\epsilon\right) & \epsilon & 0 & -2\epsilon\\
-\left(1-\epsilon\right) & \epsilon & 0 & -2\epsilon\\
0 & 0 & 0 & 0\\
2\left(1-\epsilon\right) & -2\epsilon & 0 & 4\epsilon
\end{array}\right]-\frac{E}{h}\left[\begin{array}{cccc}
-\left(1-\epsilon\right) & -\left(1-\epsilon\right) & 0 & 2\left(1-\epsilon\right)\\
\epsilon & \epsilon & 0 & -2\epsilon\\
0 & 0 & 0 & 0\\
-2\epsilon & -2\epsilon & 0 & 4\epsilon
\end{array}\right]\\
\mathbf{f} & = & g\alpha^{+}\left\{ \begin{array}{c}
1-\epsilon\\
\epsilon\\
2\left(1-\epsilon\right)\\
0
\end{array}\right\} -g\alpha^{-}\left\{ \begin{array}{c}
1-\epsilon\\
\epsilon\\
0\\
-2\epsilon
\end{array}\right\} +g\frac{E}{h}\left\{ \begin{array}{c}
-1\\
1\\
-2\\
0
\end{array}\right\} +g\frac{E}{h}\left\{ \begin{array}{c}
-1\\
1\\
0\\
-2
\end{array}\right\} 
\end{eqnarray*}
\[
\mathbf{K}u=\mathbf{f}
\]
Let us try to analyze the matrices \textbf{$\mathbf{K}_{b}$, $\mathbf{K}_{s}$
}and $\mathbf{K}_{n}$ by the respective domain.
\begin{eqnarray*}
\mathbf{A}^{-}=\mathbf{K}_{b}^{-}+\mathbf{K}_{s}^{-}+\mathbf{K}_{n}^{-} & = & \frac{\epsilon E}{h}\left[\begin{array}{cccc}
1 & -1 & 0 & 2\\
-1 & 1 & 0 & -2\\
0 & 0 & 0 & 0\\
2 & -2 & 0 & 4
\end{array}\right]+\alpha^{-}\left[\begin{array}{cccc}
\left(1-\epsilon\right)^{2} & \left(1-\epsilon\right)\epsilon & 0 & -2\left(1-\epsilon\right)\epsilon\\
\left(1-\epsilon\right)\epsilon & \epsilon^{2} & 0 & -2\epsilon^{2}\\
0 & 0 & 0 & 0\\
-2\left(1-\epsilon\right)\epsilon & -2\epsilon^{2} & 0 & 4\epsilon^{2}
\end{array}\right]\\
 &  & -\frac{E}{h}\left[\begin{array}{cccc}
-\left(1-\epsilon\right) & \epsilon & 0 & -2\epsilon\\
-\left(1-\epsilon\right) & \epsilon & 0 & -2\epsilon\\
0 & 0 & 0 & 0\\
2\left(1-\epsilon\right) & -2\epsilon & 0 & 4\epsilon
\end{array}\right]-\frac{E}{h}\left[\begin{array}{cccc}
-\left(1-\epsilon\right) & -\left(1-\epsilon\right) & 0 & 2\left(1-\epsilon\right)\\
\epsilon & \epsilon & 0 & -2\epsilon\\
0 & 0 & 0 & 0\\
-2\epsilon & -2\epsilon & 0 & 4\epsilon
\end{array}\right]
\end{eqnarray*}
\begin{eqnarray*}
\mathbf{A}^{+}=\mathbf{K}_{b}^{+}+\mathbf{K}_{s}^{+}+\mathbf{K}_{n}^{+} & = & \frac{\left(1-\epsilon\right)E}{h}\left[\begin{array}{cccc}
1 & -1 & 2 & 0\\
-1 & 1 & -2 & 0\\
2 & -2 & 4 & 0\\
0 & 0 & 0 & 0
\end{array}\right]+\alpha^{+}\left[\begin{array}{cccc}
\left(1-\epsilon\right)^{2} & \left(1-\epsilon\right)\epsilon & 2\left(1-\epsilon\right)^{2} & 0\\
\left(1-\epsilon\right)\epsilon & \epsilon^{2} & 2\left(1-\epsilon\right)\epsilon & 0\\
2\left(1-\epsilon\right)^{2} & 2\left(1-\epsilon\right)\epsilon & 4\left(1-\epsilon\right)^{2} & 0\\
0 & 0 & 0 & 0
\end{array}\right]\\
 &  & +\frac{E}{h}\left[\begin{array}{cccc}
-\left(1-\epsilon\right) & \epsilon & -2\left(1-\epsilon\right) & 0\\
-\left(1-\epsilon\right) & \epsilon & -2\left(1-\epsilon\right) & 0\\
-2\left(1-\epsilon\right) & 2\epsilon & -4\left(1-\epsilon\right) & 0\\
0 & 0 & 0 & 0
\end{array}\right]+\frac{E}{h}\left[\begin{array}{cccc}
-\left(1-\epsilon\right) & -\left(1-\epsilon\right) & -2\left(1-\epsilon\right) & 0\\
\epsilon & \epsilon & 2\epsilon & 0\\
-2\left(1-\epsilon\right) & -2\left(1-\epsilon\right) & -4\left(1-\epsilon\right) & 0\\
0 & 0 & 0 & 0
\end{array}\right]
\end{eqnarray*}
To maintain the coercivity of the system, we try to find $\alpha^{+}$
and $\alpha^{-}$ such that:
\[
\left|\mathbf{A}^{-}\right|>0
\]
We ignore the row and column corresponding to the diagonal with zero
value. 
\[
\left|\begin{array}{ccc}
\frac{3\left(1-\epsilon\right)E}{h}+\alpha^{-}\left(1-\epsilon\right)^{2} & \frac{E}{h}+\alpha^{-}\epsilon\left(1-\epsilon\right) & \frac{2\epsilon E}{h}-2\alpha^{-}\epsilon\left(1-\epsilon\right)\\
\frac{E}{h}+\alpha^{-}\epsilon\left(1-\epsilon\right) & \frac{\left(1-3\epsilon\right)E}{h}+\alpha^{-}\epsilon^{2} & \frac{-2\left(1-3\epsilon\right)E}{h}-2\alpha^{-}\epsilon^{2}\\
\frac{2\epsilon E}{h}-2\alpha^{-}\epsilon\left(1-\epsilon\right) & \frac{-2\left(1-3\epsilon\right)E}{h}-2\alpha^{-}\epsilon^{2} & \frac{4\left(1-3\epsilon\right)E}{h}+4\alpha^{-}\epsilon^{2}
\end{array}\right|>0
\]
\[
\alpha^{-}>\frac{E}{2\epsilon h}
\]
and:
\[
\left|\mathbf{A}^{+}\right|>0
\]
\[
\left|\begin{array}{ccc}
\frac{\left(3\epsilon-2\right)E}{h}+\alpha^{+}\left(1-\epsilon\right)^{2} & \frac{E\left(-1+\epsilon\right)}{h}+\alpha^{+}\epsilon\left(1-\epsilon\right) & \frac{2\left(3\epsilon-2\right)E}{h}-2\alpha^{+}\left(1-\epsilon\right)^{2}\\
\frac{E\left(-1+\epsilon\right)}{h}+\alpha^{+}\epsilon\left(1-\epsilon\right) & \frac{3\epsilon E}{h}+\alpha^{+}\epsilon^{2} & \frac{2\left(-1+\epsilon\right)E}{h}+2\alpha^{+}\epsilon\left(1-\epsilon\right)\\
\frac{2\left(3\epsilon-2\right)E}{h}-2\alpha^{+}\left(1-\epsilon\right)^{2} & \frac{2\left(-1+\epsilon\right)E}{h}+2\alpha^{+}\epsilon\left(1-\epsilon\right) & \frac{4\left(3\epsilon-2\right)E}{h}+4\alpha^{+}\left(1-\epsilon\right)^{2}
\end{array}\right|>0
\]
\[
\alpha^{+}>\frac{E}{2\left(1-\epsilon\right)h}
\]
Solving the system with $h=1$, $\epsilon=\frac{1}{2}$ we have $\alpha^{+}>\frac{E}{h}$
and $\alpha^{-}>\frac{E}{h}$. With these values, we get the solution:
\[
u=\left\{ \begin{array}{c}
-\frac{2g}{3}\\
\frac{2g}{3}\\
g\\
g
\end{array}\right\} 
\]
which is in-fact the solution obtained by assuming jump conditions.
\newpage{}

\section{Implementation in Code-Aster}

\subsection{Options and routines}

\subsubsection*{RIGI\_NITS}

This option calculates the elementary matrices necessary for Nitsche's
method. It receives material, geometry and mesh information. It takes
all the level-sets into consideration and one of the input parameter
is whether the question is of type 'Jump', 'Dirichlet+' or 'Dirichlet-'.
It also receives the stabilization parameter $\alpha$.

\subsubsection*{CHAR\_MECA\_NITS\_R}

This option calculates the elementary vectors associated with Nitsche's
method. Along with all the information that RIGI\_NITS obtains, it
also receives the value of the parameter associated to 'Jump' or 'Dirichlet'
boundary condition.

\subsubsection*{RAPH\_MECA\_NITS\_R}

This option is similar to CHAR\_MECA\_NITS\_R and calculates the elementary
vectors associated to Nitsche's method for a non-linear Newton iteration.
It utilizes the parameter associated to 'Jump' or 'Dirichlet' and
calculates the total contribution of the displacement field in the
given iteration.

\subsubsection*{PARA\_NITS}

This option handles the geometry and material information required
to calculate the paramters $\gamma$, the weighting parameter, and
$\alpha$, the stabilization parameter for the given element. This
helps us solve the given problem using the wighted discretization
discussed in section \ref{subsec:Weighted-Discretization}.

\subsubsection*{xmnits1}

This subroutine calculates the matrix $\boldsymbol{n}^{T}\boldsymbol{DB}$
by taking the geometrical, material and mesh information along with
the shape-functions', sub-elements' and gauss-points' informtion for
the particular element involved.

\subsubsection*{xmnits2}

This subroutine uses the heavy side function information and calculates
$\boldsymbol{N\tilde{H}}$ matrix. This subroutine can be used also
to calculate $\boldsymbol{n}^{T}\boldsymbol{DB}\boldsymbol{\tilde{H}}$
since $\boldsymbol{N}$ and $\boldsymbol{n}^{T}\boldsymbol{DB}$ are
of the same form for a given sub-element.

\subsubsection*{xgamma}

This routine is called by the option PARA\_NITS to compute the parameters
$\gamma$ and $\alpha$.

\subsubsection*{xmnits\_cote}

This subroutine calculates the part of the elementary matrix that
is contributed by the penalization and Nitsche's part for the 'Dirichlet'
type of problems. This single subroutine is capable of computing both
the '+' and the '-' parts of the interfacial boundary conditions.

\subsubsection*{xvnits\_cote}

This subroutine similar to xmnits\_cote computes the elementary vector
part contributed by the penalization and Nitsche's part for the 'Dirichlet'
type of problems. This subroutine in addition receives the value of
the boundary condition.

\subsubsection*{xmnits\_saut}

This subroutine calculates the part of the elementary matrix that
is contributed by the penalization and Nitsche's part for the 'Jump'
type of problems. 

\subsubsection*{xvnits\_saut}

This subroutine calculates the part of the elementary vector that
is contributed by the penalization and Nitsche's part for the 'Jump'
type of problems.

\subsubsection*{te0567}

This routine does all the computations necessary to use the option
RIGI\_NITS to compute the part of elementary matrice corresponding
to Nitsche's method.

\subsubsection*{te0568}

This routine does all the computations necessary to use the option
CHAR\_MECA\_NITS\_R and RAPH\_MECA\_NITS\_R to compute the part of
elementary vector corresponding to Nitsche's method.

\subsection{Simple test case}

\begin{figure}
\begin{centering}
\includegraphics[height=8cm]{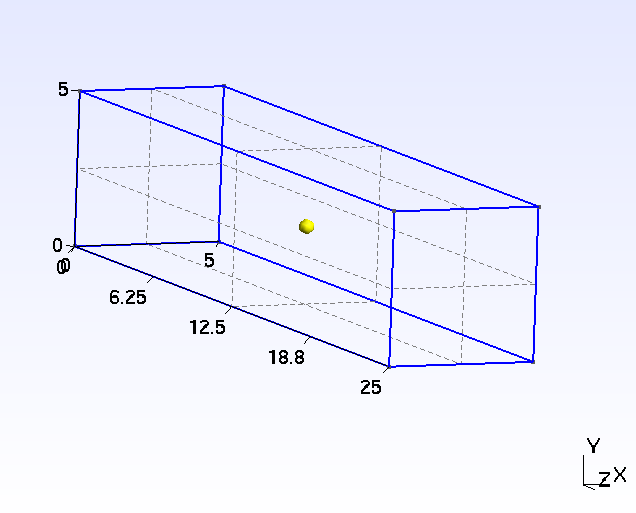}
\par\end{centering}
\caption{Rectangular block of size $5\times5\times25$ with a plane interface
at $z=12.5$}
\end{figure}
In order to test the implementation of Nitsche's method, a simple
rectangular block with an interface in between was used.The block
was tested under conditions of jump in displacement at the interface
as well as a prescribed Dirichlet conditions on both sides of the
interface. The following conditions were assumed :

\[
E=205\times10{}^{3}\mbox{ Pa}
\]
\[
\nu=0.3
\]
\[
L=25\mbox{ mm}
\]
\[
f=0\mbox{ N}\mbox{ in }\Omega
\]
\[
u_{z}(z=0)=u_{z}(z=25)=0\mbox{ mm}
\]
\[
u_{x}(x=0)=0\mbox{ mm}
\]
\[
u_{y}(y=0)=0\mbox{ mm}
\]

\begin{enumerate}
\item For the case of jump in displacement
\[
\bar{\mbox{i}}=\left[\left[u_{z}(z=12.5)\right]\right]=3\times10^{-6}\mbox{ m}
\]
\[
\bar{\mbox{j}}=\left[\left[\sigma(z=12.5)\right]\right].\mbox{n}=0\mbox{ Pa}
\]
\item For the case of Dirichlet condition
\[
g^{+}=u_{z}^{+}(z=12.5)=1.5\times10^{-6}\mbox{ m}
\]
\[
g^{-}=u_{z}^{-}(z=12.5)=-1.5\times10^{-6}\mbox{ m}
\]
\end{enumerate}
It can be noted that both conditions will generate the same output.
Due to the symmetric nature of the problem, we consider $g^{+}=g^{-}=g$.
Analytically, the solution for the above problem is:
\begin{equation}
\sigma_{zz}=-\frac{2Eg}{L}\mbox{ in }\Omega
\end{equation}
\begin{equation}
u_{x}=\frac{2\nu gx}{L}
\end{equation}
\begin{equation}
u_{y}=\frac{2\nu gy}{L}
\end{equation}
\begin{equation}
u_{z}=-\frac{2gz}{L}\mbox{ in }\Omega^{-}
\end{equation}
\begin{equation}
u_{z}=-\frac{2gz}{L}+2g\mbox{ in }\Omega^{+}
\end{equation}

The above problem was tested with both Nitsche's method and the penalization
method with varying penalization parameter. The solutions obtained
were then compared to the analytical solution in terms of error in
energy norm and displacement norm in $L_{2}$.

\begin{figure}
\begin{centering}
\includegraphics[height=8cm]{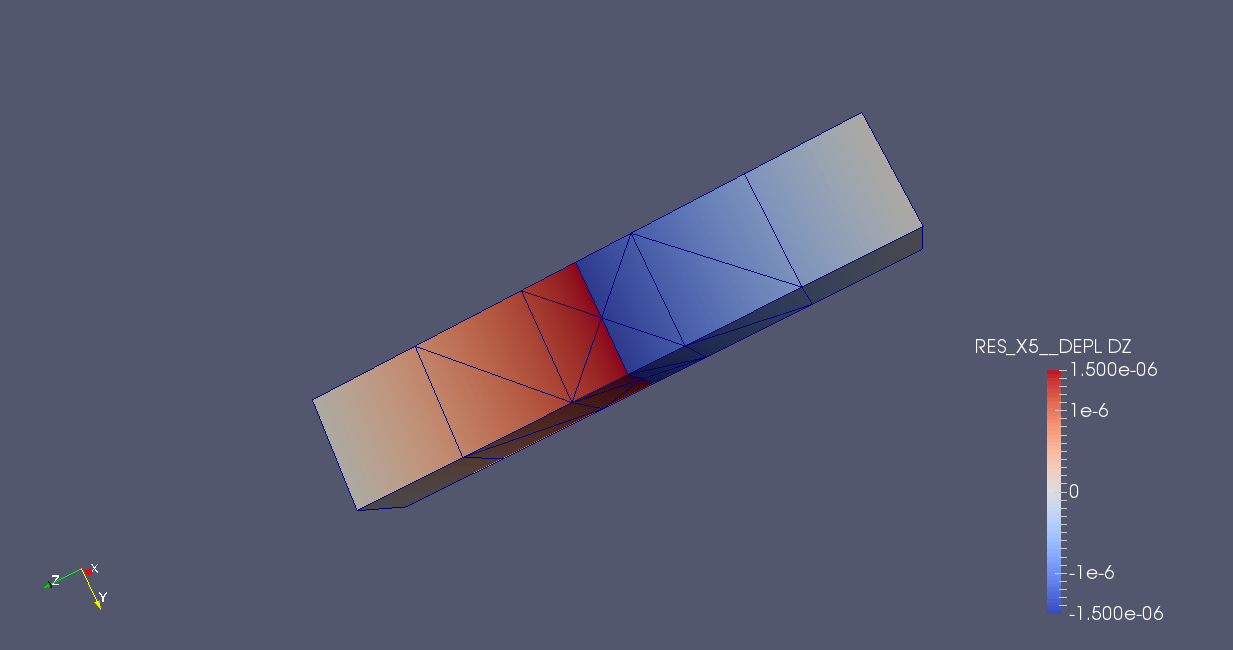}
\par\end{centering}
\caption{\label{fig:The-solution-as}the solution as obtained in Code Aster
after analyzing a problem of jump type using non-linear Newton iterations.}
\end{figure}

The solution is shown in Figure \ref{fig:The-solution-as} imprinted
on the hexahedral mesh. If we look at the error norms, both Dirichlet
and jump conditions have similar behavior, thus signifying the equivalence
of the two methods for equivalent boundary conditions. We see that
for penalization method, the solution tends to converge at higher
penalization parameter while for Nitsche's method, we obtain good
results even at very low stabilization parameter and in fact, the
machine error adds up at higher parameter (Figures \ref{fig:Relative-energy-norm}
and \ref{fig:Relative-energy-norm-1}), resulting in higher error,
even for Nitsche's method!

\begin{center}
\begin{figure}
\begin{centering}
\includegraphics[height=8cm]{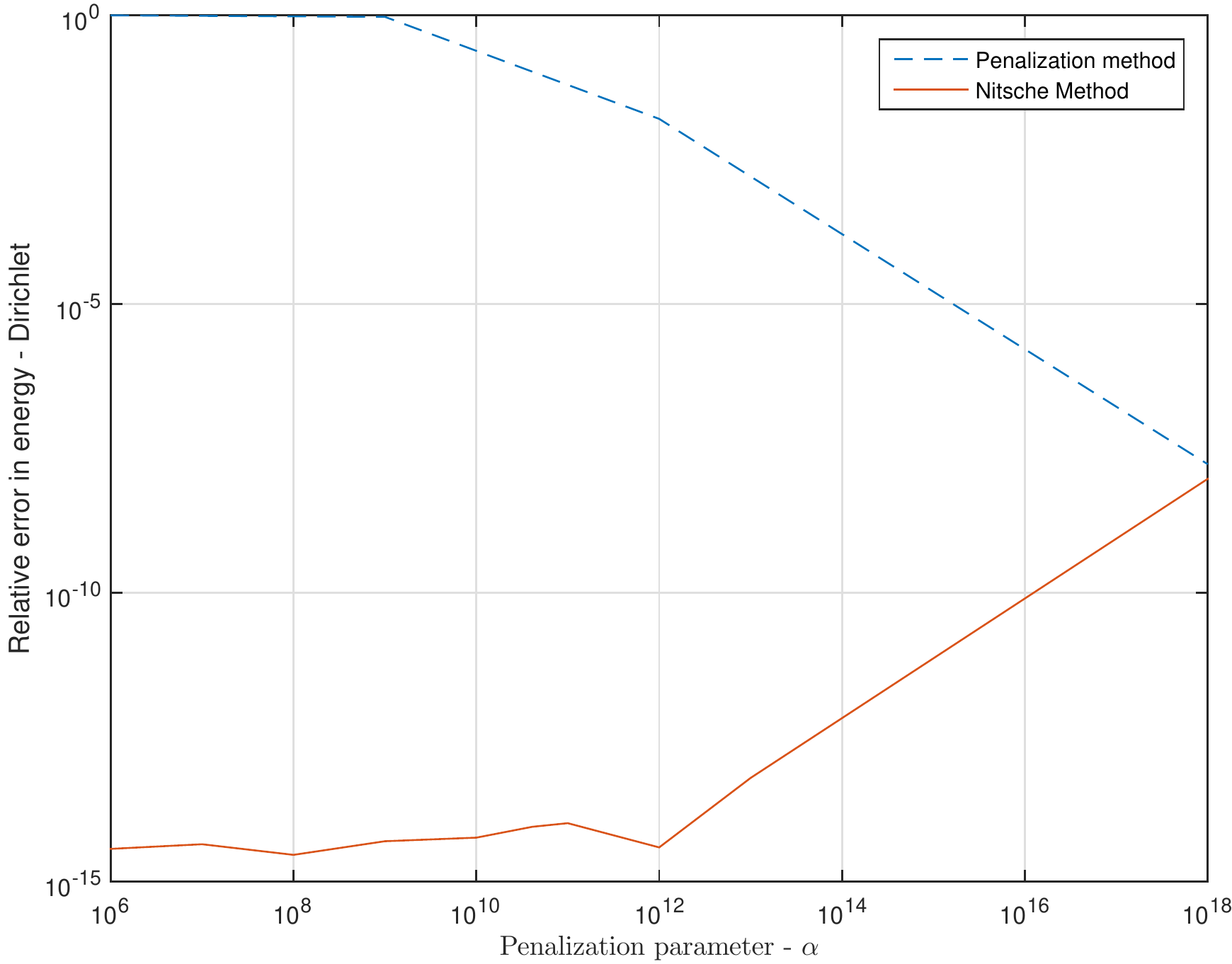}
\par\end{centering}
\caption{\label{fig:Relative-energy-norm}Relative energy norm error for Dirichlet
conditions and the study of their variation with change in stabilization
parameter in Nitsche's and penalization methods.}
\end{figure}
\begin{figure}
\begin{centering}
\includegraphics[height=8cm]{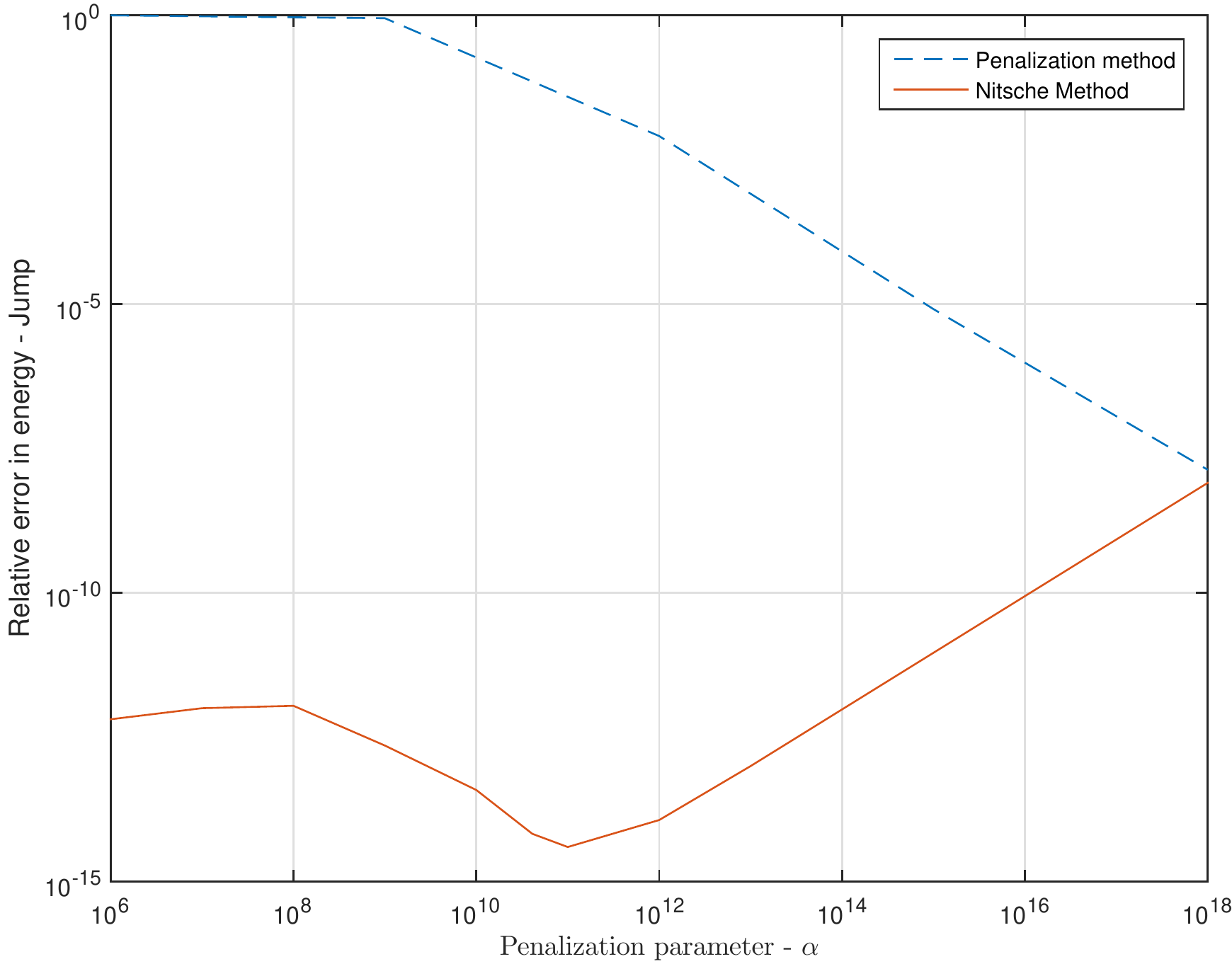}
\par\end{centering}
\caption{\label{fig:Relative-energy-norm-1}Relative energy norm error for
Jump conditions and the study of their variation with change in stabilization
parameter in Nitsche's and penalization methods.}
\end{figure}
\begin{figure}
\begin{centering}
\includegraphics[height=8cm]{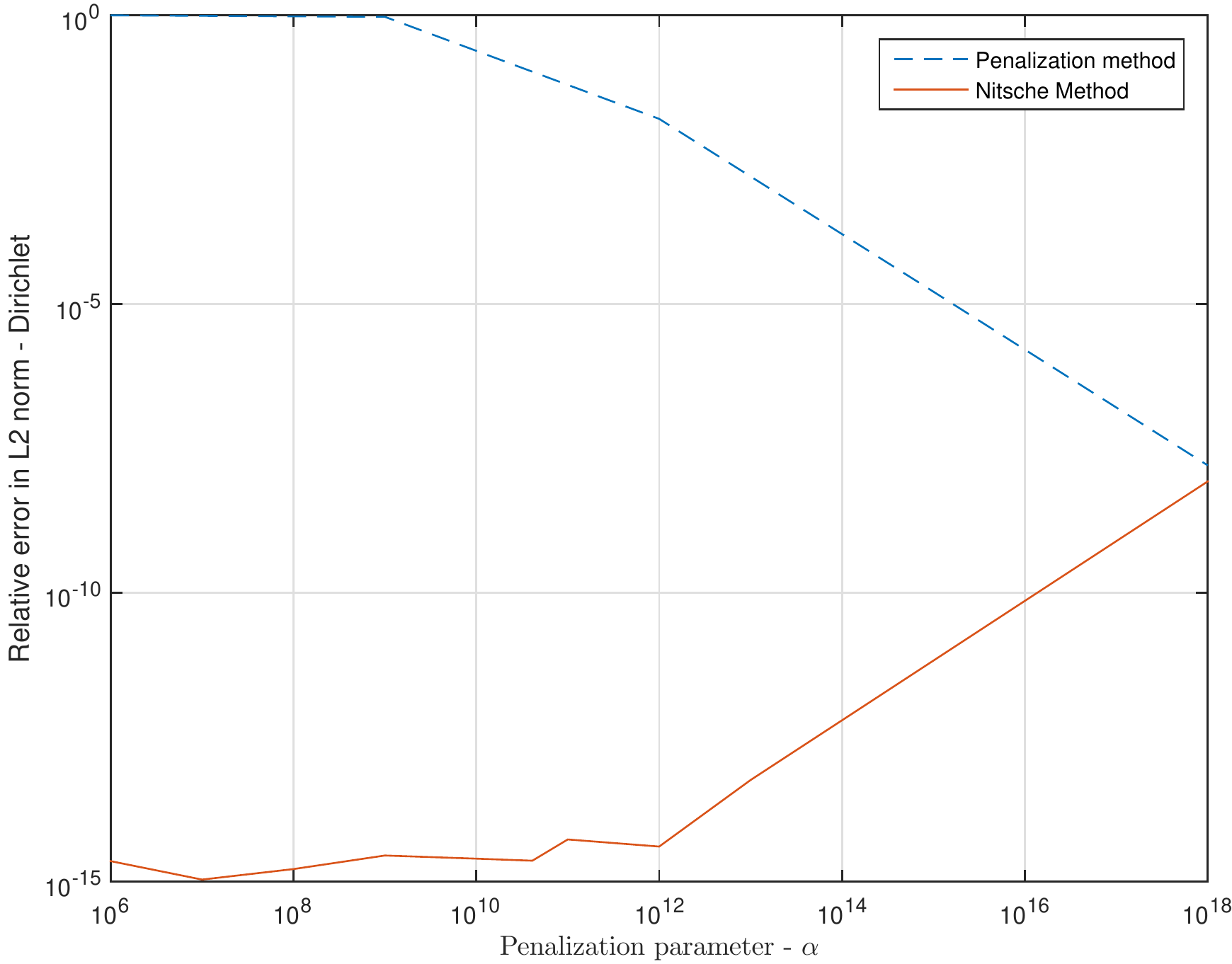}
\par\end{centering}
\caption{\label{fig:Relative-energy-norm-1-2}Relative L2 norm error for Dirichlet
conditions and the study of their variation with change in stabilization
parameter in Nitsche's and penalization methods.}
\end{figure}
\begin{figure}
\begin{centering}
\includegraphics[height=8cm]{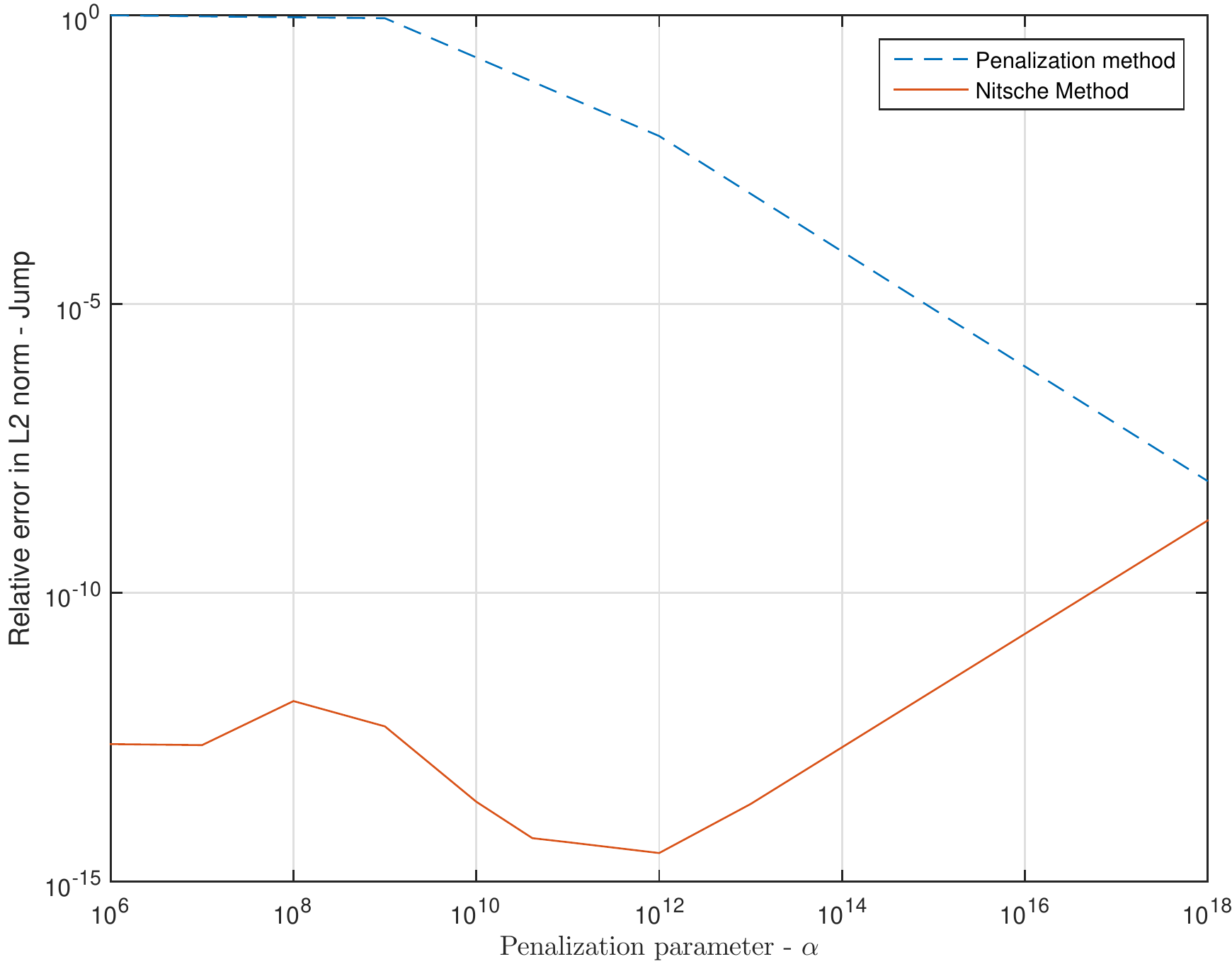}
\par\end{centering}
\caption{\label{fig:Relative-energy-norm-1-1}Relative L2 norm error for Jump
conditions and the study of their variation with change in stabilization
parameter in Nitsche's and penalization methods.}
\end{figure}
\par\end{center}

\subsection{Circular Inclusion problem}

In this two-dimensional bi-material test case, a weak discontinuity
is present, and the displacement field is continuous with discontinuous
stresses and strains.\cite{key-7,key-11} Inside a circular plate
of radius b, whose material is defines by $E_{1}=1$ and $\nu_{1}=0.25$,
a circular inclusion with radius $a$ of different material with $E_{2}=10$
and $\nu=0.3$ is considered. The loading of the structure results
from a linear displacement of the outer boundary: $u_{r}\left(b,\theta\right)=r$
and $u_{\theta}\left(b,\theta\right)=0$. The situation is depicted
in figure (\ref{fig:Problem-statement-of}). The exact solution can
be found in \cite{key-16}.

\begin{center}
\begin{figure}
\begin{centering}
\includegraphics[bb=170bp 470bp 430bp 715bp,clip]{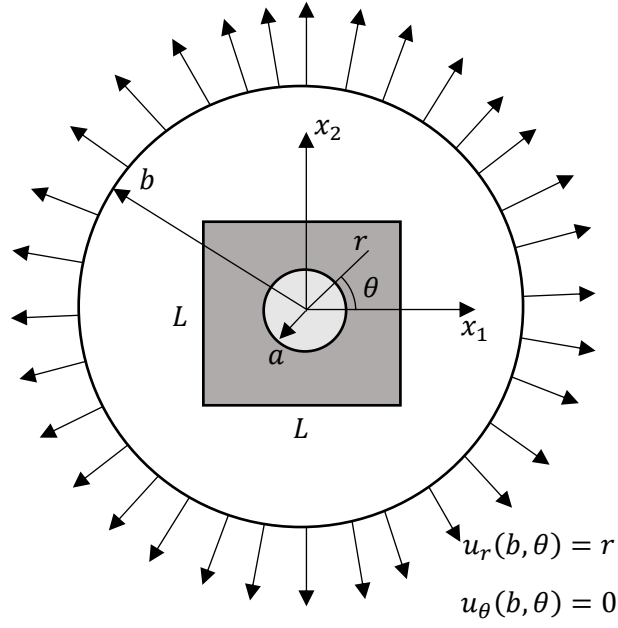}
\par\end{centering}
\caption{\label{fig:Problem-statement-of}Problem statement of the circular
inclusion problem (the gray area is the numerical domain)}
\end{figure}
\par\end{center}

The stresses are given as:
\begin{equation}
\sigma_{rr}\left(r,\theta\right)=2\mu\varepsilon_{rr}+\lambda\left(\varepsilon_{rr}+\varepsilon_{\theta\theta}\right)
\end{equation}
\begin{equation}
\sigma_{\theta\theta}\left(r,\theta\right)=2\mu\varepsilon_{\theta\theta}+\lambda\left(\varepsilon_{rr}+\varepsilon_{\theta\theta}\right)
\end{equation}
where the Lamé constants $\lambda$ and $\mu$ have to be replaced
by the appropriate values for the corresponding area, respectively.
The strains are:
\begin{equation}
\varepsilon_{rr}\left(r,\theta\right)=\begin{cases}
\left(1-\frac{b^{2}}{a^{2}}\right)\alpha+\frac{b^{2}}{a^{2}}, & 0\le r\le a\\
\left(1+\frac{b^{2}}{r^{2}}\right)\alpha-\frac{b^{2}}{r^{2}}, & a\le r\le b
\end{cases}
\end{equation}
\begin{equation}
\varepsilon_{\theta\theta}\left(r,\theta\right)=\begin{cases}
\left(1-\frac{b^{2}}{a^{2}}\right)\alpha+\frac{b^{2}}{a^{2}}, & 0\le r\le a\\
\left(1+\frac{b^{2}}{r^{2}}\right)\alpha-\frac{b^{2}}{r^{2}}, & a\le r\le b
\end{cases}
\end{equation}
and the displacements:
\begin{equation}
u_{r}\left(r,\theta\right)=\begin{cases}
\left[\left(1-\frac{b^{2}}{a^{2}}\right)\alpha+\frac{b^{2}}{a^{2}}\right]r, & 0\le r\le a\\
\left(r-\frac{b^{2}}{r^{2}}\right)\alpha+\frac{b^{2}}{r}, & a\le r\le b
\end{cases}
\end{equation}
\begin{equation}
u_{\theta}\left(r,\theta\right)=0
\end{equation}
The parameter $\alpha$ involved in these definitions is:
\begin{equation}
\alpha=\frac{\left(\lambda_{1}+\mu_{1}+\mu_{2}\right)b^{2}}{\left(\lambda_{2}+\mu_{2}\right)a^{2}+\left(\lambda_{1}+\mu_{1}\right)\left(b^{2}-a^{2}\right)+\mu_{2}b^{2}}
\end{equation}

\begin{figure}
\begin{centering}
\includegraphics[height=8cm]{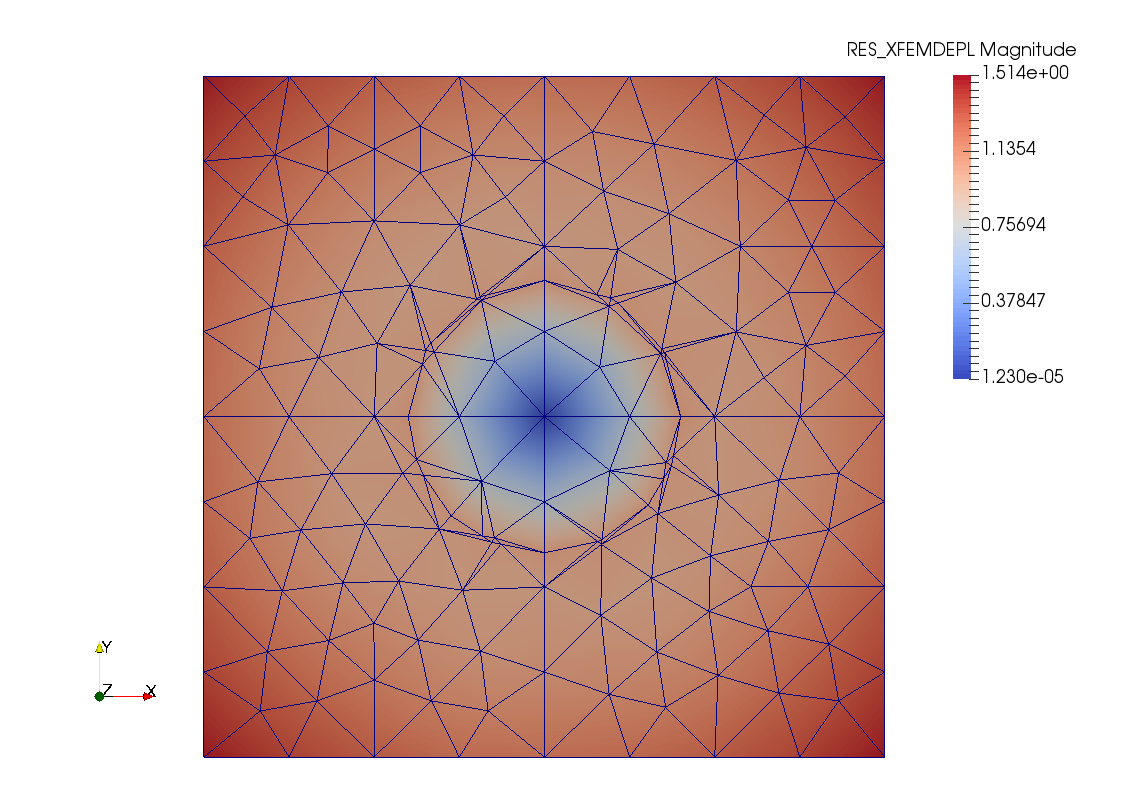}
\par\end{centering}
\caption{\label{fig:Circular-inclusion-problem}Circular inclusion problem
solution on an irregular grid of size $h=0.2$. Notice how the element
is divided into subelements at the interface.}
\end{figure}

\begin{figure}
\begin{centering}
\includegraphics[height=8cm]{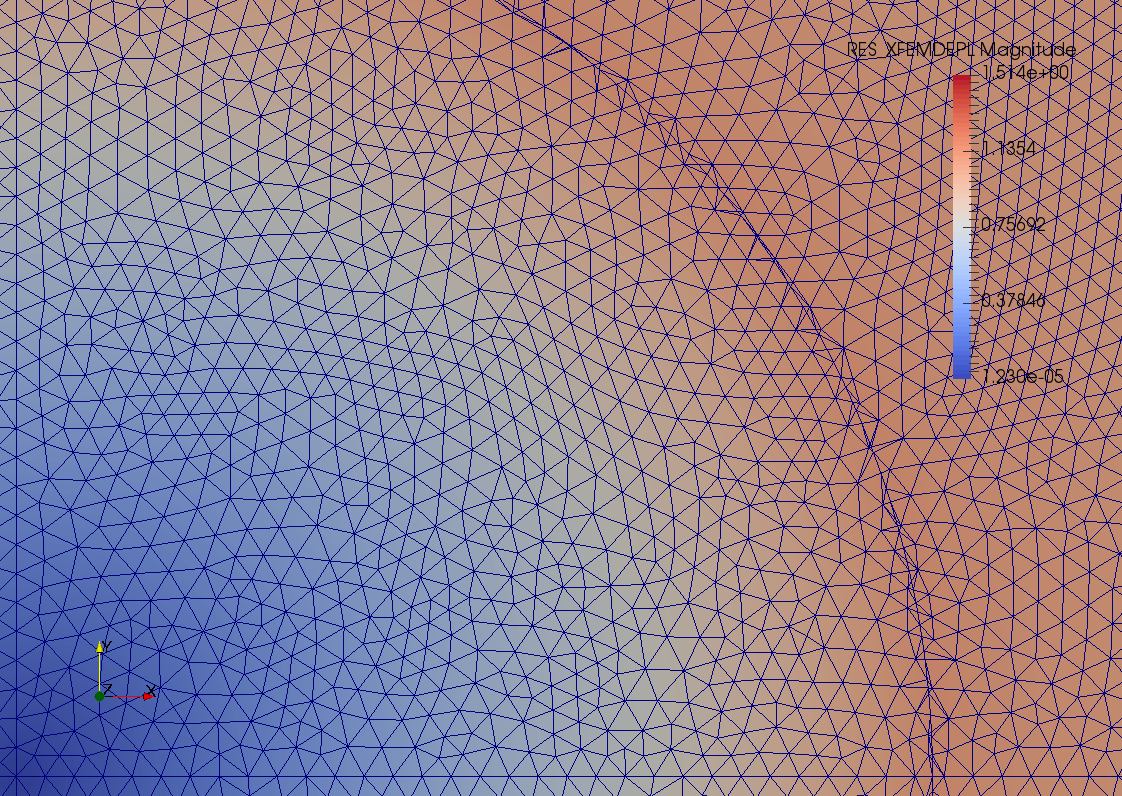}
\par\end{centering}
\caption{\label{fig:Circular-inclusion-problem-1}Circular inclusion problem
solution on an irregular grid of size $h=0.0125$. In this close up
view the effect of level set on the mesh and the solution can be seen.}
\end{figure}

\begin{figure}
\begin{centering}
\includegraphics[height=8cm]{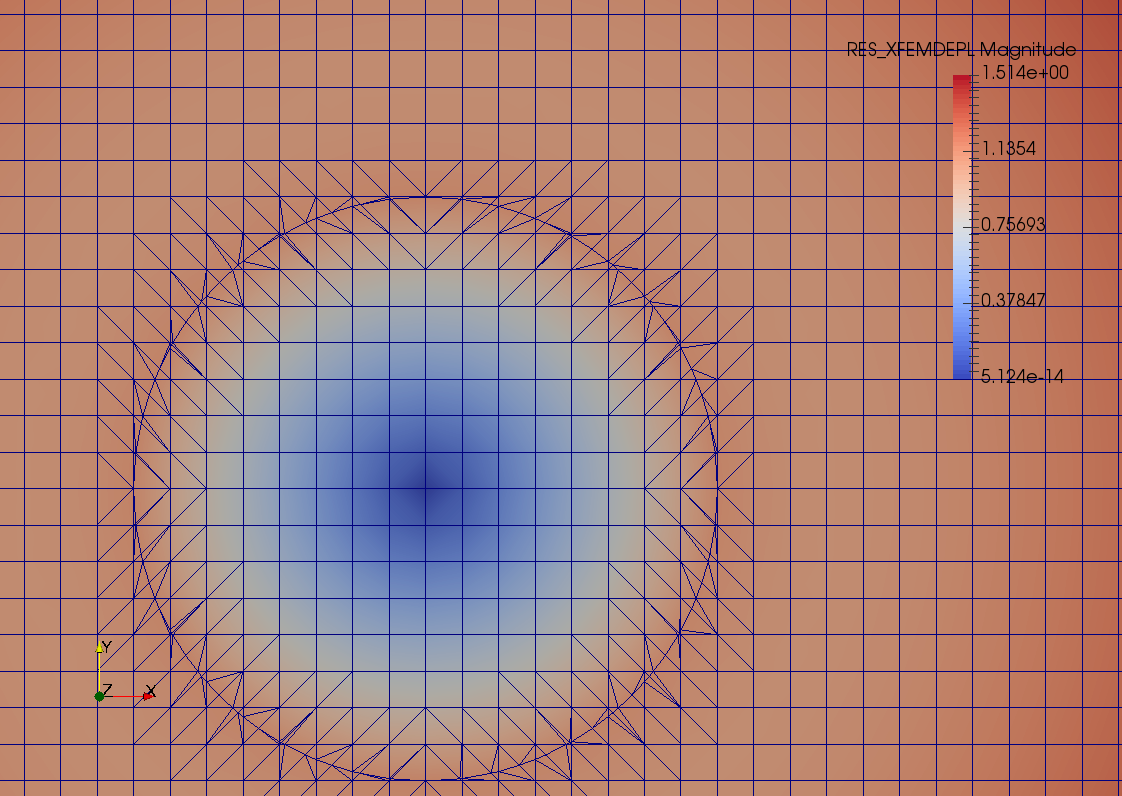}
\par\end{centering}
\caption{\label{fig:Circular-inclusion-problem-2}Circular inclusion problem
solution on an regular grid of size $h=0.05$. The enrichment is extended
to an element beyond the interface in both sides.}
\end{figure}

\begin{figure}
\begin{centering}
\includegraphics[height=5cm]{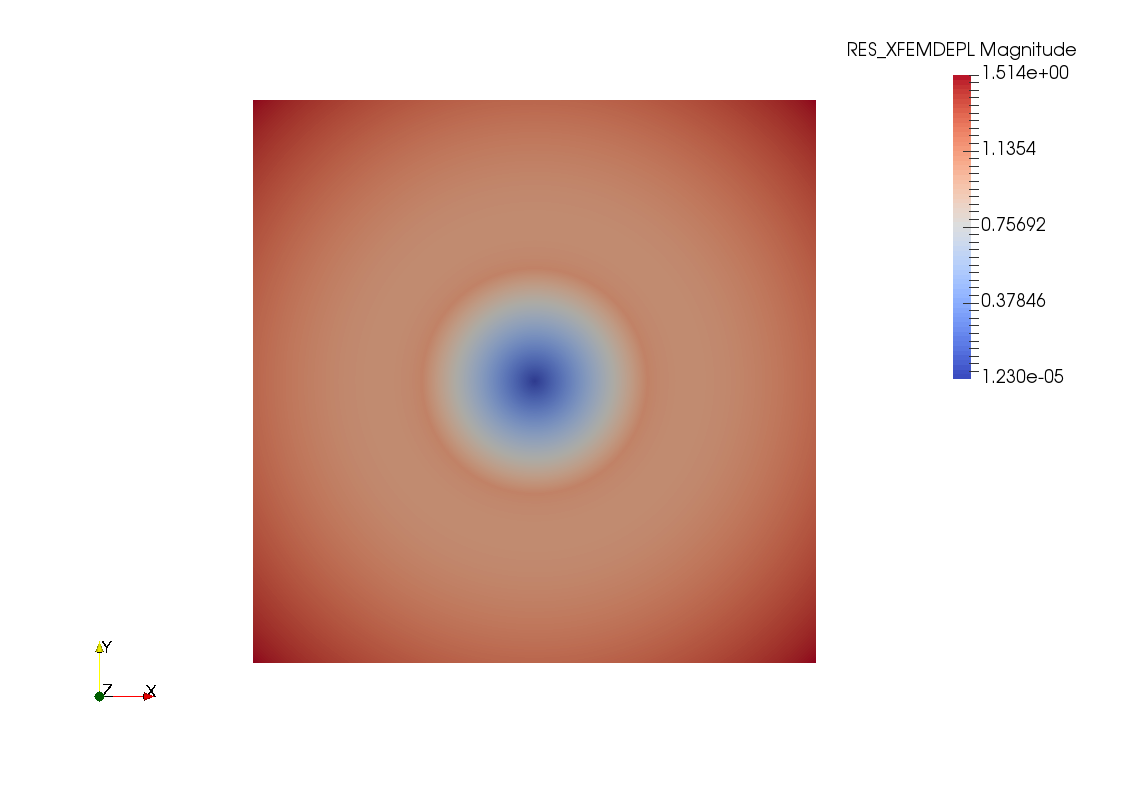}\includegraphics[height=5cm]{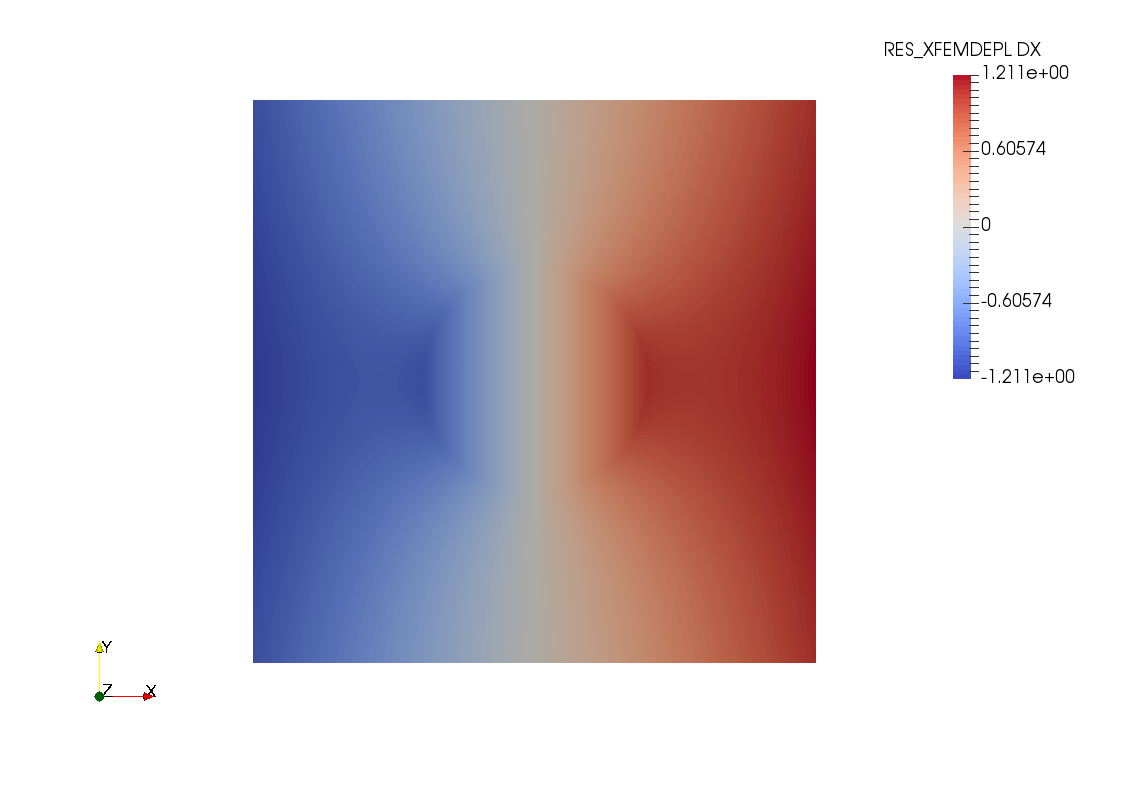}
\par\end{centering}
\begin{centering}
\includegraphics[height=5cm]{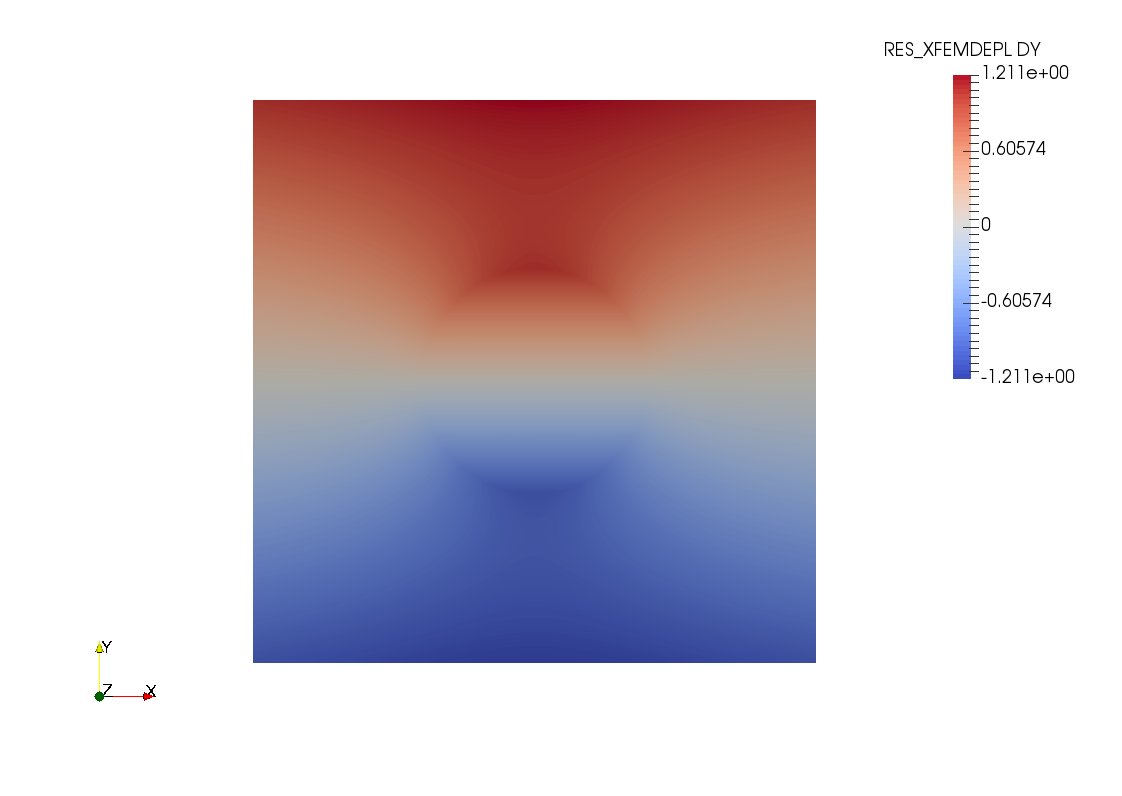}
\par\end{centering}
\caption{Circular inclusion problem solution on an regular grid of size $h=0.0125$.
\label{fig:Circular-inclusion-problem-3}Magnitude of the displacement
solution on the entire domain and its variation with respect to x
and y direction.}
\end{figure}

\begin{figure}
\begin{centering}
\includegraphics[height=5cm]{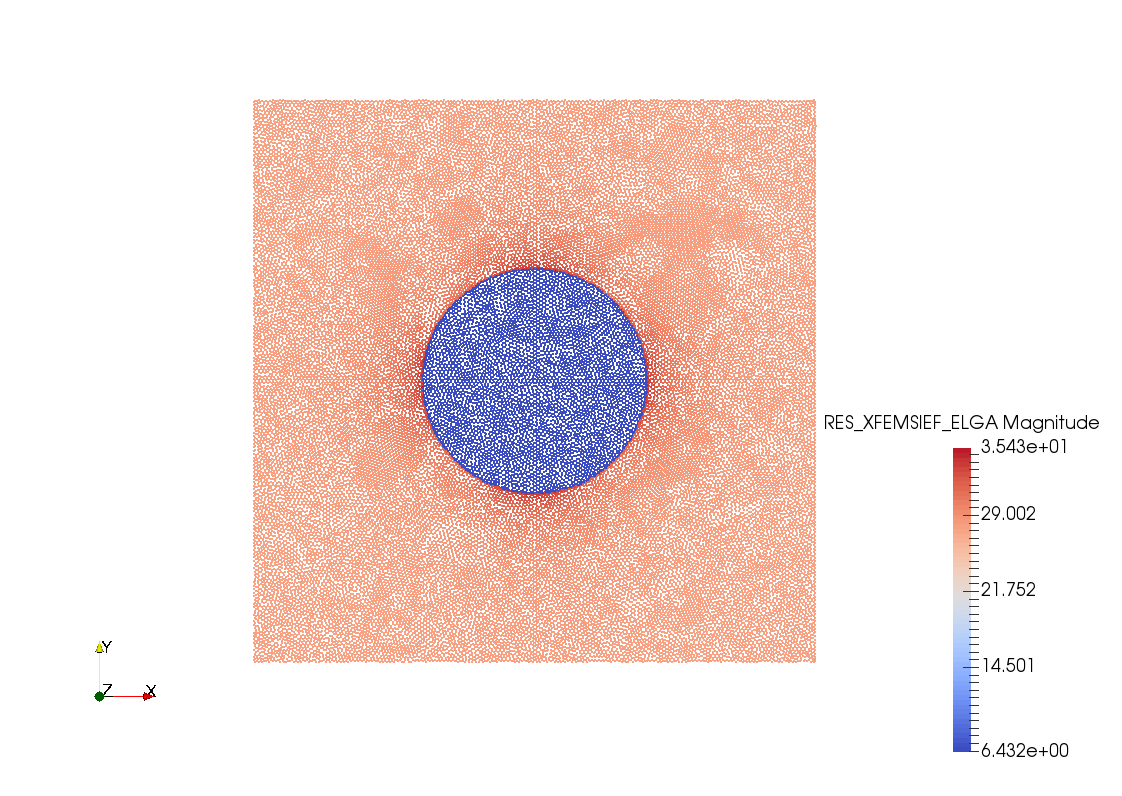}\includegraphics[height=5cm]{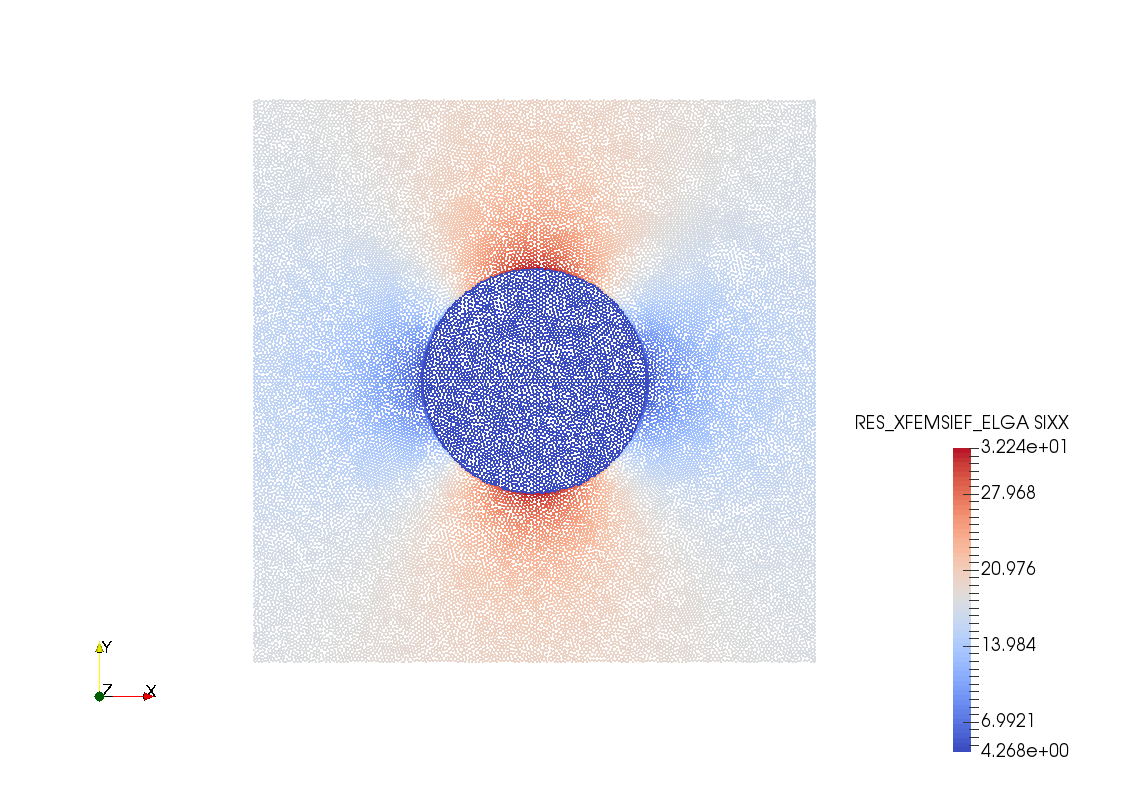}
\par\end{centering}
\begin{centering}
\includegraphics[height=5cm]{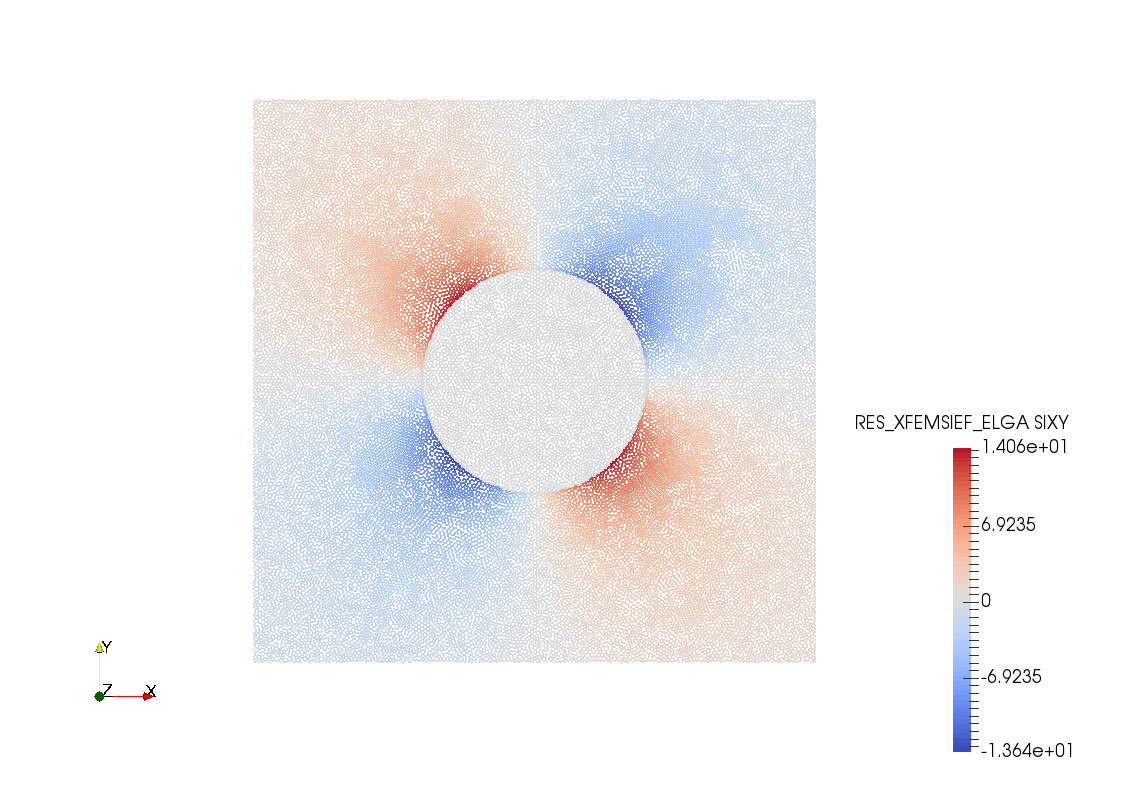}\includegraphics[height=5cm]{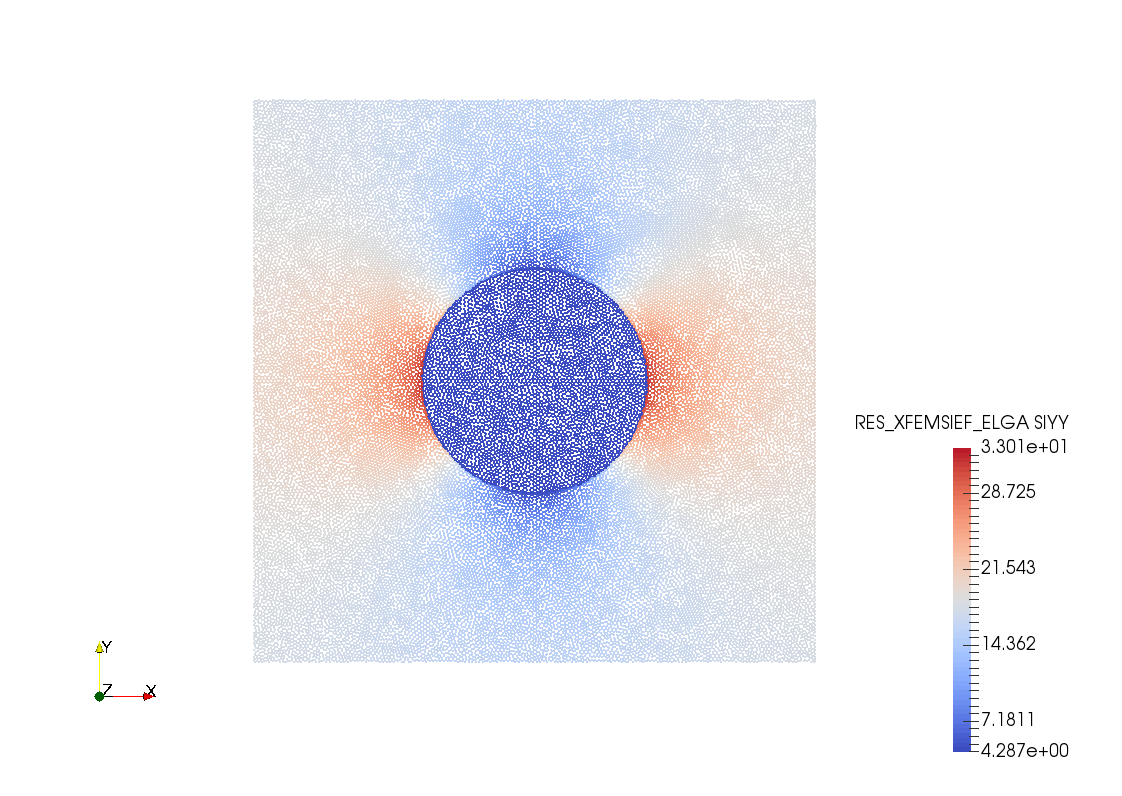}
\par\end{centering}
\caption{\label{fig:Circular-inclusion-problem-4}Circular inclusion problem
solution on an regular grid of size $h=0.0125$. Solution of the discontinuous
stresses in the domain}
\end{figure}

For the numerical model, the domain is a square of size $L\times L$
with $L=2$, the outer radius is chosen to be $b=2$ and the inner
radius $a=0.4$. The exact stresses are prescribed along the boundaries
of the square domain, and displacements are prescribed as:
\[
u_{1}\left(0,\pm1\right)=0
\]
and:
\[
u_{2}\left(\pm1,0\right)=0
\]
Plane strain conditions are assumed. Results are obtained for different
methods, Nitsche's method, penalty method and non-linear method with
Lagrange multipliers. A set of displacement and stress on gauss points
plots have been presented in figures \ref{fig:Circular-inclusion-problem}
to \ref{fig:Circular-inclusion-problem-4}. The interface is embedded
into the mesh and Code Aster divides every element cut by the interface
into sub-elements. The elements adjacent to the ones cut by interface
are also enriched. (Figure \ref{fig:Circular-inclusion-problem-2}).
The displacement plot shows a continuous change of the magnitude in
both the axis(figure \ref{fig:Circular-inclusion-problem-3}), while
the stress plot shows the discontinuous stress in the two domains.

\begin{figure}
\begin{centering}
\includegraphics[height=8cm]{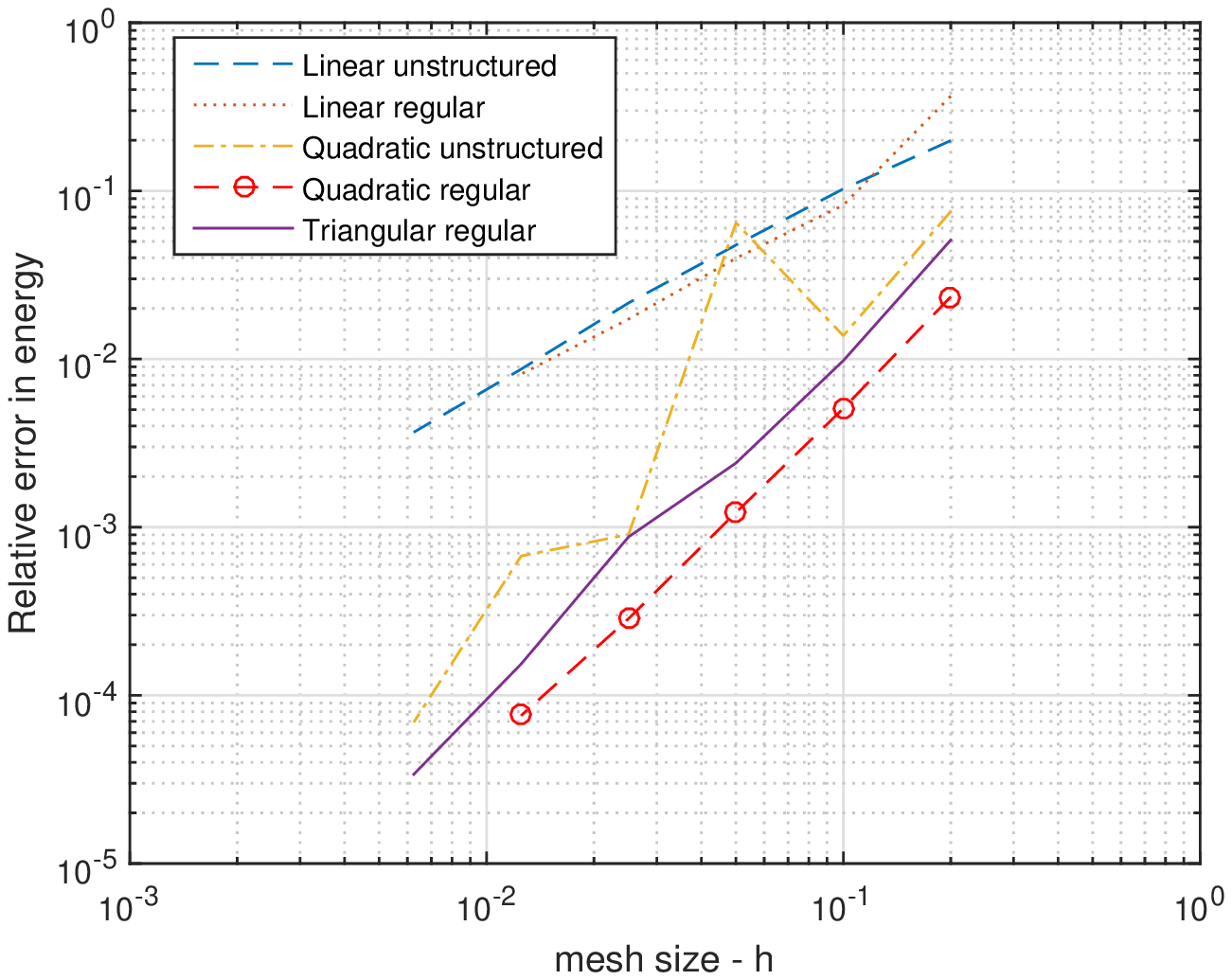}
\par\end{centering}
\caption{\label{fig:Relative-energy-norm-1-1-1}Relative energy norm error
for the circular inclusion problem - comparison with various grids
and order of the shape functions.}
\end{figure}

\begin{figure}
\begin{centering}
\includegraphics[height=8cm]{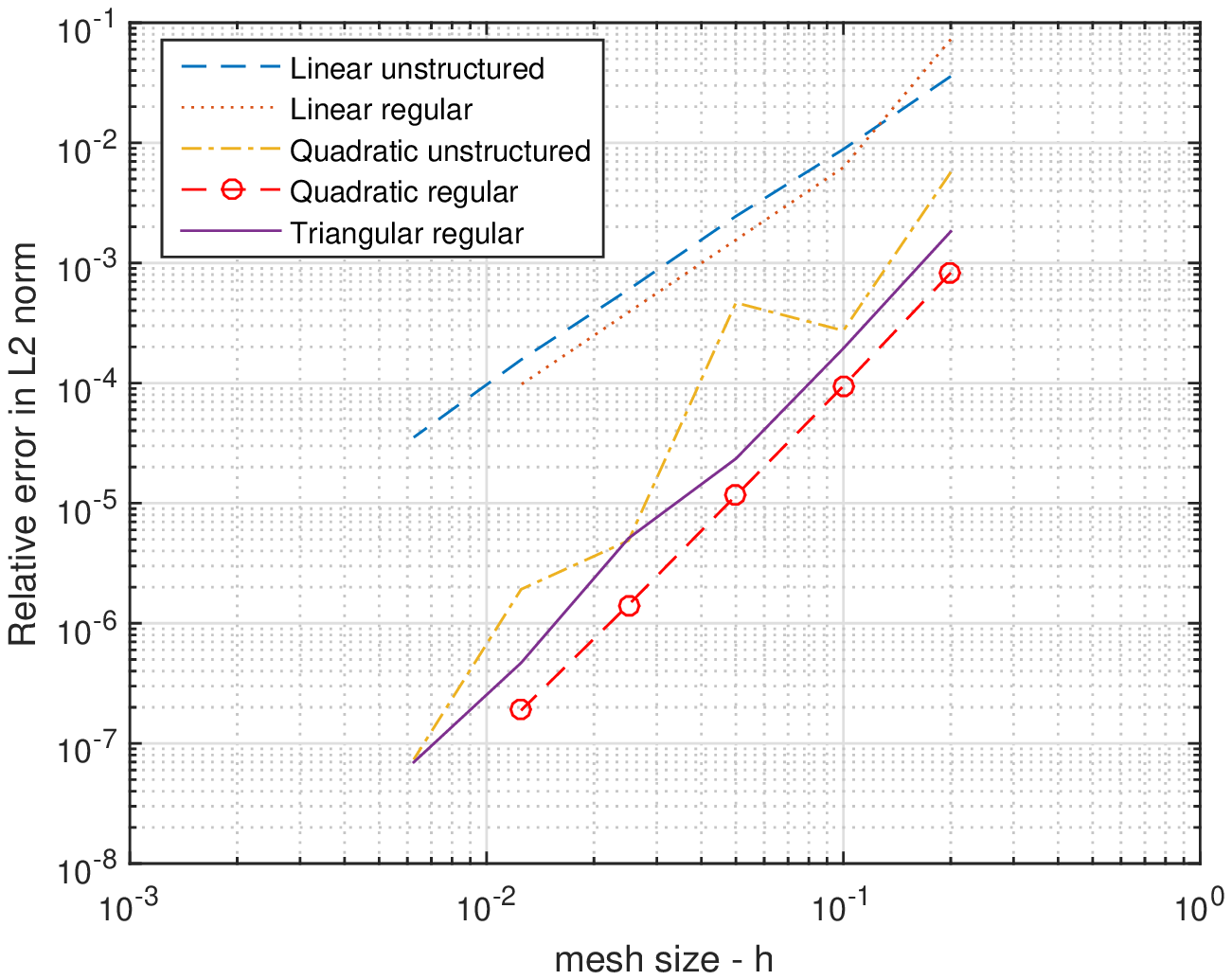}
\par\end{centering}
\caption{\label{fig:Relative-energy-norm-1-1-1-1}Relative L2 norm error for
the circular inclusion problem - comparison with various grids and
order of the shape functions.}
\end{figure}

\begin{figure}
\begin{centering}
\includegraphics[height=8cm]{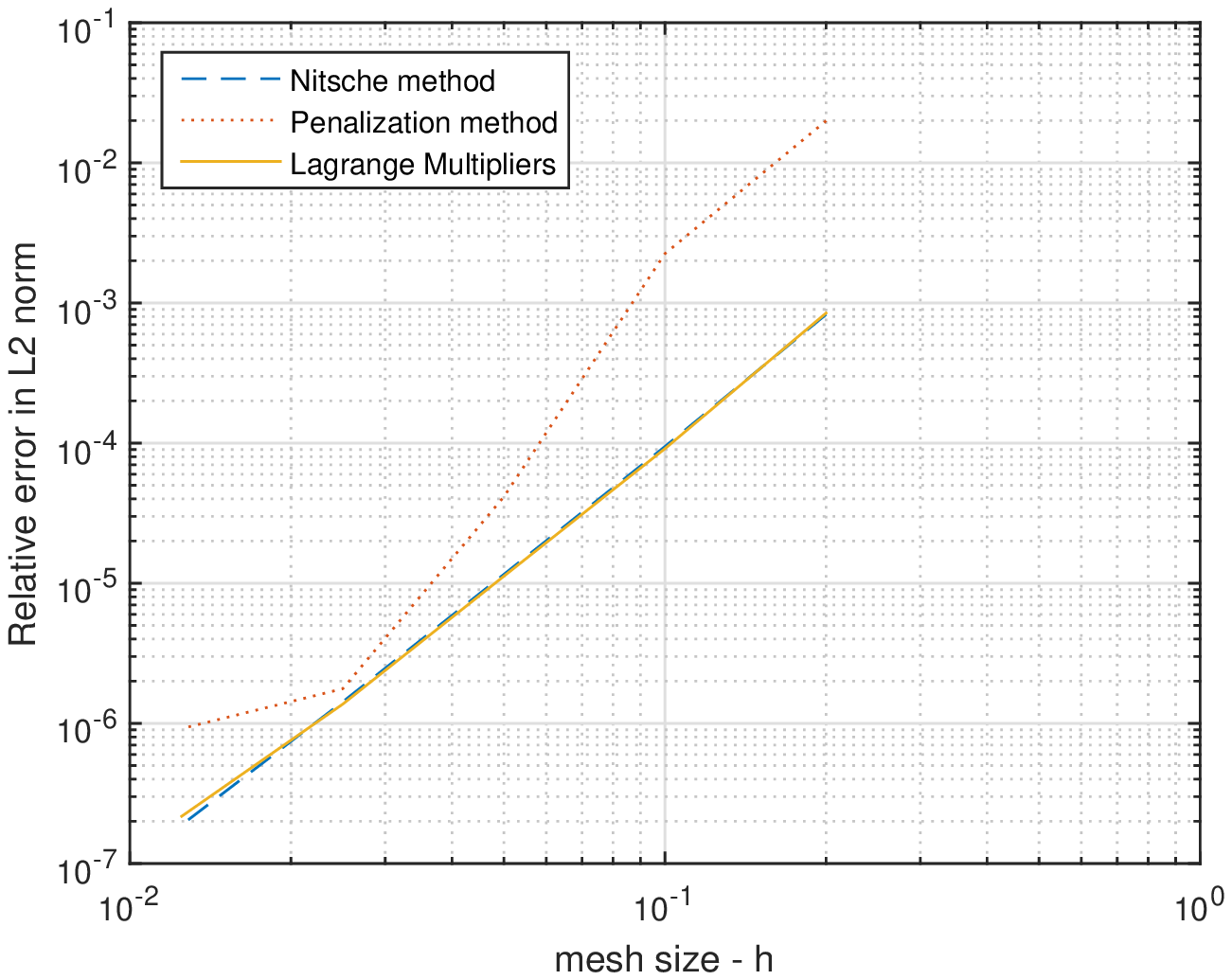}
\par\end{centering}
\caption{\label{fig:Relative-energy-norm-1-1-1-3}Relative L2 norm error for
the circular inclusion problem - comparison between the three methods
on a regular quadratic grid.}
\end{figure}

\begin{figure}
\begin{centering}
\includegraphics[height=8cm]{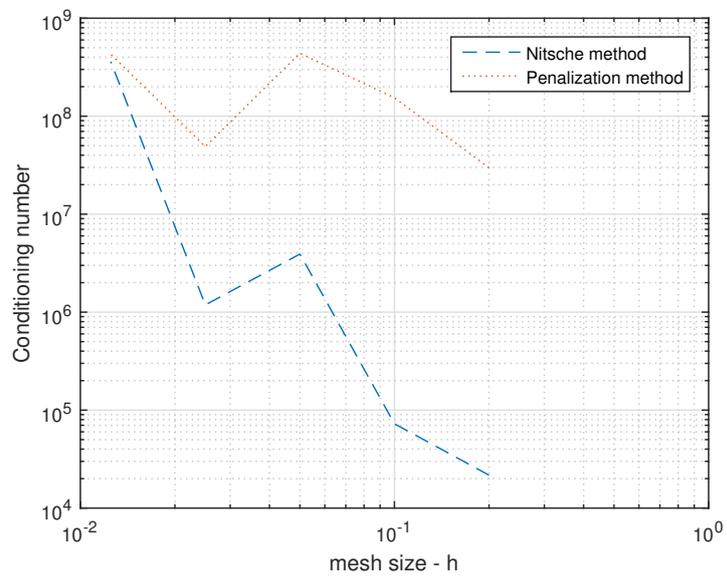}
\par\end{centering}
\caption{\label{fig:Relative-energy-norm-1-1-1-4}Condition number comparison
for the circular inclusion problem between Nitsche's method and penalty
method.}
\end{figure}

We compare the error norms by considering linear irregular, linear
regular, quadratic irregular, qudratic regular and triangular irregular
elements. Note that in quadratic elements, the shape functions assume
a quadratic nature. (Figures \ref{fig:Relative-energy-norm-1-1-1}
and \ref{fig:Relative-energy-norm-1-1-1-1}) The convergence order
is the highest for a structured grid with quadratic shape functions.
This is obvious as better approximation is obtained with higher order
of polynomial. On comparison with other methods like penalty and lagrange
multipliers (figures \ref{fig:Relative-energy-norm-1-1-1-3} and\ref{fig:Relative-energy-norm-1-1-1-4}),
we see that Nitsche's method has slightly better convergence than
penalty method. But we have seen from the previous example that we
need to prescribe very high penalization parameter to obtain this
resulting in higher conditioning of the system. Nitsche's method is
competitive on comparison to Lagrange method but comes free of additional
degree of freedom associated with Lagrange multipliers.

\begin{figure}
\begin{centering}
\includegraphics[height=8cm]{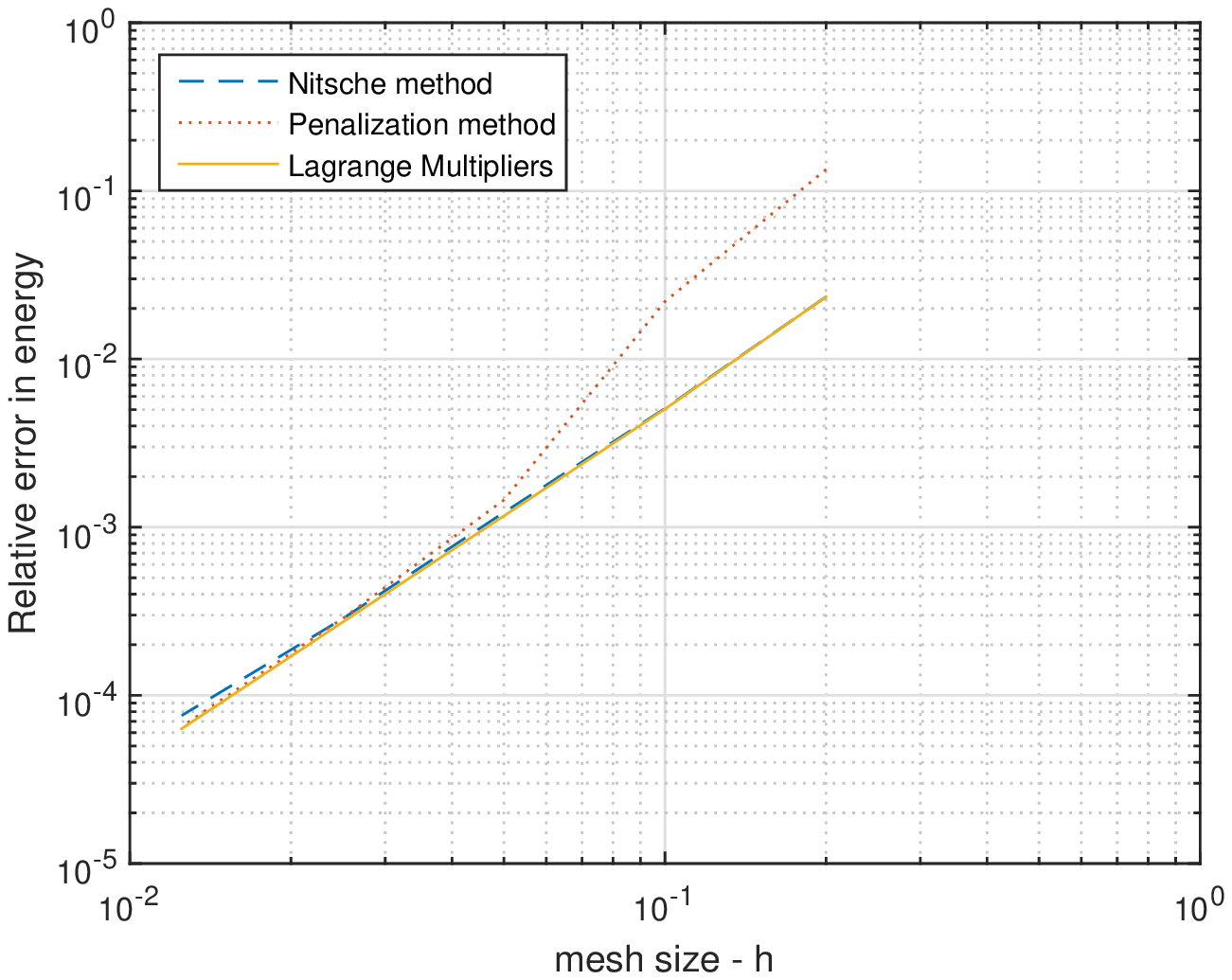}
\par\end{centering}
\caption{\label{fig:Relative-energy-norm-1-1-1-2}Relative energy norm error
for the circular inclusion problem - comparison between the three
methods on a regular quadratic grid.}
\end{figure}

\clearpage{}\newpage{}

\section{Conclusion and future scope}

As can be seen from this work, Nitsche's method was succesfully implemented
in Code-Aster using X-FEM discretization with shifted basis enrichment.
We have seen that this method is better in terms of relative error
as compared to standard penalty methods as well as relieves us from
calculating a 'free' detrimental penalization or stabilization parameter,
and instead focuses on domain dependant parameters. We calculate these
parameters based on numerical analysis, insisting on the coercivity
of the bilinear form. Also the method is better than Lagrange multipliers
method as it can calculate the solution in the same computational
range but without any extra degree of freedom that is associated with
Lagrange method. 

With the help of simple yet effective problems we have brought forth
the advantages of a method that can capture solutions whiile enforcing
internal constraints. Coupled with X-FEM, discontinuous enrichment
of finite element basis functions allows the construction of a solution
space that takes into account discontinuities at interfaces without
the disadvantage of needing to grid the interfaces.The method is straightforward
to implement, requiring only the modifications of element stiffness
routines of elements intersected by the interface. With the introduction
of a shifted basis enrichment we get increased convergence.

One of the next step would be to extend the method on tips and cracks
with internal endings. This can come in handy especially when the
two opposing lips of a crack are constrained. Also, this can be clubbed
with Nitsche's method in contact and can help analyse models with
constaint on one side and contact friction on the other with vastly
differing material properties. Also this method can be extended to
dynamic problems, like fluid flow or the seismic activities.

\newpage{}

\appendix

\section*{\addcontentsline{toc}{part}{Appendix}}

\part*{Appendix}

\section{Stabilization parameter}

The stabilization or penalization parameter, $\alpha$, is defined
such that the appropriate bilinear form is coercive. We use the local
approach here, which offers added simplicity and efficiency.\cite{key-2,key-5,key-6}
Referring back to the discrete bilinear form containing both bulk
and interfacial components, $a_{i}\left(\boldsymbol{u}^{h},\boldsymbol{v}^{h}\right)$
we further define an 'energy' norm:
\begin{equation}
\left|\left|\boldsymbol{v}^{h}\right|\right|_{E}^{2}=\left(\varepsilon\left(\boldsymbol{v}\right),\sigma\left(\boldsymbol{v}\right)\right)\label{eq:-37}
\end{equation}
Here, $\left(\centerdot,\centerdot\right)$ is the $L_{2}$ inner
product. The duality pairing $\left\langle \centerdot,\centerdot\right\rangle $
denotes integration along the interface. 

\paragraph{Dirichlet condition.}

Considering this problem as two 'one-sided' problems, we consider
one domain for simplicity, $\Omega^{-}$. We make use of the generalized
inverse estimate and there exists a configuration dependent constant
$C_{1}$, to assert coercivity, such that
\begin{equation}
\left|\left|\sigma\left(\boldsymbol{v}^{h}\right)_{,n}\right|\right|_{\Gamma_{e}^{*}}\le C_{1}\left|\left|\boldsymbol{v}^{h}\right|\right|_{\Omega_{e}^{-},E}\label{eq:-36}
\end{equation}

The gradient is constant within the element and the normal derivative
is constant along the interface, helping us obtain a lower obtain
for $C_{1}$. For the case of linear triangular element on a linear
isotropic element
\begin{equation}
\left|\left|\sigma\left(\boldsymbol{v}^{h}\right)_{,n}\right|\right|_{\Gamma_{e}^{*}}^{2}=L_{s}\left(\sigma\left(\boldsymbol{v}^{h}\right).n\right)^{2}\le L_{s}E^{2}\left|\varepsilon\left(\boldsymbol{v}^{h}\right)\right|^{2}
\end{equation}
\begin{equation}
\left|\left|\boldsymbol{v}^{h}\right|\right|_{\Omega_{e}^{-},E}^{2}=A^{-}E\left|\varepsilon\left(\boldsymbol{v}^{h}\right)\right|^{2}
\end{equation}
with $L_{s}=\mbox{meas}\left(\Gamma_{e}^{*}\right)$ and $A^{-}=\mbox{meas}\left(\Omega_{e}^{-}\right)$.
Thus we have
\[
C_{1}^{2}\ge EL_{s}/A^{-}
\]
Similarly for the case of linear tetrahedron, we can write
\[
C_{1}^{2}\ge EA_{s}/V^{-}
\]
with $A_{s}=\mbox{meas}\left(\Gamma_{e}^{*}\right)$ and $V^{-}=\mbox{meas}\left(\Omega_{e}^{-}\right)$.
Utilizing the lowest estimate for $C_{1}$ we can use, for a given
element
\begin{equation}
\alpha_{e}=2C_{1}^{2}\label{eq:-35}
\end{equation}
which provides coercivity of the bilinear form on (\ref{eq:-33}),
\begin{eqnarray}
a\left(v^{h},v^{h}\right)_{e} & = & \int_{\Omega}\varepsilon(v)\sigma\mathrm{d\Omega}-\int_{\mathcal{S}}v^{-}\left(\sigma^{-}\right).\mathbf{n}\mathrm{d}\mathrm{\Gamma}-\int_{\mathcal{S}}v^{-}\left(\sigma{}^{-}\right).\mathbf{n}\mathrm{d}\mathrm{\Gamma}+\int_{\mathcal{S}}v^{-}\alpha_{e}^{-}v^{-}\mathrm{d}\mathrm{\Gamma}\nonumber \\
 & = & \left|\left|v^{h}\right|\right|_{\Omega^{-},E}^{2}-2\left\langle \sigma^{-}\left(v^{h}\right){}_{,n},v^{h}\right\rangle _{\Gamma^{*}}+\alpha_{e}\left|\left|v^{h}\right|\right|_{\Gamma^{*}}^{2}\label{eq:-38}
\end{eqnarray}
Young's inequality (also Peter-Paul inequality) with $\epsilon>0$,
gives
\begin{equation}
2\left\langle \sigma^{-}\left(v^{h}\right){}_{,n},v^{h}\right\rangle _{\Gamma^{*}}\le\epsilon\left|\left|\sigma^{-}\left(v^{h}\right){}_{,n}\right|\right|_{\Gamma^{*}}^{2}+\frac{1}{\epsilon}\left|\left|v^{h}\right|\right|_{\Gamma^{*}}^{2}
\end{equation}
Thus, from the definition of unit vector, we have the inequality
\begin{eqnarray}
a\left(v^{h},v^{h}\right)_{e} & \ge & \left|\left|v^{h}\right|\right|_{\Omega^{-},E}^{2}-\epsilon\left|\left|\sigma^{-}\left(v^{h}\right){}_{,n}\right|\right|_{\Gamma^{*}}^{2}+\left(\alpha_{e}-\frac{1}{\epsilon}\right)\left|\left|v^{h}\right|\right|_{\Gamma^{*}}^{2}\\
 & \ge & \left(1-\epsilon C_{1}^{2}\right)\left|\left|v^{h}\right|\right|_{\Omega^{-},E}^{2}+\left(\alpha_{e}-\frac{1}{\epsilon}\right)\left|\left|v^{h}\right|\right|_{\Gamma^{*}}^{2}
\end{eqnarray}
By using $\epsilon=1/\alpha_{e}$ and $\alpha_{e}=2C_{1}^{2}$ we
get
\begin{equation}
a\left(v^{h},v^{h}\right)_{e}\ge\frac{1}{2}\left|\left|v^{h}\right|\right|_{\Omega^{-},E}^{2}\label{eq:-39}
\end{equation}
We can see that coercivity is ensured with any choice of $\alpha_{e}\ge1/\epsilon\ge C_{1}^{2}$
while (\ref{eq:-35}) provides good performance in computation.

\paragraph{Jump condition.\label{sec:Stabilization-parameter} }

The generalized inverse estimate (\ref{eq:-36}) is extended to account
for the average flux
\begin{equation}
\left|\left|\left\langle \sigma\left(\boldsymbol{v}^{h}\right)_{,n}\right\rangle \right|\right|_{\Gamma^{*}}\le C_{1}\left|\left|\boldsymbol{v}^{h}\right|\right|_{\Omega,E}
\end{equation}
in terms of energy norm (\ref{eq:-37}).

For a linear triangular element, the gradient is piecewise constant
within the element. Assuming isotropic material, $E$ is also piecewise
constant within each element, the mean flux is constant along the
interface; thus,
\begin{equation}
\left|\left|\left\langle \sigma\left(\boldsymbol{v}^{h}\right)_{,n}\right\rangle \right|\right|_{\Gamma^{*}}^{2}=L_{s}\left\langle \sigma\left(\boldsymbol{v}^{h}\right)_{,n}\right\rangle ^{2}
\end{equation}

\begin{equation}
\left|\left|\boldsymbol{v}^{h}\right|\right|_{\Omega,E}^{2}=A^{-}E^{-}\left|\varepsilon\left(\boldsymbol{v}^{h-}\right)\right|^{2}+A^{+}E^{+}\left|\varepsilon\left(\boldsymbol{v}^{h+}\right)\right|^{2}
\end{equation}
For the average flux
\begin{eqnarray}
\left\langle \sigma\left(\boldsymbol{v}^{h}\right)_{,n}\right\rangle ^{2} & = & \frac{1}{4}\left(\sigma\left(\boldsymbol{v}^{h^{+}}\right).n+\sigma\left(\boldsymbol{v}^{h^{-}}\right).n\right)^{2}\nonumber \\
\left(\sigma\left(\boldsymbol{v}^{h^{+}}\right).n+\sigma\left(\boldsymbol{v}^{h^{-}}\right).n\right) & \le & \left(1+\epsilon\right)\left(\sigma\left(\boldsymbol{v}^{h^{-}}\right).n\right)^{2}+\left(1+\frac{1}{\epsilon}\right)\left(\sigma\left(\boldsymbol{v}^{h+}\right).n\right)^{2}\mbox{ }\forall\epsilon>0
\end{eqnarray}
This follows from Young's inequality. By selecting $\epsilon=E^{+}A^{-}/E^{-}A^{+}$
we get
\begin{eqnarray}
\left\langle \sigma\left(\boldsymbol{v}^{h}\right)_{,n}\right\rangle ^{2} & \le & \left(1+\frac{E^{+}A^{-}}{E^{-}A^{+}}\right)\left(E^{-}\right)^{2}\left|\varepsilon\left(\boldsymbol{v}^{h^{-}}\right)\right|^{2}+\left(1+\frac{E^{-}A^{+}}{E^{+}A^{-}}\right)\left(E^{+}\right)^{2}\left|\varepsilon\left(\boldsymbol{v}^{h+}\right)\right|^{2}\\
 & = & \left(\frac{E^{-}}{A^{-}}+\frac{E^{+}}{A^{+}}\right)\left(A^{-}E^{-}\left|\varepsilon\left(\boldsymbol{v}^{h^{+}}\right)\right|^{2}+A^{+}E^{+}\left|\varepsilon\left(\boldsymbol{v}^{h^{-}}\right)\right|^{2}\right)
\end{eqnarray}
As can be seen, the generalized inverse estimate is satisfied for
\[
C_{1}^{2}\ge\frac{L_{s}}{4}\left(\frac{E^{-}}{A^{-}}+\frac{E^{+}}{A^{+}}\right)
\]
Similarly, for the linear tetrahedron,
\[
C_{1}^{2}\ge\frac{A_{s}}{4}\left(\frac{E^{-}}{V^{-}}+\frac{E^{+}}{V^{+}}\right)
\]

Following what has been done from (\ref{eq:-38}) to (\ref{eq:-39}),
\begin{eqnarray}
a\left(v^{h},v^{h}\right)_{e} & = & \int_{\Omega}\varepsilon(v)\sigma\mathrm{d\Omega}-\int_{\mathcal{S}}[[v]]<\sigma>.\mathbf{n}\mathrm{d\Gamma}-\int_{\mathcal{S}}[[v]]<\sigma>.\mathbf{n}\mathrm{d\Gamma}+\int_{\mathcal{S}}[[v]]\alpha[[v]]\mathrm{d\Gamma}\nonumber \\
 & = & \left|\left|v^{h}\right|\right|_{\Omega,E}^{2}-2\left\langle \left\langle \sigma\left(v^{h}\right){}_{,n}\right\rangle ,\left[\left[v^{h}\right]\right]\right\rangle _{\Gamma^{*}}+\alpha_{e}\left|\left|\left[\left[v^{h}\right]\right]\right|\right|_{\Gamma^{*}}^{2}\\
 & \ge & \left|\left|v^{h}\right|\right|_{\Omega,E}^{2}-\epsilon\left|\left|\left\langle \sigma\left(v^{h}\right){}_{,n}\right\rangle \right|\right|+\left(\alpha_{e}-\frac{1}{\epsilon}\right)\left|\left|\left[\left[v^{h}\right]\right]\right|\right|_{\Gamma^{*}}^{2}\mbox{ }\forall\epsilon>0\\
 & \ge & \left(1-\epsilon C_{1}^{2}\right)\left|\left|v^{h}\right|\right|_{\Omega,E}^{2}+\left(\alpha_{e}-\frac{1}{\epsilon}\right)\left|\left|\left[\left[v^{h}\right]\right]\right|\right|_{\Gamma^{*}}^{2}\\
 & \ge & \frac{1}{2}\left|\left|v^{h}\right|\right|_{\Omega,E}^{2}
\end{eqnarray}
This gives the same coercivity assurance as in (\ref{eq:-39}). Thus
the stabilization parameter can be chosen according to (\ref{eq:-35}).

\paragraph{Weighted parameters}

Similar to what was done in the previous section, we have, using Cauchy-Schwartz
inequality
\begin{eqnarray*}
a\left(v^{h},v^{h}\right)_{e} & = & \left|\left|v^{h}\right|\right|_{\Omega,E}^{2}-2\left\langle \left\langle \sigma\left(v^{h}\right){}_{,n}\right\rangle _{\gamma},\left[\left[v^{h}\right]\right]\right\rangle _{\Gamma^{*}}+\alpha_{e}\left|\left|\left[\left[v^{h}\right]\right]\right|\right|_{\Gamma^{*}}^{2}\\
 & \ge & \left|\left|v^{h}\right|\right|_{\Omega,E}^{2}+\alpha_{e}\left|\left|\left[\left[v^{h}\right]\right]\right|\right|_{\Gamma^{*}}^{2}-2\left|\left|\left[\left[v^{h}\right]\right]\right|\right|_{\Gamma^{*}}\left|\left|\left\langle \sigma\left(v^{h}\right)\right\rangle _{\gamma}.n\right|\right|_{\Gamma^{*}}\\
 & \ge & \left(\left|\left|v^{h}\right|\right|_{\Omega,E}-C_{1}\left|\left|\left[\left[v^{h}\right]\right]\right|\right|_{\Gamma^{*}}\right)^{2}+\left(\alpha_{e}-C_{1}^{2}\right)\left|\left|\left[\left[v^{h}\right]\right]\right|\right|_{\Gamma^{*}}
\end{eqnarray*}
with $\left|\left|\left\langle \sigma\left(v^{h}\right)\right\rangle _{\gamma}.n\right|\right|_{\Gamma^{*}}\le C_{1}\left|\left|v^{h}\right|\right|_{\Omega,E}$.
We use generalized inverse estimate, to find the lower bound for $C_{1}^{2}$.

For a linear triangular element, the gradient is piecewise constant
within the element. Assuming isotropic material, $E$ is also piecewise
constant within each element, the mean flux is constant along the
interface; thus,

\begin{equation}
\left|\left|\boldsymbol{v}^{h}\right|\right|_{\Omega,E}^{2}=A^{-}E^{-}\left|\varepsilon\left(\boldsymbol{v}^{h-}\right)\right|^{2}+A^{+}E^{+}\left|\varepsilon\left(\boldsymbol{v}^{h+}\right)\right|^{2}
\end{equation}
\begin{eqnarray}
\left|\left|\left\langle \sigma\left(\boldsymbol{v}^{h}\right)_{,n}\right\rangle _{\gamma}\right|\right|_{\Gamma^{*}}^{2} & = & L_{s}\left\langle \sigma\left(\boldsymbol{v}^{h}\right)_{,n}\right\rangle _{\gamma}^{2}\\
 & = & L_{s}\left(\gamma_{e}E\left|\varepsilon\left(\boldsymbol{v}^{h+}\right)\right|.n+\left(1-\gamma_{e}\right)E\left|\varepsilon\left(\boldsymbol{v}^{h-}\right)\right|.n\right)^{2}\\
 & \le & L_{s}\left(\left(\gamma_{e}E\left|\varepsilon\left(\boldsymbol{v}^{h+}\right)\right|\right)^{2}\left(1+\epsilon\right)+\left(\left(1-\gamma_{e}\right)E\left|\varepsilon\left(\boldsymbol{v}^{h-}\right)\right|\right)^{2}\left(1+\frac{1}{\epsilon}\right)\right)
\end{eqnarray}
from Young's inequality. By selecting $\epsilon=E^{+}A^{-}\gamma_{e}^{2}/E^{-}A^{+}\left(1-\gamma_{e}\right)^{2}$
we get
\begin{eqnarray}
\left\langle \sigma\left(\boldsymbol{v}^{h}\right).n\right\rangle _{\Gamma^{*}}^{2} & \le & \left(1+\frac{E^{+}A^{-}\gamma_{e}^{2}}{E^{-}A^{+}\left(1-\gamma_{e}\right)^{2}}\right)\left(\left(1-\gamma_{e}\right)E^{-}\left|\varepsilon\left(\boldsymbol{v}^{h^{-}}\right)\right|\right)^{2}+\left(1+\frac{E^{-}A^{+}\left(1-\gamma_{e}\right)^{2}}{E^{+}A^{-}\gamma_{e}^{2}}\right)\left(\gamma_{e}E^{+}\left|\varepsilon\left(\boldsymbol{v}^{h+}\right)\right|\right)^{2}\nonumber \\
 & = & \left(\frac{\left(1-\gamma_{e}\right)^{2}E^{-}}{A^{-}}+\frac{\gamma_{e}^{2}E^{+}}{A^{+}}\right)\left(A^{-}E^{-}\left|\varepsilon\left(\boldsymbol{v}^{h^{-}}\right)\right|^{2}+A^{+}E^{+}\left|\varepsilon\left(\boldsymbol{v}^{h^{+}}\right)\right|^{2}\right)
\end{eqnarray}
As can be seen, the generalized inverse estimate is satisfied for
\[
C_{1}^{2}\ge L_{s}\left(\frac{E^{-}\left(1-\gamma_{e}\right)^{2}}{A^{-}}+\frac{E^{+}\gamma_{e}^{2}}{A^{+}}\right)
\]
Similarly, for the linear tetrahedron,
\[
C_{1}^{2}\ge A_{s}\left(\frac{E^{-}\left(1-\gamma_{e}\right)^{2}}{V^{-}}+\frac{E^{+}\gamma_{e}^{2}}{V^{+}}\right)
\]

If we use the smart choice for $\gamma_{e}$ from (\ref{eq:-40}),
then we get 
\begin{equation}
C_{1}^{2}=\frac{L_{s}}{\left(\frac{A^{-}}{E^{-}}+\frac{A^{+}}{E^{+}}\right)}
\end{equation}
which helps us avoid numerical issues that creep up due to conforming
meshes from classical Nitsche's values for $C_{1}$. Recalling (\ref{eq:-35}),
the stabilization parameter for weighted algorithm is 
\[
\alpha_{e}=2C_{1}^{2}
\]

\newpage{}

\section*{\addcontentsline{toc}{section}{References}}

\end{document}